\documentclass[xypic,amssymb,rawfonts,amsthm,amsmath,amsfonts,syntonly,11pt]{amsart}
\usepackage[notref,notcite,final]{showkeys}
\usepackage{amssymb}

\newcounter{thecounter}

\numberwithin{thecounter}{section}
\newtheorem{lemma}[thecounter]{Lemma}
\newtheorem{prop}[thecounter]{Proposition}
\newtheorem{thm}[thecounter]{Theorem}
\newtheorem{cor}[thecounter]{Corollary}

\newtheorem{conj}[thecounter]{Conjecture}

\theoremstyle{definition}

\newtheorem{example}[thecounter]{Example}

\newtheorem{rem}[thecounter]{Remark}
\newtheorem{construction}[thecounter]{Construction}
\newtheorem{con}[thecounter]{Elaboration}
\newtheorem{notation}[thecounter]{Notation}

\numberwithin{equation}{section}

\newcommand{\onecov}{\langle 1 \rangle}

\newcommand{\map}{\operatorname{map}}

\newcommand{\K}{\operatorname{\bf K}}

\newcommand{\Top}{ \operatorname{Spaces}}

\newcommand{\ch}{\operatorname{char}}

\newcommand{\Aut}{\operatorname{Aut}}
\newcommand{\Rep}{\operatorname{Rep}}
\newcommand{\Out}{\operatorname{Out}}

\newcommand{\rk}{\operatorname{rk}}

\newcommand{\Hom}{\operatorname{Hom}}

\newcommand{\Tor}{\operatorname{Tor}}

\newcommand{\res}{\operatorname{res}}
\newcommand{\Der}{\operatorname{Der}}
\newcommand{\diag}{\operatorname{diag}}

\newcommand{\Zppmod}{\Z_{p}\mbox{-}\operatorname{mod}}

\newcommand{\covpp}[1]{\operatorname{Cov}^{p'}({#1})}
\newcommand{\covtwop}[1]{\operatorname{Cov}^{2'}({#1})}

\newcommand{\holim}{\operatorname{holim}}
\newcommand{\hocolim}{\operatorname{hocolim}}

\newcommand{\cF}{{\mathcal F}}

\newcommand{\cC}{{\mathcal{C}}}

\newcommand{\calD}{{\mathcal{D}}}

\newcommand{\coker}{\operatorname{coker}}
\newcommand{\op}{\operatorname{op}}

\newcommand{\pcom}{\hat{{}_p}}
\newcommand{\lcom}{\hat{{}_l}}
\newcommand{\threecom}{\hat{{}_3}}
\newcommand{\fivecom}{\hat{{}_5}}
\newcommand{\twocom}{\hat{{}_2}}

\newcommand{\St}{\operatorname{St}}

\newcommand{\GL}{\operatorname{GL}}
\newcommand{\SO}{\operatorname{SO}}

\newcommand{\SU}{\operatorname{SU}}
\newcommand{\PU}{\operatorname{PU}}
\newcommand{\Sp}{\operatorname{Sp}}

\newcommand{\PGL}{\operatorname{PGL}}
\newcommand{\oO}{\operatorname{O}}
\newcommand{\SL}{\operatorname{SL}}
\newcommand{\PSL}{\operatorname{PSL}}

\newcommand{\Spin}{\operatorname{Spin}}

\newcommand{\luft}{\medskip\par\noindent}
\newcommand{\QED}{\hfill {\bf $\square$}\luft}

\newcommand{\Q}{{\mathbf {Q}}}
\newcommand{\R}{{\mathbf {R}}}
\newcommand{\F}{{\mathbf {F}}}
\newcommand{\Z}{{\mathbf {Z}}}

\newcommand{\T}{{\mathbf {T}}}
\newcommand{\D}{{\mathbf D}}

\newcommand{\B}{{\mathcal{B}} }

\newcommand{\cE}{{\mathcal{E}}}
\newcommand{\A}{{\mathbf A}}

\newcommand{\W}{{\mathcal{W}}}
\newcommand{\N}{{\mathcal{N}}}

\newcommand{\C}{{\mathbf{C}}}
\newcommand{\cZ}{{\mathcal{Z}}}

\newcommand{\into}{\hookrightarrow}

\newcommand{\func}[3]{\mbox{$#1\colon #2\rightarrow #3$}}
\newcommand{\beq}{\begin{eqnarray*}}
\newcommand{\eeq}{\end{eqnarray*}}

\newcommand{\tuborg}{\left\{\begin{array}{ll}}
\newcommand{\sluttuborg}{\end{array}\right.}

\newcommand{\cd}{\operatorname{cd}}
\newfont{\bm}{msbm10}
\newcommand{\semi}{\rtimes}
\newcommand{\m}{morphism}
\newcommand{\Ob}{\operatorname{Ob}}
\newcommand{\U}{\mathrm{U}}

\input xypic
\xyoption{all}
\CompileMatrices
\newdir{ >}{{}*!/-5pt/\dir{>}} %\ar@{ >->} gives a monomorphism arrow

\newcommand{\dN}{\breve \N}
\newcommand{\dZ}{\breve \cZ}
\newcommand{\dT}{\breve T}

\newcommand{\AM}{\Phi}

\newcommand{\cB}{{\mathcal B} }

\newcommand{\longline}{\bigskip\centerline{\hbox to 5cm{\hrulefill}}\bigskip}

\begin{document}
\title{The classification of $p$-compact groups for $p$ odd} 
\author[K. Andersen, J. Grodal, J. M{\o}ller, and A. Viruel]{Kasper K. S. Andersen, Jesper Grodal, Jesper M. M{\o}ller, and Antonio Viruel}

%\date{\today}

\thanks{The first named author was supported by EU grant EEC
  HPRN-CT-1999-00119. The second named author was supported by NSF
  grant DMS-0104318, a Clay liftoff fellowship, and the Institute for
  Advanced Study for different parts of the time this research was
  carried out. The fourth named author was supported by EU grant EEC
  HPRN-CT-1999-00119, MCYT grant BFM2001-1825, and CEC-JA grant
  FQM0213. 
}
\address{Department of Mathematics, University of Aarhus, DK-8000 {{\AA}rhus C}, Denmark}
\email{kksa@math.ku.dk}
\address{Department of Mathematics, University of Chicago, Chicago, IL 60637, USA}
\email{jg@math.uchicago.edu}
\address{Department of Mathematics, University of Copenhagen, 2100
  Copenhagen OE, Denmark}
\email{moller@math.ku.dk}
\address{Departamento de \'Algebra, Geometr{\'\i}a y Topolog{\'\i}a,
  Universidad de M\'alaga, Ap. 59, \newline \indent 29080 M\'alaga, Spain}
\email{viruel@agt.cie.uma.es}

\begin{abstract}
A $p$-compact group, as defined by Dwyer and Wilkerson, is a purely homotopically defined $p$-local analog of a compact Lie group. It has long been the hope, and later the conjecture, that these objects should have a classification similar to the classification of compact Lie groups.
In this paper we finish the proof of this conjecture, for $p$ an odd prime, proving that there is a one-to-one correspondence between connected $p$-compact groups and finite reflection groups over the $p$-adic integers. We do this by providing the last, and rather intricate, piece, namely that the exceptional compact Lie groups are uniquely determined as $p$-compact groups by their Weyl groups seen as finite reflection groups over the $p$-adic integers. Our approach in fact gives a largely self-contained proof of the entire classification theorem.
\end{abstract}
\maketitle
\tableofcontents

\section{Introduction}
It has been a central goal in homotopy theory for about half a century
to single out the homotopy theoretical properties
characterizing compact Lie groups, and obtain a corresponding
classification, starting with the work of Hopf \cite{hopf41} and
Serre \cite[Ch.~IV]{serre51} on $H$-spaces and loop
spaces. Materializing old dreams of Sullivan
\cite[p.~5.96]{sullivan70} and Rector \cite{rector71}, Dwyer and
Wilkerson, in their seminal paper \cite{DW94}, introduced the notion
of a $p$-compact group, as a $p$-complete loop
space with finite mod $p$ cohomology, and proved that $p$-compact
groups have many Lie-like properties.
Even before their introduction it has been the hope \cite{RS74}, and
later the conjecture \cite{DW96,lannes96,dwyer98}, that these
objects should admit a classification much like the classification of
compact connected Lie groups, and the work toward this has
been carried out by many authors. The goal of this paper is to complete
the proof of the classification theorem for $p$ an odd prime, showing
that there is a one-to-one correspondence between connected
$p$-compact groups and finite reflection groups over the $p$-adic
integers $\Z_p$. We do this by providing the last---and rather
intricate---piece, namely that the
$p$-completions of the exceptional compact connected Lie groups are uniquely
determined as $p$-compact groups by their Weyl groups, seen as
$\Z_p$-reflection groups. In fact our method of proof gives an essentially
self-contained proof of the entire classification theorem.

We start by very briefly introducing $p$-compact groups and some
objects associated to them necessary to state the
classification theorem---we will later in the introduction return to
the history behind the various steps of the proof. We refer the reader
to \cite{DW94} for more details on $p$-compact groups and also
recommend the overview articles \cite{dwyer98,lannes96,moller95}. We point out that it is the technical advances on homotopy fixed
points by Miller \cite{miller84}, Lannes \cite{lannes92}, and others
which makes this theory possible.

A space $X$ with a loop space structure, for short a {\em loop space},
is a triple $(X,BX,e)$ where $BX$ is a pointed connected space, called
the {\em classifying space} of $X$, and $e: X \to \Omega BX$ is a
homotopy equivalence. A {\em $p$-compact group} is a loop space with
the two additional properties that $H^*(X;\F_p)$ is finite dimensional
over $\F_p$ (to be thought of as `compactness') and that $BX$ is
$\F_p$-local \cite{bousfield75}\cite[\S 11]{DW94} (or, in this
context, equivalently $\F_p$-complete \cite[Def.~I.5.1]{bk}). Often we refer to
a loop space simply as $X$. When working with a loop space we shall
only be concerned with its classifying space $BX$, since this
determines the rest of the structure---indeed, we could instead have
defined a $p$-compact group to be a space $BX$ with the above
properties. The loop space $(G\pcom,BG\pcom,e)$, corresponding to a
pair $(G,p)$  (where $p$ is a prime, $G$ a compact Lie group with
component group a finite $p$-group, and $(\cdot)\pcom$ denotes
$\F_p$-completion \cite[Def.~I.4.2]{bk} \cite[\S 11]{DW94})
is a $p$-compact group. (Note however that a compact Lie group $G$ is not
uniquely determined by $BG\pcom$, since we are only focusing on the
structure `visible at the prime $p$'; e.g., $B\SO(2n+1)\pcom \simeq
B\Sp(n)\pcom$ if $p \neq 2$, as originally proved by Friedlander
\cite{friedlander75}; see Theorem~\ref{liegroup-lattice} for a
complete analysis.)

A {\em morphism} $X \to Y$ between loop spaces is a pointed map of spaces $BX \to BY$. We say that two morphisms are {\em conjugate} if the corresponding maps of classifying spaces are freely homotopic. A morphism is called an {\em isomorphism} (or {\em equivalence}) if it has an inverse up to conjugation, or in other words if $BX \to BY$ is a homotopy equivalence. If $X$ and $Y$ are $p$-compact groups, we call a morphism a {\em monomorphism} if the homotopy fiber $Y/X$ of the map $BX \to BY$ is $\F_p$-finite.

The loop space corresponding to the $\F_p$-completed classifying space
$BT=(B\U(1)^r)\pcom$ is called a {\em $p$-compact torus} of rank $r$.
A {\em maximal torus} in $X$ is a monomorphism $i: T \to X$ such that the
homotopy fiber of $BT \to BX$ has non-zero Euler characteristic. (We
define the Euler characteristic as the alternating sum of the
$\F_p$-dimensions of the $\F_p$-homology groups.) Fundamental to the
theory of $p$-compact groups is the theorem of Dwyer-Wilkerson
\cite[Thm.~8.13]{DW94} that, analogously to the classical situation, any
$p$-compact group admits a maximal torus. It is unique in the sense
that for any other maximal torus $i': T' \to X$, there exists an
isomorphism $\varphi: T \to T'$ such that $i'\varphi$ and
$i$ are
conjugate. Note the slight difference from the classical formulation
due to the fact that a maximal torus is defined to be a
{\em map} and not a subgroup.

Fix a $p$-compact group $X$ with maximal torus $i: T \to X$ of rank
$r$. Replace the map $Bi: BT \to BX$ by an equivalent fibration, and
define the {\em Weyl space} ${\W}_X(T)$ as the topological monoid of
self-maps $BT \to BT$ over $BX$. The {\em Weyl group} is defined as
$W_X(T) = \pi_0({\W}_X(T))$ \cite[Def.~9.6]{DW94}. By
\cite[Prop.~9.5]{DW94} $W_X(T)$ is a finite group of order
$\chi(X/T)$. Furthermore, by \cite[Pf.~of~Thm.~9.7]{DW94}, if $X$ is
connected then $W_X(T)$ identifies with the set of conjugacy classes
of self-equivalences $\varphi$ of $T$ such that $i$ and $i\varphi$ are
conjugate. In other words, the canonical homomorphism $W_X(T) \to
\Aut(\pi_1(T))$ is injective, so we can view $W_X(T)$ as a subgroup of
$\GL_r(\Z_p)$, and this subgroup is independent of $T$ up to
conjugation in $\GL_r(\Z_p)$. We will therefore suppress $T$ from the
notation.

Now, by \cite[Thm.~9.7]{DW94} this exhibits $(W_X,\pi_1(T))$ as a
finite reflection group over $\Z_p$. Finite reflection groups over
$\Z_p$ have been classified for $p$ odd by Notbohm \cite{notbohm96}
extending the classification over $\Q_p$ by Clark-Ewing \cite{CE74}
 and Dwyer-Miller-Wilkerson \cite{DMW92}
(which again builds on the classification over ${\mathbf C}$ by
Shephard-Todd \cite{ST54}); we recall this classification in
Section~\ref{latticesection} and extend Notbohm's result to all
primes. Recall that a finite $\Z_p$-reflection group is a pair
$(W,L)$ where $L$ is a finitely generated free $\Z_p$-module, and $W$
is a finite subgroup of $\Aut(L)$ generated by elements $\alpha$ such
that $1 - \alpha$ has rank one. We say that two finite $\Z_p$-reflection groups $(W,L)$ and $(W',L')$
are {\em isomorphic}, if we can find a $\Z_p$-linear isomorphism $\varphi: L
\to L'$ such that the group $\varphi W \varphi^{-1}$ equals
$W'$.

Given any self-homotopy equivalence $Bf: BX \to BX$, there exists, by
the uniqueness of maximal tori, a map $B\tilde f: BT \to BT$ such
that $Bf \circ Bi$ is homotopy equivalent to $Bi \circ B\tilde
f$. Furthermore, $B\tilde f$ is unique up to the action of the Weyl
group, as is easily seen from the definitions
(cf.\ Lemma~\ref{themap}). This sets up a homomorphism $\AM:
\pi_0(\Aut(BX)) \to N_{\GL(L_X)}(W_X)/W_X$, where $\Aut(BX)$ is the
space of self-homotopy equivalences of $BX$. (This map has precursors
going back to Adams-Mahmud \cite{AM76}; see Lemma~\ref{themap} and
Theorem~\ref{ndetthm} for a more elaborate version.) The group
$N_{\GL(L_X)}(W_X)/W_X$ can be completely calculated; see
Section~\ref{autosection}.

The main classification theorem which we complete in this paper, is
the following.

\begin{thm} \label{classthm}
Let $p$ be an odd prime. The assignment that to each connected
$p$-compact group $X$ associates the pair $(W_X,L_X)$ via the canonical action of $W_X$ on $L_X =
\pi_1(T)$ defines a bijection between the set of isomorphism classes
of connected $p$-compact groups and the set of isomorphism classes of
finite $\Z_p$-reflection groups.

Furthermore, for each connected $p$-compact group $X$ the map $\AM: \pi_0(\Aut(BX)) \to N_{\GL(L_X)}(W_X)/W_X$ is an isomorphism, i.e., the group
of outer automorphisms of $X$ is canonically isomorphic to the group
of outer automorphisms of $(W_X,L_X)$.
\end{thm}

In particular this proves, for $p$ odd, Conjecture $5.3$ in
\cite{dwyer98} (see Theorem~\ref{ndetthm}). The self-map part of the
statement can be viewed as an extension to $p$-compact groups, $p$
odd, of the main result of Jackowski-McClure-Oliver
\cite{JMO92,JMO95}. Our method of proof via centralizers is `dual',
but logically independent, of the one in \cite{JMO92,JMO95} (see e.g.\
\cite{dwyer97,grodal02}).

By \cite{dw:center} the identity component of $\Aut(BX)$ is the
classifying space of a $p$-compact group $\cZ X$, which is defined to be
the {\em center} of $X$. We call $X$ {\em center-free} if $\cZ X$ is
trivial. For $p$ odd this is equivalent to $(W_X,L_X)$ being {\em center-free}, i.e., $(L_X \otimes \Z/p^\infty)^{W_X}= 0$, by \cite[Thm.~7.5]{dw:center}. Furthermore recall that a connected $p$-compact group $X$ is
called {\em simple} if $L_X\otimes \Q$ is an irreducible
$W$-representation and $X$ is called {\em exotic} if it is simple and
$(W_X,L_X)$ does not come from a $\Z$-reflection group (see
Section~\ref{latticesection}). By inspection of the classification of
finite $\Z_p$-reflection groups, Theorem~\ref{classthm} has as a
corollary that the theory of $p$-compact groups on the level of
objects splits in two parts, as has been conjectured (Conjectures
$5.1$ and $5.2$ in \cite{dwyer98}).

\begin{thm}\label{splittheory}
Let $X$ be a connected $p$-compact group, $p$ odd. Then $X$ can be
written as a product of $p$-compact groups
$$ X \cong G\pcom \times X'$$
where $G$ is a compact connected Lie group, and $X'$ is a direct
product of exotic $p$-compact groups. Any exotic $p$-compact group
is simply connected,
center-free, and has torsion free $\Z_p$-cohomology.
\end{thm}

Theorem~\ref{classthm} has both an existence and a uniqueness part to
it, the existence part being that all finite $\Z_p$-reflection groups are
realized as Weyl groups of a connected $p$-compact
group. The finite $\Z_p$-reflection groups which come from compact connected
Lie groups are of course realizable, and the finite $\Z_p$-reflection
groups where $p$ does not divide the order of the group can also
relatively easily be dealt with, as done by Sullivan
\cite[p.~5.96]{sullivan70} and Clark-Ewing \cite{CE74} long before
$p$-compact groups were officially defined. The remaining cases were
realized by Quillen \cite[\S 10]{quillen72}, Zabrodsky
\cite[4.3]{zabrodsky84}, Aguad\'e \cite{aguade89} and Notbohm-Oliver
\cite{notbohm98} \cite[Thm.~1.4]{notbohm99}. The classification of
finite $\Z_p$-reflection groups, Theorem~\ref{lattice-split},
guarantees that the construction of these examples actually enables
one to realize all finite $\Z_p$-reflection groups as Weyl groups of
connected $p$-compact groups.

The work toward the uniqueness part, to show that a connected $p$-compact group
is uniquely determined by its Weyl group, also predates the
introduction of $p$-compact groups. The quest was initiated by
Dwyer-Miller-Wilkerson \cite{DMW86,DMW92} (building on \cite{AW80}) who
proved the statement, using slightly different language, in the case
where $p$ is prime to the order of $W_X$ as well as for
$\SU(2)\twocom$ and $\SO(3)\twocom$.  Notbohm \cite{notbohm94} and
M{\o}ller-Notbohm \cite[Thm.~1.9]{MN98} extended this to a uniqueness
statement for all
$p$-compact groups $X$ where  $\Z_p[L_X]^{W_X}$ (the ring of $W_X$-invariant
polynomial functions on $L_X$) is a polynomial algebra and $(W_X,L_X)$
comes from a finite $\Z$-reflection group. Notbohm
\cite{notbohm98,notbohm99} subsequently also handled the cases where
$(W_X,L_X)$ does not come from a finite $\Z$-reflection group.  
It is worth mentioning that if $X$ has torsion free $\Z_p$-cohomology
(or equivalently, if $H^*(BX;\Z_p)$ is a polynomial algebra), then it
is straightforward to see that $\Z_p[L_X]^{W_X}$ is a polynomial
algebra (see Theorem~\ref{classicalborel}). The reverse implication is
also true, but the argument is more elaborate (see
Remark~\ref{polyvspoly} and also Theorem~\ref{borelelementaryabelian}
and Remark~\ref{torsionremark}); some of the papers quoted
above in fact operates with the a priori more restrictive assumption on $X$.

To get general statements beyond the case where $\Z_p[L_X]^{W_X}$ is a
polynomial algebra, i.e., to attack the cases where there exists
$p$-torsion in the cohomology ring, the first step is to reduce the
classification to the case of simple, center-free $p$-compact
groups. The results necessary to obtain this reduction were
achieved by the splitting theorem of Dwyer-Wilkerson \cite{dw:split}
and Notbohm \cite{notbohm00} along with properties of the center of a
$p$-compact group established by Dwyer-Wilkerson \cite{dw:center} and
M{\o}ller-Notbohm \cite{MN94}. We explain this reduction in
Section~\ref{reductionsection}; most of this reduction was already explained
by the third-named author in \cite{moller98} via different arguments.

An analysis of the classification of finite $\Z_p$-reflection groups
together with explicit calculations (see \cite{notbohm99lattice} and
Theorem~\ref{polyinvthm}) shows that, for $p$ odd, $\Z_p[L]^{W}$
is a polynomial algebra for all irreducible finite
$\Z_p$-reflection groups $(W,L)$ that are center-free, except the reflection groups coming from the
$p$-compact groups $\PU(n)\pcom$, $(E_8)\fivecom$, $(F_4)\threecom$,
$(E_6)\threecom$, $(E_7)\threecom$, and $(E_8)\threecom$. For
exceptional compact connected Lie groups the notation $E_6$ etc.\ denotes their
{\em adjoint} form.

The case $\PU(n)\pcom$ was handled by Broto-Viruel \cite{BV99}, using a
Bockstein spectral sequence argument to deduce it from the result for
$\SU(n)$, generalizing earlier partial results of Broto-Viruel
\cite{BV98} and M{\o}ller \cite{moller97puppreprint}. The remaining
step in the classification is therefore to handle the exceptional compact
connected Lie groups, in particular the problematic $E$-family at the
prime $3$, and this is what is carried out in this paper. (The fourth named
author has also given alternative proofs for $(F_4)\threecom$ and
$(E_8)\fivecom$ in \cite{viruel01} and \cite{viruel00}.)

\begin{thm}
Let $X$ be a connected $p$-compact group, for $p$ odd, with Weyl group
equal to $(W_G,L_G \otimes \Z_p)$ for $(G,p) = (F_4,3)$,
$(E_8,5)$, $(E_6,3)$, $(E_7,3)$, or $(E_8,3)$. Then $X$ is
isomorphic, as a $p$-compact group, to the $\F_p$-completion of the
corresponding exceptional group $G$.
\end{thm}

We will in fact give an essentially self-contained proof of the entire
classification Theorem~\ref{classthm}, since this comes rather
naturally out of our inductive approach to the exceptional
cases. We however still rely on the classification of finite
$\Z_p$-reflection groups (see \cite{notbohm96,notbohm99lattice} and
Sections~\ref{latticesection} and \ref{invariantsection}) as well as
the above mentioned structural results from
\cite{DW94,dw:center,MN94,dw:split,notbohm00}. We remark that
we also need not assume known a priori that `unstable Adams
operations' \cite{sullivan70,wilkerson74selfmaps,friedlander75}
exist.

The main ingredient in handling the exceptional groups, once the right inductive setup is in place, is to get
sufficiently detailed information about their many conjugacy classes
of elementary abelian $p$-subgroups, and then to use this information to show that the relevant
obstruction groups are trivial, using properties of Steinberg modules combined with formulas of Oliver
\cite{bob:steinberg} (see also \cite{grodal02}); we elaborate on this at the end of this introduction and in Section~\ref{section:ndet}.

It is possible to formulate a more topological version of the
uniqueness part of Theorem~\ref{classthm} which holds for all
$p$-compact groups ($p$ odd), not necessarily connected, which is
however easily seen to be equivalent to the first one using
\cite[Thm.~1.2]{kksa:thesis}. It should be viewed as a topological
analog of Chevalley's isomorphism theorem for linear algebraic groups
(see \cite[\S 32]{humphreys75} \cite[Thm.~1.5]{steinberg99} and
\cite{CWW74, osse97, notbohm95}).
To state it, we define the {\em maximal torus normalizer} $\N_X(T)$ to
be the loop space such that $B\N_X(T)$ is the Borel construction of
the canonical action of $\W_X(T)$ on $BT$. Note that by construction
$\N_X(T)$ comes with a morphism $\N_X(T) \to X$.
 By \cite[Prop.~9.5]{DW94},
$\W_X(T)$ is a discrete space so $B\N_X(T)$ has only two non-trivial
homotopy groups and fits into a fibration sequence $BT \to B\N_X(T)
\to BW_X$. (Beware that in general $\N_X(T)$ will not be a $p$-compact
group since its group of components $W_X$ need not be a $p$-group.)

\begin{thm}[Topological isomorphism theorem for $p$-compact groups,
  $p$ odd]\label{ndetthm}
Let $p$ be an odd prime and let $X$ and $X'$ be $p$-compact groups
with maximal torus normalizers $\N_X$ and $\N_{X'}$.
Then $X \cong X'$ if and only if $B\N_{X} \simeq B\N_{X'}$. 

Furthermore the spaces of self-homotopy equivalences $\Aut(BX)$ and
$\Aut(B\N_X)$ are equivalent as grouplike topological
monoids. Explicitly, turn $i: B\N_X \to BX$ into a fibration which we
will again denote by $i$, and let $\Aut(i)$ denote the grouplike
topological monoid of self-homotopy equivalences of the map $i$. Then
the following canonical zig-zag, given by restrictions, is a zig-zag
of homotopy equivalences:
$$
B\Aut(BX) \xleftarrow{\simeq} B\Aut(i) \xrightarrow{\simeq} B\Aut(B\N_X).
$$
\end{thm}

In the above theorem, the fact that the evaluation map $\Aut(i) \to
\Aut(BX)$ is an equivalence follows by a short general argument
(Lemma~\ref{themap}), which gives a canonical homomorphism $\AM: \Aut(BX) \xrightarrow{\simeq} \Aut(i) \to  \Aut(B\N_X)$, whereas the equivalence $\Aut(i) \to
\Aut(B\N_X)$ requires a detailed case-by-case analysis.

We point out that the classification of course gives easy, although
somewhat unsatisfactory, proofs that many theorems from Lie theory
extend to $p$-compact groups, by using that the theorem is known to be
true in the Lie case, and then checking the exotic cases.
Since the classifying spaces of the exotic $p$-compact groups
have cohomology ring a polynomial algebra, this can turn out to be
rather straightforward. In this way one for instance sees that Bott's
celebrated result about the structure of $G/T$ \cite{bott56} still
holds true for $p$-compact groups, at least on cohomology.

\begin{thm}[Bott's theorem for $p$-compact groups] \label{homogeneous}
Let $X$ be a connected $p$-compact group, $p$ odd,  with maximal torus
$T$ and Weyl group $W_X$. Then $H^*(X/T;\Z_p)$ is a free $\Z_p$-module
of dimension $|W_X|$, concentrated in even degrees.
\end{thm}

Likewise combining the classification with a case-by-case verification
for the exotic $p$-compact groups by Castellana
\cite{castellanaexotic,castellana2a}, we obtain that the Peter-Weyl
theorem holds for connected $p$-compact groups, $p$ odd:

\begin{thm}[Peter-Weyl theorem for connected $p$-compact groups]
  \label{peterweyl}
Let $X$ be a connected $p$-compact group, $p$ odd. Then there exists a
monomorphism $X \to \U(n)\pcom$ for some $n$.
\end{thm}

We also still have the `standard' formula for the fundamental group
(the subscript denotes coinvariants).

\begin{thm} \label{fundamentalgroup}
Let $X$ be a connected $p$-compact group, $p$ odd. Then 
$$\pi_1(X) = (L_X)_{W_X}$$ 
\end{thm}

The classification also gives a verification that results of Borel,
Steinberg, Demazure, and Notbohm \cite[Prop.~1.11]{notbohm99} extend
to $p$-compact groups, $p$ odd. Recall that an elementary abelian $p$-subgroup of $X$ is just a monomorphism $\nu: E \to X$, where $E \cong (\Z/p)^r$ for some $r$.
\begin{thm} \label{borelelementaryabelian} Let $X$ be a connected
  $p$-compact group, $p$ odd. The following conditions are equivalent:
\begin{enumerate}
\item  \label{borel1} $X$ has torsion free $\Z_p$-cohomology.
\item  \label{borel4} $BX$ has torsion free $\Z_p$-cohomology.
\item  \label{borel2} $\Z_p[L_X]^{W_X}$ is a polynomial algebra over $\Z_p$.
\item  \label{borel3} All elementary abelian $p$-subgroups of $X$
  factor through a maximal torus.
\end{enumerate}
\end{thm}
(See also Theorem~\ref{classicalborel} for equivalent formulations of condition $(1)$.)
Even in the Lie case, the proof of the above theorem is still not
entirely satisfactory despite much effort---see the comments
surrounding our proof in Section~\ref{consequencessection} as well as
Borel's comments \cite[p.~775]{boreloeuvresII} and the references
\cite{borel61,demazure73,steinberg75}. The {\em
  centralizer} $\cC_X(\nu)$ of an elementary abelian $p$-subgroup
$\nu: E \to X$ is defined as $\cC_X(\nu) = \Omega \map(BE,BX)_{B\nu}$;
cf.\ Section~\ref{section:ndet}. The following related result from Lie
theory also holds true.

\begin{thm} \label{borelelementaryabelianII}
Let $X$ be a connected $p$-compact group, $p$ odd. Then the following
conditions are equivalent:
\begin{enumerate}
\item
$\pi_1(X)$ is torsion free.
\item
Every rank one elementary abelian $p$-subgroup $\nu: \Z/p \to X$ has
connected centralizer $\cC_X(\nu)$.
\item
Every rank two elementary abelian $p$-subgroup factors through a
maximal torus.
\end{enumerate}
\end{thm}

Results about $p$-compact groups can in general, via Sullivan's
arithmetic square, be translated into results about finite loop
spaces, and the last theorem in this introduction is an example of
such a translation. (For another instance see \cite{ABGP03}.)
Recall that a finite loop space is a loop space $(X,BX,e)$, where $X$
is a finite CW-complex. A maximal torus of a finite loop space is
simply a map $B\U(1)^r \to BX$ for some $r$, such that the homotopy
fiber is homotopy equivalent to a finite CW-complex of non-zero
Euler characteristic. The classical maximal torus conjecture (stated
in 1974 by Wilkerson \cite[Conj.~1]{wilkerson74} as ``a popular
conjecture toward which the author is biased''), asserts that
compact connected Lie groups are the {\em only} connected finite loop
spaces which admit maximal tori. A slightly more elaborate version
states that the classifying space functor should set up a bijection
between isomorphism classes of compact connected Lie groups and
isomorphism classes of connected finite loop spaces admitting a
maximal torus, under which the outer automorphism group of the Lie
group $G$ equals the outer automorphism group of the corresponding loop
space $(G,BG,e)$. (The last part is true by \cite[Cor.~3.7]{JMO95}.) It is well
known that a proof of the conjectured classification of $p$-compact
groups for all primes $p$ would imply the maximal torus
conjecture. Our results at least imply that the conjecture is true
after inverting the single prime $2$.

\begin{thm}\label{Tconj}
Let $X$ be a connected finite loop space with a maximal torus. Then there
exists a compact connected Lie group $G$ such that
$BX[\frac{1}{2}]$ and $BG[\frac{1}{2}]$ are homotopy equivalent
spaces, where $[\frac{1}{2}]$ indicates $\Z[\frac{1}{2}]$-localization.
\end{thm}

\subsection*{Relationship to the Lie case and the conjectural picture
  for $p=2$}

We now state a common formulation of both the
classification of compact connected Lie groups and the classification of
connected $p$-compact groups for $p$ odd, which conjecturally should also hold
for $p=2$. We have not encountered this---in our opinion quite
natural---description before in the literature (compare \cite{dwyer98}
and \cite{lannes96}).

Let $R$ be an integral domain and $W$ a finite $R$-reflection
group. For an $RW$-{\em lattice} $L$ (i.e., an $RW$-module which is finitely
generated and free as an $R$-module) define $SL$ to be the sublattice
of $L$ generated by $(1-w)x$ where $w \in W$ and $x \in L$.
Define an $R$-{\em reflection datum} to be a triple $(W,L,L_0)$ where
$(W,L)$ is a
finite $R$-reflection group and $L_0$ is an $RW$-lattice such that $SL
\subseteq L_0 \subseteq L$ and $L_0$ is isomorphic to $SL'$ for some
$RW$-lattice $L'$.  (If $R = \Z_p$, $p$ odd, then `$S$' is
idempotent, since $W$ is generated by elements of order prime to $p$ so $H_1(W;L_W) = 0$, and therefore $L_0 = SL$ in
this case.)
Two reflection data
$(W,L,L_0)$ and $(W',L',L_0')$ are said to be isomorphic if there
exists an $R$-linear isomorphism $\varphi: L \to L'$ such that
$\varphi W \varphi^{-1} = W'$ and $\varphi(L_0) = L_0'$.

If $\calD$ is either the category of compact connected Lie groups or
connected $p$-compact groups, then we can consider the assignment
which to each object $X$ in $\calD$ associates the triple $(W,L,L_0)$,
where $W$ is the Weyl group, $L = \pi_1(T)$ is the {\em integral lattice}, and 
$L_0 = \ker( \pi_1(T) \to \pi_1(X))$ is the {\em coroot lattice}.

Theorems~\ref{classthm} and \ref{fundamentalgroup} as well as the
classification of compact connected Lie groups \cite[\S 4, no.\
9]{bo9} can now be reformulated as follows:

\begin{thm}\label{combithm}
Let $\calD$ be the category of compact connected Lie groups, $R = \Z$, or
connected $p$-compact groups for $p$ odd, $R = \Z_p$.  For $X$ in $\calD$ the associated triple $(W,L,L_0)$ is an $R$-reflection
datum and this assignment sets up a
bijection between the objects of $\calD$ up to isomorphism and
$R$-reflection data up to isomorphism. Furthermore the group of outer
automorphisms of $X$ equals the group of outer automorphisms of the
corresponding $R$-reflection datum.
\end{thm}

\begin{conj} \label{ptwoconj}
Theorem~\ref{combithm} is also true if $\calD$ is the category of
connected $2$-compact groups.
\end{conj}

One can check that the conjecture on objects is equivalent to the
conjecture given in \cite{dwyer98} and \cite{lannes96}, and the
self-map statement would then follow from \cite[Cor.~3.5]{JMO95}
and \cite[Thm.~3.5]{notbohm00di4preprint}. The role of the coroot lattice $L_0$ in the above theorem and conjecture is in fact only to be
able to distinguish direct summands isomorphic to $\SO(2n+1)$ from
direct summands isomorphic to $\Sp(n)$,
cf.\ Theorem~\ref{liegroup-lattice}.
Alternatively one can use the extension class $\gamma \in H^3(W;L)$
of the maximal torus normalizer (see
Section~\ref{nautomorphismsection}) rather than $L_0$ but in that
picture it is not a priori clear which triples $(W,L,\gamma)$ are
realizable. It would be desirable to have a `topological' version of
Theorem~\ref{combithm} and Conjecture~\ref{ptwoconj}, i.e., statements
on the level of automorphism spaces like Theorem~\ref{ndetthm}, but we
do not know a general formulation which incorporates this feature.

\subsection*{Organization of the paper}
The paper is organized around Section~\ref{section:ndet} which sets up the framework of the proof and gives an inductive proof of the main theorems, referring to the later sections of the paper for many key statements. The remaining sections can be read in an almost arbitrary order. We want to briefly sketch how these sections are used.

We first say a few words about Section~4--7, before describing Section~\ref{section:ndet} and the later sections in a little more detail.
The short Sections~\ref{themapsection} and \ref{nautomorphismsection} construct the map $\Phi: \Aut(BX) \to \Aut(B\N_X)$ and give an algebraic description of the automorphisms of $B\N_X$. Section~\ref{reductionsection} contains the reduction to the case of simple, center-free, connected $p$-compact groups.
In Section~\ref{nakajimasection} we prove an integral version of a theorem of Nakajima, and show how this leads to an easy criterion for inductively constructing certain $p$-compact groups; this criterion will, in the setup of the induction, lead to a construction of the exotic $p$-compact groups and show that they have torsion free $\Z_p$-cohomology. 

Armed with this information let us now summarize Section~\ref{section:ndet}. In the inductive framework of the main theorem the results in Section~\ref{nakajimasection} guarantee that we have concrete models for conjecturally all $p$-compact groups, and that those coming from exotic finite $\Z_p$-reflection groups have torsion free $\Z_p$-cohomology. Likewise, by the reduction theorems in Section~\ref{reductionsection}, we are furthermore reduced to showing that if $X'$ is an unknown connected center-free simple $p$-compact group with associated $\Z_p$-reflection group $(W,L)$ then it agrees with our known model $X$ realizing $(W,L)$. We want, using the inductive assumption, to construct a map from the centralizers in $X$ to $X'$, and show that these maps glue together to give an isomorphism $X \to X'$. To be able to glue the maps together, we need to have a preferred choice on each centralizer and know that these agree on the intersection---this is why we also have to keep track of the automorphisms of $p$-compact groups in our inductive hypothesis. If $X$ has torsion free $\Z_p$-cohomology, then every elementary abelian $p$-subgroup factors through the maximal torus, and it follows from our construction that our maps on the different centralizers of elementary abelian $p$-subgroups in $X$ to $X'$ match up, as maps in the {\em homotopy category}. This is not obvious in the case where $X$ has torsion in its $\Z_p$-cohomology, and we develop tools in Section~\ref{section:ndetlemmas} which suffice to handle all the torsion cases, on a case-by-case basis. This step should be thought of as inductively showing that $X$ and $X'$ have the same (centralizer) fusion. We now have to rigidify our maps on the centralizers from a consistent collection of maps in the homotopy category to a consistent collection map in the category of spaces. There is an obstruction theory for dealing with this issue. Again, in the case where $X$ does not have torsion there is a general argument for showing that these obstruction groups vanish, whereas we in the case where $X$ has torsion have to show this on a case-by-case basis. To deal with this we give in the purely algebraic Section~\ref{liesection} complete information about all non-toral elementary abelian $p$-subgroups of the projective unitary groups and the exceptional compact connected Lie groups, along with their Weyl groups and centralizers. This information is needed as input in Section~\ref{obstructionssection} for showing that the obstruction groups vanish. Hence we get a map in the category of spaces from the centralizers in $X$ to $X'$, which then glues together to produce a map $X \to X'$ which then by our construction is easily seen to be an isomorphism. As a by-product of the analysis we also conclude that $X$ has the right automorphism group. This proves the main theorems. Section~\ref{consequencessection} establishes the consequences of the main theorem, listed in the introduction.

There are three appendices: In Section~\ref{latticesection} we give a concise classification
of finite $\Z_p$-reflection groups generalizing Notbohm's
classification to all primes. In
Section~\ref{invariantsection} we recall Notbohm's results on
invariant rings of finite $\Z_p$-reflection groups. These facts are all used multiple times in the proof. Finally in Section~\ref{autosection} we briefly calculate the
outer automorphism groups of the finite $\Z_p$-reflection groups
to make the automorphism statement in the main result more explicit.

\subsection*{Notation}

We have tried to introduce the definitions relating to $p$-compact
group as they are used, but it is nevertheless probably helpful
for the reader unfamiliar with
$p$-compact groups to keep copies of the excellent papers \cite{DW94}
and \cite{dw:center} of Dwyer-Wilkerson (whose terminology we follow)
within reach. As a technical term we say that a $p$-compact group $X$
is {\em determined by} $\N_X$ if any $p$-compact group
$X'$ with the same maximal torus normalizer is isomorphic to $X$
(which will be true for all $p$-compact groups, $p$ odd, by
Theorem~\ref{ndetthm}).

We tacitly assume that any space in this paper has the homotopy type
of a CW-complex, if necessary replacing a given space by the
realization of its singular complex \cite{may67}.

\subsection*{Acknowledgments}
We would like to thank
H.~H.~Andersen, D.~Benson, G.~Kemper, A.~Kleschev, G.~Malle, and
J.-P.~Serre for helpful correspondence.
We also thank J.~P.~May, H.~Miller, and the referee for their comments and suggestions. 
We would in particular like to thank W.~Dwyer,
D.~Notbohm, and C.~Wilkerson for several useful tutorials on their
beautiful work, which this paper builds upon.

%%%%%%%%%%%%%%%%%%%%%%%%%%%%%%%%%%%%%%%%%%%%%%%%%%%%%%%%%%%%%%%%%%%%%%%%%%%

\section{Skeleton of the proof of the main Theorems~\ref{classthm} and \ref{ndetthm}}\label{section:ndet}

The purpose of this section is to give the skeleton of the proof of the main
Theorems~\ref{classthm} and \ref{ndetthm}, but in the proofs referring
forward to the remaining sections in the paper for the proof of many key statements, as explained in the organizational remarks in the introduction.

We start by explaining the strategy in general terms, which is done
via a grand induction---for simplicity we focus first on the
uniqueness statement. Suppose that $X$ is a known $p$-compact group
and $X'$ is another $p$-compact group with the same maximal torus
normalizer. We want to construct an isomorphism $X \to X'$, by
decomposing $X$ in terms of centralizers of its non-trivial elementary
abelian $p$-subgroups, as we will explain below. Using an inductive assumption we can construct a homomorphism from each of these centralizers to $X'$, and we want see
that we can do this in a coherent way, so that they glue together to
give the desired map $X \to X'$.

We first explain the centralizer decomposition.
It is a theorem of Lannes \cite[Thm.~3.1.5.1]{lannes92} and
Dwyer-Zabrodsky \cite{DZ87} (see also \cite[Thm.~3.2]{JMO92}),
that for an elementary abelian $p$-group $E$ and a compact Lie group $G$ with component group a $p$-group, we have a homotopy equivalence
$$\coprod_{\nu \in \Rep(E,G)} \hspace{-10pt}BC_G(\nu(E))\pcom \xrightarrow{\simeq} \map(BE,BG\pcom)$$
induced by the adjoint of the canonical map
$BE \times
BC_G(\nu(E)) \to BG$. Here $\Rep(E,G)$ denotes the set of
homomorphisms $E \to G$, modulo conjugacy in $G$.

Generalizing this, one defines, for a $p$-compact group $X$, an 
 {\em elementary abelian $p$-subgroup} of $X$ to be a monomorphism $\nu: E \to X$, and its
 {\em centralizer} to be
the $p$-compact group $\cC_X(\nu)$ with classifying space $B\cC_X(\nu)
= \map(BE,BX)_{B\nu}$. 
By a theorem of Dwyer-Wilkerson
\cite[Props.~5.1~and~5.2]{DW94} this actually is a $p$-compact
group and the evaluation map to $X$ is a monomorphism.
 Note however that
$\cC_X(\nu)$ is not defined as a subobject of $X$, i.e., the map to $X$
is defined in terms of $\nu$, unlike in the Lie case.

For a $p$-compact group $X$, let $\A(X)$ denote the Quillen category
of $X$. The objects of $\A(X)$ are conjugacy classes of monomorphisms 
\func{\nu}{E}{X} of non-trivial elementary abelian $p$-groups $E$
into $X$. The morphisms $(\nu: E \to X) \rightarrow (\nu': E' \to X)$
of $\A(X)$ consists of all group monomorphisms \func{\rho}{E}{E'}
such that $\nu$ and $\nu'\rho$ are conjugate.

The centralizer construction gives a functor
\begin{equation*}\label{eq:BCX}
  B\cC_X: \A(X)^{\op} \to \Top
\end{equation*}
that takes the mono\m\ $(\nu: E \to X) \in\Ob(\A(X))$ to its centralizer
$B\cC_X(\nu) = \map(BE,BX)_{B\nu}$ and a morphism $\rho$  to
composition with $B\rho: BE \to BE'$.

The centralizer decomposition theorem of Dwyer-Wilkerson \cite[Thm.~8.1]{dw:center}, generalizing a
theorem for compact Lie groups by Jackowski-McClure \cite[Thm.~1.3]{JM92}, says that the
evaluation map
\begin{equation*}
  \hocolim_{\A(X)}B\cC_X \to BX
\end{equation*}
induces an isomorphism on mod $p$ homology.
If $X$ is connected and center-free, then for all $\nu$, the centralizer
$\cC_X(\nu)$ is a $p$-compact group with smaller cohomological dimension, hence
setting the stage for a proof by induction, cf.\ \cite[\S 9]{dw:center}. (The {\em cohomological
dimension} of a $p$-compact group $X$ is defined as $\cd(X) = \max\{ n
| H^n(X;\F_p) \neq 0 \}$; see \cite[Def.~6.14]{DW94} and
\cite[Lem.~3.8]{dw:split}.)

To make use of this we need a way to construct
a map from centralizers of elementary abelian $p$-subgroups in $X$
to any other $p$-compact group $X'$ with the same
maximal torus normalizer $\N$. Let $\N$ be embedded via homomorphisms $j: \N \to X$ and $j': \N \to X'$ respectively.
If $\nu: E \to X$ can be factored through a maximal torus $i: T \to X$, i.e.,
if there exists
$\mu: E \to T$ such that $i\mu = \nu$, then $\mu$ is unique up to
conjugation as a map to $\N$
by \cite[Prop.~3.4]{dw:split}, and furthermore by
\cite[Pf.~of~Thm.~7.6(1)]{dw:center}, $\cC_\N(\mu)$ is a maximal torus normalizer
in $\cC_X(\nu)$, where centralizers in $\N$ are defined in the same way as in a $p$-compact group. In this case $j'\mu$ will be an elementary abelian
$p$-subgroup of $X'$, which we have assigned without making any
choices, and $\cC_{X'}(j'\mu)$ will have maximal torus normalizer
$\cC_\N(\mu)$. Suppose that $\cC_X(\nu)$ is
determined by $\N_{\cC_X(\nu)}$ (i.e., any $p$-compact group with maximal torus normalizer isomorphic to $\N_{\cC_X(\nu)}$ is isomorphic to $\cC_X(\nu)$) and that the homomorphism $\AM: \Aut(B\cC_X(\nu)) \to \Aut(B\N_{\cC_X(\nu)})$, defined after Theorem~\ref{ndetthm}, is an equivalence. Since $\cC_X(\nu)$ is determined by its maximal torus normalizer, surjectivity of $\pi_0(\Phi)$ implies that there exists an isomorphism $h_\nu$ making the diagram
\begin{equation}\label{rankonediagram}
\xymatrix{
 & \cC_\N(\mu) \ar[dr]^-{j'} \ar[dl]_-{j} & \\
\cC_X(\nu) \ar[rr]^-{h_{\nu}}_-{\cong} & & \cC_{X'}(j'\mu)
}
\end{equation}
commute, and $h_\nu$ is unique up to conjugacy, by the injectivity of $\pi_0(\Phi)$. (In fact the space of such $h_\nu$ is contractible, since $\Phi$ is an equivalence.) This constructs the desired map $\varphi_\nu:\cC_X(\nu) \xrightarrow{h_\nu} \cC_{X'}(j'\mu) \to X'$ for elementary abelian $p$-subgroups $\nu: E \to X$ which factor
through the maximal torus. An elementary abelian $p$-subgroup is called {\em toral} if it has this property, and {\em non-toral} if not.

We want to construct maps also for non-toral elementary abelian $p$-subgroups, by utilizing
the centralizers of rank one elementary abelian $p$-subgroups, which are always toral by
\cite[Prop.~5.6]{DW94} if $X$ is connected. For this we need to recall
the construction of adjoint maps.

\begin{construction}[Adjoint maps] \label{adjoint}
Let $A$ be an abelian $p$-compact group (i.e., a $p$-compact group such that $\cZ A \to A$ is an isomorphism), $X$ a $p$-compact group,
and $\nu: A \to X$ a homomorphism. Suppose that $E$ is a elementary abelian $p$-subgroup of $A$ and note that we have a canonical map
$$BA \times BE
\xrightarrow{\text{mult}} BA \to BX$$
whose homotopy class only depends on the conjugacy class of $\nu$. Since furthermore
$$\pi_0(\map(BA \times BE,BX)) \cong \coprod_{\xi \in [BE,BX]} \pi_0( \map(BA,\map(BE,BX)_\xi))$$
every element $\nu: A \to X$ gives rise to an element
$\tilde \nu: A \to \cC_X(\nu|_{E})$ making the diagram
$$
\xymatrix{
& \cC_X(\nu|_{E}) \ar[d]^-{\text{ev}} \\
 A \ar[ur]^-{\tilde \nu} \ar[r]^-\nu & X
}
$$
commutative. Here $\tilde \nu$ is well-defined up to conjugacy in terms of
the conjugacy class of $\nu$. We will always use the notation
$\widetilde{(\cdot)}$ for this construction.
\end{construction}

Let $\nu: E \to X$ be an {\em arbitrary} non-trivial elementary
abelian $p$-subgroup of a connected $p$-compact group $X$ and let $V$
be a rank one subgroup of $E$. Then $\nu|_{V}$ is toral by
\cite[Prop.~5.6]{DW94}, i.e., it factors through $T$ and the map $\mu: V \to T \to \N$ is unique up to conjugation in $\N$. Furthermore if $\cC_X(\nu|_{V})$ is determined by $\cC_\N(\mu)$ and
$\AM: \Aut(B\cC_X(\nu|_{V})) \xrightarrow{\cong} \Aut(B\N_{\cC_X(\nu|_{V})})$ then
 $h_{\nu|_{V}}$ is defined as before, and we can look at the composite
$$
\varphi_{\nu,V}:\cC_X(\nu) \longrightarrow \cC_X(\nu|_{V})
\xrightarrow[\cong]{h_{\nu|_V}} \cC_{X'}(j'\mu) \longrightarrow X'.
$$
This is the definition we will use in general.
It is easy to see using adjoint maps that this construction
generalizes the previous one in the case where $\nu$ is toral, under suitable
inductive assumptions (cf.\ the proof of
Theorem~\ref{collectlemma} below). However if $\nu$ is non-toral it is not
obvious that this map is independent of the choice of subgroup $V$
of $E$, which is needed in order to get a map (in the homotopy category) from the centralizer diagram of $BX$ to $BX'$.
Checking that this is the case basically amounts to inductively
establishing that elementary abelian $p$-subgroups and their
centralizers are conjugate in the same way in $X$ and $X'$, i.e., that
they have the same fusion.
Furthermore we want see that
this diagram can be rigidified to a diagram in the category of spaces,
to get an induced map from the homotopy colimit of the
centralizer diagram. The next theorem states precisely what needs to be
checked---the calculations to verify that these conditions are indeed
satisfied for all simple center-free $p$-compact groups is essentially the contents of the
rest of the paper.

\begin{thm}\label{collectlemma}
Let $X$ and $X'$ be two connected
$p$-compact groups with the same maximal torus normalizer $\N$
embedded via $j$ and $j'$ respectively. Assume that $X$ satisfies the following
inductive assumption:
\begin{itemize}
\item[($\star$)]\label{inductive} For all rank one
elementary abelian $p$-subgroups $\nu: E \to X$ of $X$ the centralizer
$\cC_X(\nu)$ is determined by $\N_{\cC_X(\nu)}$ and $\AM: \Aut(B\cC_X(\nu)) \xrightarrow{\cong} \Aut(B\N_{\cC_X(\nu)})$ when $\nu$ is of
rank one or two.
\end{itemize}
Then we have the following.
\begin{enumerate}
\item \label{del1}
Assume that for every rank two non-toral elementary abelian $p$-subgroup $\nu: E \to X$ 
the induced map $\varphi_{\nu,V}$ is independent of the choice of the rank one
subgroup $V$ of $E$.
Then there exists a map in the homotopy category
of spaces from the centralizer
diagram of $BX$ to $BX'$ (seen as a constant diagram), i.e., an
element in $\lim^0_{\nu \in \A(X)}[B\cC_X(\nu),BX']$, given via the maps
$\varphi_{\nu,V}$ described above.
\item \label{del2}
Assume furthermore that $\lim^i_{\nu \in \A(X)}\pi_j(B\cZ\cC_X(\nu))=0$ for
$j=1,2$ and $i=j, j+1$. Then there is a lift of this element in
$\lim^0$ to a map in the (diagram) category of spaces. This produces an
isomorphism \func{f}{X}{X'} under $\N$, unique up to conjugacy, and $\AM:
\Aut(BX) \xrightarrow{\cong} \Aut(B\N)$.
\end{enumerate}
\end{thm}

\begin{proof}
As explained before the theorem if $\nu: E \to X$ has rank one, then
$\nu$ factors through $T$ to give a map $\mu: E \to \N$, unique up to
conjugation in $\N$, so the inductive assumption ($\star$) guarantees that we can construct a map $\cC_X(\nu)
\to X'$, under $C_\N(\mu)$, and this map is well defined up to
conjugation in $X'$.

We now want to see that in the case where $E$ has rank two, the map
$\varphi_{\nu,V}$ is in fact independent of the choice of the rank
one subgroup $V$.
Assume first that $\nu : E \to X$ is toral and let $\mu: E \to \N$ be a factorization of $\nu$ through $T$. By adjointness we have the following commutative diagram.
\begin{equation}\label{cooldiag}
\xymatrix{
 & & \cC_\N(\mu) \ar[ddll]_-j \ar[ddrr]^-{j'} & & \\
 & & \cC_{\cC_\N(\mu|_{V})}(\tilde \mu) \ar[u]_-\cong \ar[dl] \ar[dr] \\
\cC_X(\nu) & \cC_{\cC_X(\nu|_{V})}(\tilde \nu) \ar[l]_-\cong
\ar[rr]_-{\cong}^-{\widetilde{h_{\nu|_{V}}}} & & \cC_{\cC_{X'}(j'\mu|_{V})}(\widetilde{j'\mu}) \ar[r]^-\cong & \cC_{X'}(j'\mu)
}
\end{equation}
where $\widetilde{h_{\nu|_{V}}}$ is the map induced from $h_{\nu|_{V}}$ on
the centralizers. The rank two uniqueness assumption ($\star$) now guarantees that the bottom left-to-right composite $\psi$ is independent of the choice of $V$.

However, for any particular choice of $V$ we have a commutative diagram
$$
\xymatrix{
\cC_X(\nu) \ar[r]^-\psi_-\cong \ar[d] & \cC_{X'}(j'\mu) \ar[d] \ar[dr] \\
\cC_X(\nu|_{V}) \ar[r]^-{h_{\nu|_{V}}}_-\cong & \cC_{X'}(j'\mu|_{V}) \ar[r] & X'
}
$$
and since $\psi$ is independent of $V$ this shows that
$\varphi_{\nu,V}$ is independent of $V$ as wanted. This handles the
rank two toral case. For the rank two non-toral case we are simply
assuming that $\varphi_{\nu,V}$ is independent of $V$. (Note that the problem which prevents the toral argument to carry over to the non-toral case is that we cannot choose a uniform $\mu: E \to \N$ such that $\mu|_{V}$ factors through $T$ for all $V$, since this would imply that $E$ itself was toral.)

The fact that $\varphi_{\nu,V}$ is independent of $V$ when $E$ is of
rank two implies the statement in general: Let $\nu : E \to X$
be an elementary abelian $p$-subgroup of rank at least three. If $V_1$ and $V_2$ are two different rank
one subgroups of $E$, we set $U = V_1 \oplus V_2$ and consider the
following diagram
$$
\xymatrix{
 & \cC_X(\nu|_{V_1})  \ar[dr]^-{\varphi_{\nu|_{V_1}}} & \\
\cC_X(\nu) \ar[ur] \ar[r] \ar[dr] & \cC_X(\nu|_{U}) \ar[u] \ar[d] & X' \\
 & \cC_X(\nu|_{V_2}) \ar[ur]_-{\varphi_{\nu|_{V_2}}}
}
$$
Here the left-hand side of the diagram is constructed by adjointness and hence commutes, and the right-hand side of the diagram commutes up to conjugation by the rank two assumption. This shows that the top left-to-right composite $\varphi_{\nu,V_1}$ is conjugate to the bottom left-to-right composite $\varphi_{\nu,V_2}$, i.e., the map $\varphi_{\nu,V}$ is independent of the choice of rank one subgroup $V$ in general. We hence drop the subscript $V$ and denote this map by $\varphi_\nu$.

With these preparations we can now easily finish the proof of part \eqref{del1} of the theorem. We need to see that for an arbitrary morphism  $\rho: (\nu: E \to X) \to (\nu': E' \to X)$ in $\A(X)$ the diagram
$$
\xymatrix{
\cC_X(\nu') \ar[rr]^-{\cC_X(\rho)} \ar[rd]_-{\varphi_{\nu'}} &
&\cC_X(\nu) \ar[ld]^-{\varphi_{\nu}} \\
& X'
}
$$
commutes.
Suppose first that $E'$ has rank one, and let $\mu: E' \to T \to \N$
be the factorization of $\nu'$ through $T$. The statement follows
here since the diagram
$$
\xymatrix{
 \cC_X(\nu') \ar[r]_-{\cong}^-{h_{\nu'}} \ar[d]^-{\cong}_-{\cC_X(\rho)} & \cC_{X'}(j'\mu) \ar[d]^-{\cong}_-{\cC_{X'}(\rho)} \ar[r] & X'\\
         \cC_X(\nu'\rho) \ar[r]^-{h_{\nu'\rho}}_-\cong  &
         \cC_{X'}(j'\mu\rho) \ar[ur]
}
$$
commutes up to conjugation, by the uniqueness in the rank one case, since we can view the diagram of
isomorphisms as taking place under $\cC_\N(\mu) \to
\cC_\N(\mu\rho)$. The general case follows from the rank one case, by
the independence of choice of rank one subgroup: If $V$ is a rank one
subgroup of $E$ set $V' = \rho(V)$ and observe that by adjointness
the diagram
$$
\xymatrix{
\cC_X(\nu') \ar[r] \ar[d]^-{\cC_X(\rho)} & \cC_X({\nu'}|_{V'}) \ar[d]^-{\cC_X(\rho)} \\
\cC_X(\nu) \ar[r] & \cC_X(\nu|_{V})
}
$$
commutes. Hence we have constructed a map up to conjugacy from the
centralizer diagram of $X$ to $X'$ (seen as a constant diagram), or in
other words we have given an element
$$
[\vartheta] \in \lim_{\nu \in
  \A(X)}\hspace{-10pt}{}^0\pi_0(\map(B\cC_X(\nu),BX'))
$$
This concludes the proof of part \eqref{del1}.

Using \cite[Rem.~after Def.~6.3]{DW96} \cite[Lem.~11.15]{dw:center}
(which says that the centralizer diagram of a $p$-compact group is
`centric') it is easy to see that the map $\varphi_\nu: \cC_X(\nu) \to X'$ induces a homotopy equivalence $$\map(B\cC_X(\nu),B\cC_X(\nu))_1 \xrightarrow{\simeq} \map(B\cC_X(\nu),BX')_{\varphi_\nu}$$
where the first term equals the classifying space of the center $B\cZ\cC_X(\nu)$ by definition \cite{dw:center}.
Since this is natural it gives a canonical identification of
the functor $\nu \mapsto \pi_i(\map(B\cC_X(\nu),BX')_{[\vartheta]})$ with
$\nu \mapsto \pi_i(B\cZ\cC_X(\nu))$.

By obstruction theory (see \cite[Prop.~3]{wojtkowiak87}
\cite[Prop.~1.4]{JO96}) the existence obstructions for lifting this to
an element in $\pi_0(\holim_{\A(X)}\map(B\cC_X(\nu),X')) \cong
\pi_0(\map(BX,BX'))$ lie in
$$\lim_{\nu \in
  \A(X)}\hspace{-10pt}{}^{i+1} \pi_i(\map(B\cC_X(\nu),BX')_{[\vartheta]})
\cong  \lim_{\nu \in \A(X)}\hspace{-10pt}{}^{i+1}\pi_i(B\cZ\cC_X(\nu))
\mbox{, } i \geq 1$$
But by assumption all these groups are identically zero, so our element $[\vartheta]$
lifts to a map $Bf: BX \to BX'$.

We now want to see that the construction of $f$ forces it to be an
isomorphism. Let $\N_p$ denote a $p$-normalizer of $T$, i.e., the
union of components in $\N$ corresponding to a Sylow $p$-subgroup of
$W$. Since $\N_p$ has non-trivial center (by standard facts about $p$-groups), we can find
a central rank one elementary abelian $p$-subgroup $\mu: V \to T \to \N_p$, so we can view $\N_p$ as sitting inside $\cC_\N(\mu)$. Hence by construction the diagram
$$
\xymatrix{
 & \N_p \ar[dl]_-j \ar[dr]^-{j'} \\
X \ar[rr]^-{f} & & X'
}
$$
commutes up to conjugation, and in particular $fj: \N_p \to X'$ is a
monomorphism. This implies that $f$ is a monomorphism as well:
$H^*(B\N_p;\F_p)$ is finitely generated over $H^*(BX';\F_p)$ via
$H^*(Bf \circ Bj;\F_p)$ by \cite[Prop.~9.11]{DW94}. By an application of the
transfer \cite[Thm.~9.13]{DW94} the map $H^*(Bj;\F_p): H^*(BX;\F_p) \to
H^*(B\N_p;\F_p)$ is a monomorphism, and since $H^*(BX';\F_p)$ is noetherian by
\cite[Thm.~2.4]{DW94} we conclude that $H^*(BX;\F_p)$ is finitely generated
over $H^*(BX';\F_p)$ as well.
Hence $f: X\to X'$ is a monomorphism by another application of
\cite[Prop.~9.11]{DW94}. Since we can identify the maximal tori of
$X'$ and $X$,
the definition of the Weyl group defines a map between the Weyl groups
$W_X \to W_{X'}$, which has to be injective since the Weyl groups act
faithfully on $T$ (by \cite[Thm.~9.7]{DW94}).  But since we know that
$X$ and $X'$ have the same maximal torus normalizer, the above map of
Weyl groups is an isomorphism. By \cite[Thm.~4.7]{dw:center}
(or \cite[Prop.~3.7]{MN94} using \cite[Thm.~9.7]{DW94}) this means
that $f$ is indeed an isomorphism.

We now want to argue that this $f$ is a map under $\N$. By
Lemma~\ref{themap} we know that there exists $Bg\in \Aut(B\N)$, unique up
  to conjugation, such that 
$$
\xymatrix{
\N \ar[r]^-g \ar[d]_-{j} & \N \ar[d]^-{j'}\\
X \ar[r]^-f & X'
}
$$
commutes up to conjugation. By covering space theory and Sylow's
theorem we can restrict $g$ to a self-map $g'$ making the diagram 
$$
\xymatrix{
\N_p \ar[r]^-{g'} \ar[d]_-{j} & \N_p \ar[d]^-{j'}\\
            X \ar[r]^-f & X'
}
$$
commute. Furthermore any other map $\N_p \to \N_p$ fitting in this
diagram will be conjugate to $g'$ in $\N$, by the proof of Lemma~\ref{themap}.
However by construction $f$ is a map under $\N_p$ so $g'$ is conjugate
in $\N$ to the canonical inclusion. But by Propositions~\ref{discreten}
and \ref{autn} we see that if an outer automorphism of $\N$ restricts
to the identity on $\N_p$ then it is in fact the identity. This shows
that $g = 1$, i.e., that $f$ is a map under $\N$. This also shows that
$\AM: \pi_0(\Aut(BX)) \to \pi_0(\Aut(B\N))$ is surjective, since for any
automorphism $g: \N \to \N$, $jg$ is also a maximal torus normalizer
in $X$ by \cite[Thm.~1.2(3)]{moller99}.

Note that if the component of $\Aut(B\N)$ of the identity map,
$\Aut_1(B\N)$, is not contractible we can find a rank
one elementary abelian
$p$-subgroup $\nu:V \to T$ such that $\cC_\N(\nu) \xrightarrow{\cong} \N$
which by assumption means that $\AM: \Aut(BX) \xrightarrow{\cong}
\Aut(B\N)$. So we can assume that $\Aut_1(B\N)$ is contractible in
which case $\Aut_1(BX)$ is as well by \cite[Thms.~1.3~and~7.5]{dw:center}.

The only remaining claim in the theorem is that the map $\AM:
\pi_0(\Aut(BX)) \to \pi_0(\Aut(B\N))$ is injective under the additional assumption that $\lim^{i}_{\nu \in
  \A(X)}\pi_i(B\cZ\cC_X(\nu)) =0 $, $i \geq 1$. In other words we have to
see that any self-equivalence $f$ of $X$ which, up to conjugacy,
induces the identity on $\N$ is in fact conjugate to the identity. But
if we examine the above argument with $X'=X$, the map on centralizers
of rank one elements induced by $f$ has to be the identity by the rank
one uniqueness assumption. The maps for higher rank are centralizers
of maps of rank one, so they as well have to be the identity. Hence
$f$ maps to the same element as the identity in $\lim_{\nu \in
  \A(X)}^0\pi_0(\map(B\cC_X(\nu),BX))$, which means that $f$ actually
is the identity by the vanishing of the obstruction groups (again, e.g., by
\cite[Prop.~4]{wojtkowiak87} or \cite[Prop.~1.4]{JO96}).
\end{proof}

\begin{rem}
Note how the assumptions of the theorem fail (as they should) for the group
$\SO(3)$ at the prime $2$ which is not determined by its maximal torus
normalizer. In this case the element $\diag(-1,-1,1)$ in the maximal
torus $\SO(2) \times 1$ is fixed under the Weyl group action and
has centralizer equal to the maximal torus normalizer $\oO(2)$.
\end{rem}

Define the {\em cohomological dimension} $\cd(W,L)$ of a finite
$\Z_p$-reflection group $(W,L)$ to be $2\cdot (\mbox{the number of
  reflections in }W) + \rk L$, and note that it follows easily from
\cite[Lem.~3.8]{dw:split} and \cite[Thm.~7.2.1]{benson93} that for $X$ a
connected $p$-compact group $\cd(X) = \cd(W_X,L_X)$.
We are now ready to give the proof of the main Theorems~\ref{classthm} and \ref{ndetthm}, referring forward to the rest of the  paper---the statements we refer to can however easily be taken at face value and returned to later.

\begin{proof}[Proof of Theorems~\ref{classthm}, and \ref{ndetthm} using Section~3--9 and 11--12]
We simultaneously show that Theorems~\ref{classthm} and \ref{ndetthm} hold by an induction on the cohomological dimension of
$X$ and $(W,L)$. We will furthermore add to the induction hypothesis the statement
that if $X$ is connected and $\Z_p[L_X]^{W_X}$ is a polynomial ring,
then $H^*(BX;\Z_p) \cong H^*(BT;\Z_p)^{W_X}$.

By the Component Reduction Lemma~\ref{ndetcomponent},
Theorem~\ref{ndetthm} holds for a $p$-compact group $X$ if it holds
for its identity component $X_1$, so we can assume that $X$ is
connected.

By a result of the first-named author \cite[Thm.~1.2]{kksa:thesis}, if $(W,L)$ is realized as the Weyl
group of a $p$-compact group $X$, then $\N_X$ will be split, i.e., the unique possible $k$-invariant of $B\N_X$ is zero and $B\N_X \simeq (BT)_{hW}$. (See also \cite{tits66}, \cite{DW01}, and
\cite{neumann99} for the Lie case.) Furthermore, by the Component Group Formula (Lemma~\ref{pprimelemma}) we can read off the fundamental group of $X$ from $\N_X$. So, to prove Theorems~\ref{classthm}
and \ref{ndetthm} we have to show that given any finite
$\Z_p$-reflection group $(W,L)$ there exists a unique connected
$p$-compact group $X$ realizing $(W,L)$, with self-maps satisfying
$\AM: \Aut(BX) \xrightarrow{\cong} \Aut(B\N_X)$, since this implies
$\AM: \pi_0(\Aut(BX)) \xrightarrow{\cong} N_{\GL(L)}(W)/W$ by
Propositions~\ref{discreten} and \ref{autn}.

We first deal with the existence part. By the classification of finite
$\Z_p$-reflection groups (Theorem~\ref{lattice-split}), $(W,L)$ can be
written as a product of exotic finite $\Z_p$-reflection groups
and a finite $\Z_p$-reflection group of the form $(W_G,L_G\otimes
\Z_p)$ for some compact connected Lie group $G$. The factor $(W_G,L_G
\otimes \Z_p)$ can of course be realized by $G\pcom$, and in this case
$H^*(BG\pcom;\Z_p) \cong H^*(BT;\Z_p)^{W_G}$ if and only if
$\Z_p[L_G\otimes \Z_p]^{W_G}$ is a
polynomial algebra by the invariant theory appendix (Theorems~\ref{polyinvthm} and
\ref{classicalborel}). If $(W,L)$ is an exotic finite $\Z_p$-reflection
group then $\Z_p[L]^W$ is a polynomial algebra by
Theorem~\ref{polyinvthm} and $(W,L)$ satisfies $\dT^W = 0$ by
the classification of finite $\Z_p$-reflection groups
Theorem~\ref{lattice-split}, where $\dT \cong L \otimes \Z/p^\infty$ is a discrete approximation to $T$. By our integral version of a theorem of
Nakajima (Theorem~\ref{nakajima}), the subgroup $W_V$ of $W$ fixing a
non-trivial elementary abelian $p$-subgroup $V$ in $\dT$ is again a
$\Z_p$-reflection group, and since reflections in $W_V$ are also
reflections in $W$ (and $W_V$ is a proper subgroup of $W$), we see
that $(W_V,L)$ has smaller cohomological dimension than $(W,L)$. 
Hence by the induction hypothesis, the assumptions of Inductive Polynomial Realization Theorem~\ref{polynomialrealizationthm} are
satisfied. So, by this theorem there exists
a (unique) connected $p$-compact group $X$ with Weyl group $(W,L)$ and this
satisfies $H^*(BX;\Z_p) \cong H^*(BT;\Z_p)^{W_X}$.

We now want to show that $X$ is uniquely determined by $(W,L) = (W_X,L_X)$ and
that $X$ satisfies $\AM: \Aut(BX) \xrightarrow{\cong} \Aut(B\N)$, i.e., that
$X$ satisfies the conclusion of Theorem~\ref{ndetthm} (and hence that of
Theorem~\ref{classthm}).
By the Center Reduction Lemma~\ref{ndetcenter} we can assume that $X$ is
center-free. Likewise by the splitting theorem
\cite[Thms.~1.4~and~1.5]{dw:split} together with the Product Automorphism Lemma~\ref{ndetprod} we can assume that $X$ is simple. By the classification of finite $\Z_p$-reflection groups (Theorem~\ref{lattice-split}) and the invariant theory appendix 
(Theorem~\ref{polyinvthm}) either $(W,L)$ has the property that $\Z_p[L]^W$ is a
polynomial algebra, or $(W,L)$ is one of the reflection groups
$(W_{\PU(n)},L_{\PU(n)} \otimes \Z_p)$ (with $p\!\mid\!n$),
$(W_{E_8},L_{E_8} \otimes \Z_5)$, $(W_{F_4},L_{F_4} \otimes \Z_3)$,
$(W_{E_6},L_{E_6} \otimes \Z_3)$, $(W_{E_7},L_{E_7} \otimes \Z_3)$, or
$(W_{E_8},L_{E_8} \otimes \Z_3)$.

We will go through these cases individually. We can assume that $X$ is either constructed via the Inductive Polynomial Realization Theorem, or $X=G\pcom$ for the relevant compact connected Lie group $G$. Let $X'$ be a connected $p$-compact group with Weyl group $(W,L)$. We want to see that the assumptions of Theorem~\ref{collectlemma} are satisfied.
For this we use the calculation of the elementary abelian $p$-subgroups in Section~\ref{liesection} sometimes together with a specialized lemma from Section~\ref{section:ndetlemmas} to see that the assumptions of Theorem~\ref{collectlemma}\eqref{del1} are satisfied. The assumptions of Theorem~\ref{collectlemma}\eqref{del2} follow from the Obstruction Vanishing Theorem~\ref{obstructionsmainthm}.

If $\Z_p[L]^W$ is a polynomial algebra, then by the Inductive Polynomial Realization Theorem, $X$ satisfies $H^*(BX;\Z_p) \cong
H^*(B^2L;\Z_p)^W$. Hence all elementary abelian $p$-subgroups of $X$
are toral by an application of Lannes $T$-functor (cf.\ Lemma~\ref{torsionfreetoral}).
In particular $X$ has no rank two non-toral elementary abelian
$p$-subgroups so the assumptions of
Theorem~\ref{collectlemma}\eqref{del1} are satisfied. By the Obstruction Vanishing Theorem~\ref{obstructionsmainthm}
the assumptions of Theorem~\ref{collectlemma}\eqref{del2} also hold, and
hence Theorem~\ref{collectlemma} implies that there exists an isomorphism of
$p$-compact groups $X \to X'$, and that $X$ satisfies the conclusion
of Theorem~\ref{ndetthm}.

Now consider $(W,L) = (W_{\PU(n)},L_{\PU(n)} \otimes \Z_p)$ where
$p\!\mid\!n$. Theorem~\ref{punelementary} says that $\PU(n)$ has exactly one conjugacy class of 
rank two non-toral elementary abelian $p$-subgroups $E$ and gives its Weyl group and centralizer. We divide into two cases. If $n\neq p$, Lemma~\ref{equivariancelemma} implies that the
assumptions of Theorem~\ref{collectlemma}\eqref{del1} are satisfied. If $n=p$, Lemma~\ref{selfcentralizing2} implies that again the assumptions of Theorem~\ref{collectlemma}\eqref{del1} are satisfied.
 In both cases the assumptions of Theorem~\ref{collectlemma}\eqref{del2} are satisfied by the Obstruction Vanishing Theorem~\ref{obstructionsmainthm}, so Theorem~\ref{ndetthm} holds for $X$.

If $(W,L) = (W_G,L_G\otimes \Z_p)$ for $(G,p) = (E_8,5)$, $(F_4,3)$,
$(2E_7,3)$, or $(E_8,3)$ then $G$ (and hence $X$) does not have any rank two non-toral
elementary abelian $p$-subgroups by Theorem~\ref{alggrpfacts}\eqref{simply-con}, so 
the assumptions of Theorem~\ref{collectlemma}\eqref{del1} are vacuously satisfied and the assumptions of Theorem~\ref{collectlemma}\eqref{del2}
hold by the Obstruction Vanishing Theorem~\ref{obstructionsmainthm}, so Theorem~\ref{ndetthm} holds also in this case.

Finally, if $(W,L) = (W_G,L_G\otimes \Z_p)$ for $(G,p) = (E_6,3)$ there
are by Theorem~\ref{nontorals:E6-3} two isomorphism classes of rank two non-toral elementary abelian $3$-subgroups $E^{2a}_{E_6}$ and $E^{2b}_{E_6}$ in $\A(X)$, $X=G\pcom$. These both satisfy the assumptions of Theorem~\ref{collectlemma}\eqref{del1} by Lemma~\ref{equivariancelemma}, using the information about the centralizers in Theorem~\ref{nontorals:E6-3}.
Since the assumptions of Theorem~\ref{collectlemma}\eqref{del2} as usual 
are satisfied by the Obstruction Vanishing Theorem~\ref{obstructionsmainthm} 
we conclude by Theorem~\ref{collectlemma} that Theorem~\ref{ndetthm} holds for $X$
as well. This concludes the proof of the main theorem.
\end{proof}

\begin{rem}
Note that taking the case $(W_{E_6},L_{E_6} \otimes \Z_3)$ last in the
above theorem is a bit misleading, since groups with adjoint form
$E_6$ appear as centralizers in $E_7$ and $E_8$, so a separate
inductive proof of uniqueness in those cases would require knowing
uniqueness of $E_6$.
\end{rem}

\begin{rem}
The very careful reader might have noticed that the proof of the splitting result
in \cite{kksa:thesis}, which we use in the above proof,
refers to a uniqueness result in \cite{BV98} in the case of $(W_{\PU(3)},L_{\PU(3)}\otimes \Z_3)$.
We now quickly sketch a more direct way to get the splitting in this case, which we were told by Dwyer-Wilkerson: We need to see that a
$3$-compact group with Weyl group $(W_{\PU(3)},L_{\PU(3)} \otimes \Z_3)$
has to have split maximal torus normalizer $\N$. So, suppose that $X$
is a hypothetical $3$-compact group as above but with non-split
maximal torus normalizer. By a transfer argument (cf.\
\cite[Thm.~9.13]{DW94}) $\N_3$ has to be non-split as well. Since every
elementary abelian $3$-subgroup in $X$ can be conjugated into $\N_3$
(since $\chi(X/\N_p)$ is prime to $p$), this means that all elementary
abelian $3$-subgroups in $X$ are toral. Furthermore by
\cite[Prop.~3.4]{dw:split} conjugation between toral elementary
abelian $p$-subgroups is controlled by the Weyl group, so the Quillen
category of $X$ in fact agrees with the Quillen category of $\N$. The
category has up to isomorphism one object of rank two and two objects of rank
one. The centralizers $\cC_\N(V)$ of these are respectively $T$, $T:
\Z/2$, and $T \cdot \Z/3$. The unique $3$-compact groups corresponding
to these centralizers are in fact given by
$B\cC_\N(V)\threecom$. Hence the map $B\N \to BX$ is an equivalence by
the centralizer cohomology decomposition theorem \cite[Thm.~8.1]{dw:center}. But since $\N$ is non-split, we can find a map $\Z/3^2
\to \N$, which is not conjugate in $\N$ to a map into $T$. Hence the
corresponding map $\Z/3^2 \to X$ cannot be conjugated into $T$ either,
contradicting \cite[Prop.~5.6]{DW94}.
\end{rem}

%%%%%%%%%%%%%%%%%%%%%%%%%%%%%%%%%%%%%%%%%%%%%%%%%%%%%%%%%%%%%%%%%%%%%%%%%%%

\section{Two lemmas used in Section~\ref{section:ndet}} \label{section:ndetlemmas}
In this section we prove two lemmas which are used to verify the assumption in
Theorem~\ref{collectlemma}\eqref{del1} for a non-toral elementary abelian
$p$-subgroup of rank two---see the text preceding Theorem~\ref{collectlemma} for an explanation of these assumptions; we continue with the notation of Section~\ref{section:ndet}.
We first need a proposition which establishes a bound on the Weyl
group of a self-centralizing rank two non-toral elementary abelian
$p$-subgroup of a connected $p$-compact group. (The {\em Weyl group}
of an elementary abelian $p$-subgroup $\nu:E\rightarrow X$ of a
$p$-compact group $X$ is the subgroup of $\GL(E)$ consisting of
elements $\alpha$ such that $\nu\alpha$ is homotopic to $\nu$.)
Let $\dN_X$ and $\dT$ denote
discrete approximations to $\N_X$ and $T$ respectively, i.e., $\dT
\cong L \otimes \Z/p^\infty$ and $\dN_X$ is an extension of
$W_X$ by $\dT$ such that $B\dN_X \to B\N_X$ is an $\F_p$-equivalence---we refer to
 \cite[\S 3]{dw:center} for facts about discrete
approximations.

\begin{prop} \label{selfcentralizingranktwo}
Let $X$ be a connected $p$-compact group, and let $\nu: E \to X$ be a rank two
elementary abelian $p$-subgroup with $\cC_X(\nu) \cong E$. Then $\SL(E) \subseteq W(\nu)$, where
$W(\nu)$ denotes the Weyl group of $\nu$.
\end{prop}

\begin{proof}
Let $V$ be an arbitrary rank one subgroup of $E$ and consider the adjoint map
$\tilde \nu: E \to \cC_X(\nu|_{V})$.
Let $\dN_p$ denote a discrete approximation
to the $p$-normalizer $\N_p$ of a maximal torus in $\cC_X(\nu|_{V})$, which has positive rank
since $X$ is assumed connected. Since $\chi(\cC_X(\nu|_V)/\N_p)$ is not divisible by $p$ we can
factor $\tilde \nu$ through $\dN_p$ (see \cite[Prop.~2.14(1)]{dw:center}),
and by elementary theory of $p$-groups $N_{\dN_p}(E)$
contains a $p$-group strictly larger than $E$. By assumption
$\cC_X(\nu) \cong E$ so $\cC_{\cC_X(\nu|_{V})}(\tilde \nu) \cong E$, and hence
$C_{\dN}(E)=E$. Thus $N_{\dN}(E)/C_{\dN}(E)\subseteq W(\nu)\subseteq \GL(E)$
contains a subgroup of order $p$ stabilizing $V$. Since $V$ was arbitrary, this
shows that $W(\nu)$ contains all Sylow $p$-subgroups in $\SL(E)$, and hence
$\SL(E)$ itself, cf.\ \cite[Satz~II.6.7]{huppert67}.
\end{proof}

\begin{lemma} \label{selfcentralizing2}
Let $X$ and $X'$ be two connected $p$-compact groups with the same maximal
torus normalizer $\N$ embedded via $j$ and $j'$ respectively. Assume that for
all elementary abelian $p$-subgroups $\eta: E \to X$ of $X$ of rank one the
centralizer $\cC_X(\eta)$ is determined by $\N_{\cC_X(\eta)}$ and
$\AM: \Aut(B\cC_X(\eta)) \xrightarrow{\cong} \Aut(B\N_{\cC_X(\eta)})$.

If $\nu: E \to X$ is a rank two non-toral elementary abelian $p$-subgroup of
$X$ such that $\cC_X(\nu) \cong E$ then the map $\varphi_{\nu,V}:
\cC_X(\nu) \to X'$ is independent of the choice of the rank one subgroup $V$ of $E$ (i.e.,
the assumption of Theorem~\ref{collectlemma}\eqref{del1} is satisfied for
$\nu$).
\end{lemma}

\begin{proof}
Fix a rank one subgroup $V\subseteq E$ and let $\mu:V \to T \to \N$ be the factorization of the toral elementary abelian $p$-subgroup $\nu|_V : V\to X$ through $T$, unique as a map to $\N$. Then $\varphi_{\nu,V} : E \cong \cC_X(\nu) \to X'$ is an elementary abelian $p$-subgroup of $X'$ and since we have an isomorphism $h_{\nu|_V} : \cC_X(\nu|_V) \xrightarrow{\cong} \cC_{X'}(j'\mu)$ by assumption it follows by adjointness that $\cC_{X'}(\varphi_{\nu,V}) \cong E$. By Proposition~\ref{selfcentralizingranktwo} we get $\SL(E) \subseteq W_X(\nu)$ and $\SL(E) \subseteq W_{X'}(\varphi_{\nu,V})$.

Now let $\alpha\in \SL(E) \subseteq W_X(\nu)$. Then $\alpha(V) \xrightarrow{\alpha^{-1}} V \xrightarrow{\mu} \N$ is the factorization of $(\nu \circ \alpha^{-1})|_{\alpha(V)} \cong \nu|_{\alpha(V)}$ through $T$, unique as a map to $\N$. Now consider the diagram
$$
\xymatrix@C=40pt{
E \ar[r]^-{\cong} \ar[d]^-{\alpha^{-1}} & \cC_X(\nu) \ar[r] & \cC_X(\nu|_V) \ar[d] \ar[r]^-{h_{\nu|_V}}_-{\cong} & \cC_{X'}(j'\mu) \ar[d] \ar[r] & X' \ar@{=}[d] \\
E \ar[r]^-{\cong} & \cC_X(\nu) \ar[r] & \cC_X(\nu|_{\alpha(V)}) \ar[r]^-{h_{\nu|_{\alpha(V)}}}_-{\cong} & \cC_{X'}(j'\mu\circ(\alpha^{-1}|_{\alpha(V)})) \ar[r] & X'. \\
}
$$
The left-hand and right-hand squares obviously commute and the middle square commutes by our assumption on rank one subgroups. We thus conclude that $\varphi_{\nu,\alpha(V)} = \varphi_{\nu,V}\circ\alpha$ for all $\alpha\in \SL(E)$. Since $W_{X'}(\varphi_{\nu,V})$ contains $\SL(E)$ and $\SL(E)$ acts transitively on the rank one subgroups of $E$ it follows that $\varphi_{\nu,V'} = \varphi_{\nu,V}$ for any rank one subgroup $V'$ of $E$ as desired.
\end{proof}

\begin{lemma}\label{equivariancelemma}
Let $X$ and $X'$ be two connected $p$-compact groups with the same maximal
torus normalizer $\N$ embedded via $j$ and $j'$ respectively.
Assume the inductive hypothesis $(\star)$ of Theorem~\ref{collectlemma},
i.e., that for all elementary abelian $p$-subgroups $\eta: E \to X$ of
$X$ the centralizer $\cC_X(\eta)$ is determined by $\N_{\cC_X(\eta)}$ when $\eta$ has rank one and that $\AM: \Aut(B\cC_X(\eta)) \xrightarrow{\cong} \Aut(B\N_{\cC_X(\eta)})$ when $\eta$ has rank one or two.

If $\nu: E \to X$ is a rank two non-toral elementary abelian $p$-subgroup of
$X$ such that $\cC_X(\nu)_1$ is non-trivial then the map
$\varphi_{\nu,V}: \cC_X(\nu) \to X'$ is independent of the choice of the rank
one subgroup $V$ of $E$ (i.e., the assumption of
Theorem~\ref{collectlemma}\eqref{del1} is satisfied for $\nu$).
\end{lemma}

\begin{proof}
Choose a rank one elementary abelian $p$-subgroup $\xi:U = \Z/p \to
\cC_X(\nu)_1$ in the center of the $p$-normalizer of a maximal torus
in $\cC_X(\nu)$, which is always possible since the action of a finite
$p$-group on a non-trivial $p$-discrete torus has a non-trivial fixed point. Let $\xi \times \nu: U \times E \to X$ be the map defined by adjointness.
For any rank one subgroup $V$ of $E$, consider the map
$\xi \times \nu|_{V} : U \times V \to X$ obtained by restriction. By construction $\xi \times \nu|_{V}$ is the adjoint of the composite
$U \xrightarrow{\xi} \cC_X(\nu)_1 \xrightarrow{\res} \cC_X(\nu|_{V})_1$, so
$\xi \times \nu|_{V} : U \times V \to X$ factors through a maximal
torus in $X$, since every rank one elementary abelian $p$-subgroup in
a connected $p$-compact group factors through a maximal torus by \cite[Prop.~5.6]{DW94}.
We want to see that $\xi \times \nu$ is a monomorphism, using the
theory of kernels \cite[\S 7]{DW94}: If $\xi \times \nu$ was not a
monomorphism then it would have a rank one kernel $K$, which by the
choice of $\xi$ cannot be equal to $U$. But this would mean that, for
some rank one subgroup $V'$ of $E$, both $\nu$ and $\xi \times
\nu|_{V'}$ would be monomorphisms of rank two and factor through the
monomorphism $\overline{(\xi \times \nu)}: (U\times E)/K \to X$ of
rank two. But this is a contradiction since $\xi \times \nu|_{V'}$ is
toral and $\nu$ is not.

Now consider the following diagram
$$
\xymatrix{
 & \cC_X(\nu|_V) \ar[rd]^-{\varphi_{\nu|_V}} \\
U \times E \ar[r] \ar[ru] \ar[rd] & \cC_X(\xi\times \nu|_V) \ar[u] \ar[d] & X' \\
 & \cC_X(\xi) \ar[ru]_-{\varphi_\xi}
}
$$
Here the left-hand side of the diagram is constructed by taking adjoints of
$\xi \times \nu$ and hence it commutes. The right-hand side is also forced to
commute by our inductive assumption, as explained in the beginning of the
proof of Theorem~\ref{collectlemma}, since $\xi \times \nu|_{V}$ is toral of
rank two.
We can hence without ambiguity define $(\xi \times \nu)'$ as either the top
left-to-right composite (for some rank one subgroup $V \subseteq E$) or the bottom
left-to-right composite. We let $\nu'$ denote the restriction of
$(\xi \times \nu)'$ to $E$.

Finally consider the diagram
$$
\xymatrix{
\cC_X(\xi \times \nu) \ar[r]^-{\widetilde{h_{(\xi \times \nu)|_V}}}_-\cong \ar[d] & \cC_{X'}((\xi \times \nu)') \ar[d] \\
\cC_X(\nu) \ar[r]^-{\widetilde{h_{\nu|_V}}}_-\cong & \cC_{X'}(\nu')
}
$$
and note that as before the inductive assumption guarantees that
${\widetilde{h_{(\xi \times \nu)|_V}}} = {\widetilde{h_{(\xi \times \nu)|_U}}}$, since $\xi \times \nu|_{V}$ is toral, and in particular it is independent
of the choice of $V$. We claim that this forces the same to be true for the bottom map.
Note that by our induction hypothesis, automorphisms of $\cC_X(\nu)$ are determined by the restriction to a maximal torus normalizer, which again are determined by their restriction to a $p$-normalizer of a maximal torus, by the description of automorphisms of maximal torus normalizers in
Propositions~\ref{discreten} and \ref{autn}. But by
our choice of $\xi$, the centralizer $\cC_X(\xi \times \nu )$ contains a
$p$-normalizer of a maximal torus in $\cC_X(\nu)$. This shows that also the bottom map
 $\widetilde{h_{\nu|_V}}$ is
independent of $V$. Hence $\varphi_{\nu,V}: \cC_X(\nu) \xrightarrow{\widetilde{h_{\nu|_V}}} \cC_{X'}(\nu') \xrightarrow{\text{ev}} X'$ is also independent of $V$ as wanted.
\end{proof}

%%%%%%%%%%%%%%%%%%%%%%%%%%%%%%%%%%%%%%%%%%%%%%%%%%%%%%%%%%%%%%%%%%%%%%%%%%%

\section{The map $\AM: \Aut(BX) \to \Aut(B\N_X)$} \label{themapsection}
The purpose of this very short section is to construct the map
$\AM: \Aut(BX) \to \Aut(B\N_X)$ which we will prove is an
equivalence. We have been unable to find this description in the
literature.

For a fibration $f: \cE \to \B$ we let $\Aut(f)$ denote the space of
commutative diagrams
$$
\xymatrix{
\cE \ar[d]_-f \ar[r] & \cE \ar[d]^-f\\
\B \ar[r]  & \B
}
$$
such that the horizontal maps are homotopy equivalences. (This is a
subspace of $\Aut(\cE) \times \Aut(\B)$.)

\begin{lemma}[Adams-Mahmud lifting] \label{themap}
Let $X$ be a $p$-compact group with maximal torus normalizer
$\N_X$. Turn the inclusion of the maximal torus normalizer into a
fibration $i: B\N_X \to BX$. Then the restriction map $\Aut(i) \to
\Aut(BX)$ is an equivalence of grouplike topological monoids.

In particular any self-homotopy equivalence of $BX$ lifts to a
self-homotopy equivalence of $B\N_X$, which is unique in the strong
sense that the space of lifts is contractible. Choosing a homotopy
inverse to the homotopy equivalence $B\Aut(i) \to B\Aut(BX)$, we get a
canonical map
$$\AM: B\Aut(BX) \xrightarrow{\cong} B\Aut(i) \to B\Aut(B\N_X).$$
\end{lemma}

\begin{proof}
For any $\varphi \in \Aut(BX)$, there exists, e.g.\ by
\cite[Thm.~1.2(3)]{moller99}, a map $\psi \in \Aut(B\N_X)$ such that
$\varphi i$ is homotopic to $i\psi$. Since $i$ is assumed to be a
fibration, $\psi$ can furthermore be modified such that the equality
is strict. This shows that the evaluation map $\Aut(i) \to \Aut(BX)$
is surjective on components. This map of grouplike topological monoids
is furthermore easily seen to have the homotopy lifting property. To
see that it is a homotopy equivalence we hence just have to verify
that the fiber $\Aut_{BX}(B\N_X)$ over the identity map is
contractible. First observe that, by \cite[Prop.~8.11]{DW94} and the definitions, there is a unique map $BT \to B\N_X$ over $BX$, up to homotopy. This shows that the homotopy fixed point space $(X/\N_X)^{hT}$ is at least connected. (We refer to \cite[\S 10]{DW94} for basic facts and
definitions about homotopy actions.) Consider the following diagram in which the rows and columns are fibrations 
$$
\xymatrix{
\W_X \ar[r] \ar@{=}[d] & X/T \ar[r] \ar[d] & X/\N_X \ar[d] \\
\W_X \ar[r] \ar[d] & BT \ar[r] \ar[d] & B\N_X \ar[d]\\
\ast \ar[r] & BX \ar@{=}[r] & BX.
}
$$
Take homotopy $T$-fixed points of the top row, which by definition equals taking spaces of liftings of the map from $BT$ to the bottom row to the corresponding term in the middle row.
This produces an induced fibration sequence
$\W_X \to (X/T)^{hT} \to (X/\N_X)^{hT}$ since $(X/\N_X)^{hT}$ is connected. However the map $\W_X \to (X/T)^{hT}$ is the identity, so $(X/\N_X)^{hT}$ is in fact contractible.
By \cite[Lem.~10.5 and Rem.~10.9]{DW94} we can rewrite
$(X/\N_X)^{h\N_X} \simeq ((X/\N_X)^{hT})^{hW_X}$, which shows that
$(X/\N_X)^{h\N_X}$ is contractible as well. Hence any self-map of
$B\N_X$ over $BX$ is an equivalence, and $\Aut_{BX}(B\N_X)$ is contractible as
wanted.
\end{proof}

%%%%%%%%%%%%%%%%%%%%%%%%%%%%%%%%%%%%%%%%%%%%%%%%%%%%%%%%%%%%%%%%%%%%%%%%%%%

\section{Automorphisms of maximal torus normalizers}\label{nautomorphismsection}
The aim of this short section is to establish some easy facts about
automorphisms of maximal torus normalizers which are needed to carry
out the reduction to connected, center-free simple $p$-compact groups
in Section~\ref{reductionsection}. At the same time the section serves
to make the automorphism statement of Theorem~\ref{classthm} more
explicit.

Recall that an {\em extended $p$-compact torus} is a loop space $\N$
such that $W = \pi_0(\N)$ is a finite group and the identity component
$\N_1$ of $\N$ is a $p$-compact torus $T$. Let $\dN$ be the discrete
approximation to $\N$ (see \cite[3.12]{dw:center}), and recall that
$\dN$ will have a unique largest $p$-divisible subgroup $\dT$, which
will be a discrete approximation to $T$.

\begin{prop}  \label{discreten}
For an extended $p$-compact torus $\N$, the obvious map associating to
a self-homotopy equivalence of $B\dN$ a self-homotopy equivalence of
$B\N$ via fiber-wise $\F_p$-completion \cite[Ch.~I, \S 8]{bk}
induces an
equivalence of aspherical grouplike topological monoids
$$\Aut(B\dN)\pcom \xrightarrow{\cong} \Aut(B\N).$$

If $\pi_0(\N)$ acts faithfully on $\pi_1(\N_1)$ then $\Aut_1(B\dN)$,
the component of $\Aut(B\dN)$ of the identity map, has the homotopy
type of $B(\dT^W)$ where $\dT$ is a discrete approximation to $T$.
\end{prop}

\begin{proof}[Sketch of proof]
The statement on the level of component groups follows directly from
\cite[Prop.~3.1]{dw:center}. (The point is that the homotopy fiber of $B\dN
\to B\N$ will have homotopy type $K(V,1)$ for a $\Q_p$-vector space
$V$, and hence the existence and uniqueness obstructions to lifting a
map $B\dN \to B\dN$ to $B\N$ lie in $H^n(\dN;V)$ where $n = 2, 1$
which are easily seen to be zero.) It is likewise easy to see that
both spaces are aspherical and that we get a homotopy equivalence of
the identity components. The last statement is also obvious.
\end{proof}

Let $L$ be a finitely generated free $\Z_p$-module and suppose that $W
\subseteq \Aut(\dT)$, where we set $\dT = L \otimes
\Z/p^\infty$. Consider the second cohomology group $H^2(W;\dT)$ which
classifies extensions of $W$ by $\dT$ with the fixed action of $W$ on
$\dT$. Given an isomorphism $\alpha: L \to L'$ sending $W \subseteq
\GL(L)$ to $W' \subseteq \GL(L')$ we get an isomorphism of cohomology
groups $H^2(W;\dT) \to H^2(W';\dT')$ by sending an extension $\dT
\xrightarrow{i} \dN \xrightarrow{\pi} W$ to the extension $\dT' \xrightarrow{i \circ
  \alpha^{-1}} \dN \xrightarrow{c_{\alpha}\circ \pi} W'$, where $c_\alpha$
denotes conjugation by $\alpha$. An isomorphism between two triples
$(W,L,\gamma)$ and $(W',L',\gamma')$, where $\gamma$ and $\gamma'$ are
extension classes, is an isomorphism $L \to L'$ sending $W$ to $W'$
and $\gamma$ to $\gamma'$. The automorphism group of a triple
$(W,L,\gamma)$ thus identifies with
$$
{}^\gamma N_{\Aut(\dT)}(W) = \{ \alpha \in N_{\Aut(\dT)}(W) \,  | \, \alpha(\gamma) = \gamma \in H^2(W;\dT)\}.
$$
It follows directly from the definition (and using that $\dT$ is
characteristic in $\dN$) that two triples as above are isomorphic if
and only if the associated groups $\dN$ and $\dN'$ are isomorphic,
where $\dN$ is obtained from the extension $1 \to \dT \to \dN \to W
\to 1$ given by $\gamma$, and analogously for $\gamma'$. However,
$\dN$ and $(W,L,\gamma)$ in general have slightly different
automorphism groups, as described in the following lemma (see also
\cite{wells71}):

\begin{prop} \label{autn}
In the notation above, for any exact sequence $1 \to \dT \to \dN
\xrightarrow{\pi} W \to 1$ with extension class $\gamma$ we have a canonical
exact sequence
\begin{equation} \label{normseq}
1 \to \Der(W,\dT) \to \Aut(\dN) \to {}^\gamma N_{\Aut(\dT)}(W) \to 1
\end{equation}
where we embed the derivations $\Der(W,\dT)$ in $\Aut(\dN)$ by sending
a derivation $s$ to the automorphism given by $x \mapsto s(\pi(x))x$,
and the map $\Aut(\dN) \to {}^\gamma N_{\Aut(\dT)}(W)$ is given by
restricting an automorphism $\varphi \in \Aut(\dN)$ to $\dT$.

This exact sequence has an exact subsequence $1 \to \dT/\dT^W
\to\dN/Z\dN \to W \to 1$ and the quotient exact sequence is
$$1 \to H^1(W;\dT) \to \Out(\dN) \to {}^\gamma N_{\Aut(\dT)}(W)/W \to 1.$$ 

In particular if $(W,L)$ is a finite $\Z_p$-reflection group and $p$
is odd then by \cite[Thm.~3.3]{kksa:thesis}\cite[Pf.~of Prop.~3.5]{JMO92}
$H^1(W;\dT) = 0$, so we get an isomorphism $\Out(\dN) \xrightarrow{\cong}
{}^\gamma N_{\Aut(\dT)}(W)/W$.
\end{prop}

\begin{proof} 
Let $\varphi \in \Aut(\dN)$, and consider the restriction map $\varphi
\mapsto  \varphi|_{\dT} \in \Aut(\dT)$.
Note that for all $x \in \dN$, $l \in \dT$ we have $$(\varphi \circ
c_x)(l) = \varphi(xlx^{-1}) = \varphi(x) \varphi(l) \varphi(x)^{-1} =
(c_{\varphi(x)} \circ \varphi)(l)$$ so $\varphi|_{\dT} \in
N_{\Aut(\dT)}(W)$. That the image is the elements which fix the
extension class follows easily from the definitions: The diagram
$$
\xymatrix{
\dT \ar[r]^-{i\circ \varphi^{-1}} \ar@{=}[d] & \dN \ar[r]^-{c_\varphi \circ \pi} \ar[d]^-\varphi & W \ar@{=}[d] \\
\dT \ar[r]^-i & \dN \ar[r]^-\pi & W
}
$$
shows that $\varphi$ leaves $\gamma$ invariant.
Likewise, to see that the right map in \eqref{normseq} is surjective
let $\psi \in {}^\gamma N_{\Aut(\dT)}(W)$ and let $\dT \to \tilde N \to
W$ be the extension obtained by first pushing forward along $\psi: \dT
\to \dT$ and then pulling back along $\psi^{-1}( - )\psi: W \to
W$. Since $\psi$ fixes $\gamma$ there exists an isomorphism $\tilde N
\to \dN$ making the following diagram commute
$$
\xymatrix@C=35pt{
\dT \ar[r]^-{\psi} \ar[d] & \dT  \ar[d] \ar@{=}[r] & \dT \ar[d]\\
 \dN \ar[r] \ar[d] & \tilde N \ar[r] \ar[d] & \dN \ar[d] \\
 W \ar[r]^-{\psi(-)\psi^{-1}} & W \ar@{=}[r] & W
}
$$
which shows that $\Aut(\dN) \to {}^\gamma N_{\Aut(\dT)}(W)$ is
surjective.

Now suppose $\varphi \in \Aut(\dN)$ restricts to the identity on
$\dT$. For $x \in \dN$ and $l\in \dT$ we have
$$\varphi(x) l \varphi(x^{-1})
= \varphi(x) \varphi(l) \varphi(x^{-1}) = \varphi(xlx^{-1}) =
xlx^{-1}$$
so the induced map $\varphi: W \to W$ is the identity since $W$ acts
faithfully on $\dT$. This means that we can define a map $s: W \to
\dT$ by $s(w) = \varphi(\tilde w) \tilde w^{-1}$ where $\tilde w$ is a
lift of $w$, and this is easily seen to be a derivation. Furthermore
taking the automorphism of $\dN$ associated to $s$ gives back
$\varphi$, which establishes exactness in the middle, and we have
proved the existence of the first exact sequence.

The existence of the short exact subsequence is clear, noting that
$Z\dN = \dT^W$ (since $W$ acts faithfully on $\dT$) and that
$\dT/\dT^W$ embeds in $\Der(W,\dT)$ as the principal derivations by
sending $l$ to the derivation $w \mapsto l (w\cdot l)^{-1}$. The last exact
sequence is now obvious.
\end{proof}

\begin{rem}
See \cite[Thm.~1.2]{hammerli02} for a related exact sequence for compact
connected Lie groups, fitting with the conjectured classification of
connected $p$-compact groups for $p=2$.
\end{rem}

\begin{prop} \label{algprod}
Suppose $\{(W_i,L_i,\gamma_i)\}_{i=0}^k$ is a collection of pair-wise
non-isomorphic triples where $L_i$ is a finitely generated free
$\Z_p$-module,  $W_i$ is a finite subgroup of $\GL(L_i)$ such that
$L_i \otimes \Q$ is an irreducible $W_i$-module, and $\gamma_i \in
H^2(W_i;\dT_i)$. Assume that $W_0 = 1$ but that $W_i$ is non-trivial
for $i \geq 1$. Let $(W,L,\gamma) = \prod_{i=0}^k
(W_i,L_i,\gamma_i)^{m_i}$ denote the product. Then
$$
\GL_{m_0}(\Z_p) \times \left(\prod_{i=1}^k \left({}^{\gamma_i}
N_{\GL(L_i)}(W_i)/W_i\right) \wr \Sigma_{m_i}\right) \xrightarrow{\cong} {}^\gamma
N_{\GL(L)}(W)/W.
$$
\end{prop}

\begin{proof}
The map of the proposition is injective by definition, and we have to
see that it is surjective.
To lessen confusion write
$$
L = (L_{0,1} \oplus \cdots \oplus L_{0,m_0}) \oplus (L_{1,1} \oplus
\cdots \oplus L_{1,m_1}) \oplus \cdots \oplus (L_{k,1} \oplus \cdots
\oplus L_{k,m_k}),
$$
which we consider as a $W = 1 \times (W_{1,1} \times \cdots \times
W_{1,m_1)} \times \cdots \times (W_{k,1} \times \cdots
W_{k,m_k})$-module, where $(W_{i,j},L_{i,j})$ is isomorphic to
$(W_i,L_i)$ as $\Z_p$-reflection groups.

Consider $\varphi \in {}^\gamma N_{\GL(L)}(W)$; we need to see that
this has the prescribed form. First note that for every $w \in W$
there exists a unique $\tilde w \in W$ such that
$$
\varphi(wx) = \tilde w \varphi(x) \mbox{ for all } x \in L.
$$
Let $\alpha$ denote the corresponding element in $\Aut(W)$ given by
$w\mapsto \tilde w$.

Note that the  above splitting of $L$ induces a splitting of
${}^\alpha L$, where the superscript means that we consider $L$ as a
$W$-module through $\alpha$.

Let $M$ and $N$ be indecomposable summands of $L$.
By definition of $\alpha$ the canonical map
$$\varphi_{MN} : M \to L \xrightarrow{\varphi} {}^\alpha L \to {}^\alpha N$$
is $W$-equivariant. Therefore this map, after tensoring
with $\Q$, has to be either an isomorphism or zero.  Since all the
non-trivial summands of $L \otimes \Q$ and ${}^\alpha L \otimes \Q$
occur with multiplicity one there is for each non-trivial $M$ at most
one $N$ for which the map can be non-zero, and this $N$ is necessarily
non-trivial. Since $\varphi$ is an isomorphism there is hence exactly
one such $N$, and the map $\varphi_{MN}$ has to be an
isomorphism. Note furthermore that since $\varphi_{MN}$ gives an
isomorphism between $M = L_{i,j}$ and ${}^\alpha N = {}^\alpha
L_{k,l}$ as $W_{i,j}$-modules then $(\alpha, \varphi)$ induces an
isomorphism between the reflection groups $(W_{i,j},L_{i,j})$ and $(W
\cap \GL(L_{k,l}),L_{k,l})$, which by assumption has to send
$\gamma_i$ to $\gamma_k$, so $i=k$. This shows that $\varphi$ is of
the required form.
\end{proof}

%%%%%%%%%%%%%%%%%%%%%%%%%%%%%%%%%%%%%%%%%%%%%%%%%%%%%%%%%%%%%%%%%%%%%%%%%%%

\section{Reduction to connected, center-free simple $p$-compact groups}\label{reductionsection}

In this section we prove some lemmas, which, together with the
splitting theorems of Dwyer-Wilkerson \cite{dw:split} and Notbohm
\cite{notbohm00}, reduce the proof of Theorem~\ref{ndetthm} to the
case of connected, center-free simple $p$-compact groups. This reduction is
known and most of it appears in \cite{moller98} (relying on earlier
work of that author). We here provide a self-contained and a bit
more direct proof using \cite{dw:center}.

\begin{lemma}[Product Automorphism Lemma] \label{ndetprod}
Let $X$ and $X'$ be $p$-compact groups with maximal torus normalizers
$\N$ and $\N'$. Then $\N \times \N'$ is a maximal torus normalizer for
$X \times X'$ and the following statements hold:
\begin{enumerate}

\item \label{prod1}
$\Aut_1(BX) \times \Aut_1(BX') \xrightarrow{\cong} \Aut_1(BX \times BX')$ and
$\Aut_1(B\N) \times \Aut_1(B\N') \xrightarrow{\cong} \Aut_1(B\N \times B\N')$,
where $\Aut_1$ denotes the set of homotopy equivalences homotopic to the
identity.

\item \label{prod2}
If $\AM: \Aut(BX) \to \Aut(B\N)$ and $\AM: \Aut(BX') \to \Aut(B\N')$ are
injective on $\pi_0$, then so is $\AM: \Aut(B(X \times X')) \to \Aut(B(\N
\times \N'))$.

\item \label{prod3}
Suppose that $p$ is odd and that $X$ is connected with $X = X_1 \times
\cdots \times X_k$ such that each $X_i$ is simple and determined by
its maximal torus normalizer. If, for each $i$, $\AM: \Aut(BX_i) \to
\Aut(B\N_{X_i})$ is surjective on $\pi_0$ then so is $\AM: \Aut(BX) \to
\Aut(B\N)$.
\end{enumerate}
\end{lemma}

\begin{proof} 
Recall that the map $\AM: \Aut(BX) \to \Aut(B\N)$ was described in
Section~\ref{themapsection}.
To see \eqref{prod1} first note that
\begin{equation}\label{diag1}
\map(BX \times BX',BX \times BX') \simeq \map(BX,\map(BX',BX)) \times
\map(BX',\map(BX,BX')).
\end{equation}
The evaluation map $\map(BX',BX)_{0} \to BX$ is an
equivalence by the Sullivan conjecture for $p$-compact groups
\cite[Thm.~9.3 and Prop.~10.1]{dw:center}, where the subscript $0$ denotes the component of the constant map.
Since the component of the identity map on the left-hand side of
\eqref{diag1} lands in the component of the constant map in
$\map(BX',BX)$ this shows that $\map(BX \times BX',BX \times BX')_1
\simeq \map(BX,BX)_1 \times \map(BX',BX')_1$ as wanted. (The statement
just says that the center of a product of $p$-compact groups is the
product of the centers, which of course also follows from the
equivalence of the different definitions of the center from
\cite{dw:center}.)

To see \eqref{prod2} suppose that $\varphi$ is a self-equivalence of
$BX \times BX'$ such that its restriction to a self-equivalence of
$B(\N \times \N')$ becomes homotopic to the identity. The restriction
$\varphi|_{BX \times \ast}$ composed with the projection onto $BX'$
becomes null homotopic upon restriction to $B\N$, which, e.g.\ by
\cite[Cor.~6.6]{moller96}, implies that it is null homotopic. Likewise
the projection of $\varphi|_{\ast \times BX'}$ onto $BX'$ becomes
homotopic to the identity map upon restriction to $B\N$, which by
assumption means that the projection of $\varphi_{\ast \times BX'}$
onto $BX'$ is the identity. But by adjointness, repeating the argument
of the first claim, this implies that $\varphi$ composed with the
projection onto $BX'$ is homotopic to the projection map onto $BX'$
(this is \cite[Lem.~5.3]{dw:center}). By symmetry this holds for the
projection onto $BX$ as well, and we conclude that $\varphi$ is
homotopic to the identity as wanted.

Finally, to see $\eqref{prod3}$, note that Propositions~\ref{discreten}, \ref{autn}, and \ref{algprod} give a complete description of
$\pi_0(\Aut(B\N)) \cong \Out(\dN)$. The assumption that each $X_i$ is
determined by its maximal torus normalizer means by definition that if
$\N_{X_i}$ is isomorphic to $\N_{X_j}$ then $X_i$ is isomorphic to
$X_j$. It is now clear from the description of $\Out(\dN)$ and the
assumptions on the $X_i$'s, that all elements in $\Out(\dN)$ can be
realized by self-equivalences of $BX$.
\end{proof}

\begin{rem} \label{ptwoproducts}
Part \eqref{prod3} of the above lemma is in general false for
$p=2$. For instance if $X = \SO(3)\twocom$ then it is easy to
calculate directly (or appeal to \cite[Cor.~3.5]{JMO95}) that for both $Y = X$ and
$Y = X \times X$ we have $\AM: \pi_0(\Aut(BY)) \xrightarrow{\cong}
N_{\GL(L_Y)}(W_Y)/W_Y$. But for $Y = X \times X$ we have
$H^1(W_{Y};\dT_{Y}) \cong \Z/2 \times \Z/2$, so $B\N_{Y}$ has
non-trivial automorphisms which restrict to the identity on $BT_Y$
(see Proposition~\ref{autn}).
\end{rem}

\begin{rem}
Note that the assumption that the factors $X_i$ in \eqref{prod3} are
determined by $\N_i$ of course appears for a good reason. The
statement that $\AM: \pi_0(\Aut(BX)) \to \pi_0(\Aut(B\N))$ is
surjective for all $p$-compact groups $X$ implies that all $p$-compact
groups are determined by their maximal torus normalizer, as is seen by
taking products. Hence the first part of Theorem~\ref{ndetthm} in fact
follows from the second part.
\end{rem}

Recall the observation that for $p$ odd the component group of $X$ is
determined by $W_X$:

\begin{lemma}[Component Group Formula] \label{pprimelemma}
Let $X$ be a $p$-compact group for $p$ odd, with maximal torus
normalizer $j:\N \to X$. By definition $W_X = \pi_0(\N)$. The map
$\pi_0(j): W_X = \pi_0(\N) \to \pi_0(X)$ is surjective.  The kernel is
$O^{p}(W_X)$, the subgroup generated by elements of order prime to
$p$. It can also be identified with the Weyl group of the identity
component $X_1$ of $X$, and is the largest $\Z_p$-reflection subgroup
of $W_X$.
\end{lemma}

\begin{proof}
By \cite[Rem.~2.11]{dw:center} $\pi_0(j)$ is surjective with kernel
the Weyl group of the identity component of $X$. Since $\pi_0(X)$ is
a $p$-group, $O^{p}(\pi_0(\N))$ is contained in the kernel. On the
other hand, since $p$ is odd, the Weyl group of $X_1$ is generated by
elements of order prime to $p$, since it is a $\Z_p$-reflection group,
so equality has to hold.
\end{proof}

\begin{rem} \label{ptwocomponents}
For $p=2$ the component group of $X$ cannot be read off from $\N_X$,
and one would have to remember $\pi_0(X)$ as part of the data. For
instance the $2$-compact groups $\SO(3)\twocom$ and $\oO(2)\twocom$
have the same maximal torus normalizers, namely $\oO(2)\twocom$. Note
however that if $X$ is the centralizer of a toral abelian subgroup $A$ of a
connected $p$-compact group $Y$, then the component group of $X$ {\em
  can} be read off from $A$ and $\N_Y$ (see
\cite[Thm.~7.6]{dw:center}); a case of frequent interest.
\end{rem}

Before proceeding recall that by \cite{DKS89} (see also
\cite[Prop.~11.9]{dw:center}) we have, for a fibration  $\cF \to \cE
\xrightarrow{f} \cB$, a fibration sequence
$$
\map(\cB,B\Aut(\cF))_{C(f)} \to B\Aut(f) \to B\Aut(\cB).
$$
Here $C(f)$ denotes the components corresponding to the orbit of the
$\pi_0(\Aut(\cB))$-action on the class in $[\cB,B\Aut(\cF)]$
classifying the fibration.

We are interested in when the map of grouplike topological monoids
$\Aut(f) \to \Aut(\cE)$ is a homotopy equivalence. This will follow if
we can see that $\Aut_1(f) \to  \Aut_1(\cE)$ and $\pi_0(\Aut(f)) \to
\pi_0(\Aut(\cE))$ are equivalences. By an easy general argument given
in \cite[Prop.~11.10]{dw:center} the statement about the identity
components follows if $\cB \to \map(\cF,\cB)_0$ is an
equivalence, where the subscript $0$ denotes the component of the constant map.

\begin{lemma}[Component Reduction Lemma] \label{ndetcomponent}
Let $X$ be a $p$-compact group with maximal torus normalizer $\N$, and
assume that $p$ is odd (so that $\pi_0(X)$ can be read off from
$\N$). Let $\N_1$ denote the kernel of the map $\N \to \pi_0(X)$,
which is a maximal torus normalizer for $X_1$.

If $\AM: \Aut(BX_1) \xrightarrow{\cong} \Aut(B\N_1)$, then $\AM: \Aut(BX) \xrightarrow{\cong}
\Aut(B\N)$. If furthermore $BX_1$ is determined by $B\N_1$ then $BX$
is determined by $B\N$.
\end{lemma}

\begin{proof}
First note that by an inspection of Euler characteristics and using
\cite[Thm.~1.2(3)]{moller99}, $\N_1$ is indeed a maximal torus
normalizer in $X_1$. Set $\pi = \pi_0(X)$ for short. We
want to apply the setup described before the lemma to the fibrations
$BX_1 \to BX \to B\pi$ and $B\N_1 \to B\N \to B\pi$ and to see that in
both cases the map of monoids $\Aut(f) \to \Aut(\cE)$ are homotopy
equivalences. By the remarks above this follows if it is an
isomorphism on $\pi_0$ and $\B \to \map(\cF,\B)_0$ is an
equivalence. The statement about $\pi_0$ is true in both cases since a
self-map of $\cE$ determines a unique self-map of $B\pi$. Likewise it
is easy to see that $B\pi \xrightarrow{\simeq} \map(BX_1,B\pi)_0$ and that
$B\pi \xrightarrow{\simeq} \map(B\N_1,B\pi)_0$.
This means that our map $B\Aut(BX) \to B\Aut(B\N)$ (from
Lemma~\ref{themap}) fits in a map of fibration sequences
$$
\xymatrix{
\map(B\pi,B\Aut(BX_1))_{C(f)} \ar[d] \ar[r] & B\Aut(BX) \ar[d] \ar[r] & B\Aut(B\pi) \ar@{=}[d] \\
\map(B\pi,B\Aut(B\N_1))_{C(f)} \ar[r] & B\Aut(B\N) \ar[r] &
B\Aut(B\pi).
}
$$
Here the maps between the fibers and base spaces are homotopy equivalences
by assumption, and we conclude that we get a homotopy equivalence
between the total spaces as well.

Now assume furthermore that $X_1$ is determined by $\N_1$, and let
$X'$ be another $p$-compact group with maximal torus normalizer
$\N$. By Lemma~\ref{pprimelemma},  $\pi = \pi_0(X) \cong
\pi_0(X')$ and $\N_1$ is also a maximal torus normalizer in
$X'_1$.

We want to show that the two fibrations $BX \to B\pi$ and $BX' \to
B\pi$ are equivalent as fibrations over $B\pi$, or equivalently that
the $\pi$-spaces $BX_1$ and $BX_1'$ are $h\pi$-equivalent, i.e., that
we can find a zig-zag of $\pi$-maps which are non-equivariant
equivalences connecting the two (see  e.g., \cite{DDK80equivariant}
where this equivalence relation is called equivariant weak homotopy
equivalence).

By the assumptions on $X_1$ we can choose a homotopy equivalence $Bf:
BX_1 \to BX'_1$ such that
$$
\xymatrix{
 & B\N_1  \ar[dr]^-{Bj'} \ar[dl]_-{Bj} &\\
BX_1 \ar[rr]^-{Bf} & & BX_1'
}
$$
commutes up to homotopy, and $Bf$ is unique up to homotopy.

We now want to see that we can change $Bf$ so that it becomes a
$\pi$-map. For this, consider the restriction $\pi$-map
$$\map(BX_1,BX_1') \to \map(B\N_1,BX_1').$$
By the assumptions on $\Aut(BX_1)$ this map sends distinct components of
$\map(BX_1,BX_1')$ corresponding to homotopy equivalences to distinct
components of $\map(B\N_1,BX_1')$. Moreover, by the proof of
Lemma~\ref{themap}, we have a homotopy equivalence
$\map(BX_1,BX_1')_{Bf} \simeq \map(B\N_1,BX_1')_{Bf \circ Bj}$. In
particular the component $\map(BX_1,BX_1')_{Bf}$ is preserved under
the $\pi$-action, since this obviously is so for
$\map(B\N_1,BX_1')_{Bj'}$. Furthermore since
$\map(B\N_1,BX_1')_{Bj'}^{\pi}$ contains $Bj'$ we see that
$\map(BX_1,BX_1')_{Bf}^{h\pi} \simeq \map(B\N_1,BX_1')_{Bj'}^{h\pi}$
is non-empty, and so there exists a $\pi$-map $E\pi \times BX_1 \to
BX_1'$ which is a homotopy equivalence. This shows that $BX_1$ and
$BX_1'$ are $h\pi$-homotopy equivalent as wanted.
\end{proof}

\begin{rem} \label{centralizerassumptionremark}
If $X$ is a connected $p$-compact group, and $p$ is odd, then it
follows from \cite[Thm.~7.5]{dw:center} that $Z(\dN)$ is a discrete
approximation to the center of $X$. The proof of the above lemma
extends this to $X$ non-connected {\em provided} we know that the
self-maps of $X_1$ are detected by their restriction to $\N_1$, which
will be a consequence of Theorem~\ref{ndetthm}. Having to appeal to
this is a bit unfortunate but seems unavoidable. The point is that if
there for a connected $p$-compact group $X$, existed a
self-equivalence $\sigma$ of $X$ of finite $p$-power order not
detected by $\N$, then we could form $X \semi \langle \sigma\rangle$,
where $\sigma$ would be central in the normalizer but not in the whole
group. (See also Lemma~\ref{centralizercenter}.)
\end{rem}

\begin{lemma}[Center Reduction Lemma] \label{ndetcenter}
Let $X$ be a connected $p$-compact group with center $\cZ$. 
\begin{enumerate}
\item \label{ndetcenterobj}
If $\AM: \pi_0(\Aut(BX/\cZ)) \to \pi_0(\Aut(B\N/\cZ))$ is surjective and
$X/\cZ$ is determined by $\N/\cZ$ then $X$ is determined by $\N$.
\item \label{ndetcentermor}
If $p$ is odd and $\AM: \Aut(BX/\cZ) \to \Aut(B\N/\cZ)$ is a homotopy
equivalence then $\AM: \Aut(BX) \to \Aut(B\N)$ is as well.
\end{enumerate}
\end{lemma}

\begin{proof}
To prove the first statement, suppose that $X$ and $X'$ have the same
maximal torus normalizer $\N$, choose fixed inclusions $j: \N \to X$
and $j': \N \to X'$, and let $\cZ$ be the center of $X$, which we can
view as a subgroup $\N$ via an inclusion $i: \cZ \to \N$. We claim
that $\cZ$ is also central in $X'$. It is central in the connected
component $X'_1$ by the formula for the center in
\cite[Thm.~7.5]{dw:center}. Furthermore, $\pi = \pi_0(X') \cong
W_{X}/W_{X'_1}$ acts trivially on $\cZ$ so the lift of $j'i$ to a map
$k: \cZ \to X_1'$, is unique up to conjugacy. Therefore we have
fibration sequences
$$\xymatrix{\map(B\cZ,BX'_1)_k \ar[r] \ar[d]^\simeq & \map(B\cZ,BX')_{j'i} \ar[r] \ar[d] & \map(B\cZ,B\pi)_0 \ar[d]^\simeq \\
BX_1' \ar[r] & BX' \ar[r] & B\pi}$$
where the left vertical map is an equivalence since $\cZ$ is central in $X'_1$. Hence the middle map in the above diagram is an equivalence as well, and $\cZ$ is central in $X'$ as claimed. Now assume that $X/\cZ$ is isomorphic to $X'/\cZ$. If
$\AM: \pi_0(\Aut(BX/\cZ)) \to \pi_0(\Aut(B\N/\cZ))$ is surjective we can
furthermore choose the homotopy equivalence $BX/\cZ
\to BX'/\cZ$ in such a way that
$$
\xymatrix{
 & B\N/\cZ \ar[dl]_-{j/\cZ} \ar[dr]^-{j'/\cZ} & \\
BX/\cZ \ar[rr] & & BX'/\cZ
}
$$
commutes up to homotopy.

We have canonical maps $BX/\cZ \to B^2\cZ$ and $BX'/\cZ \to B^2\cZ$
classifying the extensions, and we claim that in fact the bottom
triangle in the diagram
$$
\xymatrix{
 & B\N/\cZ \ar[dr] \ar[dl] &\\
BX/\cZ \ar[dr] \ar[rr] & & BX'/\cZ \ar[dl]\\
& B^2\cZ
}
$$ 
commutes up to homotopy.
By construction the outer square commutes up to homotopy (since both
composites agree with the classifying map $B\N/\cZ \to B^2\cZ$
since $j$ and $j'$ are fixed). Since the top triangle also
commutes up to homotopy, an application of the transfer
\cite[Thm.~9.13]{DW94}, using that $B^2\cZ$ is a product of
Eilenberg-Mac~Lane spaces and that $\chi((X/\cZ)/(\N/\cZ))=1$, shows
that the bottom triangle commutes up to homotopy as well. Since we
have constructed a map $BX/\cZ \to BX'/\cZ$ over $B^2\cZ$ we get an
induced homotopy equivalence $BX \to BX'$. (Note that this
construction does not a priori give this map as a map under $B\N$.)

We now want to get the second statement about automorphism
groups. Consider the homotopy commutative diagram
$$
\xymatrix{
B\N \ar[r]^-{f'} \ar[d]& B\N/\cZ \ar[d]\\
BX \ar[r]^-{f} &BX/\cZ
}
$$
where we can suppose that the two horizontal maps $f'$ and $f$ are fibrations.

We first claim that we can replace $B\Aut(f)$ with $B\Aut(BX)$ and
$B\Aut(f')$ with $B\Aut(B\N)$.
As in the case of the component group (see the proof of
Lemma~\ref{ndetcomponent}) we just have to justify that in the
appropriate fibration sequences we have equivalences $\cB \to
\map(\cF,\cB)_0$ and $\pi_0(\Aut(f)) \to \pi_0(\Aut(\cE))$.  The map
$BX/\cZ \to \map(B\cZ,BX/\cZ)_0$ is a homotopy equivalence since the
trivial map is central \cite[Prop.~10.1]{dw:center}. That $B\N/\cZ \to
\map(B\cZ,B\N/\cZ)_0$ is an equivalence follows by a similar (but easier)
argument.

By Lemma~\ref{themap} a self-equivalence of $BX$ induces a unique
self-equivalence of $B\N$, and hence a canonical self-equivalence of
$B\cZ$. Now, by the description of $X/\cZ$ as a Borel construction
(given in \cite[Pf.~of Prop.~8.3]{DW94}) we get a canonical
self-equivalence of $BX/\cZ$. This self-equivalence is furthermore
unique, in the sense that given a diagram
$$
\xymatrix{
BX \ar[r]^-{g} \ar[d] &  BX \ar[d] \\
BX/\cZ \ar[r]^-{g'} & BX/\cZ
}
$$ 
the homotopy type of $g'$ is uniquely given by that of $g$. To see
this note that by Lemma~\ref{themap} the diagram restricts to a unique
diagram
$$
\xymatrix{
B\N \ar[r]^-{\tilde g} \ar[d] & B\N \ar[d] \\
B\N/\cZ \ar[r]^-{\tilde g'} & B\N/\cZ.
}
$$
By looking at discrete approximations we see that the  homotopy class
of $\tilde g'$ is determined by $\tilde g$. Since by assumption the
homotopy class of $g'$ is determined by $\tilde g'$, we conclude that
a self-equivalence of $BX$ induces a unique self-equivalence of
$BX/\cZ$, and so $\pi_0(\Aut(f)) \cong \pi_0(\Aut(BX))$. The last part
of the argument furthermore shows that also $\pi_0(\Aut(f')) \cong
\pi_0(\Aut(B\N))$.

We hence have the following diagram where the rows are fibration sequences
$$
\xymatrix{
\map(BX/\cZ,B\Aut(B\cZ))_{C(f)} \ar[r] \ar[d] & 
B\Aut(BX) \ar[r] \ar[d] & B\Aut(BX/\cZ) \ar[d]\\
        \map(B\N/\cZ,B\Aut(B\cZ))_{C(f')} \ar[r] & B\Aut(B\N) \ar[r] & B\Aut(B\N/\cZ).
}
$$
Examining when the middle vertical arrow is a homotopy equivalence
reduces to finding out when the restriction map
$\map(BX/\cZ,B\Aut(B\cZ))_{C(f)} \to
\map(B\N/\cZ,B\Aut(B\cZ))_{C(f')}$ is a homotopy equivalence, which we
now analyze.

Note that since $B\cZ$ is a product of Eilenberg-Mac~Lane spaces
we have a fibration sequence
$$
B^2\cZ \to B\Aut(B\cZ) \to B\Aut(\dZ)
$$
where $\dZ$ is the discrete approximation to $\cZ$ and $\Aut(\dZ)$ is
the discrete group of automorphisms. Since our extensions are central
this gives a diagram of fibration sequences
$$
\xymatrix@C=18pt{
\map(BX/\cZ,B^2\cZ)_{C(f)} \ar[r] \ar[d] &  \map(BX/\cZ,B\Aut(B\cZ))_{C(f)} \ar[r] \ar[d] & \map(BX/\cZ,B\Aut(\dZ))_0 \ar[d] \\
\map(B\N/\cZ,B^2\cZ)_{C(f')} \ar[r] &  \map(B\N/\cZ,B\Aut(B\cZ))_{C(f')} \ar[r] & \map(B\N/\cZ,B\Aut(\dZ))_0.
}
$$
Again, in this diagram the map between the base spaces is obviously an
equivalence, so we are reduced to studying
\begin{equation}\label{someeqn}
\map(BX/\cZ,B^2\cZ)_{C(f)} \to \map(B\N/\cZ,B^2\cZ)_{C(f')}.
\end{equation}
Since $B^2\cZ$ is a product of Eilenberg-Mac~Lane spaces a transfer
argument (cf.\ \cite[Thm.~9.13]{DW94}) shows that this gives an embedding
as a retract. Since we assume $\AM: \pi_0(\Aut(BX/\cZ)) \cong
\pi_0(B\Aut(B\N/\cZ))$ we furthermore get that this is an isomorphism
on $\pi_0$ by the definition of $C(f)$ and $C(f')$. Let
$(W,L')$ denote the Weyl group of $X/\cZ$.
Write $B\cZ \simeq B^2A \times BA'$, where $A$ is a finite sum of
copies of $\Z_p$
and $A'$ is finite (cf.\ \cite[Thm.~1.1]{dw:center}). On $\pi_1$ the
map \eqref{someeqn} identifies with
$$
H^1(BX/\cZ;A') \oplus H^2(BX/\cZ;A) \to H^1(B\N/\cZ;A') \oplus H^2(B\N/\cZ;A).
$$
The group $H^1(B\N/\cZ;A')$ is zero since $\pi_1(B\N/\cZ) = W$ is
generated by elements of order prime to $p$, since $p$ is assumed to
be odd.

Furthermore, $H^2(B\N/\cZ;A)$ is related via the Serre spectral sequence
to the groups
$$H^2(BW;H^0(B^2L';A)),\ H^1(BW;H^1(B^2L';A)), \mbox{ and
}H^0(BW;H^2(B^2L';A)).$$
The first of these groups is zero since $W$
is generated by elements of order prime to $p$ by the assumption that $p$ is
odd. The second is obviously zero, and the last group is zero since
$H^0(W;\Hom(L',\Z_p)) = \Hom((L')_W,\Z_p) = 0$ because $(L')_W$ is finite.

Hence we get an isomorphism on $\pi_1$, since we already know that the
map is injective. On $\pi_2$ and $\pi_3$ the map identifies with
$$
H^0(BX/\cZ;A') \oplus H^1(BX/\cZ;A) \to H^0(B\N/\cZ;A') \oplus
H^1(B\N/\cZ;A)
$$
and $H^0(BX/\cZ;A) \to H^0(B\N/\cZ;A)$ respectively, and these maps are
obviously isomorphisms.
Hence $\map(BX/\cZ,B^2\cZ)_{C(f)} \to \map(B\N/\cZ,B^2\cZ)_{C(f')}$ is
a homotopy equivalence, which via the fibration sequences above implies
that $B\Aut(BX) \to B\Aut(B\N)$ is a homotopy equivalence as wanted.
\end{proof}

\begin{rem}
Consider $BX = B(\SU(3) \times S^1)\twocom$. This has center $\cZ =
(S^1)\twocom$ and $X/\cZ = \SU(3)\twocom$. By direct calculation (or appeal
to \cite[Cor.~3.5]{JMO95}) we have $B\Aut(BX/\cZ) \xrightarrow{\simeq} B\Aut(B\N/\cZ)$.
However $\AM: \pi_0(\Aut(BX)) \to \pi_0(\Aut(B\N))$ is not surjective by
Proposition~\ref{autn}, since $\Hom(W_{\SU(3)},\Z/2^\infty) = \Z/2$. This
shows that the assumption that $p$ is odd is necessary in
the last part of the above lemma.
\end{rem}

\begin{rem}
Suppose that $X$ is a connected $p$-compact
group. Fibration sequences with base space $B^2\pi_1(X)$ and fiber
$B(X\onecov)$ are in one-to-one correspondence with the set of maps
$[B^2\pi_1(X),B\Aut(B(X\onecov))]$. Likewise self-equivalences of
$BX$ can be expressed in terms of self-equivalences of $B(X\onecov)$
and $\pi_1(X)$, analogously to the lemmas above. Hence if we a
priori knew that Theorem~\ref{fundamentalgroup}  held true, i.e., if
we could read off $\pi_1(X)$ from $\N_X$ then the above methods
would reduce the proof of the main theorems to the simply connected
case, which could be used advantageously in the proofs. (See also
Remark~\ref{genpioneremark}.)
\end{rem}

\begin{rem}
The assumption in Lemma~\ref{ndetcenter}\eqref{ndetcenterobj} that
$\AM: \pi_0(\Aut(BX/\cZ)) \to \pi_0(\Aut(B\N/\cZ))$ is surjective has the
following origin. We have a canonical restriction map $H^2(BX/\cZ;\dZ)
\to H^2(B\N/\cZ;\dZ)$, which is injective by a transfer argument. Two
extension classes in $H^2(BX/\cZ;\dZ)$ give rise to isomorphic total
spaces if the extension classes are conjugate via the actions of
$\Aut(BX/\cZ)$ and $\Aut(\dZ)$ on $H^2(BX/\cZ;\dZ)$. The total spaces have
isomorphic maximal torus normalizers if the extension classes have
images in $H^2(B\N/Z;\dZ)$ which are conjugate under the actions of
$\Aut(B\N/\cZ)$ and $\Aut(\dZ)$, which could a priori be a weaker
notion.
\end{rem}

%%%%%%%%%%%%%%%%%%%%%%%%%%%%%%%%%%%%%%%%%%%%%%%%%%%%%%%%%%%%%%%%%%%%%%%%%%%

\section{An integral version of a theorem of Nakajima and realization
of $p$-compact groups} \label{nakajimasection}

The goal of this section is to prove an integral version of an algebraic
result of Nakajima (Theorem~\ref{nakajima}) and use this to prove
Theorem~\ref{polynomialrealizationthm} which, as part of our
inductive proof of Theorem~\ref{classthm}, will allow us to
construct the center-free $p$-compact groups corresponding to
$\Z_p$-reflection groups $(W,L)$ such that $\Z_p[L]^W$ is a
polynomial algebra. This will provide the existence part of
Theorem~\ref{classthm}. We feel that this way of showing existence, is
perhaps more straightforward than previous approaches; compare for instance
\cite{notbohm99}. (We refer to the introduction for the history behind
this result.)

\begin{thm} \label{nakajima}
Let $p$ be an odd prime and let $(W,L)$ be a finite
$\Z_p$-reflection group. For a subspace $V$ of $L\otimes \F_p$
we let $W_V$ denote the pointwise stabilizer of $V$ in $W$.
Then the following conditions are equivalent:
\begin{enumerate}
\item \label{nak:Lpoly}
$\Z_p[L]^W$ is a polynomial algebra.
\item \label{nak:stab-poly}
$\Z_p[L]^{W_V}$ is a polynomial algebra for all non-trivial subspaces
$V\subseteq L\otimes \F_p$.
\item \label{nak:stab-refl}
$(W_V,L)$ is a $\Z_p$-reflection group for all non-trivial subspaces
$V\subseteq L\otimes \F_p$.
\end{enumerate}
\end{thm}

\begin{rem}
An analog of the implication $\eqref{nak:Lpoly}\Rightarrow
\eqref{nak:stab-poly}$ where the ring $\Z_p$ is replaced by a field
was proven by Nakajima \cite[Lem.~1.4]{nakajima79} (in the case of
finite fields see also \cite[Thm.~1.4]{DW98} and \cite[Cor.~10.6.1]{ns02}). For fields of positive
characteristic the implication $\eqref{nak:stab-refl}\Rightarrow
\eqref{nak:Lpoly}$ does not hold; see \cite[Ex.~2.2]{km97} for more information
about this case.
Our proof unfortunately involves the classification of finite
$\Z_p$-reflection groups and some case-by-case checking. (See the
discussion following the proof of Theorem~\ref{borelelementaryabelian}
for related information.)
\end{rem}

\begin{proof}[Proof of Theorem~\ref{nakajima}]
To start, note that the implication $\eqref{nak:stab-poly}\Rightarrow
\eqref{nak:stab-refl}$ follows from the fact that if $\Z_p[L]^{W_V}$
is a polynomial algebra then $\Q_p[L \otimes \Q]^{W_V}$ is as well, so
$(W_V,L)$ is a $\Z_p$-reflection group by the Shephard-Todd-Chevalley
theorem (\cite[Thm.~7.2.1]{benson93} or \cite[Thm.~7.4.1]{smith95}).

To go further we want to see that the theorem is well behaved under
products, i.e., that if $(W,L) = (W',L') \times (W'',L'')$, then the
theorem holds for $(W,L)$ if it holds for $(W',L')$ and
$(W'',L'')$. This follows from the fact that the stabilizer in $W'
\times W''$ of an arbitrary subgroup of $(L' \otimes \F_p) \oplus (L''
\otimes \F_p)$ equals the stabilizer of the smallest product subgroup
containing it, combined with the fact that the tensor product of two
algebras is a polynomial algebra if and only if each of the factors
is. Hence to prove the remaining implications it follows from
Theorem~\ref{lattice-split} that it suffices to consider separately
the cases where $(W,L)$ comes from a compact connected Lie group and
the cases where $(W,L)$ is one of the exotic $\Z_p$-reflection groups.

Assume first that $(W,L) = (W_G,L_G \otimes \Z_p)$ for a compact
connected Lie group $G$.
If $\Z_p[L]^W$ is a polynomial algebra then by
Theorem~\ref{polyinvthm} (which involves case-by-case considerations
and $p$ odd) $BX = BG\pcom$ satisfies
$H^*(BX;\Z_p) \cong H^*(B^2L;\Z_p)^W$. We can identify $V \subseteq L
\otimes \F_p$ with a toral elementary abelian $p$-subgroup in $X$ and
by \cite[Thm.~1.3]{DW98}
$H^*(B\cC_X(V);\Z_p)$ is again a polynomial algebra concentrated in even
degrees. In particular
$\cC_X(V)$ is connected and by \cite[Thm.~7.6(1)]{dw:center} $W_{\cC_X(V)} =
W_V$. Hence, by Theorem~\ref{classicalborel}, $H^*(B\cC_X(V);\Z_p) \cong
H^*(B^2L;\Z_p)^{W_V}$, so $\Z_p[L]^{W_V}$ is a polynomial
algebra. This shows that
$\eqref{nak:Lpoly}\Rightarrow \eqref{nak:stab-poly}$ when $(W,L)$ comes
from a compact connected Lie group. To prove  $\eqref{nak:stab-refl}\Rightarrow
\eqref{nak:Lpoly}$ for Lie groups suppose that $(W,L)$ is a finite
$\Z_p$-reflection group corresponding to a $p$-compact group $X =
G\pcom$ such that $(W_V,L)$ is a $\Z_p$-reflection group for all
non-trivial  $V\subseteq L\otimes \F_p$. Since $p$ is
odd it follows by \cite[Thm.~7.6]{dw:center} that $\cC_X(V)$ is
connected for all non-trivial $V \subseteq L \otimes \F_p$. Hence,
since $X$ is assumed to come from a compact connected Lie group
\cite[Thm.~B]{borel61} (or \cite[Thm.~2.28]{steinberg75}) implies that
$H^*(BX;\Z_p)$ does not have torsion and hence $H^*(BX;\Z_p) \cong
H^*(B^2L;\Z_p)^W$ (cf.\ Theorem~\ref{classicalborel}). So $\Z_p[L]^W$
is a polynomial algebra as wanted.

Next we assume that $(W,L)$ is one of the exotic $\Z_p$-reflection
groups. By Theorem~\ref{polyinvthm}, $\Z_p[L]^W$ is a polynomial
algebra, so we only need to prove that
$\Z_p[L]^{W_V}$ is a polynomial algebra for any non-trivial
$V\subseteq L \otimes \F_p$. Furthermore, by
Theorem~\ref{polyinvthm}\eqref{lattice-simplycon}, $\F_p[L \otimes
\F_p]^W$ is a polynomial algebra.  Nakajima's result
\cite[Lem.~1.4]{nakajima79} shows that $\F_p[L\otimes \F_p]^{W_V}$ is
a polynomial algebra as well. Thus we are done if $p\nmid |W_V|$ by
Lemma~\ref{nonmodular}, which in particular covers the cases where
$p\nmid |W|$.

If $(W,L)$ belongs to family number $2$ on the Clark-Ewing list,
then since $p$ is odd, it is easily seen from the form of the
representing matrices (see Section~\ref{latticesection} for a concrete
description) that reduction mod $p$ gives a bijection between
reflections in $(W,L)$ and $(W,L \otimes \F_p)$. As
$\F_p[L\otimes \F_p]^{W_V}$ is a polynomial algebra it follows by the
Shephard-Todd-Chevalley theorem \cite[Thm.~7.2.1]{benson93} that
$W_V\subseteq \GL(L\otimes \F_p)$ is a
reflection group. Thus $(W_V,L)$ is a $\Z_p$-reflection group. Since
the representing matrices are monomial, it follows by
\cite[Thm.~2.4]{nakajima79} that $\Z_p[L]^{W_V}$ is a polynomial algebra.

By Theorem~\ref{lattice-split} only four
cases remain, namely the Zabrodsky-Aguad\'e cases
$(W_{12},p=3)$, $(W_{29},p=5)$, $(W_{31},p=5)$ and $(W_{34},p=7)$.
For each of these a direct computation, for instance easily done
with the aid of a computer, shows that if $S$
is a Sylow $p$-subgroup of $W$, then $U=(L\otimes \F_p)^S$ is
$1$-dimensional and $(W_U,L)$ is isomorphic to
$(\Sigma_p,L_{\SU(p)}\otimes \Z_p)$. (A more ad hoc construction of this reflection subgroup can also be found in Aguad\'e
\cite{aguade89}.)
Hence we see
that if $V\subseteq L\otimes \F_p$ is non-trivial then either $p\nmid |W_V|$
or $V$ is $W$-conjugate to $U$. But in these cases we already know
that $\Z_p[L]^{W_V}$ is a polynomial algebra.
\end{proof}

\begin{thm}[Inductive Polynomial Realization Theorem]\label{polynomialrealizationthm}
Let $p$ be an odd prime and let $(W,L)$ be a finite
$\Z_p$-reflection group with the property that $\Z_p[L]^W$ is a
polynomial algebra over $\Z_p$. Note that for any non-trivial elementary abelian $p$-subgroup $V$ of $\dT = L \otimes \Z/p^\infty$, the subgroup $W_V$ of $W$ fixing $V$ pointwise is again a finite $\Z_p$-reflection group by Theorem~\ref{nakajima}.

Assume that for all such $V$ there exists a $p$-compact
group $F(V)$  with discrete approximation to its maximal torus
normalizer given by $\dT \semi W_V$ such that $F(V)$ is determined by
$\N_{F(V)}$, $\AM: \Aut(BF(V)) \xrightarrow{\cong} \Aut(B\N_{F(V)})$, and
$H^*(BF(V);\Z_p) \xrightarrow{\cong} H^*(B^2L;\Z_p)^{W_V}$. Then there exists
a connected $p$-compact group $X$ with discrete approximation to its maximal
torus normalizer given by $\dT \semi W$ satisfying the same properties
as listed for $F(V)$.
\end{thm}

\begin{proof}
First note that by Theorem~\ref{nakajima}
$(W_V,L)$ is again a $\Z_p$-reflection group and $\Z_p[L]^{W_V}$ is a
polynomial algebra, so the assumptions make sense.
Set $\dN = \dT \semi W$. We want to construct a
candidate `centralizer decomposition' diagram. Let $\A$ be the
category with objects the non-trivial elementary abelian $p$-subgroups
$V$ of $\dT$ and morphisms the homomorphisms between them induced by
inclusions of subgroups and conjugation by elements in $W$. We now
define a functor $F$ from $\A^{\op}$ to $p$-compact groups and
conjugation classes of morphisms. On objects we send $V$ to $F(V)$. By
assumption  $j_V: C_{\dN}(V) \to F(V)$ is a discrete approximation
to the maximal torus normalizer in $F(V)$.
Now let $\varphi: V \to V'$ be a morphism in $\A$, induced by
conjugation by an element $x \in
W$ and consider the diagram
$$
\xymatrix{
V' \ar[r] & C_{\dN}(V') \ar[r]^-{c_{x^{-1}}}
\ar[d]^-{j_{V'}} & C_{\dN}(V) \ar[d]^-{j_V}\\
& F(V') & F(V)
}
$$
Taking the centralizer of the composite map $x^{-1}: V' \to F(V)$
we get a space $\cC_{F(V)}(x^{-1}) = \Omega
\map(BV',BF(V))_{Bx^{-1}}$, which has discrete approximation to its
maximal torus normalizer equal to $C_{\dN}(V')$. By assumption we get
a unique (up to conjugacy) isomorphism $F(V') \to \cC_{F(V)}(x^{-1})$
under $C_{\dN}(V')$. By composing with the evaluation
$\cC_{F(V)}(x^{-1}) \to F(V)$, we get a morphism $F(\varphi): F(V') \to
F(V)$. We need to check that this gives us a well-defined functor
from $\A^{\op}$ to the homotopy category of spaces, i.e., that for $V
\xrightarrow{\varphi} V' \xrightarrow{\psi} V''$, $F(\psi \varphi)$ is conjugate to
$F(\varphi) F(\psi)$. To see this suppose that $\psi$ is induced by
conjugation by $y \in W$ and consider the following diagram with obvious maps
$$
\xymatrix{
F(V) &  \cC_{F(V)}(x^{-1}) \ar[l]_-{\text{ev}} & F(V')
\ar[l]_-{\cong} & \cC_{F(V')}(y^{-1}) \ar[l]_-{\text{ev}} \ar[dl]^-\cong & \ar[l]_-{\cong} F(V'') \ar[ddll]^-{\cong}\\
& & \cC_{\cC_{F(V)}(x^{-1})}(\widetilde{x^{-1}y^{-1}}) \ar[ul]^-{\text{ev}} \ar[d]^-\cong & & \\
& & \cC_{F(V)}(x^{-1}y^{-1}) \ar[uull]^-{\text{ev}}
}
$$
(Here $\tilde{(\cdot)}$ denotes the adjoint map which is explained
in Construction~\ref{adjoint}.) Note that the bottom composite from
$F(V'')$ to $F(V)$ is $F(\psi \varphi)$ and the top composite is
$F(\varphi) F(\psi)$. The top triangle is commutative, since the
lower isomorphism in that triangle is just the map obtained by
taking centralizers of the upper one. The rightmost square is
homotopy commutative, since the corresponding square of
isomorphisms between centralizers in $\dN$ is commutative, using our
assumptions that maps are detected here. Finally, the leftmost
square is homotopy commutative, by definition of the adjoint
construction.

We hence get a well-defined functor $BF: \A^{\op} \to
Ho(\Top)$, where $Ho(\Top)$ denotes the homotopy category of spaces,
on objects given by $V \mapsto BF(V)$. By construction the functor
obtained when taking cohomology of this diagram, can be identified
with the canonical functor which on objects is given by $V \mapsto
H^*(B\dT;\Z_p)^{W_V}$.

We want to lift this to a diagram in the category of spaces. The
obstruction theory for doing this is described in
\cite[Thm.~1.1]{DK92}, noting that by \cite[Lem.~11.15]{dw:center}
our diagram is a so-called centric diagram so the assumptions of that
theorem are satisfied.

By looking at their cohomology we see that all the spaces $F(V)$
are connected and hence by \cite[Thm.~7.5]{dw:center} have center
given by ${\dT}^{W_V}$, since $p$ is odd. In particular (see e.g.\
Lemma~\ref{centralizercenter}) the homotopy groups of $\cZ F(V)$ are
given by $\pi_0(\cZ F(V)) = H^1(W_V;L)$ and $\pi_1(\cZ F(V)) =
L^{W_V}$. By \cite[\S 8]{DW93} (for details see
Section~\ref{obstructionssection}) $\lim_{V \in \A}^*\pi_*(F(-)) = 0$,
so by \cite[Thm.~1.1]{DK92} there exists a (unique) lift of our
functor $BF$ to a functor $\widetilde{BF}$ landing in $\Top$. Set
$BX = (\hocolim_\A\widetilde{BF})\pcom $.

The spectral sequence for calculating the cohomology of a
homotopy colimit \cite[XII.4.5]{bk} has $E_2$-term given by
$E_2^{i,j} = \lim^i_{V \in \A} H^j(B\dT;\Z_p)^{W_V}$. But again by
\cite[\S 8]{DW93} these groups vanish for $i>0$ and for $i=0$ give
$\lim^0_\A H^*(B\dT;\Z_p)^{W_V} \cong H^*(B\dT;\Z_p)^W$. Hence the
spectral sequence collapses onto the vertical axis, and we get
$H^*(BX;\Z_p) \cong H^*(B\dT;\Z_p)^W$.

Since $H^*(BX;\Z_p)$ is a polynomial algebra $H^*(X;\Z_p)$ will be an
exterior algebra on odd generators (cf.\
Theorem~\ref{classicalborel}), so $X$ is indeed a connected $p$-compact
group. The fact that $X$ is determined by $\N$ and satisfies $\AM:
\Aut(BX) \xrightarrow{\cong} \Aut(B\N)$, also follows easily from the
above---the details are given in the proof of
Theorem~\ref{collectlemma}.
\end{proof}

\begin{rem}
Note that Theorem~\ref{polynomialrealizationthm} in itself does
not quite give a stand-alone proof of the realization and uniqueness
of all center-free $p$-compact groups with Weyl group satisfying that
$\Z_p[L]^W$ is a polynomial algebra, since $(W_V,L)$ is not
center-free which prevents the obvious induction from working. Compare to the proof 
of Theorem~\ref{classthm}; the main problem is that unitary groups will occur in most 
decompositions, but their adjoint forms do not have cohomology rings which are polynomial algebras.
\end{rem}

%%%%%%%%%%%%%%%%%%%%%%%%%%%%%%%%%%%%%%%%%%%%%%%%%%%%%%%%%%%%%%%%%%%%%%%%%%%

\section{Non-toral elementary abelian $p$-subgroups of simple center-free Lie groups}\label{liesection}

In this section we determine, for an odd prime $p$, all conjugacy
classes of non-toral elementary abelian $p$-subgroups $E$, of any
simple center-free compact Lie group $G$, as well as their
centralizers $C_G(E)$ and {\em Weyl groups} $W(E) =
N_G(E)/C_G(E)$. (Recall that a subgroup of $G$ is called {\em toral}
if it is contained in a torus in $G$ and {\em non-toral} otherwise.)

Our strategy is as follows. Since $p$ is odd, the groups $G$ we need
to consider are the projective unitary groups $\PU(n)$ and the
exceptional groups. The groups $\PU(n)$ are easy to deal with and we
only expand slightly on the work of Griess \cite{griess}. For the
exceptional groups the maximal non-toral elementary abelian $p$-subgroups are
also determined by Griess \cite{griess}. We first find these subgroups explicitly and then get lower bounds for their Weyl
groups by producing explicit elements in their normalizers. From this
we are able to identify the non-maximal non-toral elementary abelian
$p$-subgroups and get lower bounds for their Weyl groups. Finally we
get exact results on the Weyl groups by computing centralizers.

In accordance with the standard literature we will in this section
state and prove all theorems in the context of linear algebraic groups
over the complex numbers $\C$---we state
in Proposition~\ref{algebraicvslie} why this is equivalent to
considering compact Lie groups. (The results for $G(\C)$ can
furthermore be translated into results for $G(F)$ for any
algebraically closed field $F$ of characteristic prime to $p$, see
\cite[Thm.~1.22]{GR98} and \cite{FM88}.)

This section is divided into five subsections. The first recalls some
results from the theory of linear algebraic groups and discusses the
relationship with compact Lie groups. In the second subsection we
determine the elementary abelian $p$-subgroups of the projective
unitary groups and the final subsections deal with the
elementary abelian $3$-subgroups of the groups of type $E_6$, $E_7$
and $E_8$ respectively. (The remaining non-trivial cases $E_8(\C),
p=5$ and $F_4(\C), p=3$ are treated completely in \cite[Lem.~10.3 and
Thm.~7.4]{griess}.)

For some of our computations for the groups $3E_6(\C)$ and $E_8(\C)$
we have used the computer algebra system {\sc Magma} \cite{magma1},
although this reliance on computers could if needed be replaced by
some rather tedious hand calculations.

\begin{notation}
We now collect some notation which will be used multiple times throughout the computations in this section.
We use standard names for the linear algebraic groups we consider,
e.g.\ $3E_6(\C)$ denotes the simply connected group of type $E_6$ over
$\C$ and $E_6(\C)$ denotes its adjoint version. We let $\T_n$ denote
an $n$-dimensional torus, i.e., $\T_n=(\C^\times)^n$.

To describe centralizers we follow
standard notation for extensions of groups, cf.\ the $\mathbb{ATLAS}$
\cite[p.~xx]{CCNPW85}. Thus $A:B$ denotes a semidirect product,
$A\cdot B$ denotes a non-split extension and $A\circ_C B$ denotes a
central product.

Whenever $E$ is a concrete elementary abelian $p$-group of rank $n$
we will always fix an ordered basis of $E$, so that $\GL(E)$
identifies with $\GL_n(\F_p)$. We make the
standing convention that all matrices acts on columns.

We identify a permutation $\sigma$ in the symmetric group $\Sigma_n$
with its permutation matrix $A=[a_{ij}]$ given by
$a_{ij}=\delta_{i,\sigma(j)}$ where $\delta$ is the Kronecker delta.

If $K$ is a field, we let $M_n(K)$ denote the set of $n\times n$-matrices
over $K$. For $a_1,\ldots, a_n\in K$ we let $\diag(a_1,\ldots, a_n)\in
M_n(K)$ denote the diagonal matrix with the $a_i$'s in the
diagonal. For $1\leq i,j\leq n$, $e_{ij}\in M_n(K)$ denotes the
matrix whose only non-zero entry is $1$ in position $(i,j)$. Given matrices
$A_1\in M_{n_1}(K)$, \ldots, $A_m\in M_{n_m}(K)$ we let $A_1\oplus
\ldots\oplus A_m$ denote the $n\times n$-block matrix with the $A_i$'s
in the diagonal, $n=n_1+\ldots +n_m$. We also need the `blowup'
homomorphism $\Delta_{n,m}: M_{n}(K) \longrightarrow M_{mn}(K)$
defined by replacing each entry $a_{ij}$ by $a_{ij} I_m$, where
$I_m\in M_m(K)$ is the identity matrix.

As $p=3$ for all the exceptional groups we consider, we use some special
notation. An arbitrary element of $\F_3$ is denoted by $*$, and
$\varepsilon$ denotes an element of the multiplicative group
$\F_3^\times$. We let $\omega=e^{2\pi i/3}$ and $\eta =
e^{2\pi i/9}$ and define elements $\beta, \gamma, \tau_1, \tau_2\in
\SL_3(\C)$ by $\beta = \diag(1,\omega,\omega^2)$,
$$
\gamma = (1,2,3) =
\begin{bmatrix}
0 & 0 & 1\\
1 & 0 & 0\\
0 & 1 & 0
\end{bmatrix}
,\qquad
\tau_1 = \frac{e^{-\pi i/18}}{\sqrt{3}}
\begin{bmatrix}
1 & \omega^2 & 1 \\
1 & 1 & \omega^2 \\
\omega^2 & 1 & 1
\end{bmatrix}
$$
and $\tau_2 = \diag(\eta,\eta^{-2},\eta)$. Note that
$\beta^{\tau_1}=\beta\gamma$, $\gamma^{\tau_1}=\gamma$,
$\beta^{\tau_2}=\beta$ and $\gamma^{\tau_2}=\beta\gamma$.

To distinguish subgroups we use {\em class distributions}. As an example the group $3E_6(\C)$ contains $7$ conjugacy classes of elements of order $3$ labeled $\mathbf{3A}$, $\mathbf{3B}$, $\mathbf{3B'}$, $\mathbf{3C}$, $\mathbf{3D}$, $\mathbf{3E}$ and $\mathbf{3E'}$, cf. \ref{classes:3E6}, and the fact that the subgroup $E^4_{3E_6}$ from Theorem~\ref{nontorals:3E6-3} below has class distribution $\mathbf{3C^{78}3E^13{E'}^1}$ means that it (apart from the identity) contains  $78$ elements from the conjugacy class $\mathbf{3C}$ and one element from each of the conjugacy classes $\mathbf{3E}$ and $\mathbf{3E'}$.
\end{notation}

\subsection{Recollection of some results on linear algebraic groups} \label{lietheory}

Recall that a (not necessarily connected) linear
algebraic group $G$ is called {\em reductive} if its unipotent
radical, i.e., the largest normal connected unipotent subgroup of $G$,
is trivial.

\begin{thm} \label{alggrpfacts}
Let $G$ be a linear algebraic group over an algebraically closed field $K$.
\begin{enumerate}
\item \label{centr-toral}
If $A$ is a subgroup of $G$ and $S$ is some subset of $A$, then
$A$ is toral in $G$ if and only if $A$ is toral in $C_G(S)$.
\item \label{weyl-section}
If $H$ is a maximal torus of $G$, then two subsets of $H$ are conjugate in $G$
if and only if they are conjugate in $N_G(H)$. If $A$ is a toral subgroup of
$G$, then $W(A)=N_G(A)/C_G(A)$ is isomorphic to a subquotient of the Weyl group
$W=N_G(H)/H$ of $G$.
\item \label{simply-con}
Assume that $G$ is a connected reductive group such that $G'$ is simply
connected. Then the centralizer of any semisimple element in $G$
is connected. In particular, if $A$ is an abelian subgroup of $G$ consisting of
semisimple elements generated by at most two elements, then $A$ is toral.
\item \label{semisimpleauto}
If $G$ is reductive and $\sigma$ is a semisimple automorphism of
$G$, then the fixed point subgroup $G^\sigma$ is reductive and
contains a regular element of $G$.
\item \label{covering}
Assume that $G$ is a connected reductive group, let $Z\subseteq G$ be a
central subgroup, and let $\pi:G\rightarrow G/Z$ be the quotient
homomorphism. If $A$ is a subgroup of $G$, then $A$ is toral in $G$
if and only if $\pi(A)$ is toral in $G/Z$.
\item \label{dim-centr}
Assume $\ch K=0$ and let $\mathfrak{g}$ be the Lie algebra of
$G$. If $S\subseteq G$ is a finite subset of $G$, then the Lie algebra of
$C_G(S)$ is given by
$$
\mathfrak{c}_\mathfrak{g}(S)=\{x\in \mathfrak{g} |
\operatorname{Ad}(s)(x)=x \mbox{ for all $s\in S$} \}.
$$
In particular, if $S\subseteq G$ is a finite subgroup, then
$$
\dim C_G(S) = \frac{1}{|S|} \sum_{s\in S} \operatorname{tr}_{\mathfrak{g}}\;
\operatorname{Ad}(s).
$$
\end{enumerate}
\end{thm}

\begin{proof}
\eqref{centr-toral}:
Obviously, if $A$ is toral in $C_G(S)$ then $A$ is toral in
$G$. Conversely, if $A$ is toral in $G$, then $A\subseteq H$ for a
torus $H$ in $G$. Since $S\subseteq A$ we get $H\subseteq C_G(S)$ and thus
$A$ is toral in $C_G(S)$.

\eqref{weyl-section}:
The first part follows by the Frattini argument: Assume that $A,A^g\subseteq
H$ are conjugate subsets of $H$. Then $H$ and $H^{g^{-1}}$ are maximal tori of
$C_G(A)$ and thus conjugate in $C_G(A)$ (cf.\
\cite[Cor.~21.3.A]{humphreys75}). Thus we may write $H=H^{g^{-1}c}$
for some $c\in C_G(A)$ and we conclude that $n=g^{-1}c\in
N_G(H)$. Then $A^{n^{-1}}=A^{c^{-1}g}=A^g$, which proves the first
part. The second part follows similarly,
cf.\ \cite[Prop.~1.1(i)]{Liebeck-Seitz}.

\eqref{simply-con}:
The first part which is due to Steinberg is proved in
\cite[Thm.~3.5.6]{carter:rep}. The second part follows from the first,
cf.\ \cite[II.5.1]{SS70}.

\eqref{semisimpleauto}:
We can assume $G$ to be connected. In case $G$
is semisimple and simply connected the first claim is proved in
\cite[Thm.~8.1]{steinberg68} and the general case reduces to this
one. Indeed we can find a finite cover $\tilde{G}$ of $G$ which is a
direct product of a semisimple simply connected group and a torus, and
$\sigma$ lifts to a semisimple automorphism of $\tilde{G}$ by
\cite[9.16]{steinberg68}. For the second claim see
\cite[Thms.~2~and~3]{winter66} or \cite[Pf.~of Thm.~II.5.16]{SS70} in case
$G$ is semisimple; the general case clearly reduces to this one.

\eqref{covering}:
By \cite[Cor.~21.3.C]{humphreys75} we know that if $H$ is a maximal
torus of $G$, then $\pi(H)$ is a maximal torus of $G/Z$, and all maximal tori
of $G/Z$ are of this form. Thus if $A$ is toral in $G$, then $\pi(A)$ is toral
in $G/Z$. Conversely, if $H'$ is a maximal torus of $G/Z$ containing $\pi(A)$,
then by the above we have $H'=\pi(H)$ for some maximal torus $H$ of $G$. Thus
we get $A\subseteq \left<H,Z\right>$. However since $G$ is connected and
reductive, we get $Z\subseteq H$ by \cite[Cor.~26.2.A(b)]{humphreys75}.
Thus $A\subseteq H$ and we are done.

\eqref{dim-centr}:
In case $S$ consists of a single element, the first part follows from
\cite[Thm.~13.4(a)]{humphreys75} (note that the connectivity assumption
in  \cite[Thm.~13.4]{humphreys75} is only used in \cite[Thm.~13.4(b)]{humphreys75}). The general
case follows from this by applying \cite[Thm.~12.5]{humphreys75} to
the centralizers $C_G(s)$, $s\in S$.

Now assume that $S\subseteq G$ is a finite subgroup, and let $\chi$
denote the character of the adjoint representation of $G$ restricted
to $S$. Then the dimension of
$$
\mathfrak{c}_\mathfrak{g}(S)=\{x\in \mathfrak{g} |
\operatorname{Ad}(s)(x)=x \mbox{ for all $s\in S$} \}.
$$
equals the multiplicity of the trivial character in $\chi$. By the
orthogonality relations this is given by
$$
\left(\chi\!\mid\! 1\right) = \frac{1}{|S|} \sum_{s\in S} \chi(s),
$$
and we are done.
\end{proof}

We also need the following result whose proof is extracted from
\cite{serre99}.

\begin{thm} \label{alggroups:preferredlifts}
Let $G$ be a reductive linear algebraic group, $H$ a maximal torus of $G$ and
let $N=N_G(H)$. Let $U\subseteq N$ be a subgroup consisting of
semisimple elements such that $U/(U\cap H)$ is cyclic. Let $S$ be the
identity component of $H^U$ (the subgroup of $H$ fixed by $U$), and
assume that $S$ is a maximal torus of $C_G(U)$. Then
$C_N(U)=N_{C_G(U)}(S)$ and in particular $C_N(U)$ is a maximal torus
normalizer in $C_G(U)$.
\end{thm}

\begin{proof}
As any element of $C_N(U)$ normalizes $H^U$ and hence also its
identity component $S$, the inclusion $C_N(U)\subseteq N_{C_G(U)}(S)$
is clear. Suppose conversely that $x\in N_{C_G(U)}(S)$. Let $C=U\cap
H$. From \cite[2.15(d)]{BT65} it follows that $G^C$ is reductive.
By assumption the cyclic group $U/C$ acts by semisimple automorphisms
on $G^C$. It now follows from Theorem~\ref{alggrpfacts}\eqref{semisimpleauto}
that $G^U = (G^C)^{U/C}$ is reductive and that every maximal torus of
$G^U$ is contained in a unique maximal torus of $G^C$. Since
$C\subseteq H$, we see that $H$ is the maximal torus of $G^C$
containing $S$. As $H^x$ is also a maximal torus of $G^C$ and
$H^x\supseteq S^x= S$ we conclude that $H^x=H$. Thus $x\in C_N(U)$
proving the result.
\end{proof}

We now explain the relationship between reductive complex linear
algebraic groups and compact Lie groups. If $G$ is a complex linear algebraic
group then the underlying variety of $G$ is an affine complex variety. By
endowing this variety with the usual Euclidean topology instead of the
Zariski topology we may view $G$ as a complex Lie group since the group
operations are given by polynomial maps.

\begin{prop}\label{algebraicvslie}
Let $G$ be a complex linear algebraic group.
\begin{enumerate}
\item
Viewed as a Lie group, $G$ contains a maximal compact subgroup which is unique
up to conjugacy, and for any such subgroup $K$ we have a diffeomorphism
$G\cong K\times \R^s$ for some $s$.
\item
Let $K$ be a maximal compact subgroup of $G$, and let $S,S'\subseteq
K$ be two subsets. If $S'=S^g$ for some $g\in G$, then there exists
$k\in K$ such that $x^k=x^g$ for all $x\in S$.
\item
Assume that $G$ is reductive. If $S$ is a finite subgroup of $G$, then
$C_G(S)$ is also reductive. If $K$ is a maximal compact subgroup of $G$
containing $S$, then $C_K(S)$ is a maximal compact subgroup of $C_G(S)$.
\item
If $G$ is reductive and $K$ is a maximal compact subgroup of $G$, then we have
a diffeomorphism $Z(G)\cong Z(K)\times \R^s$ for some $s$.
\end{enumerate}
\end{prop}

\begin{proof}
Note first that the identity component $G_1$ of $G$ seen as a Lie group
coincides with the identity component of $G$ seen as a linear algebraic group
\cite[Ch.\ 3, \S 3, no.\ 1]{onishchik-vinberg}. Thus $G/G_1$ is finite
by \cite[Prop.~7.3(a)]{humphreys75}. The first claim is now part of the
Cartan-Chevalley-Iwasawa-Malcev-Mostow theorem
\cite[Ch.~XV, Thm.~3.1]{hochschild} and the second claim also follows
from this, cf.\ \cite[Ch.~V, \S 24.7, Prop.~2]{borel91}.

In case $G$ is reductive it is possible to give a more explicit
form of the decomposition above. By \cite[Thm.~8.6]{humphreys75} we
may assume that $G$ is a closed subgroup of
$\GL(V)$ for some complex vector space $V$. From
\cite[Thm.~5.2.8]{onishchik-vinberg} it follows that $G$ has a compact
real form $K$ and we may thus choose a non-degenerate Hermitian inner
product on $V$ which is invariant under $K$ (e.g.\ by
\cite[Thm.~3.4.2]{onishchik-vinberg}). Let $\U(V)$ denote the
set of operators in $\GL(V)$ which are unitary with respect to the
chosen inner product. Using \cite[Problems~5.2.3 and
5.2.4]{onishchik-vinberg} we see that $G\subseteq \GL(V)$ is
self-adjoint and that $K=G\cap \U(V)$. The last part now follows by
combining \cite[Cor.~2 of Thm.~5.2.2]{onishchik-vinberg}
with \cite[Cor.~2 of Thm.~5.2.1]{onishchik-vinberg}.

If $S$ is a subgroup of $K$, then $S$ is self-adjoint since $K$ consists
of unitary operators. In particular $C_G(S)$ is also a self-adjoint subgroup of
$\GL(V)$, so by \cite[Cor.~3 of Thm.~5.2.1]{onishchik-vinberg} it
follows that $C_K(S)=C_G(S)\cap \U(V)$ is a maximal compact subgroup of
$C_G(S)$.

It only remains to prove that $C_G(S)$ is reductive for a finite subgroup $S$
of $G$. However by \cite[Problem~6.11]{onishchik-vinberg} and
\cite[Ch.~4,\,\S 1, no.~2]{onishchik-vinberg} we see that a complex
linear algebraic group is reductive if and only if its Lie algebra is
reductive. Thus it suffices to prove that the Lie algebra of $C_G(S)$ is
reductive. However by Theorem~\ref{alggrpfacts}\eqref{dim-centr} this Lie
algebra equals
$$
\mathfrak{c}_\mathfrak{g}(S)=\{x\in \mathfrak{g} |
\operatorname{Ad}(s)(x)=x \mbox{ for all $s\in S$} \}.
$$
where $\mathfrak{g}$ denotes the Lie algebra of $G$. The claim now follows from
\cite[Ch.~V,\,\S 2, no.~2, Prop.~8]{chevalley:3}.
\end{proof}

\subsection{The projective unitary groups} \label{punsection}
The purpose of this short subsection is to describe the non-toral
elementary abelian subgroups of $\PGL_n(\C)$, which by
Proposition~\ref{algebraicvslie} is equivalent to finding them for its
maximal compact subgroup $\PU(n)$, as well as to give information about centralizers
and Weyl groups. The subgroups are easily determined and are described
in \cite[Thm.~3.1]{griess}---we here just add some extra information about
centralizers and Weyl groups which we need in our proof of
Theorem~\ref{classthm}.

We first introduce a useful subgroup. If $p^r$ divides $n$ write $n
=p^r k$ and consider the extraspecial group $p^{1+2r}_+$ embedded in
$\GL_n(\C)$ by taking $k$ copies of one of the $p-1$ faithful
irreducible $p^r$-dimensional representations. (They all have the same
image; see \cite[Satz~V.16.14]{huppert67}.) Note that this embedding
maps the center of $p^{1+2r}_+$ to the elements of order $p$ in the
center of $\GL_n(\C)$. Let $\Gamma_r$ denote the subgroup of
$\GL_n(\C)$ given by the subgroup generated by the image of
$p^{1+2r}_+$ and the center of $\GL_n(\C)$. Note that as an abstract
group $\Gamma_r$ fits into an extension sequence
$$ 1 \to \C^\times \to \Gamma_r \to \bar \Gamma_r \to 1$$
where $\C^\times$ identifies with the center of $\GL_n(\C)$ and $\bar
\Gamma_r \cong (\Z/p)^{2r}$ identifies with the image of $\Gamma_r$ in
$\PGL_n(\C)$.  (The matrices for $\Gamma_r$ are written explicitly for
$k=1$ in \cite[p.~56--57]{oliver94} where it is called $\Gamma^U_{p^r}$.)

\begin{thm}\label{punelementary}
Suppose $E$ is a non-toral elementary abelian $p$-subgroup of
$\PGL_n(\C)$ for an arbitrary prime $p$. Then, up to conjugacy, $E$
can be written as $E = \bar \Gamma_r \times \bar A$, for some $r\geq 1$ with
$n=p^r k$ and some abelian subgroup $A$ of $C_{\GL_n(\C)}(\Gamma_r) \cong
\GL_k(\C)$.

For a given $r$, the conjugacy classes of such subgroups $E$
are in one-to-one correspondence with the conjugacy classes of toral
elementary abelian $p$-subgroups $\bar A$ of $\PGL_k(\C) \cong
C_{\PGL_n(\C)}(\bar \Gamma_r)_1$ (allowing the trivial subgroup), and
the centralizer of $E$ is given by $C_{\PGL_n(\C)}(E) \cong \bar
\Gamma_r \times C_{\PGL_k(\C)}(\bar A)$.

The Weyl group equals
$$
W_{\PGL_n(\C)}(E) =
\begin{bmatrix}
\Sp(\bar \Gamma_r) & 0 \\
* & W_{\PGL_k(\C)}(\bar A)
\end{bmatrix}
$$
Here  $\Sp(\bar \Gamma_r)$ is the symplectic group relative to the
symplectic product coming from the commutator product $[\cdot,\cdot]:
\bar \Gamma_r \times \bar \Gamma_r \to \Z/p \subseteq \C^\times$ and
the symbol $*$ denotes a $\operatorname{rank} \bar A \times 2r$ matrix with
arbitrary entries.

An element $ \alpha \in \Sp(\bar \Gamma_r) \subseteq
W_{\PGL_n(\C)}(E)$ acts as up to conjugacy as $\alpha \times 1$ on
$C_{\PGL_n(\C)}(E) \cong \bar \Gamma_r \times C_{\PGL_k(\C)}(\bar
A)$.
\end{thm}

\begin{proof}[Sketch of proof:]
The existence of the decomposition $E = \bar \Gamma_r \times \bar A$
follows from Griess \cite[Thm.~3.1]{griess} and the statements about
uniqueness follow by representation theory of the extraspecial
$p$-groups (cf.\ \cite[Ch.~5.5]{gorenstein} or
\cite[Satz~V.16.14]{huppert67}). Since the image of $p^{1+2r}_+$ is
the sum of $k$ identical irreducible representations we have
$C_{\GL_n(\C)}(\Gamma_r) \cong \GL_k(\C)$ by Schur's lemma (see also
\cite[Prop.~4]{oliver94}). From this the centralizer in $\PGL_n(\C)$
can easily be worked out.

In the case where $\bar A$ is trivial the statement about Weyl groups
is given in \cite[Thm.~6]{oliver94} (and just uses elementary
character theory). The general case follows similarly, again using
character theory. 

For the statement about the Weyl group action, first note that $\Out(
\bar \Gamma_r \times \PGL_k(\C)) \cong \Aut(\bar \Gamma_r) \times
\Out(\PGL_k(\C))$. An element $\alpha \in \Sp(\bar \Gamma_r) =
W_{\PGL_n(\C)}(\bar \Gamma_r)$ acts as an inner automorphism on $\PGL_k(\C)$
since this is true for the action on $C_{\GL_n(\C)}(\Gamma_r) \cong
\GL_k(\C)$ by character theory. Hence we can choose a representative
$g \in N_{\PGL_n(\C)}(\bar \Gamma_r)$ of $\alpha$ which acts as
$\alpha \times 1$ on $C_{\PGL_n(\C)}(\bar \Gamma_r) \cong \bar \Gamma_r \times
\PGL_k(\C)$. Hence $g$ is also a representative of $\alpha \in
\Sp(\bar \Gamma_r) \subseteq W_{\PGL_n(\C)}(\bar \Gamma_r \times \bar
A)$. The claim now follows.
\end{proof}

\subsection{The groups $E_6(\C)$ and $3E_6(\C)$, $p=3$}\label{type:E6}
In this subsection we consider the elementary abelian $3$-subgroups of
the groups of type $E_6$ over $\C$. The group $3E_6(\C)$ has two
non-isomorphic faithful irreducible $27$-dimensional
representations. These have highest weight $\lambda_1$ and $\lambda_6$
respectively and are dual to each other. An explicit construction of
$3E_6(\C)$ based on one of these representations was originally given
by Freudenthal \cite{freudenthal}. This construction is described in
more detail in \cite[\S 2]{CoWa} from which we take most of our
notation. In particular we let $\K$ be the $27$-dimensional complex
vector space consisting of triples $m=(m_1,m_2,m_3)$ of complex $3\times
3$-matrices $m_i$, $1\leq i\leq 3$, where addition and scalar
multiplication is defined coordinatewise. We define a cubic form
$\left<\cdot\right>$ on $\K$ by
$$
\left<m\right> = \det(m_1)+\det(m_2)+\det(m_3)-
\operatorname{tr}(m_1 m_2 m_3).
$$
Then $3E_6(\C)$ is the subgroup of $\GL(\K)$ preserving
the form $\left<\cdot\right>$. Moreover the stabilizer in
$3E_6(\C)$ of the element $(I_3,0,0)\in\K$ is the group
$F_4(\C)$. For $g_1, g_2, g_3\in \SL_3(\C)$ we
have the element $s_{g_1, g_2, g_3}$ of $3E_6(\C)$ given by
$$
s_{g_1, g_2, g_3}\left(m_1,m_2,m_3\right) = \left(g_1 m_1 g_2^{-1}, g_2 m_2
g_3^{-1}, g_3 m_3 g_1^{-1}\right)
$$
for $m=(m_1,m_2,m_3)\in\K$. This gives a representation of
$\SL_3(\C)^3$ which has kernel $C_3$ generated by $(\omega I_3,\omega
I_3,\omega I_3)$, and we thus get an embedding of $\SL_3(\C)^3/C_3$ in
$3E_6(\C)$. We will denote the element $s_{g_1, g_2, g_3}$ by
$\left[g_1,g_2,g_3\right]$.

We let $\{e_{j,k}^i\}$, $1\leq i,j,k\leq 3$ be the natural basis of
$\K$ consisting of the elements $e_{j,k}^i$ whose entries are all $0$
except for the $(j,k)$-entry of the $i$th matrix which equals $1$. The
elements of $3E_6(\C)$ which act diagonally with respect to this
basis of $\K$ form a maximal torus $H$ in
$3E_6(\C)$. Let $m_i^{j,k}$ denote the $(j,k)$-entry of the matrix
$m_i$. We then have $H$-invariant subgroups
$$
u_{\alpha_1}(t) = \left[
I_3, I_3+t e_{1,3}, I_3
\right],\qquad
u_{-\alpha_1}(t) = \left[
I_3, I_3+t e_{3,1}, I_3
\right],
$$
$$
u_{\alpha_2}(t) = \left[
I_3+t e_{2,1}, I_3, I_3
\right],\qquad
u_{-\alpha_2}(t) = \left[
I_3+t e_{1,2}, I_3, I_3
\right],
$$
$$
u_{\alpha_3}(t) = \left[
I_3, I_3+t e_{2,1}, I_3
\right],\qquad
u_{-\alpha_3}(t) = \left[
I_3, I_3+t e_{1,2}, I_3
\right],
$$
\begin{align*}
u_{\alpha_4}(t): \left( m_i \right)_{i=1, 2, 3} &\mapsto
\left( m_i + t\cdot
\begin{bmatrix}
0 & -m_{i+2}^{2,3} & 0 \\
0 & 0 & 0 \\
0 & m_{i+2}^{2,1} & 0 \\
\end{bmatrix}
\right)_{i=1, 2, 3}
\end{align*}
\begin{align*}
u_{-\alpha_4}(t): \left( m_i \right)_{i=1, 2, 3} &\mapsto
\left( m_i + t\cdot
\begin{bmatrix}
0 & 0 & 0 \\
m_{i+1}^{3,2} & 0 & -m_{i+1}^{1,2} \\
0 & 0 & 0 \\
\end{bmatrix}
\right)_{i=1, 2, 3}
\end{align*}
$$
u_{\alpha_5}(t) = \left[
I_3, I_3, I_3+t e_{2,1}
\right],\qquad
u_{-\alpha_5}(t) = \left[
I_3, I_3, I_3+t e_{1,2}
\right],
$$
$$
u_{\alpha_6}(t) = \left[
I_3, I_3, I_3+t e_{1,3}
\right],\qquad
u_{-\alpha_6}(t) = \left[
I_3, I_3, I_3+t e_{3,1}
\right].
$$
Here, in the description of $u_{\pm\alpha_4}(t)$, the $m_i$'s should
be counted cyclicly mod $3$, e.g.\ $m_{i+2}=m_1$ for $i=2$.

The associated roots $\alpha_i$, $1\leq i\leq 6$, of these root
subgroups form a simple system in the root system $\Phi(E_6)$ of
$3E_6(\C)$ (our numbering agrees with
\cite[Planche~V]{bourbaki-lie4:6}). For this simple system, the
highest weight of $\K$ is $\lambda_1$.
Furthermore the root subgroups
$u_{\pm\alpha_i}$, $1\leq i\leq 6$, have been chosen so that they
satisfy the conditions in \cite[Prop.~8.1.1(i) and Lem.~8.1.4(i)]{spr:linalggrp}, i.e., they form
part of a {\it realization} (\cite[p.~133]{spr:linalggrp}) of
$\Phi(E_6)$ in $3E_6(\C)$. For $\alpha=\pm\alpha_i$, $1\leq
i\leq 6$, and $t\in\C^\times$, we define the elements
$$
n_\alpha(t)=u_\alpha(t) u_{-\alpha}(-1/t) u_\alpha(t),\qquad
h_\alpha(t)=n_\alpha(t) n_\alpha(1)^{-1}.
$$
Then the maximal torus consists of the elements
$h(t_1, t_2, t_3, t_4, t_5, t_6)=\prod_{i=1}^6 h_{\alpha_i}(t_i)$ and
the normalizer $N(H)$ of the maximal torus is generated by $H$
and the elements $n_i=n_{\alpha_i}(1)$, $1\leq i\leq 6$.
It should be noted that this notation differs from the one used in
\cite{CoWa}. More precisely, the element
$h(\alpha,\beta,\gamma,\delta,\varepsilon,\zeta)$ in \cite[p.~109]{CoWa} equals
$h(\delta,\alpha^{-1},\gamma^{-1},\beta,\varepsilon^{-1},\zeta)$ in our
notation, and the elements $n_1$, $n_2$, $n_3$, $n_4$, $n_5$ and $n_6$
in \cite[p.~109]{CoWa} equal
$n_1 h_{\alpha_1}(-1) h_{\alpha_3}(-1)$, $n_2 h(-1,1,1,-1,1,-1)$,
$n_3 h_{\alpha_1}(-1)$, $n_4$, $n_5 h_{\alpha_6}(-1)$ and
$n_6 h_{\alpha_5}(-1) h_{\alpha_6}(-1)$ respectively in our notation.

From the description of the root system of
type $E_6$ in \cite[Planche~V]{bourbaki-lie4:6} we see that the
center $Z$ of $3E_6(\C)$ is cyclic of order $3$ and is
generated by the element $z=\left[I_3,\omega^2 I_3,\omega I_3\right]$. We
also consider the element $a=\left[\omega I_3,I_3,I_3\right]$. A
straightforward computation shows that the roots of the centralizer
$C_{3E_6(\C)}(a)$ are
$$
\{\pm\alpha_1,\pm\alpha_2,\pm\alpha_3,\pm\alpha_5,\pm\alpha_6,
\pm\widetilde{\alpha}, \pm(\alpha_1+\alpha_3),\pm(\alpha_5+\alpha_6),
\pm(\alpha_2-\widetilde{\alpha})\},
$$
where $\widetilde{\alpha}$ is the longest root. The Dynkin diagram for
this centralizer is the same as the extended Dynkin diagram for $E_6$
with the node $\alpha_4$ removed. In particular it has type
$A_2A_2A_2$ and a simple system of roots is given by
$\{\alpha_1,\alpha_3,\alpha_5,\alpha_6,\alpha_2,
-\widetilde{\alpha}\}$. Since $3E_6(\C)$ is simply connected,
Theorem~\ref{alggrpfacts}\eqref{simply-con} implies that the centralizer
$C_{3E_6(\C)}(a)$ is connected, and thus it is generated by
the maximal torus $H$ and the root subgroups $u_{\pm\alpha}(t)$ where
$\alpha$ runs through the simple roots
$\{\alpha_1, \alpha_3, \alpha_5, \alpha_6, \alpha_2, -\widetilde{\alpha}\}.$
Now note that
$u_{\widetilde{\alpha}}(t) = \left[ I_3+t e_{3,1}, I_3, I_3\right]$ and
$u_{-\widetilde{\alpha}}(t) = \left[ I_3+t e_{1,3}, I_3, I_3\right]$
are root subgroups with associated roots $\widetilde{\alpha}$ and
$-\widetilde{\alpha}$ respectively. Since these along with $H$ and the root
subgroups $u_{\pm\alpha_1}$, $u_{\pm\alpha_2}$, $u_{\pm\alpha_3}$,
$u_{\pm\alpha_5}$ and $u_{\pm\alpha_6}$ generate the subgroup
$\SL_3(\C)^3/C_3$ of $3E_6(\C)$ from above, we
conclude that $C_{3E_6(\C)}(a)=\SL_3(\C)^3/C_3$.

To describe the conjugacy classes of elementary abelian $3$-subgroups
we introduce the following elements in $\SL_3(\C)^3/C_3 \subseteq 3E_6(\C)$:
$$
x_1=\left[I_3,\beta,\beta\right],\qquad
x_2=\left[\beta,\beta,\beta\right],\qquad
y_1=\left[I_3,\gamma,\gamma^2\right],\qquad
y_2=\left[\gamma,\gamma,\gamma\right].
$$
We also need the following elements in $N(H)$:
\begin{equation*}
\begin{split}
s_1 =\; & n_1 n_3 n_4 n_2 n_5 n_4 n_3 n_1 n_6 n_5 n_4 n_2 n_3 n_4 n_5 n_6,\\
s_2 =\; & n_1 n_2 n_3 n_1 n_4 n_2 n_3 n_1 n_4 n_3 n_5 n_4 n_2 n_3 n_1 n_4
n_3 n_5 n_4 n_2 n_6 n_5 n_4 n_2 n_3 n_1 n_4\cdot\\
& n_3 n_5 n_4 n_2 n_6 n_5 n_4 n_3 n_1.
\end{split}
\end{equation*}
The actions of these elements are as follows:
$$
s_1(m_1,m_2,m_3) = (m_3,m_1,m_2),\qquad
s_2(m_1,m_2,m_3) = (m_3^{\mathtt T},m_2^{\mathtt T},m_1^{\mathtt T}),
$$
where $m_i^{\mathtt T}$ denotes the transpose of $m_i$. Thus these
elements act by conjugation on the subgroup $\SL_3(\C)^3/C_3$ as
follows
$$
\left[g_1,g_2,g_3\right]^{s_1} = \left[g_2, g_3, g_1\right],\qquad
\left[g_1,g_2,g_3\right]^{s_2} = \left[
\left(g_1^{-1}\right)^{\mathtt T},\left(g_3^{-1}\right)^{\mathtt T},
\left(g_2^{-1}\right)^{\mathtt T}\right].
$$

\begin{lemma} \label{E2a-conjugation}
We have
\begin{eqnarray*}
z&=&h(\omega,1,\omega^2,1,\omega,\omega^2),\quad
a=h(\omega,1,\omega^2,1,\omega^2,\omega),\quad
x_1=h(\omega,1,\omega,1,\omega,\omega),\\
x_2&=&h(1,\omega^2,\omega^2,1,\omega^2,1),\quad
y_1=n_1 n_3 n_5 n_6\cdot h_{\alpha_5}(-1),\\
y_2&=&n_1 n_2 n_3 n_4 n_3 n_1 n_5 n_4 n_2 n_3 n_4 n_5 n_6 n_5 n_4 n_2
n_3 n_1 n_4 n_3 n_5 n_4 n_6 n_5\cdot h_{\alpha_2}(-1).
\end{eqnarray*}
Moreover conjugation by the element
$$
n_1 n_4 n_2 n_3 n_1 n_4 n_5 n_4 n_6 n_5 n_4 n_2 n_3 n_1 n_4\cdot
h_{\alpha_2}(-1) h_{\alpha_4}(-1)
$$
acts as follows:
\begin{alignat*}{5}
a &\mapsto x_2,\qquad & x_2 &\mapsto a,\qquad & y_1 &\mapsto s_1,\qquad & y_2
&\mapsto y_2^2,\qquad & x_2 x_1^{-1} &\mapsto
h_{\alpha_4}(\omega)=\left[\tau_2,\tau_2,\tau_2\right].
\end{alignat*}
\end{lemma}

\begin{proof}
Both parts of the lemma may be checked by direct
computation. The second part also follows from the
first by using the following relations in $N(H)$: The element $n_i$
has image $s_{\alpha_i}$ in $W$ (\cite[Lem.~8.1.4(i)]{spr:linalggrp}),
we have $n_i^2=h_{\alpha_i}(-1)$ (\cite[Lem.~8.1.4(ii)]{spr:linalggrp}) and
$$
n_i n_j n_i\ldots = n_j n_i n_j\ldots
$$
for $1\leq i,j\leq 6$, where the number of factors on both sides
equals the order of $s_{\alpha_i} s_{\alpha_j}$ in $W$
(\cite[Prop.~9.3.2]{spr:linalggrp}).
\end{proof}

\begin{notation} \label{classes:3E6}
For our calculations, we need some information on the conjugacy
classes of elements of order $3$ in $3E_6(\C)$. These are given in
\cite[Table~2]{CoWa}: There are $7$ such conjugacy classes, which we
label $\mathbf{3A}$, $\mathbf{3B}$, $\mathbf{3B'}$, $\mathbf{3C}$,
$\mathbf{3D}$, $\mathbf{3E}$ and $\mathbf{3E'}$, where $\mathbf{3B'}$ and
$\mathbf{3E'}$ denote the inverses of the classes $\mathbf{3B}$ and
$\mathbf{3E}$. This notation is almost identical to the notation in
\cite[Table~2]{CoWa}, but differs from \cite[Table~VI]{griess}.
We will need the following, which follows
quickly from \cite[Table~2]{CoWa} using the action of $W$ on $H$: We
have $z\in\mathbf{3E}$, $a, x_2, y_2\in\mathbf{3C}$, $x_1,
y_1\in\mathbf{3D}$ and $x_2 x_1^{-1}\in\mathbf{3A}$. Multiplication by
$z$ acts as follows on the conjugacy classes:
\begin{alignat*}{7}
\mathbf{3A} &\mapsto \mathbf{3B},\; &
\mathbf{3B} &\mapsto \mathbf{3B'},\; &
\mathbf{3B'} &\mapsto \mathbf{3A},\; &
\mathbf{3C} &\mapsto \mathbf{3C},\; &
\mathbf{3D} &\mapsto \mathbf{3D},\; &
\mathbf{3E} &\mapsto \mathbf{3E'},\; &
\mathbf{3E'} &\mapsto \mathbf{1},
\end{alignat*}
where $\mathbf{1}$ denotes the conjugacy class consisting of the
identity element.
\end{notation}

\begin{thm}\label{nontorals:3E6-3}
The conjugacy classes of non-toral elementary abelian $3$-subgroups of
$3E_6(\C)$ are given by the following table.
$$\begin{array}{|c|c|c|c|c|}
\hline
\mbox{rank} & \mbox{name} & \mbox{ordered basis} &
\mbox{$3E_6(\C)$-class distribution} &
C_{3E_6(\C)}(E) \\
\hline
\hline
3 & E^3_{3E_6} & \{a,x_2,y_2\} & \mathbf{3C^{26}} & E^4_{3E_6} \\
\hline
\hline
4 & E^4_{3E_6} & \{z,a,x_2,y_2\} &
\mathbf{3C^{78}3E^{1}{3E'}^{1}} & E^4_{3E_6} \\
\hline
\end{array}$$
Their Weyl groups with respect to the given ordered bases are as follows:
$$
W(E^3_{3E_6}) =  \SL_3(\F_3),\:\:
W(E^4_{3E_6}) = \left[
\begin{array}{c|c}
1 &
\begin{array}{ccc}
* & * & *
\end{array} \\
\hline
\begin{array}{c}
0 \\
0 \\
0
\end{array}
& \SL_3(\F_3)
\end{array}
\right]
$$
\end{thm}

\begin{proof}
{\em Non-toral subgroups:}
By \cite[Thm.~11.13]{griess}, there
are two conjugacy classes of non-toral elementary abelian
$3$-subgroups in $3E_6(\C)$, one non-maximal of rank $3$ and one
maximal of rank $4$. We may concretely realize these as follows.
Consider the subgroups
$$
E^3_{3E_6} = \left<a,x_2,y_2\right> \mbox{ and }
E^4_{3E_6} =\left<z,a,x_2,y_2\right>,
$$
which are readily seen to be elementary abelian $3$-subgroups of rank $3$ and
$4$ respectively. Both subgroups are contained in
$C_{3E_6(\C)}(a)=\SL_3(\C)^3/C_3$, and since
$\beta,\gamma\in\SL_3(\C)$ do not commute, we see that the
preimages of $E^3_{3E_6}$ and $E^4_{3E_6}$ under the projection
$\SL_3(\C)^3 \rightarrow \SL_3(\C)^3/C_3$ are
non-abelian. Thus by Theorem~\ref{alggrpfacts}\eqref{covering} $E^3_{3E_6}$
and $E^4_{3E_6}$ are non-toral in $\SL_3(\C)^3/C_3=
C_{3E_6(\C)}(a)$ and hence also non-toral in
$3E_6(\C)$ by Theorem~\ref{alggrpfacts}\eqref{centr-toral}.
Hence these two subgroups represent the two conjugacy
classes of non-toral elementary abelian $3$-subgroups in $3E_6(\C)$.

{\em Lower bounds for Weyl groups:}
By \cite[Thm.~7.4]{griess} there is a unique non-toral
elementary abelian $3$-subgroup $E$ of $F_4(\C)$ of rank $3$ whose
Weyl group in $F_4(\C)$ equals $\SL_3(\F_3)$. Since we have an
inclusion $F_4(\C)\subseteq 3E_6(\C)$ this subgroup may also be
considered as a subgroup of $3E_6(\C)$ and its Weyl group in
$3E_6(\C)$ must contain $\SL_3(\F_3)$. In particular it has order
divisible by $13$ and since $13\nmid |W(E_6)|$, we conclude by
Theorem~\ref{alggrpfacts}\eqref{weyl-section} that $E$ is non-toral
in $3E_6(\C)$ as well. Thus by the above $E$ must be conjugate to
$E^3_{3E_6}$, and hence $W(E^3_{3E_6})$ contains
$\SL_3(\F_3)$.
From this we immediately see that $W(E^4_{3E_6})$
contains the group $1\times \SL_3(\F_3)$.

Note that the element $\left[I_3,\beta,\beta^2\right]$ commutes with $z$,
$a$ and $x_2$ and conjugates $y_2$ to $y_2 z$. Thus it normalizes
$E^4_{3E_6}$ and produces the element $I_4+e_{1,4}$ in
$W(E^4_{3E_6})$. As a result we see that $W(E^4_{3E_6})$ contains the group
$$
\left[
\begin{array}{c|c}
1 &
\begin{array}{ccc}
* & * & *
\end{array} \\
\hline
\begin{array}{c}
0 \\
0 \\
0
\end{array}
& \SL_3(\F_3)
\end{array}
\right]
$$

{\em Class distributions:} Since $a\in\mathbf{3C}$ by
\ref{classes:3E6} and $W(E^3_{3E_6})$ contains $\SL_3(\F_3)$ which
acts transitively on $E^3_{3E_6}-\{1\}$, the class
distribution of $E^3_{3E_6}$ follows immediately. Using this and
the information given in \ref{classes:3E6} about multiplication by
$z$, the class distribution of $E^4_{3E_6}$ follows.

{\em Centralizers:}
Since $C_{3E_6(\C)}(a)=\SL_3(\C)^3/C_3$ we directly get
\begin{xalignat*}{1}
C_{3E_6(\C)}(a,x_2) &= C_{\SL_3(\C)^3/C_3}(x_2)=
\left<y_2,\left(\T_2\times\T_2\times \T_2\right)/C_3\right>,\\
C_{3E_6(\C)}(a,x_2,y_2) &= \left<x_2,y_2,\left(\left<\omega I_3\right>
\times \left<\omega I_3\right> \times \left<\omega
I_3\right>\right)/C_3\right> = E^4_{3E_6},
\end{xalignat*}
proving that $C_{3E_6(\C)}(E^3_{3E_6}) = C_{3E_6(\C)}(E^4_{3E_6}) =
E^4_{3E_6}$.

{\em Exact Weyl groups:}
From the lower bounds above and the fact that $z$ is central we get
$\SL_3(\F_3)\subseteq W(E^3_{3E_6}) \subseteq \GL_3(\F_3)$ and
$$
\left[
\begin{array}{c|c}
1 &
\begin{array}{ccc}
* & * & *
\end{array} \\
\hline
\begin{array}{c}
0 \\
0 \\
0
\end{array}
& \SL_3(\F_3)
\end{array}
\right]
\subseteq W(E^4_{3E_6})\subseteq
\left[
\begin{array}{c|c}
1 &
\begin{array}{ccc}
* & * & *
\end{array} \\
\hline
\begin{array}{c}
0 \\
0 \\
0
\end{array}
& \GL_3(\F_3)
\end{array}
\right]
$$
As $C_{3E_6(\C)}(a,x_2) = \left<y_2,\left(\T_2\times\T_2\times
\T_2\right)/C_3\right>$, we see that no element in
$C_{3E_6(\C)}(a,x_2)$ conjugates $y_2$ to $y_2^{-1}$. Hence
$\diag(1,1,2)\notin W(E^3_{3E_6})$ and
$\diag(1,1,1,2)\notin W(E^4_{3E_6})$ which shows that the Weyl groups are
the ones given in the theorem.
\end{proof}

We now turn to the group $E_6(\C)$. As above, let $Z$ be the center
of $3E_6(\C)$ and let $\pi:3E_6(\C)\rightarrow E_6(\C)= 3E_6(\C)/Z$
denote the quotient homomorphism. For $g\in 3E_6(\C)$ we write
$\overline{g}$ instead of $\pi(g)$ and similarly we let
$\overline{S}=\pi(S)$ for a subset $S\subseteq 3E_6(\C)$.

\begin{lemma} \label{weyl:rank2}
Let $E$ be a rank two non-toral elementary abelian $3$-subgroup of
$E_6(\C)$. Then the Weyl group $W(E)$ is a subgroup of
$\SL_2(\F_3)$.
\end{lemma}

\begin{proof}
Let $\{\overline{g_1},\overline{g_2}\}$ be an ordered basis of $E$.
By Theorem~\ref{alggrpfacts} parts~\eqref{covering}
and \eqref{simply-con} the subgroup
$\left<g_1,g_2\right>\subseteq 3E_6(\C)$ is non-abelian.
Thus setting
$z'=\left[g_1,g_2\right]\in Z$ we
have $z'\neq 1$. Assume that $\sigma\in W(E)$ is represented by the matrix
$
\begin{bmatrix}
a_{11} & a_{12} \\
a_{21} & a_{22}
\end{bmatrix}
$,
i.e., we have $\sigma(\overline{g_1})=(\overline{g_1})^{a_{11}}
(\overline{g_2})^{a_{21}}$ and
$\sigma(\overline{g_2})=(\overline{g_1})^{a_{12}}
(\overline{g_2})^{a_{22}}$. Since $\sigma$ is given by a conjugation in
$E_6(\C)$, it lifts to a conjugation in $3E_6(\C)$. Now the
relation $\left[g_1,g_2\right]=z'\in Z$ shows that
${(z')}^{a_{11}\cdot a_{22}-a_{12}\cdot a_{21}}=z'$, so
$\sigma\in \SL_2(\F_3)$ since $z'\neq 1$.
\end{proof}

\begin{thm} \label{nontorals:E6-3}
The conjugacy classes of non-toral elementary abelian $3$-subgroups of
$E_6(\C)$ are given by the following table:
{\Small
$$\begin{array}{|c|c|c|c|c|c|}
\hline
\mbox{rank} & \mbox{name} & \mbox{ordered basis} &
\mbox{$3E_6(\C)$-class distribution} &
C_{E_6(\C)}(E) & Z(C_{E_6(\C)}(E))\\
\hline
\hline
2 & E^{2a}_{E_6} & \{\overline{y_1},\overline{x_2}\} &
\mathbf{3C^{18}3D^{6}3E^{1} {3E'}^{1}} & E^{2a}_{E_6}\times
\PSL_3(\C)& E^{2a}_{E_6}\\
\hline
2 & E^{2b}_{E_6} & \{\overline{y_1},\overline{x_1}\} &
\mathbf{3D^{24}3E^{1} {3E'}^{1}} &  E^{2b}_{E_6}\times G_2(\C) &
E^{2b}_{E_6} \\
\hline
\hline
3 & E^{3a}_{E_6} & \{\overline{a},\overline{y_1},\overline{x_2}\} &
\mathbf{3C^{60}3D^{18}3E^{1} {3E'}^{1}} &
E^{3a}_{E_6}\circ_{\left<\overline{a} \right>}
(\T_2:\left<\overline{y_2} \right>) & E^{3a}_{E_6}\\
\hline
3 & E^{3b}_{E_6} & \{\overline{a},\overline{x_2},\overline{y_2}\} &
\mathbf{3C^{78}3E^{1} {3E'}^{1}} & E^{3b}_{E_6}\cdot {(C_3)}^3 & E^{3b}_{E_6}\\
\hline
3 & E^{3c}_{E_6} & \{\overline{a},\overline{y_1},\overline{x_1}\}
& \mathbf{3C^{6}3D^{72}3E^{1} {3E'}^{1}} & E^{3c}_{E_6}
\circ_{\left<\overline{a}\right>} \SL_3(\C) & E^{3c}_{E_6}\\
\hline
3 & E^{3d}_{E_6} & \{\overline{x_2x_1^{-1}},\overline{y_1},
\overline{x_1}\} &
\mathbf{3A^{2}3B^{2}{3B'}^{2}3C^{48}3D^{24}3E^{1} {3E'}^{1}} &
E^{3d}_{E_6} \circ_{\left<\overline{x_2x_1^{-1}} \right>} \GL_2(\C) &
E^{3d}_{E_6} \circ_{\left<\overline{x_2x_1^{-1}} \right>} \T_1\\
\hline
\hline
4 & E^{4a}_{E_6} &
\{\overline{a},\overline{y_2},\overline{y_1},\overline{x_2}\} &
\mathbf{3C^{186}3D^{54}3E^{1} {3E'}^{1}} & E^{4a}_{E_6} & E^{4a}_{E_6}\\
\hline
4 & E^{4b}_{E_6} &
\{\overline{a},\overline{x_2x_1^{-1}},\overline{y_1},\overline{x_1}\}
& \mathbf{3A^{6}3B^{6}{3B'}^{6}3C^{150}3D^{72}3E^{1} {3E'}^{1}} &
E^{4b}_{E_6}\circ_{\left<\overline{a},\overline{x_2x_1^{-1}} \right>} \T_2 &
E^{4b}_{E_6}\circ_{\left<\overline{a},\overline{x_2x_1^{-1}} \right>}
\T_2\\
\hline
\end{array}$$
} 
In  particular ${}_3 Z(C_{E_6(\C)}(E))=E$ for all non-toral elementary
abelian $3$-subgroups $E$ of $E_6(\C)$. (In the table the $3E_6(\C)$-class
distribution of $E\subseteq E_6(\C)$ denotes the class distribution of
$\pi^{-1}(E)\subseteq 3E_6(\C)$.)

The Weyl groups of these subgroups with
respect to the given ordered bases are as follows:
$$
W(E^{2a}_{E_6})=
\begin{bmatrix}
\varepsilon_1 & * \\
0 & \varepsilon_1
\end{bmatrix}
,\:\:
W(E^{2b}_{E_6})=\SL_2(\F_3),\:\:
W(E^{3a}_{E_6})=
\begin{bmatrix}
\varepsilon_1 & * & * \\
0 & \varepsilon_2 & * \\
0 & 0 & \varepsilon_2
\end{bmatrix}
$$
$$
W(E^{3b}_{E_6})= \SL_3(\F_3),\:\:
W(E^{3c}_{E_6})=\left[
\begin{array}{c|c}
\varepsilon &
\begin{array}{cc}
* & *
\end{array} \\
\hline
\begin{array}{c}
0 \\
0
\end{array} &
\SL_2(\F_3)
\end{array}
\right],\:\:
W(E^{3d}_{E_6})= \GL_1(\F_3)\times \SL_2(\F_3)
$$
$$
W(E^{4a}_{E_6}) = \left[
\begin{array}{c|cc}
\GL_2(\F_3) &
\begin{array}{c}
* \\
*
\end{array} &
\begin{array}{c}
* \\
*
\end{array} \\
\hline
\begin{array}{cc}
0 & 0
\end{array} &
\operatorname{det} & * \\
\begin{array}{cc}
0 & 0
\end{array} &
0 & \operatorname{det}
\end{array}
\right],\:\:
W(E^{4b}_{E_6})=\left[
\begin{array}{cc|c}
\varepsilon_1 & * &
\begin{array}{cc}
* & *
\end{array} \\
0 & \varepsilon_2 &
\begin{array}{cc}
0 & 0
\end{array} \\
\hline
\begin{array}{c}
0 \\
0
\end{array} &
\begin{array}{c}
0 \\
0
\end{array} &
\SL_2(\F_3)
\end{array}
\right]
$$
where $\operatorname{det}$ denotes the determinant of the matrix from
$\GL_2(\F_3)$ in the description of $W(E^{4a}_{E_6})$.
\end{thm}

\begin{proof}
{\em Maximal non-toral subgroups:} By \cite[Thm.~11.14]{griess}, there
are two conjugacy classes of maximal non-toral elementary abelian
$3$-subgroups in $E_6(\C)$, both of rank $4$. We may
concretely realize these as follows.
Consider the subgroups
$$
E_a = \left<z,a,y_1,y_2,x_2\right> \mbox{ and }
E_b = \left<z,a,x_2x_1^{-1},y_1,x_1\right>
$$
of $C_{3E_6(\C)}(a)=\SL_3(\C)^3/C_3$. Since the commutator
subgroup of both of these equals $Z$, we see that
$E^{4a}_{E_6}=\pi(E_a)$ and $E^{4b}_{E_6}=\pi(E_b)$ are elementary abelian
$3$-subgroups of rank $4$ in $E_6(\C)$. It follows from
Theorem~\ref{alggrpfacts}\eqref{covering} that both $E^{4a}_{E_6}$ and
$E^{4b}_{E_6}$ are non-toral in $E_6(\C)$.
We will see below that their class distributions are the ones given in
the table. From this it follows that they are not conjugate and thus
represent the two conjugacy classes of maximal elementary abelian
$3$-subgroups in $E_6(\C)$.

{\em Lower bounds for Weyl groups of maximal non-toral subgroups:}
We now find lower bounds for the Weyl groups of the maximal
non-toral elementary abelian $3$-subgroups by conjugating with
elements coming from the centralizer
$C_{3E_6(\C)}(a)=\SL_3(\C)^3/C_3$ and the normalizer $N(H)$ of the
maximal torus.

The elements $\overline{[\beta^2,I_3,I_3]}$,
$\overline{[I_3,\tau_1,\tau_1^2]}$, $\overline{s_1}$ and
$\overline{s_2}$ normalize $E_{E_6}^{4a}$ and conjugation by these
elements induces the automorphisms on $E^{4a}_{E_6}$ given by the matrices
$I_4+e_{1,2}$, $I_4+e_{3,4}$, $I_4+e_{2,3}$ and $\diag(2,1,2,2)$.
Moreover, by Lemma~\ref{E2a-conjugation} we may
conjugate the ordered basis of $E^{4a}_{E_6}$ into the ordered
basis
$\{\overline{x_2},\overline{y_2^2},\overline{s_1},\overline{a}\}.$
Noting that the element $\overline{[\tau_1,\tau_1,\tau_1]}$
commutes with $\overline{y_2}$, $\overline{s_1}$ and
$\overline{a}$ and conjugates $\overline{x_2}$ into $\overline{x_2
y_2}$, we see that $W(E^{4a}_{E_6})$ contains the element $I_4+e_{2,1}$.
The above matrices are easily seen to generate the group
$$
W'(E^{4a}_{E_6})=\left[
\begin{array}{c|cc}
\GL_2(\F_3) &
\begin{array}{c}
* \\
*
\end{array} &
\begin{array}{c}
* \\
*
\end{array} \\
\hline
\begin{array}{cc}
0 & 0
\end{array} &
\operatorname{det} & * \\
\begin{array}{cc}
0 & 0
\end{array} &
0 & \operatorname{det}
\end{array}
\right]
$$
and thus $W^{4a}_{E_6}$ contains this group.

Now consider $E^{4b}_{E_6}$ and let $\sigma=-(2,3)\in \SL_3(\C)$.
We then see that the elements $\overline{[I_3,\tau_1,\tau_1^2]}$,
$\overline{[I_3,\tau_2\beta,\tau_2^2]}$, $\overline{[\sigma,I_3,I_3]}$,
$\overline{[\gamma,I_3,I_3]}$, $\overline{[I_3,\beta^2,I_3]}$ and
$\overline{s_2}$ normalize $E^{4b}_{E_6}$, and conjugation by these
elements induces the automorphisms on $E^{4b}_{E_6}$ given by the matrices
$I_4+e_{3,4}$, $I_4+e_{4,3}$, $\diag(1,2,1,1)$, $I_4+e_{1,2}$,
$I_4+e_{1,3}$ and $-I_4$. These matrices generate the group
$$
W'(E^{4b}_{E_6})=\left[
\begin{array}{cc|c}
\varepsilon_1 & * &
\begin{array}{cc}
* & *
\end{array} \\
0 & \varepsilon_2 &
\begin{array}{cc}
0 & 0
\end{array} \\
\hline
\begin{array}{c}
0 \\
0
\end{array} &
\begin{array}{c}
0 \\
0
\end{array} &
\SL_2(\F_3)
\end{array}
\right]
$$
and thus $W^{4b}_{E_6}$ contains this group.

{\em Orbit computation:} Any elementary abelian $3$-subgroup of
rank one is toral since $E_6(\C)$ is connected. Since
$E^{4a}_{E_6}$ and $E^{4b}_{E_6}$ are representatives of the
maximal non-toral elementary abelian $3$-subgroups, we may find the
conjugacy classes of non-toral elementary abelian $3$-subgroups of
ranks $2$ and $3$ by studying subgroups of these.

Under the action of $W'(E^{4a}_{E_6})$, the set of rank $2$ subgroups of
$E^{4a}_{E_6}$ has orbit representatives
$$
E_{E_6}^{2a}=\left<\overline{y_1},\overline{x_2}\right>,
\left<\overline{a},\overline{x_2}\right>,
\left<\overline{a},\overline{y_1}\right>
\mbox{and}
\left<\overline{a},\overline{y_2}\right>,
$$
and under the action of $W'(E^{4b}_{E_6})$, the set of rank $2$ subgroups of
$E^{4b}_{E_6}$ has orbit representatives
$$
E_{E_6}^{2a}=\left<\overline{y_1},\overline{x_2}\right>,
E_{E_6}^{2b}=\left<\overline{y_1},\overline{x_1}\right>,
\left<\overline{a},\overline{x_2}\right>,
\left<\overline{a},\overline{y_1}\right>,
\left<\overline{a},\overline{x_2x_1^{-1}}\right>
\mbox{and}
\left<\overline{x_2x_1^{-1}},\overline{x_1}\right>.
$$
Similarly we find that under the action of $W'(E^{4a}_{E_6})$, the set of rank
$3$ subgroups of $E^{4a}_{E_6}$ has orbit representatives
$$
E_{E_6}^{3a}=\left<\overline{a},\overline{y_1},\overline{x_2}\right>,
E_{E_6}^{3b}=\left<\overline{a},\overline{y_2},\overline{x_2}\right>
\mbox{and}
\left<\overline{a},\overline{y_1},\overline{y_2}\right>,
$$
and that under the action of $W'(E^{4b}_{E_6})$, the set of rank $3$ subgroups
of $E^{4b}_{E_6}$ has orbit representatives
$$
E_{E_6}^{3a}=\left<\overline{a},\overline{y_1},\overline{x_2}\right>,
E_{E_6}^{3c}=\left<\overline{a},\overline{y_1},\overline{x_1}\right>,
E_{E_6}^{3d}=\left<\overline{x_2x_1^{-1}},\overline{y_1},\overline{x_1}\right>
\mbox{and}
\left<\overline{a},\overline{x_2x_1^{-1}},\overline{x_1}\right>.
$$

{\em Other non-toral subgroups:} The subgroups
$\left<\overline{a},\overline{x_2}\right>$,
$\left<\overline{a},\overline{x_2 x_1^{-1}}\right>$,
$\left<\overline{x_2 x_1^{-1}},\overline{x_1}\right>$ and
$\left<\overline{a},\overline{x_2 x_1^{-1}},\overline{x_1}\right>$
are visibly toral. Since the elements $\beta$ and $\gamma$ are
conjugate in $\SL_3(\C)$, the subgroup
$\left<\overline{a},\overline{y_1},\overline{y_2}\right>$ is
conjugate to the subgroup
$\left<\overline{a},\overline{\left[I_3,\beta,\beta^2\right]},
\overline{x_2}\right>$ which is obviously toral. Thus the subgroups
$\left<\overline{a},\overline{y_1},\overline{y_2}\right>$,
$\left<\overline{a},\overline{y_1}\right>$ and
$\left<\overline{a},\overline{y_2}\right>$ are also toral. Using
the fact that $\left[y_1,x_1\right]=\left[y_1,x_2\right]=z$ we see
from Theorem~\ref{alggrpfacts}\eqref{covering} that both
$E^{2a}_{E_6}$ and $E^{2b}_{E_6}$ are non-toral in $E_6(\C)$. Since
the subgroups $E^{3a}_{E_6}$, $E^{3c}_{E_6}$ and $E^{3d}_{E_6}$ all
contain either $E^{2a}_{E_6}$ or $E^{2b}_{E_6}$ they are also
non-toral. Using Theorem~\ref{alggrpfacts}\eqref{covering} we see
that the subgroup $E^{3b}_{E_6}$ is non-toral in $E_6(\C)$, since
$\pi^{-1}(E^{3b}_{E_6})=E^{4}_{3E_6}$ is non-toral in
$3E_6(\C)$ by Theorem~\ref{nontorals:3E6-3}.

{\em Class distributions:} Using \ref{classes:3E6} and the actions
of the groups $W'(E^{4a}_{E_6})$ and $W'(E^{4b}_{E_6})$ it is not
hard to verify the class distributions in the table. As an example
consider the subgroup $E^{4b}_{E_6}$. From the action of
$W'(E^{4b}_{E_6})$ we see that $E^{4b}_{E_6}-\{1\}$
contains $2$ elements conjugate to $\overline{a}$, $6$ elements
conjugate to $\overline{x_2x_1^{-1}}$, $24$ elements conjugate to
$\overline{x_1}$ and $48$ elements conjugate to $\overline{x_2}$.
Thus by \ref{classes:3E6}, the set
$\pi^{-1}(E^{4b}_{E_6}-\{1\})$ contains $6$ elements from
each of the classes $\mathbf{3A}$, $\mathbf{3B}$ and
$\mathbf{3B'}$, $3\cdot(2+48)=150$ elements from the class
$\mathbf{3C}$ and $3\cdot 24=72$ elements from the class
$\mathbf{3D}$. Including the elements $z$ and $z^2$ from the
classes $\mathbf{3E}$ and $\mathbf{3E'}$ respectively, we get the
class distribution of $\pi^{-1}(E^{4b}_{E_6})-\{1\}$ given
in the table. Similar computations give the remaining entries in
the table. Since these distributions are different we see that the
subgroups in the table are not conjugate and thus they provide a set
of representatives for the conjugacy classes of non-toral
elementary abelian $3$-subgroups of $E_6(\C)$.

{\em Lower bounds for other Weyl groups:} We now show that the
other matrix groups in the theorem are all lower bounds for the
remaining Weyl groups. To do this consider one of the non-maximal
subgroups $E$ from the table. We then have $E\subseteq E_{E_6}^{4a}$
or $E\subseteq E_{E_6}^{4b}$, and we get a lower bound on $W(E)$
by considering the action on $E$ of the subgroup of
$W'(E_{E_6}^{4a})$ or $W'(E_{E_6}^{4b})$ fixing $E$. As an example
we see that $E_{E_6}^{2a}\subseteq E_{E_6}^{4a}$ and that the
stabilizer of $E_{E_6}^{2a}$ inside $W'(E_{E_6}^{4a})$ is
$$
\left[
\begin{array}{c|cc}
\GL_2(\F_3) &
\begin{array}{c}
0 \\
0
\end{array} &
\begin{array}{c}
0 \\
0
\end{array} \\
\hline
\begin{array}{cc}
0 & 0
\end{array} &
\operatorname{det} & x \\
\begin{array}{cc}
0 & 0
\end{array} &
0 & \operatorname{det}
\end{array}
\right]
$$
where $\operatorname{det}$ is the determinant of the matrix from
$\GL_2(\F_3)$. The action of
such a matrix on $E_{E_6}^{2a}$ is given by
\begin{alignat*}{2}
\overline{y_1} &\mapsto (\overline{y_1})^{\operatorname{det}},\qquad &
\overline{x_2} &\mapsto (\overline{y_1})^{x} (\overline{x_2})^{\operatorname{det}}.
\end{alignat*}
Thus $W(E_{E_6}^{2a})$ contains the group
$$
W'(E^{2a}_{E_6})=
\begin{bmatrix}
\varepsilon_1 & * \\
0 & \varepsilon_1
\end{bmatrix}
$$
as claimed. Similar computations show that for the subgroups
$E=E_{E_6}^{2b}$, $E_{E_6}^{3a}$, $E_{E_6}^{3c}$ and $E_{E_6}^{3d}$, the group
$W'(E)$ occurring in the theorem is a lower bound for the Weyl group $W(E)$.

For the subgroup
$E_{E_6}^{3b}=\left<\overline{a},\overline{x_2},\overline{y_2}\right>$
we know the structure of $W(\pi^{-1}(E_{E_6}^{3b}))=W(E_{3E_6}^4)$ by
Theorem~\ref{nontorals:3E6-3}. From this we immediately get
$W(E_{E_6}^{3b})=\SL_3(\F_3)$.

{\em Exact Weyl groups:} We now prove that the lower bounds on the
Weyl groups established above are in fact equalities.
By Lemma~\ref{weyl:rank2} the Weyl groups $W(E_{E_6}^{2a})$ and
$W(E_{E_6}^{2b})$ are subgroups of $\SL_2(\F_3)$. From this we see
that $W(E_{E_6}^{2b})=\SL_2(\F_3)$ and that $W(E_{E_6}^{2a})$ is
equal to either $W'(E_{E_6}^{2a})$ or $\SL_2(\F_3)$, since these
are the only subgroups of $\SL_2(\F_3)$ containing
$W'(E_{E_6}^{2a})$. We have
$E_{E_6}^{2a}=\left<\overline{y_1},\overline{x_2}\right>$, and by
\ref{classes:3E6} the elements $\overline{y_1}$ and
$\overline{x_2}$ are not conjugate in $E_6(\C)$. In particular
$W(E_{E_6}^{2a})$ cannot act transitively on the
non-trivial elements of $E_{E_6}^{2a}$, and we conclude that
$W(E_{E_6}^{2a})=W'(E_{E_6}^{2a})$ is the group from above.

For each of the remaining non-toral subgroups we now show that a strictly
larger Weyl group would contradict the Weyl group results already
established. The subgroups $E = E_{E_6}^{3a}$, $E_{E_6}^{3d}$, and
$E_{E_6}^{4b}$ all contain $E_{E_6}^{2a}$. A direct computation shows
that any proper overgroup of $W'(E)$ in $\GL(E)$ contains an element
which normalizes the subgroup  $E_{E_6}^{2a}$ and induces an
automorphism which does not lie in $W(E_{E_6}^{2a})$. Hence $W(E) =
W'(E)$. If $E =  E_{E_6}^{3c}$ a similar argument, using the subgroup
$E_{E_6}^{2b}$, again shows that $W(E) = W'(E)$. Consider finally
$E = E_{E_6}^{4a}$. Each proper overgroup of $W'(E)$ contains an element which
normalizes one of the subgroups  $E_{E_6}^{2a}$  or $E_{E_6}^{3b}$ and
induces an automorphism on it not contained in its Weyl group. Hence
$W(E) = W'(E)$. This concludes the proof that the Weyl groups listed
in the theorem are the correct ones.

{\em Centralizers:}
Let $\Theta:\SL_3(\C) \longrightarrow \SL_3(\C)^3/C_3\subseteq 3E_6(\C)$
denote the homomorphism given by $\Theta(g)=\left[g,g,g\right]$ for $g\in
\SL_3(\C)$. By Lemma~\ref{E2a-conjugation} the subgroup
$E^{2a}_{E_6}=\left<\overline{x_2},\overline{y_1}\right>$ is
conjugate to the subgroup $\left<\overline{a},\overline{s_1}\right>$.
Since $a^{s_1}=a z^2$ we obtain $C_{E_6(\C)}(\overline{a})=
\overline{\left<s_1,\SL_3(\C)^3/C_3\right>}$, and hence
$$
C_{E_6(\C)}(\overline{a},\overline{s_1})=
\overline{\left<a,s_1,z,\Theta(\SL_3(\C))\right>}=
\overline{\left<a,s_1,z\right> \times \Theta(\SL_3(\C))}=
\left<\overline{a},\overline{s_1}\right> \times \PSL_3(\C),
$$
proving the claims for $E^{2a}_{E_6}$. By a slight abuse of notation,
we let $\overline{g}$ denote the image of
$g\in\SL_3(\C)$ in the quotient $\PSL_3(\C)$. From
Lemma~\ref{E2a-conjugation} we then see that the elements
$\overline{a}$, $\overline{y_2}$ and $\overline{x_2x_1^{-1}}$ in
$C_{E_6(\C)}(E^{2a}_{E_6})$ correspond to the elements
$\overline{\beta}$, $\overline{\gamma^2}$ and $\overline{\tau_2}$
in the $\PSL_3(\C)$-component of
$C_{E_6(\C)}(E^{2a}_{E_6})$. Thus we immediately get
\begin{xalignat*}{2}
C_{E_6(\C)}(E^{3a}_{E_6}) &= E^{2a}_{E_6}\times
C_{\PSL_3(\C)}(\overline{\beta}),\qquad &
C_{E_6(\C)}(E^{3d}_{E_6}) &= E^{2a}_{E_6}\times
C_{\PSL_3(\C)}(\overline{\tau_2}),\\
C_{E_6(\C)}(E^{4a}_{E_6}) &= E^{2a}_{E_6}\times
C_{\PSL_3(\C)}(\overline{\beta},\overline{\gamma^2}), \qquad &
C_{E_6(\C)}(E^{4b}_{E_6}) &= E^{2a}_{E_6}\times
C_{\PSL_3(\C)}(\overline{\beta},\overline{\tau_2}).
\end{xalignat*}
Note that $C_{\PSL_3(\C)}(\overline{\beta})=
\overline{\T_2:\left<\gamma\right>}$, giving
$C_{\PSL_3(\C)}(\overline{\beta},\overline{\gamma^2})=
\left<\overline{\beta},\overline{\gamma}\right>$ and
$C_{\PSL_3(\C)}(\overline{\beta},\overline{\tau_2})=
\overline{\T_2}$. From this the results on $E^{3a}_{E_6}$, $E^{4a}_{E_6}$ and
$E^{4b}_{E_6}$ follow directly. Note also that
$C_{\PSL_3(\C)}(\overline{\tau_2})\cong \GL_2(\C)$ from which we
deduce the claims about $E^{3d}_{E_6}$.

Now consider the subgroup $E^{3b}_{E_6}$. Since $C_{E_6(\C)}(\overline{a})=
\overline{\left<s_1,\SL_3(\C)^3/C_3\right>}$ we get
\begin{xalignat*}{1}
C_{E_6(\C)}(\overline{a},\overline{x_2}) &= \overline{ \left< s_1,y_1,y_2,
\left(\T_2\times \T_2\times \T_2\right)/C_3
\right>}, \\
C_{E_6(\C)}(\overline{a},\overline{x_2},\overline{y_2}) &= \overline{ \left<
s_1,y_1,y_2,\left[I_3,\beta,\beta^2\right],x_2, \left(\left<\omega I_3\right>
\times \left<\omega I_3\right> \times \left<\omega
  I_3\right>\right)/C_3 \right>}
\end{xalignat*}
and thus
$C_{E_6(\C)}(E^{3b}_{E_6})=\left<E^{3b}_{E_6},\overline{s_1},\overline{y_1},
\overline{\left[I_3,\beta,\beta^2\right]}\right>$. It is now easy to
check that $C_{E_6(\C)}(E^{3b}_{E_6})$ has the structure
$E^{3b}_{E_6}\cdot {(C_3)}^3$ and that $Z(C_{E_6(\C)}(E^{3b}_{E_6}))=E^{3b}_{E_6}$.
For the subgroup $E^{3c}_{E_6}$ we obtain
\begin{xalignat*}{1}
C_{E_6(\C)}(\overline{a},\overline{x_1}) &= \overline{ \left<
y_1,\left(\SL_3(\C)\times \T_2\times \T_2\right)/C_3
\right>}, \\
C_{E_6(\C)}(\overline{a},\overline{x_1},\overline{y_1}) &= \overline{
\left< y_1,x_1,\left(\SL_3(\C)\times \left<\omega I_3\right>\times
\left<\omega I_3\right>\right)/C_3 \right>}.
\end{xalignat*}
Thus $C_{E_6(\C)}(E^{3c}_{E_6})$ equals the
central product $E^{3c}_{E_6}\circ_{\left<\overline{a}\right>}
\SL_3(\C)$ and we obtain the claims about $E^{3c}_{E_6}$.

Finally consider the subgroup
$E^{2b}_{E_6}=\left<\overline{y_1},\overline{x_1}\right>$. If $g\in
\pi^{-1}(C_{E_6(\C)}(E^{2b}_{E_6}))$ then $\left[g,y_1\right],
\left[g,x_2\right]\in Z$, and since $\left[y_1,x_2\right]=z$ it
follows that $g\in \pi^{-1}(E^{2b}_{E_6}) \circ_Z
C_{3E_6(\C)}(\pi^{-1}(E^{2b}_{E_6}))$. Thus we have
$C_{E_6(\C)}(E^{2b}_{E_6}) = E^{2b}_{E_6} \times
\overline{C_{3E_6(\C)}(\pi^{-1}(E^{2b}_{E_6}))}$.

A direct computation shows that $C_{3E_6(\C)}(x_1)$ has type
$T_2D_4$ and a system of simple roots of the centralizer is
given by $\{\alpha_1+\alpha_3+\alpha_4, \alpha_2,
\alpha_4+\alpha_5+\alpha_6, \alpha_3+\alpha_4+\alpha_5\}$. From
this we see that the $2$-dimensional torus $\T_2$ consists of the
elements $h(\alpha,1,\gamma,1,\alpha,\gamma)$ where
$\alpha,\gamma\in\C^\times$. Moreover we see that
$C_{3E_6(\C)}(x_1)=\T_2\circ_C \Spin(8,\C)$, where
$C=Z(\Spin(8,\C))=C_2\times C_2$ consists of the elements
$h(\alpha,1,\gamma,1,\alpha,\gamma)$, $\alpha,\gamma=\pm 1$.

Let $\sigma$ denote the automorphism of $C_{3E_6(\C)}(x_1)$ given by
conjugation with $y_1$. A direct check shows that the map from
$C$ to $C$ given by $x\mapsto x^{-1} x^\sigma$ is surjective. It then
follows that 
$$
C_{3E_6(\C)}(\pi^{-1}(E^{2b}_{E_6})) =
\left(\T_2\circ_C \Spin(8,\C) \right)^\sigma =
\T_2^\sigma \circ_{C^\sigma} \Spin(8,\C)^\sigma.
$$
We have $\T_2^\sigma=\left<z\right>$, so
$C_{3E_6(\C)}(\pi^{-1}(E^{2b}_{E_6})) = \left<z\right> \times
\Spin(8,\C)^\sigma$. Using the class distribution of
$\pi^{-1}(E^{2b}_{E_6})$ found above together with
\cite[Table~2]{CoWa} and Theorem~\ref{alggrpfacts}\eqref{dim-centr} we
find
$$
\dim C_{3E_6(\C)}(\pi^{-1}(E^{2b}_{E_6})) =
\frac{1}{3^3}\cdot \left(3\cdot 78 + 24\cdot
(30+24\omega+24\omega^2)\right) = 14.
$$
Thus $\Spin(8,\C)^\sigma$ has dimension $14$ and since
$Z(\Spin(8,\C))^\sigma=1$ we also see that $\Spin(8,\C)^\sigma$ has
rank less than $4$. From this it follows that the identity component
of $\Spin(8,\C)^\sigma$ must have type $G_2$. By
\cite[Thm.~8.1]{steinberg68} $\Spin(8,\C)^\sigma$ is
connected, so $\Spin(8,\C)^\sigma = G_2(\C)$ and hence
$C_{3E_6(\C)}(\pi^{-1}(E^{2b}_{E_6})) = \left<z\right> \times
G_2(\C)$. Combining this with the computation from above we conclude
$C_{E_6(\C)}(E^{2b}_{E_6}) = E^{2b}_{E_6} \times G_2(\C)$.
\end{proof}

\subsection{The group $E_8(\C)$, $p=3$} \label{type:E8}
In this section we consider the elementary abelian $3$-subgroups of
the group $E_8(\C)$. By using \cite[Table~2, p.~214]{bourbaki-lie7:8}
we see that the smallest faithful representation of $E_8(\C)$ is the
adjoint representation, i.e., the representation given by the action of
$E_8(\C)$ on its Lie algebra $\mathfrak{e}_8$, which has dimension
$248$. For our computations, we explicitly construct this
representation on a computer by following the recipe in
\cite[Ch.~4]{carter:simple}. As explained in
\cite[Ch.~4]{carter:simple} there is some ambiguity in choosing a
Chevalley basis of $\mathfrak{e}_8$; we return to this problem below.

Letting $\Phi(E_8)$ denote the root
system of type $E_8$ (we use the notation of
\cite[Planche~VII]{bourbaki-lie4:6}), we have in particular a maximal
torus $H$ generated by the elements $h_{\alpha_i}(t)$, $1\leq i\leq
8$, $t\in\C^\times$ (\cite[p.~92, p.~97]{carter:simple})
and root subgroups $u_\alpha(t)$, $\alpha\in\Phi(E_8)$,
$t\in\C$. The normalizer $N(H)$ of the maximal torus, is
generated by $H$ and the elements $n_i=n_{\alpha_i}$, $1\leq i\leq 8$
(\cite[p.~93, p.~101]{carter:simple}). We let
$$
h(t_1,t_2,t_3,t_4,t_5,t_6,t_7,t_8) = \prod_{i=1}^8 h_{\alpha_i}(t_i).
$$

Note that by \cite[p.~100 and Lem.~6.4.4]{carter:simple} the
root subgroups $u_{\alpha}$ form a {\it realization}
(\cite[p.~133]{spr:linalggrp}) of $\Phi(E_8)$ in $E_8(\C)$. In
particular we have the following relations: The element $n_i$
has image $s_{\alpha_i}$ in $W=W(E_8)$ (\cite[Lem.~8.1.4(i)]{spr:linalggrp}),
we have $n_i^2=h_{\alpha_i}(-1)$ (\cite[Lem.~8.1.4(ii)]{spr:linalggrp}) and
$$
n_i n_j n_i\ldots = n_j n_i n_j\ldots
$$
for $1\leq i,j\leq 8$, where the number of factors on both sides
equals the order of $s_{\alpha_i} s_{\alpha_j}$ in $W$
(\cite[Prop.~9.3.2]{spr:linalggrp}).

Now let $\overline{a}=h_{\alpha_1}(\omega) h_{\alpha_2}(\omega)
h_{\alpha_3}(\omega^2)\in E_8(\C)$.
Direct computation shows that for any root $\alpha\in\Phi(E_8)$ we have
$\alpha(\overline{a})=\omega^{2 (\lambda_2,\alpha)}$.
From this we see that the Dynkin diagram of the centralizer
$C_{E_8(\C)}(\overline{a})$ is the
same as the extended Dynkin diagram of $E_8$ with the node $\alpha_2$
removed. Thus it has type $A_8$ and a simple system of roots is given by
$$
\{\alpha_1, \alpha_3, \alpha_4, \alpha_5, \alpha_6, \alpha_7,
\alpha_8, -\widetilde{\alpha}\},
$$
where $\widetilde{\alpha}$ is the longest root. As in
\cite[Planche~I]{bourbaki-lie4:6} we identify $\Phi(A_8)$ with the
set of elements in $\R^9$ of the form $e_i-e_j$ with
$i\neq j$ and $1\leq i,j\leq 9$, where $e_i$ denotes the $i$th
canonical basis vector in $\R^9$. We now consider $\SL_9(\C)$,
which is the simply connected group of type $A_8$ over
$\C$. Given a root $\alpha'=e_i-e_j\in \Phi(A_8)$ we let
$u'_{\alpha'}(t)=I_9+t e_{i,j}$ for $t\in\C$.
With respect to the maximal torus consisting of
diagonal matrices, this is a root subgroup of $\SL_9(\C)$
corresponding to the root $\alpha'$. The roots $\alpha'_i=e_i-e_{i+1}$,
$1\leq i\leq 8$, form a simple system in $\Phi(A_8)$.
It now follows that we can choose the Chevalley basis of
$\mathfrak{e}_8$ in such a way that
\begin{alignat*}{4}
u'_{\pm\alpha'_1}(t) & \mapsto u_{\pm\alpha_1}(t), & \qquad
u'_{\pm\alpha'_2}(t) & \mapsto u_{\pm\alpha_3}(t), & \qquad
u'_{\pm\alpha'_3}(t) & \mapsto u_{\pm\alpha_4}(t), & \qquad
u'_{\pm\alpha'_4}(t) & \mapsto u_{\pm\alpha_5}(t), \\
u'_{\pm\alpha'_5}(t) & \mapsto u_{\pm\alpha_6}(t), & \qquad
u'_{\pm\alpha'_6}(t) & \mapsto u_{\pm\alpha_7}(t), & \qquad
u'_{\pm\alpha'_7}(t) & \mapsto u_{\pm\alpha_8}(t), & \qquad
u'_{\pm\alpha'_8}(t) & \mapsto u_{\mp\widetilde{\alpha}}(t)
\end{alignat*}
defines a homomorphism $\SL_9(\C)\rightarrow E_8(\C)$
onto the centralizer $C_{E_8(\C)}(\overline{a})$, and we fix a certain
such choice. It is easy to check
that this homomorphism has kernel $C_3=\left<\omega I_9\right>$
and thus we may make the identification
$C_{E_8(\C)}(\overline{a})= \SL_9(\C)/C_3$. For any
$g\in \SL_9(\C)$ we denote by $\overline{g}$ its image in
$\SL_9(\C)/C_3= C_{E_8(\C)}(\overline{a}) \subseteq
E_8(\C)$. In particular we see
that $a=\eta I_9$ corresponds to the element $\overline{a}$ from above.
Define the following elements in $\SL_9(\C):$
\begin{alignat*}{2}
x_1 =\; & \diag(1,\omega,\omega^2,1,\omega,\omega^2,1,\omega,\omega^2),
& \qquad
x_2 =\; & \diag(1,1,1,\omega,\omega,\omega,\omega^2,\omega^2,\omega^2), \\
x_3 =\; & \diag(1,1,1,1,1,1,\omega,\omega,\omega), & \qquad
y_1 =\; & (1,2,3)(4,5,6)(7,8,9), \\
y_2 =\; & (1,4,7)(2,5,8)(3,6,9). & &
\end{alignat*}
From the explicit homomorphism above we easily find
\begin{alignat*}{2}
\overline{a} =\; & h_{\alpha_1}(\omega) h_{\alpha_2}(\omega)
h_{\alpha_3}(\omega^2),\qquad &
\overline{x_1} =\; & h_{\alpha_1}(\omega) h_{\alpha_5}(\omega)
h_{\alpha_8}(\omega), \\
\overline{x_2} =\; & h_{\alpha_1}(\omega) h_{\alpha_3}(\omega^2)
h_{\alpha_5}(\omega^2) h_{\alpha_6}(\omega),\qquad &
\overline{x_3} =\; & h_{\alpha_1}(\omega^2) h_{\alpha_3}(\omega)
h_{\alpha_5}(\omega^2) h_{\alpha_6}(\omega).
\end{alignat*}
With our particular choice of Chevalley basis
a direct computation in $E_8(\C)$ shows that
\begin{equation*}
\begin{split}
n_{-\widetilde{\alpha}} =\; & n_8 n_7 n_6 n_5 n_4 n_2 n_3 n_1 n_4 n_3 n_5
n_4 n_2 n_6 n_5 n_4 n_3 n_1 n_7 n_6 n_5 n_4 n_2 n_3 n_4 n_5 n_6 n_7
n_8\cdot \\
& n_7 n_6 n_5 n_4 n_2 n_3 n_1 n_4 n_3 n_5 n_4 n_2 n_6 n_5 n_4 n_3
n_1 n_7 n_6 n_5 n_4 n_2 n_3 n_4 n_5 n_6 n_7 n_8.
\end{split}
\end{equation*}
(A different choice of Chevalley basis may effect this expression by
an order two element in $H$. If a Chevalley basis is chosen such that
the above formula holds then all further formulas will be independent
of the choice.)

From this and the explicit homomorphism above we find, either by direct
computation or by using the relations in $N(H)$, that
\begin{equation*}
\begin{split}
\overline{y_1} =\; & n_1 n_3 n_5 n_6 n_7 n_6 n_5 n_4 n_2 n_3 n_1 n_4 n_3 n_5 n_4 n_2 n_6
n_5 n_4 n_3 n_1 n_7 n_6 n_5\cdot \\
& n_4 n_2 n_3 n_4 n_5 n_6 n_7 n_8 n_7 n_6 n_5 n_4 n_2 n_3 n_1 n_4 n_3
n_5 n_4 n_2 n_6 n_5 n_4 n_3\cdot \\
& n_1 n_7 n_6 n_5 n_4 n_2 n_3 n_4 n_5 n_6 n_7 n_8\cdot
h_{\alpha_1}(-1) h_{\alpha_2}(-1) h_{\alpha_7}(-1),
\end{split}
\end{equation*}
\begin{equation*}
\begin{split}
\overline{y_2} =\; & n_2 n_3 n_1 n_4 n_2 n_3 n_4 n_5 n_4 n_2 n_3 n_4 n_6 n_5 n_4 n_2 n_3
n_1 n_4 n_7 n_6 n_5 n_4\cdot \\
& n_2 n_3 n_1 n_4 n_3 n_5 n_6 n_7 n_8 n_7 n_6 n_5 n_4 n_2 n_3 n_1 n_4
n_3 n_5 n_4 n_2 n_6 n_5\cdot \\
& n_4 n_3 n_1 n_7 n_6 n_5 n_4 n_2 n_3 n_4 n_5 n_8 n_7 n_6\cdot
h_{\alpha_2}(-1) h_{\alpha_5}(-1).
\end{split}
\end{equation*}

\begin{notation} \label{classes:E8}
To distinguish subgroups of $E_8(\C)$, we need some
information on the conjugacy classes of elements of order $3$.
These are given in \cite[Table~VI]{griess} (which is taken from
\cite[Table~4]{CoGr}): There are
$4$ such conjugacy classes, which we label $\mathbf{3A}$,
$\mathbf{3B}$, $\mathbf{3C}$ and $\mathbf{3D}$. Moreover these classes
may be distinguished by their traces on $\mathfrak{e}_8$. Since the
trace of the element $h\in H$ is given by $8+\sum_{\alpha\in\Phi(E_8)}
\alpha(h)$ we get
$\overline{a}\in\mathbf{3A}$, $\overline{x_1}, \overline{x_2},
\overline{x_3}, \overline{y_1}, \overline{y_2}\in\mathbf{3B}$ and
$\overline{x_3 a^{-1}}\in\mathbf{3D}$.
\end{notation}

\begin{notation} \label{symplectic}
If $K$ is a field and $n$ is a natural number, we define {\it the
  group of symplectic similitudes} as $\operatorname{GSp}_{2n}(K) = \{X\in
\GL_{2n}(K) | X^t B X=cB, c\in K^{\times} \}$, where 
$$
B=\underbrace{
\begin{bmatrix}
0 & -1 \\
1 & 0
\end{bmatrix}
\oplus\ldots\oplus
\begin{bmatrix}
0 & -1 \\
1 & 0
\end{bmatrix}
}_{\text{$n$ times}}
$$
We define the homomorphism $\chi:\operatorname{GSp}_{2n}(K)\rightarrow
K^{\times}$ by $\chi(X) =c$, where $X^t B X=cB$. The kernel of
$\chi$ is {\it the symplectic group} $\Sp_{2n}(K)$. (The
notation $\operatorname{CSp}$ is also used in the litterature, cf.\ e.g.\
\cite{lusztig75}.)
\end{notation}

\begin{thm} \label{nontorals:E8-3}
The conjugacy classes of non-toral elementary abelian $3$-subgroups of
$E_8(\C)$ are given by the following table.
$$\begin{array}{|c|c|c|c|c|c|}
\hline
\mbox{rank} & \mbox{name} & \mbox{ordered basis} &
\mbox{$E_8(\C)$-class dist.} &
C_{E_8(\C)}(E) & Z(C_{E_8(\C)}(E))\\
\hline
\hline
3 & E^{3a}_{E_8} & \{\overline{x_1},\overline{y_1},\overline{a}\} &
\mathbf{3A^{18}3B^{8}} & E^{3a}_{E_8}\times
\PSL_3(\C) & E^{3a}_{E_8} \\
\hline
3 & E^{3b}_{E_8} & \{\overline{x_1},\overline{y_1},\overline{x_3}\} &
\mathbf{3B^{26}} & E^{3b}_{E_8}\times G_2(\C) & E^{3b}_{E_8}\\
\hline
\hline
4 & E^{4a}_{E_8} &
\{\overline{x_1},\overline{y_1},\overline{x_3},
\overline{x_3 a^{-1}}\} & \mathbf{3A^{52}3B^{26}3D^{2}} &
E^{4a}_{E_8}\circ_{\left<\overline{x_3 a^{-1}}\right>} \GL_2(\C) &
E^{4a}_{E_8}\circ_{\left<\overline{x_3 a^{-1}}\right>} \T_1\\
\hline
4 & E^{4b}_{E_8} &
\{\overline{x_2},\overline{x_1},\overline{y_1},\overline{a}\} &
\mathbf{3A^{54}3B^{26}} &
E^{4b}_{E_8}\circ_{\left<\overline{x_2} \right>}
(\T_2:\left<\overline{y_2} \right>) &  E^{4b}_{E_8}\\
\hline
4 & E^{4c}_{E_8} &
\{\overline{x_2},\overline{x_1},\overline{y_1},\overline{x_3}\} &
\mathbf{3B^{80}} &
E_{E_8}^{4c}\circ_{\left<\overline{x_2}\right>} \SL_3(\C) &
E^{4c}_{E_8}\\
\hline
\hline
5 & E^{5a}_{E_8} &
\{\overline{x_2},\overline{x_1},\overline{y_1},\overline{x_3},
\overline{x_3 a^{-1}}\} & \mathbf{3A^{156}3B^{80}3D^{6}} &
E^{5a}_{E_8}
\circ_{\left<\overline{x_2},\overline{x_3 a^{-1}}\right>} \T_2 &
E^{5a}_{E_8}
\circ_{\left<\overline{x_2},\overline{x_3 a^{-1}} \right>} \T_2\\
\hline
5 & E^{5b}_{E_8} &
\{\overline{x_1},\overline{y_1},\overline{x_2},\overline{y_2},
\overline{a}\} & \mathbf{3A^{162}3B^{80}} & E^{5b}_{E_8} & E^{5b}_{E_8} \\
\hline
\end{array}$$
In  particular ${}_3 Z(C_{E_8(\C)}(E))=E$ for all non-toral elementary
abelian $3$-subgroups $E$ of $E_8(\C)$.

The Weyl groups of these subgroups with respect to the given ordered
bases are as follows:
$$
W(E^{3a}_{E_8})=\left[
\begin{array}{c|c}
\GL_2(\F_3) &
\begin{array}{c}
* \\
*
\end{array} \\
\hline
\begin{array}{cc}
0 & 0
\end{array}
& \operatorname{det} \\
\end{array}
\right],\:\:
W(E^{3b}_{E_8}) = \SL_3(\F_3),\:\:
W(E^{4a}_{E_8})= \SL_3(\F_3) \times \GL_1(\F_3)
$$
$$
W(E^{4b}_{E_8})=\left[
\begin{array}{c|c|c}
\varepsilon &
\begin{array}{cc}
* & *
\end{array}
& * \\
\hline
\begin{array}{c}
0 \\
0
\end{array}
& \GL_2(\F_3) &
\begin{array}{c}
* \\
*
\end{array}
\\
\hline
0 &
\begin{array}{cc}
0 & 0
\end{array}
& \operatorname{det}
\end{array}
\right],\:\:
W(E^{4c}_{E_8}) = \left[
\begin{array}{c|c}
\varepsilon &
\begin{array}{ccc}
* & * & *
\end{array} \\
\hline
\begin{array}{c}
0 \\
0 \\
0
\end{array}
& \SL_3(\F_3)
\end{array}
\right]
$$
$$
W(E^{5a}_{E_8})=\left[
\begin{array}{c|c|c}
\varepsilon_1 &
\begin{array}{ccc}
* & * & *
\end{array}
& * \\
\hline
\begin{array}{c}
0 \\
0 \\
0
\end{array}
& \SL_3(\F_3) &
\begin{array}{c}
0 \\
0 \\
0
\end{array} \\
\hline
0 &
\begin{array}{ccc}
0 & 0 & 0
\end{array}
& \varepsilon_2
\end{array}
\right],\:\:
W(E^{5b}_{E_8})=\left[
\begin{array}{c|c}
\operatorname{GSp}_4(\F_3) &
\begin{array}{c}
* \\
* \\
* \\
*
\end{array} \\
\hline
\begin{array}{cccc}
0 & 0 & 0 & 0
\end{array}
& \chi \\
\end{array}
\right]
$$
where $\operatorname{det}$ is the determinant of the matrix from
$\GL_2(\F_3)$ in the description of $W(E^{3a}_{E_8})$ and
$W(E^{4b}_{E_8})$. In the description of $W(E_{E_8}^{5b})$, $\chi$
denotes the value of the homomorphism
$\chi:\operatorname{GSp}_4(\F_3)\rightarrow \F_3^{\times}$ defined in
\ref{symplectic} evaluated on the matrix from
$\operatorname{GSp}_4(\F_3)$.
\end{thm}

\begin{rem}
Our information on the subgroups $E^{5a}_{E_8}$ and $E^{3b}_{E_8}$
corrects \cite[Table~II~and~Lem.~11.5]{griess} and
\cite[13.2]{griess} respectively.
\end{rem}

\begin{proof}[Proof of Theorem~\ref{nontorals:E8-3}]
{\em Maximal non-toral subgroups:} By \cite[Lems.~11.7 and
11.9]{griess}, any maximal non-toral elementary abelian
$3$-subgroup of $E_8(\C)$ contains an element of type
$\mathbf{3A}$. We may thus find representatives inside
$C_{E_8(\C)}(\overline{a})=\SL_9(\C)/C_3$. From
\cite[Cor.~11.10]{griess}, it follows that there are two conjugacy
classes of maximal non-toral elementary abelian $3$-subgroups both of
rank $5$.
Moreover, by \cite[Lem.~11.5]{griess}, their preimages in
$\SL_9(\C)$ may be chosen to have the form $(3^{1+2}_+\circ_{C_3} C_9)\times
C_3\times C_3$ and $3^{1+4}_+\circ_{C_3} C_9$.
Using the representation theory of extraspecial $p$-groups
(cf.\ \cite[Ch.~5.5]{gorenstein} or \cite[Satz~V.16.14]{huppert67}) we find that they are represented by
$E^{5a}_{E_8} = \left<\overline{x_2}, \overline{x_1}, \overline{y_1},
  \overline{x_3}, \overline{x_3 a^{-1}}\right>$
and $E^{5b}_{E_8} = \left<\overline{x_1}, \overline{y_1},
  \overline{x_2}, \overline{y_2}, \overline{a}\right>$.

{\em Lower bounds for Weyl groups of maximal non-toral subgroups:}
We can find lower bounds for the Weyl groups of $E^{5a}_{E_8}$ and
$E^{5b}_{E_8}$ by conjugating with elements in the centralizer
$C_{E_8(\C)}(\overline{a})=\SL_9(\C)/C_3$ and the normalizer
$N(H)$ of the maximal torus.

Note that
\begin{xalignat*}{3}
  a &= {\eta I_3 \oplus \eta I_3 \oplus \eta I_3}, &
x_1 &= {\beta \oplus \beta \oplus \beta}, &
x_2 &= {I_3 \oplus \omega I_3 \oplus \omega^2 I_3}, \\
x_3 &= {I_3 \oplus I_3 \oplus \omega I_3}, &
y_1 &= {\gamma \oplus \gamma \oplus \gamma} & &
\end{xalignat*}
and ${(A\oplus B\oplus C)}^{y_2}={B\oplus C\oplus A}$. Conjugation
by ${\tau_1\oplus\tau_1\oplus\tau_1}$, ${\tau_2\oplus\tau_2\oplus\tau_2}$ and
${I_3 \oplus \beta^2 \oplus \beta}$ gives
\begin{xalignat}{6}
{\tau_1 \oplus \tau_1 \oplus \tau_1} :
a &\mapsto a, & x_1 &\mapsto x_1 y_1, & x_2 &\mapsto x_2, & x_3 &\mapsto x_3, &
y_1 &\mapsto y_1, & y_2 &\mapsto y_2.\label{kon1} \\
{\tau_2 \oplus \tau_2 \oplus \tau_2} :
a &\mapsto a, & x_1 &\mapsto x_1, & x_2 &\mapsto x_2, & x_3 &\mapsto x_3, & y_1
&\mapsto x_1 y_1, & y_2 &\mapsto y_2.\label{kon2} \\
{I_3 \oplus \beta^2 \oplus \beta} :
a &\mapsto a, & x_1 &\mapsto x_1, & x_2 &\mapsto x_2, & x_3 &\mapsto x_3, & y_1
&\mapsto x_2 y_1, & y_2 &\mapsto x_1 y_2.\label{kon3}
\end{xalignat}
Now consider the subgroup $E_{E_8}^{5a}$. From \eqref{kon1}--\eqref{kon3} we see
that the elements $\overline{\tau_1 \oplus \tau_1 \oplus \tau_1}$,
$\overline{\tau_2 \oplus \tau_2 \oplus \tau_2}$ and
$\overline{I_3 \oplus \beta^2 \oplus \beta}$ normalize $E_{E_8}^{5a}$ and
that conjugation by these elements induces the automorphisms on $E_{E_8}^{5a}$
given by the matrices $I_5+e_{3,2}$, $I_5+e_{2,3}$ and $I_5+e_{1,3}$.

Letting $\sigma=-(1,4)(2,5)(3,6)\in\SL_9(\C)$ we have
${(A\oplus B\oplus C)}^{\sigma}={B\oplus A\oplus C}$. Using this and the
above we obtain that $\overline{\sigma}$, $\overline{y_2}$ and
$\overline{I_3\oplus I_3\oplus \beta^2}$ normalize $E_{E_8}^{5a}$ and
that conjugation by these elements induces the automorphisms on
$E_{E_8}^{5a}$ given by the matrices
$\diag(2,1,1,1,1)$, $I_5+e_{1,4}+e_{1,5}$ and $I_5+e_{4,3}$.
By using the relations in $N(H)$ given above or by direct computation,
it may be checked that conjugation by the element
\begin{equation*}
\begin{split}
& n_1 n_2 n_4 n_2 n_3 n_5 n_4 n_2 n_3 n_1 n_4 n_3 n_5 n_4 n_6 n_5 n_4 n_2
n_3 n_4 n_7 n_6 n_5 n_4 n_8 n_7 n_6\cdot \\
& n_5 n_4 n_2 n_3 n_1 n_4 n_3 n_5 n_4 n_2 n_6 n_5 n_4 n_7\cdot h(1,1,-1,-1,-1,1,-1,-1)
\end{split}
\end{equation*}
induces the automorphism on $E_{E_8}^{5a}$ represented by the matrix
$\diag(1,1,1,1,2)+e_{2,4}$. It is easy to see that the above matrices
generate the group
$$
W'(E_{E_8}^{5a})=\left[
\begin{array}{c|c|c}
\varepsilon_1 &
\begin{array}{ccc}
* & * & *
\end{array}
& * \\
\hline
\begin{array}{c}
0 \\
0 \\
0
\end{array}
& \SL_3(\F_3) &
\begin{array}{c}
0 \\
0 \\
0
\end{array} \\
\hline
0 &
\begin{array}{ccc}
0 & 0 & 0
\end{array}
& \varepsilon_2
\end{array}
\right]
$$
and thus $W(E_{E_8}^{5a})$ contains this group.

Next consider the subgroup $E_{E_8}^{5b}$. From \eqref{kon1}--\eqref{kon3} we
see that the elements $\overline{\tau_1 \oplus \tau_1 \oplus \tau_1}$,
$\overline{\tau_2 \oplus \tau_2 \oplus \tau_2}$ and
$\overline{I_3 \oplus \beta^2 \oplus \beta}$ normalize $E_{E_8}^{5b}$ and
that conjugation by these elements induces the automorphisms on
$E_{E_8}^{5b}$ given by the matrices $I_5+e_{2,1}$, $I_5+e_{1,2}$ and
$I_5+e_{1,4}+e_{3,2}$. Now note that $a={\Delta_{3,3}(\eta I_3)}$,
$x_2={\Delta_{3,3}(\beta)}$ and $y_2={\Delta_{3,3}(\gamma)}$. Since
$\Delta_{3,3}(M_1)$ commutes with $M_2\oplus M_2\oplus M_2$
for any $M_1, M_2\in M_3(\C)$, the elements
$\overline{\Delta_{3,3}(\tau_1)}$ and
$\overline{\Delta_{3,3}(\tau_2)}$ normalize $E_{E_8}^{5b}$. The
automorphisms induced on $E_{E_8}^{5b}$ by conjugation with these
elements have the matrices $I_5+e_{4,3}$ and $I_5+e_{3,4}$.
By using the relations in $N(H)$ given above or by direct computation,
we get that conjugation by the element
\begin{equation*}
\begin{split}
& n_2 n_8 n_7 n_6 n_5 n_4 n_2 n_3 n_1 n_4 n_3 n_5 n_4 n_2 n_6 n_5 n_4 n_3 n_1
n_7 n_6\cdot \\
& n_5 n_4 n_2 n_3 n_4 n_5 n_6 n_7 n_8 \cdot h(1,-1,-1,-1,-1,1,1,1)
\end{split}
\end{equation*}
induces the automorphism on $E_{E_8}^{5b}$ represented by
the matrix $\diag(1,2,1,2,2)+e_{3,5}$. It now follows that
$W(E_{E_8}^{5b})$ contains the group
$$
W'(E_{E_8}^{5b})=
\left[
\begin{array}{c|c}
\operatorname{GSp}_4(\F_3) &
\begin{array}{c}
* \\
* \\
* \\
*
\end{array} \\
\hline
\begin{array}{cccc}
0 & 0 & 0 & 0
\end{array}
& \chi \\
\end{array}
\right]
$$

{\em Lower bounds for other Weyl groups:} We now show that the
other Weyl groups in the theorem are all lower bounds. To do this
consider one of the non-maximal subgroups $E$ from the table. We then
have $E\subseteq E_{E_8}^{5a}$, and we get a lower bound on $W(E)$
by considering the action on $E$ of the subgroup of
$W'(E_{E_8}^{5a})$ fixing $E$. As an example we find that the
stabilizer of $E_{E_8}^{3a}$ inside $W'(E_{E_8}^{5a})$ is
$$
\begin{bmatrix}
\varepsilon_1 & 0 & 0 & x & x \\
0 & a_{11} & a_{12} & a_{13} & 0 \\
0 & a_{21} & a_{22} & a_{23} & 0 \\
0 & 0 & 0 & \operatorname{det} & 0 \\
0 & 0 & 0 & 0 & \operatorname{det}
\end{bmatrix}
$$
where $\operatorname{det}=a_{11}\cdot a_{22}-a_{12}\cdot a_{21}\neq
0$. The action of such a matrix on $E_{E_8}^{3a}$ is given by
\begin{xalignat*}{3}
\overline{x_1} &\mapsto (\overline{x_1})^{a_{11}} (\overline{y_1})^{a_{21}}, &
\overline{y_1} &\mapsto (\overline{x_1})^{a_{12}} (\overline{y_1})^{a_{22}}, &
\overline{a} &\mapsto (\overline{x_1})^{a_{13}} (\overline{y_1})^{a_{23}}
(\overline{a})^{\operatorname{det}}.
\end{xalignat*}
Thus $W(E_{E_8}^{3a})$ contains the group
$$
W'(E^{3a}_{E_8})=\left[
\begin{array}{c|c}
\GL_2(\F_3) &
\begin{array}{c}
* \\
*
\end{array} \\
\hline
\begin{array}{cc}
0 & 0
\end{array}
& \operatorname{det} \\
\end{array}
\right]
$$
as claimed. Similar computations show that for the remaining subgroups
$E=E_{E_8}^{3b}$, $E_{E_8}^{4a}$, $E_{E_8}^{4b}$ and $E_{E_8}^{4c}$,
the group $W'(E)$ occurring in the theorem is a lower bound for the
Weyl group $W(E)$.

{\em Orbit computation:} Note first that all elementary abelian
$3$-subgroups of rank at most two are toral by
Theorem~\ref{alggrpfacts}\eqref{simply-con}. By using the lower
bounds on the Weyl groups of $E^{5a}_{E_8}$ and $E^{5b}_{E_8}$
established above, we may find a set of representatives for the
conjugacy classes of subgroups of $E^{5a}_{E_8}$ and
$E^{5b}_{E_8}$ of ranks $3$ and $4$.

Under the action of $W'(E_{E_8}^{5a})$, the set of rank
$3$ subgroups of $E_{E_8}^{5a}$ has orbit representatives
\begin{equation*}
\begin{split}
E_{E_8}^{3a}=\left<\overline{x_1},\overline{y_1},\overline{a}\right>,
E_{E_8}^{3b}=\left<\overline{x_1},\overline{y_1},\overline{x_3}\right>,
\left<\overline{x_1},\overline{x_2},\overline{y_1}\right>, & \\
\left<\overline{a},\overline{x_1},\overline{x_2}\right>,
\left<\overline{a},\overline{x_1},\overline{x_3}\right>
\mbox{ and }
\left<\overline{a},\overline{x_2},\overline{x_3}\right>, &
\end{split}
\end{equation*}
and under the action of $W'(E_{E_8}^{5b})$, the set of rank $3$ subgroups
of $E_{E_8}^{5b}$ has orbit representatives
$$
E_{E_8}^{3a}=\left<\overline{x_1},\overline{y_1},\overline{a}\right>,
\left<\overline{x_1},\overline{x_2},\overline{y_1}\right>
\mbox{ and }
\left<\overline{a},\overline{x_1},\overline{x_2}\right>.
$$

Similarly we find that under the action of $W'(E_{E_8}^{5a})$, the set
of rank $4$ subgroups of $E_{E_8}^{5a}$ has orbit representatives
\begin{equation*}
\begin{split}
E_{E_8}^{4a}=\left<\overline{x_1},\overline{y_1},\overline{x_3},
\overline{x_3a^{-1}}\right>,
E_{E_8}^{4b}=\left<\overline{x_2},\overline{x_1},\overline{y_1},
\overline{a}\right>, & \\
E_{E_8}^{4c}=\left<\overline{x_2},\overline{x_1},\overline{y_1},
\overline{x_3}\right>
\mbox{ and }
\left<\overline{a},\overline{x_1},\overline{x_2},
\overline{x_3}\right>, &
\end{split}
\end{equation*}
and that under the action of $W'(E_{E_8}^{5b})$, the set of rank $4$ subgroups
of $E_{E_8}^{5b}$ has orbit representatives
$$
E_{E_8}^{4b}=\left<\overline{x_2},\overline{x_1},\overline{y_1},
\overline{a}\right>
\mbox{ and }
E_0=\left<\overline{x_1},\overline{x_2},\overline{y_1},\overline{y_2}\right>.
$$

{\em Class distributions:} Recall that by \ref{classes:E8},
$\overline{a}$ is in the conjugacy class $\mathbf{3A}$,
$\overline{x_1}$ and $\overline{x_2}$ are in the class
$\mathbf{3B}$ and $\overline{x_3a^{-1}}$ belongs to the class
$\mathbf{3D}$. Using the actions of $W'(E_{E_8}^{5a})$ and
$W'(E_{E_8}^{5b})$ it is then straightforward to verify the class
distributions given in the table. As an example consider the subgroup
$E_{E_8}^{5a}$. Under the action of $W'(E_{E_8}^{5a})$ it contains
$156$ elements conjugate to $\overline{a}$, $78$ elements
conjugate to $\overline{x_1}$, $2$ elements conjugate to
$\overline{x_2}$ and $6$ elements conjugate to $\overline{x_3
a^{-1}}$, which gives the class distribution in the table. Similar
computations give the results for the remaining subgroups.

We also see that the subgroup
$E_0=\left<\overline{x_1},\overline{x_2},\overline{y_1},\overline{y_2}\right>$
has class distribution $\mathbf{3B^{80}}$ and from the class distribution of
$E_{E_8}^{5b}$ we get $E_0=(E^{5b}_{E_8}\cap\mathbf{3B})\cup\{1\}$.
It then follows from \cite[Lem.~11.5]{griess} that $E_0$ is toral.

{\em Other non-toral subgroups:} We see directly that the subgroups
$$
\left<\overline{a},\overline{x_1},\overline{x_2},\overline{x_3}\right>,
\left<\overline{a},\overline{x_1},\overline{x_2}\right>,
\left<\overline{a},\overline{x_1},\overline{x_3}\right>
\mbox{ and }
\left<\overline{a},\overline{x_2},\overline{x_3}\right>
$$
are toral. Since the subgroup
$\left<\overline{x_1},\overline{x_2},\overline{y_1}\right>$ is a subgroup of
$E_0$ it is also toral. Alternatively, from the action of
$W'(E_{E_8}^{5a})$ we see that it is conjugate to the subgroup
$\left<\overline{x_1},\overline{x_2},\overline{x_3}\right>$, which is visibly
toral. Thus any non-toral elementary abelian $3$-subgroup of
$E_8(\C)$ is conjugate to a subgroup in the table. Moreover,
since their class distributions differ, none of these subgroups are conjugate.

To see that the subgroups in the table are non-toral we may proceed as
follows. The subgroup $E_{E_8}^{3a}$ contains the element $\overline{a}$, so by
Theorem~\ref{alggrpfacts}\eqref{centr-toral} it is toral in $E_8(\C)$
if and only if it is toral in
$C_{E_8(\C)}(\overline{a})=\SL_9(\C)/C_3$. However this is
not the case by Theorem~\ref{alggrpfacts}\eqref{covering}, since its preimage
in $\SL_9(\C)$ is non-abelian. The subgroups $E_{E_8}^{4a}$ and
$E_{E_8}^{4b}$ are thus also non-toral since they contain $E_{E_8}^{3a}$. We
saw above that the Weyl group of $E_{E_8}^{3b}$ contains
$\SL_3(\F_3)$, which has order divisible by $13$. Since $13\nmid |W(E_8)|$
it follows from Theorem~\ref{alggrpfacts}\eqref{weyl-section}
that $E_{E_8}^{3b}$ is non-toral. Since $E_{E_8}^{4c}$ contains $E_{E_8}^{3b}$
it is also non-toral.

{\em Centralizers:} The subgroups $E=E_{E_8}^{3a}$, $E_{E_8}^{4a}$,
$E_{E_8}^{4b}$, $E_{E_8}^{5a}$ and $E_{E_8}^{5b}$ are easy to deal
with since they all contain $\overline{a}$, and hence we have
$C_{E_8(\C)}(E)=C_{\SL_9(\C)/C_3}(E)$ for these. It is however
notationally convenient first to change the representatives as follows.
Define $x_4=\tau_2^{-1}\oplus\tau_2^{-1}\oplus\tau_2^{-1}\in \SL_9(\C)$,
and note that conjugation by $(2,7,3,4)(5,8,9,6)\in \SL_9(\C)$ acts
as follows:
\begin{alignat*}{6}
a &\mapsto a,\qquad & x_1 &\mapsto x_2,\qquad & x_2 &\mapsto
x_1^2,\qquad & x_3a^{-1} &\mapsto x_4,\qquad & y_1 &\mapsto y_2,
\qquad & y_2 &\mapsto y_1^2.
\end{alignat*}
In particular we see that $E_{E_8}^{3a}$ is conjugate to
$\left<\overline{x_2},\overline{y_2},\overline{a}\right>$. Moreover we
have
$$
C_{E_8(\C)}(\overline{a},\overline{x_2})= \overline{\left<
y_2,\{A\oplus B\oplus C\mid \det ABC=1\} \right>}.
$$
From this we directly get
\begin{equation*}
\begin{split}
C_{E_8(\C)}(\overline{a},\overline{x_2},\overline{y_2}) =\; &
\overline{\left< x_2,y_2,\{A\oplus A\oplus A\mid \left(\det A\right)^3=1\}
\right>}\\
=\; & \overline{\left< x_2,y_2,a,\{A\oplus A\oplus A\mid \det A=1\}
\right>}\\
\cong\; & \left<\overline{x_2},\overline{y_2},\overline{a}\right> \times
\PSL_3(\C).
\end{split}
\end{equation*}
Thus $C_{E_8(\C)}(E^{3a}_{E_8})=E^{3a}_{E_8} \times
\PSL_3(\C)$ and $Z(C_{E_8(\C)}(E^{3a}_{E_8}))=E^{3a}_{E_8}$.
From the above we see that the elements
$\overline{x_2}$, $\overline{x_3a^{-1}}$ and $\overline{y_2}$ in
$C_{E_8(\C)}(E^{3a}_{E_8})$ correspond to the elements
$\overline{\beta^2}$, $\overline{\tau_2^{-1}}$ and
$\overline{\gamma^2}$ in the $\PSL_3(\C)$-component
of $C_{E_8(\C)}(E^{3a}_{E_8})$. From this we easily compute the
structure of $C_{E_8(\C)}(E)$ for the subgroups $E$ which contain
$E^{3a}_{E_8}$, cf.\ the proof of Theorem~\ref{nontorals:E6-3}.

For the computation of the centralizers of $E^{3b}_{E_8}$ and
$E^{4c}_{E_8}$ we consider the element $g = h_{\alpha_1}(\omega)
h_{\alpha_3}(\omega^2)\in E_8(\C)$. By using \cite[Tables~4 and
6]{CoGr} we see that $g$ belongs to the conjugacy class
$\mathbf{3B}$ and that the centralizer $C_{E_8(\C)}(g)$ has
type $E_6A_2$. The precise structure of this centralizer may be found
as follows. Since $E_8(\C)$ is simply connected,
Theorem~\ref{alggrpfacts}\eqref{simply-con} implies
that $C_{E_8(\C)}(g)$ is connected. Setting
$$
\alpha'=\alpha_1+\alpha_2+2\alpha_3+3\alpha_4+2\alpha_5+\alpha_6,
$$
we see that
$\{\alpha_5,\alpha_8,\alpha_6,\alpha_7,\alpha',\alpha_2\}\cup
\{\alpha_1,\alpha_3\}$ is a system of simple roots of
$C_{E_8(\C)}(g)$. (The simple systems of the components of type
$E_6$ and $A_2$ have been ordered so that the numbering is
consistent with \cite[Planches~I and
V]{bourbaki-lie4:6}.)
From this we get an explicit homomorphism $3E_6(\C)\times
\SL_3(\C)\rightarrow E_8(\C)$ onto the centralizer
$C_{E_8(\C)}(g)$. The kernel is given by
$\left<(z,\omega^2 I_3)\right>$, where $z\in 3E_6(\C)$ denotes the
central element defined in Section~\ref{type:E6}. Thus
$C_{E_8(\C)}(g) = 3E_6(\C)\circ_{C_3} \SL_3(\C)$, and we denote elements in
this central product by $A\cdot B$ where $A\in 3E_6(\C)$ and
$B\in\SL_3(\C)$. In particular we have $g=z\cdot I_3=1\cdot \omega I_3$.

Now consider the subgroup
$E=\left<z\cdot I_3,x_1\cdot\beta,y_1\cdot\gamma\right>$ which is seen
to be an elementary abelian $3$-subgroup of rank $3$ (here the elements
$x_1,y_1\in 3E_6(\C)$ from Section~\ref{type:E6} should not be
confused with the elements $x_1,y_1\in \SL_9(\C)$ from above). We have
$$
C_{E_8(\C)}(z\cdot I_3,x_1\cdot\beta) =
C_{3E_6(\C)\circ_{C_3} \SL_3(\C)}(x_1\cdot\beta) =
\left<y_1\cdot\gamma,C_{3E_6(\C)}(x_1)\circ_{C_3}
  C_{\SL_3(\C)}(\beta)\right>.
$$
We note that $y_1\cdot\gamma$ is not conjugate to its inverse in
$C_{E_8(\C)}(z\cdot I_3,x_1\cdot\beta)$ since no element in
$C_{\SL_3(\C)}(\beta)$ conjugates $\gamma$ into $\gamma^{-1}$ times a
power of $\omega I_3$. Thus $\diag(1,1,-1)\notin W(E)$ and $W(E)\neq
\GL_3(\F_3)$. From the above we also get
\begin{equation*}
\begin{split}
C_{E_8(\C)}(E) & =
\left<y_1\cdot\gamma,x_1\cdot\beta,C_{3E_6(\C)}(x_1,y_1)\circ_{C_3}
  C_{\SL_3(\C)}(\beta,\gamma)\right>\\
& = \left<y_1\cdot\gamma,x_1\cdot\beta,
C_{3E_6(\C)}(\pi^{-1}(E_{E_6}^{2b}))\circ_{C_3} Z(\SL_3(\C))\right>\\
& = \left<y_1\cdot\gamma,x_1\cdot\beta, (\left<z\right>\times
G_2(\C))\circ_{C_3} Z(\SL_3(\C))\right>\\
& = E \times G_2(\C),
\end{split}
\end{equation*}
using the computation of $C_{3E_6(\C)}(\pi^{-1}(E_{E_6}^{2b}))$ from
the last part of the proof of Theorem~\ref{nontorals:E6-3}. Since the
preimage of $E$ in $3E_6(\C)\times \SL_3(\C)$ is non-abelian it follows
from Theorem~\ref{alggrpfacts}\eqref{covering} that $E$ is non-toral in
$3E_6(\C)\circ_{C_3} \SL_3(\C)$. Now
Theorem~\ref{alggrpfacts}\eqref{centr-toral} shows that $E$ is
non-toral in $E_8(\C)$ (alternatively one can also just observe that
$C_{E_8(\C)}(E)$ has rank less than $8$).
From what we have already proved we then see that $E$ is conjugate to
either $E_{E_8}^{3a}$ or $E_{E_8}^{3b}$ in $E_8(\C)$. Since we already
know $C_{E_8(\C)}(E_{E_8}^{3a})$ we conclude that $E$ must be
conjugate to $E_{E_8}^{3b}$ (alternatively one can also compute the
class distribution of $E$ directly). In particular we have
$C_{E_8(\C)}(E_{E_8}^{3b}) = E_{E_8}^{3b}\times G_2(\C)$ and
$W(E_{E_8}^{3b})\neq \GL_3(\F_3)$.

Using the inclusion $3E_6(\C)\subseteq 3E_6(\C)\circ_{C_3}
\SL_3(\C)\subseteq E_8(\C)$ we may also consider the subgroup
$E_{3E_6}^4\subseteq 3E_6(\C)$ from Theorem~\ref{nontorals:3E6-3}
as a subgroup of $E_8(\C)$. Since $E_{3E_6}^4$ is non-toral in
$3E_6(\C)$, it is also non-toral in $3E_6(\C)\circ_{C_3}\SL_3(\C)$, and
hence also in $E_8(\C)$ by
Theorem~\ref{alggrpfacts}\eqref{centr-toral}. Thus $E_{3E_6}^4$ must be
conjugate in $E_8(\C)$ to one of the subgroups $E_{E_8}^{4a}$,
$E_{E_8}^{4b}$ or $E_{E_8}^{4c}$. Comparing the class distributions we
can rule out $E_{E_8}^{4a}$ and $E_{E_8}^{4b}$, so
$E_{3E_6}^4$ is conjugate to $E_{E_8}^{4c}$. From
Theorem~\ref{nontorals:3E6-3} we have
$C_{3E_6(\C)}(E_{3E_6}^4)= E_{3E_6}^4$. Hence
$C_{E_8(\C)}(E_{3E_6}^4)= E_{3E_6}^4\circ_{C_3} \SL_3(\C)$ and thus
$C_{E_8(\C)}(E_{E_8}^{4c})= E_{E_8}^{4c}\circ_{C_3} \SL_3(\C)$. 
We determine the precise structure of the central product below after
the computation of $W(E_{E_8}^{4c})$.

{\em Exact Weyl groups:}
Recall from above that $E_{E_8}^{3a}$ is conjugate to
$\left<\overline{x_2},\overline{y_2},\overline{a}\right>$.
If $W(E_{E_8}^{3a})$ was strictly larger than the group $W'(E_{E_8}^{3a})$
from above, $W(E_{E_8}^{3a})$ would have to contain one of the groups
$$
\left[
\begin{array}{c|c}
\GL_2(\F_3) &
\begin{array}{c}
* \\
*
\end{array} \\
\hline
\begin{array}{cc}
0 & 0
\end{array}
& \varepsilon \\
\end{array}
\right]
\mbox{ or }
\SL_3(\F_3)
$$
since these are the minimal overgroups of $W'(E_{E_8}^{3a})$ inside
$\GL_3(\F_3)$. Thus $W(E_{E_8}^{3a})$ would have to contain one of
the matrices $\diag(1,2,1)$ or $I_3+e_{3,2}$. This would
imply the existence of an element in $C_{E_8(\C)}(\overline{x_2},\overline{a})$
which conjugates $\overline{y_2}$ into either $\overline{y_2}^2$ or
$\overline{y_2 a}$. However we saw above that
$$
C_{E_8(\C)}(\overline{x_2},\overline{a})=
\overline{\left<
y_2,\{A\oplus B\oplus C\mid \det ABC=1\}
\right>},
$$
and from it follows that no such element exists. Thus
$W(E_{E_8}^{3a})=W'(E_{E_8}^{3a})$ as claimed. For the subgroup
$E_{E_8}^{3b}$ we have $\SL_3(\F_3)\subseteq W(E_{E_8}^{3b})\neq
\GL_3(\F_3)$ and hence $W(E_{E_8}^{3b})=\SL_3(\F_3)$.

As in the proof of Theorem~\ref{nontorals:E6-3} we show that
the remaining Weyl groups equal the lower bounds already established,
by looking at what a strictly larger Weyl group would imply for the
subgroups $E_{E_8}^{3a}$ and $E_{E_8}^{3b}$. For $E = E_{E_8}^{4a}$,
$E_{E_8}^{4b}$, $E_{E_8}^{5a}$, and $E_{E_8}^{5b}$, any proper
overgroup of $W(E)$ contains an element which normalizes
$E_{E_8}^{3a}$ but induces an automorphism on it not contained in its
Weyl group. For $E_{E_8}^{4c}$ the result follows by considering the
subgroup $E_{E_8}^{3b}$.

It remains only to determine the precise structure of the central
product $C_{E_8(\C)}(E_{E_8}^{4c})= E_{E_8}^{4c}\circ_{C_3}
\SL_3(\C)$. From the structure of $W(E_{E_8}^{4c})$ we see that the
subgroup $\left<\overline{x_2} \right>$ is invariant under the action
of $W(E_{E_8}^{4c})$. Thus a conjugation which sends to $E_{E_8}^{4c}$
to $E_{3E_6}^4$ must send $\left<\overline{x_2} \right>$ to a
$W(E_{3E_6}^4)$-invariant subgroup of $E_{3E_6}^4$ of rank one. From the
structure of $W(E_{3E_6}^4)$ we see that there is only one such
subgroup, namely $\left<z\cdot I_3 \right>=\left<1\cdot \omega I_3
\right>$. As this is exactly the center of the $\SL_3(\C)$-component
of $C_{E_8(\C)}(E_{3E_6}^4)= E_{3E_6}^4\circ_{C_3} \SL_3(\C)$, we see
that $C_{E_8(\C)}(E_{E_8}^{4c}) =
E_{E_8}^{4c}\circ_{\left<\overline{x_2}\right>} \SL_3(\C)$.
\end{proof}

\subsection{The group $2E_7(\C)$, $p=3$}\label{type:E7}
In this section we consider the elementary abelian $3$-subgroups of
$2E_7(\C)$. We let $H$ be a maximal torus of $2E_7(\C)$,
$\Phi(E_7)$ be the root system relative to $H$,
and choose a {\it realization} (\cite[p.~133]{spr:linalggrp})
$\left(u_\alpha\right)_{\alpha\in\Phi(E_7)}$ of $\Phi(E_7)$ in
$2E_7(\C)$. By \cite[Lem.~8.1.4(iv)]{spr:linalggrp} we may suppose that the root
subgroups $\left(u'_\alpha\right)_{\alpha\in\Phi(E_6)}$ in
$3E_6(\C)\subseteq 2E_7(\C)$ coming from the choice of root subgroups for
$3E_6(\C)$ from Section~\ref{type:E6} satisfy $u_\alpha=u'_\alpha$ for
$\alpha\in\Phi(E_6)$.

For $\alpha=\alpha_i$, $1\leq i\leq 7$, and
$t\in\C^\times$ we define the elements
$$
n_\alpha(t)=u_\alpha(t) u_{-\alpha}(-1/t) u_\alpha(t),\qquad
h_\alpha(t)=n_\alpha(t) n_\alpha(1)^{-1}.
$$
Then the maximal torus consists of the elements $\prod_{i=1}^7
h_{\alpha_i}(t_i)$ and the normalizer $N(H)$ of the maximal torus is
generated by $H$ and the elements $n_i=n_{\alpha_i}(1)$,
$1\leq i\leq 7$.

As in Section~\ref{type:E6} we define the following elements in
$3E_6(\C)\subseteq 2E_7(\C)$:
\begin{eqnarray*}
z &=& h_{\alpha_1}(\omega) h_{\alpha_3}(\omega^2) h_{\alpha_5}(\omega)
h_{\alpha_6}(\omega^2),\qquad
a=h_{\alpha_1}(\omega) h_{\alpha_3}(\omega^2) h_{\alpha_5}(\omega^2)
h_{\alpha_6}(\omega),\\
x_2 &=& h_{\alpha_2}(\omega^2) h_{\alpha_3}(\omega^2)
h_{\alpha_5}(\omega^2),\\
y_2 &=& n_1 n_2 n_3 n_4 n_3 n_1 n_5 n_4 n_2 n_3 n_4 n_5 n_6 n_5 n_4 n_2
n_3 n_1 n_4 n_3 n_5 n_4 n_6 n_5\cdot h_{\alpha_2}(-1).
\end{eqnarray*}

\begin{notation} \label{classes:2E7}
The conjugacy classes of elements of order $3$ in $2E_7(\C)$
are given in \cite[Table~VI]{griess} and \cite[Table~6]{CoGr} from
which we take our notation. There are
five such conjugacy classes, which we label $\mathbf{3A}$,
$\mathbf{3B}$, $\mathbf{3C}$, $\mathbf{3D}$ and $\mathbf{3E}$.
(We take this opportunity to note the following corrections to these references:
 In \cite[Table~VI]{griess} the eigenvalue multiplicities for the class
$\mathbf{3C}$ in $2E_7(\C)$ should be $43$, $45$, $45$; for the class
$\mathbf{3E}$ they should be $67$, $33$, $33$; the centralizer type
of the class $\mathbf{3E}$ should be $D_6T_1$; in
\cite[Table~6]{CoGr} the centralizer types for the classes $\mathbf{3A}$ and
$\mathbf{3E}$ should be $A_6T_1$ and $D_6T_1$ respectively.)
By direct computation we easily obtain the inclusions
\begin{xalignat*}{3}
\mathbf{3C}[3E_6] &\subseteq \mathbf{3C}[2E_7], &
\mathbf{3E}[3E_6] &\subseteq \mathbf{3B}[2E_7], &
\mathbf{3E'}[3E_6] &\subseteq \mathbf{3B}[2E_7],
\end{xalignat*}
corresponding to the inclusion $3E_6(\C)\subseteq 2E_7(\C)$.
\end{notation}

\begin{thm} \label{nontorals:2E7-3}
The conjugacy classes of non-toral elementary abelian $3$-subgroups of
$2E_7(\C)$ are given by the following table.
$$\begin{array}{|c|c|c|c|c|c|}
\hline
\mbox{rank} & \mbox{name} & \mbox{ordered basis} &
\mbox{$2E_7(\C)$-class dist.} &
C_{2E_7(\C)}(E) & Z(C_{2E_7(\C)}(E))\\
\hline
\hline
3 & E^3_{2E_7} & \{a,x_2,y_2\} & \mathbf{3C^{26}} &
E^{3}_{2E_7}\times \SL_2(\C) & E^{3}_{2E_7}\times Z(2E_7(\C))\\
\hline
\hline
4 & E^4_{2E_7} & \{z,a,x_2,y_2\} & \mathbf{3B^{2}3C^{78}} &
E^4_{2E_7}\circ_{\left<z\right>} \T_1 & E^4_{2E_7}\circ_{\left<z\right>} \T_1\\
\hline
\end{array}$$
In particular ${}_3 Z(C_{2E_7(\C)}(E))=E$ for all non-toral elementary
abelian $3$-subgroups $E$ of $2E_7(\C)$.

The Weyl groups of these subgroups with
respect to the given ordered bases are as follows:
$$
W(E^3_{2E_7}) = \SL_3(\F_3),\:\:
W(E^4_{2E_7}) = \left[
\begin{array}{c|c}
\varepsilon &
\begin{array}{ccc}
* & * & *
\end{array} \\
\hline
\begin{array}{c}
0 \\
0 \\
0
\end{array}
& \SL_3(\F_3)
\end{array}
\right]
$$
\end{thm}

\begin{rem}
Our information on the rank $3$ subgroup $E^3_{2E_7}$
corrects \cite[Table~II~and~Thm.~11.16]{griess}.
\end{rem}

\begin{proof}[Proof of Theorem~\ref{nontorals:2E7-3}]
{\em Non-toral subgroups:} From the way the realization
$\left(u_\alpha\right)_{\alpha\in\Phi(E_7)}$ is chosen above, it
follows from Theorem~\ref{nontorals:3E6-3} that $E^{3}_{2E_7}$ and
$E^{4}_{2E_7}$ are elementary abelian $3$-subgroups of $2E_7(\C)$
and that
$$
W(E^3_{2E_7})\supseteq\SL_3(\F_3),\:\:
W(E^4_{2E_7})\supseteq\left[
\begin{array}{c|c}
1 &
\begin{array}{ccc}
* & * & *
\end{array} \\
\hline
\begin{array}{c}
0 \\
0 \\
0
\end{array}
& \SL_3(\F_3)
\end{array}
\right]
$$
In particular both $W(E^3_{2E_7})$ and $W(E^4_{2E_7})$ have
orders divisible by $13$ and since $13\nmid \left|W(E_7)\right|$,
we conclude by Theorem~\ref{alggrpfacts}\eqref{weyl-section} that
$E^3_{2E_7}$ and $E^4_{2E_7}$ are non-toral in $2E_7(\C)$. By
\cite[Thm.~11.16]{griess} there are precisely two
conjugacy classes of non-toral elementary abelian $3$-subgroups in
$2E_7(\C)$, and thus $E^3_{2E_7}$ and $E^4_{2E_7}$ represent these two
conjugacy classes.

{\em Class distributions:} The class distributions follows
directly from the class distributions of the subgroups $E^3_{3E_6}$
and $E^4_{3E_6}$ given in Theorem~\ref{nontorals:3E6-3} and the
information in \ref{classes:2E7} about the behavior of conjugacy
classes in $3E_6(\C)$ under the inclusion $3E_6(\C)\subseteq
2E_7(\C)$.

{\em Weyl groups:} Using our realization
$\left(u_\alpha\right)_{\alpha\in\Phi(E_7)}$ we define a
canonical map $\phi:W\rightarrow N(H)$ as follows
(\cite[9.3.3]{spr:linalggrp}): If $w=s_{\alpha_{i_1}}\ldots
s_{\alpha_{i_r}}$ is a reduced expression for $w\in W$ we let
$\phi(w)=n_{i_1}\ldots n_{i_r}$ (by \cite[Props.~8.3.3 and
9.3.2]{spr:linalggrp} this does not depend on the reduced
expression for $w$). Note that the element $\phi(w)$ is a
representative of $w\in W$ in $N(H)$. Now let $w_0\in W$ be the
longest element in $W$, and let $n_0=\phi(w_0)$. From
\cite[Planche~VI]{bourbaki-lie4:6} it follows that $w_0$ equals the
scalar transformation $-1$, and so conjugation by $n_0$ acts as inversion
on $H$. Now let $w\in W$ and define $w'$ by $w w'=w_0$. Since
$w_0$ is central in our case, we have $\left(w w'\right) w^{-1} =
w^{-1} \left(w w'\right) = w'$ so we conclude that $w' w = w w' =
w_0$. Now let $\ell$ be the length function on $W$. By
\cite[p.~16]{humphreys90} we have $\ell(w)+\ell(w')=\ell(w_0)$. In
general the map $\phi$ is not a homomorphism, but we do have $\phi(w_1
w_2)=\phi(w_1) \phi(w_2)$ if $\ell(w_1 w_2)=\ell(w_1)+\ell(w_2)$
by \cite[Ex.~9.3.4(1)]{spr:linalggrp}. From this it follows that
$\phi(w)\phi(w') = \phi(w')\phi(w) = \phi(w_0) = n_0$, and we
conclude that $n_0$ commutes with $\phi(w)$ for all $w\in W$.

Now consider the element
$$
w = s_1 s_2 s_3 s_4 s_3 s_1 s_5 s_4 s_2 s_3 s_4 s_5 s_6 s_5 s_4 s_2
s_3 s_1 s_4 s_3 s_5 s_4 s_6 s_5.
$$
Since the length of an element is given by the number of
positive roots it sends to negative roots (\cite[Cor.~1.7]{humphreys90}),
the above product is a reduced expression for $w$. Thus we
have $y_2=\phi(w) h_{\alpha_2}(-1)$. From the above we then conclude
that conjugation by $n_0$ acts as follows:
\begin{alignat*}{4}
z &\mapsto z^2,\qquad & a &\mapsto a^2,\qquad & x_2 &\mapsto x_2^2,\qquad & y_2
&\mapsto y_2.
\end{alignat*}
Thus $n_0$ normalizes $E^4_{2E_7}$ and gives the element
$\diag(2,2,2,1)$ in $W(E^4_{2E_7})$. Combined with the above we
conclude that
$$
W(E^4_{2E_7})\supseteq\left[
\begin{array}{c|c}
\varepsilon &
\begin{array}{ccc}
* & * & *
\end{array} \\
\hline
\begin{array}{c}
0 \\
0 \\
0
\end{array}
& \SL_3(\F_3)
\end{array}
\right]
$$
From the inclusion $\Phi(E_7)\subseteq \Phi(E_8)$ we get the inclusion
$2E_7(\C)\subseteq E_8(\C)$, so we may consider
$E^3_{2E_7}$ and $E^4_{2E_7}$ as subgroups of $E_8(\C)$ as well. Since
the orders of their Weyl groups in $2E_7(\C)$ are divisible by $13$
and $13\nmid |W(E_8)|$, we see from
Theorem~\ref{alggrpfacts}\eqref{weyl-section} that $E^3_{2E_7}$ and
$E^4_{2E_7}$ remain non-toral in $E_8(\C)$. Using
Theorem~\ref{nontorals:E8-3} and the class distributions from above we
conclude that $E^3_{2E_7}$ and $E^4_{2E_7}$ are conjugate to
$E^{3b}_{E_8}$ and $E^{4c}_{E_8}$ respectively in
$E_8(\C)$. Theorem~\ref{nontorals:E8-3} now shows that the lower
bounds found above are indeed the Weyl groups of $E^3_{2E_7}$ and
$E^4_{2E_7}$ in $2E_7(\C)$.

{\em Centralizers:}
For the computation of the centralizer of $E^3_{2E_7}$ we consider the
element $g = h_{\alpha_1}(\omega)h_{\alpha_3}(\omega^2) \in
2E_7(\C)$. By using \cite[Table~6]{CoGr} we see that $g$ belongs to
the conjugacy class $\mathbf{3C}$ and that the centralizer
$C_{2E_7(\C)}(g)$ has type $A_5A_2$. The precise structure of this
centralizer may be found as follows. Since $2E_7(\C)$ is simply
connected, Theorem~\ref{alggrpfacts}\eqref{simply-con} implies that
$C_{2E_7(\C)}(g)$ is connected. Setting
$$
\alpha'=\alpha_1+\alpha_2+2\alpha_3+3\alpha_4+2\alpha_5+\alpha_6,
$$
we see that
$\{\alpha_5,\alpha_6,\alpha_7,\alpha',\alpha_2\}\cup
\{\alpha_1,\alpha_3\}$ is a system of simple roots of
$C_{2E_7(\C)}(g)$. (The simple systems of the components of type
$A_5$ and $A_2$ have been ordered so that the numbering is
consistent with
\cite[Planche~I]{bourbaki-lie4:6}.)
From this we get an explicit homomorphism $\SL_6(\C)\times
\SL_3(\C)\rightarrow 2E_7(\C)$ onto the centralizer
$C_{2E_7(\C)}(g)$. The kernel is given by
$\left<(\omega I_6,\omega^2 I_3)\right>$.
Thus $C_{2E_7(\C)}(g) = \SL_6(\C)\circ_{C_3} \SL_3(\C)$, and we denote
elements in this central product by $A\cdot B$ where $A\in \SL_6(\C)$ and
$B\in\SL_3(\C)$. In particular we have $g = \omega I_6\cdot I_3 =
I_6\cdot \omega I_3$.

Now consider the subgroup
$E=\left<\omega I_6\cdot I_3, (\beta\oplus\beta)\cdot\beta,
(\gamma\oplus\gamma)\cdot\gamma^2\right>$ which is seen to be an
elementary abelian $3$-subgroup of rank $3$. We have
\begin{equation*}
\begin{split}
C_{2E_7(\C)}(\omega I_6\cdot I_3, (\beta\oplus\beta)\cdot\beta) &=
C_{\SL_6(\C)\circ_{C_3} \SL_3(\C)}((\beta\oplus\beta)\cdot\beta)\\
& = \left<(\gamma\oplus\gamma)\cdot\gamma^2,
C_{\SL_6(\C)}(\beta\oplus\beta)\circ_{C_3} C_{\SL_3(\C)}(\beta)\right>.
\end{split}
\end{equation*}
From this we get
\begin{equation*}
\begin{split}
C_{2E_7(\C)}(E) & =
\left<(\gamma\oplus\gamma)\cdot\gamma^2,
(\beta\oplus\beta)\cdot\beta,
C_{\SL_6(\C)}(\beta\oplus\beta,\gamma\oplus\gamma)\circ_{C_3}
C_{\SL_3(\C)}(\beta,\gamma)\right>\\
& = \left<(\gamma\oplus\gamma)\cdot\gamma^2,
(\beta\oplus\beta)\cdot\beta,
C_{\SL_6(\C)}(\beta\oplus\beta,\gamma\oplus\gamma)\circ_{C_3}
Z(\SL_3(\C))\right>.
\end{split}
\end{equation*}
Here $C_{\SL_6(\C)}(\beta\oplus\beta,\gamma\oplus\gamma)=
\Delta_{2,3}(\{A\in\GL_2(\C)|(\det A)^3=1\})$ is generated by
$\Delta_{2,3}(\omega I_2)=\omega I_6$ and
$\Delta_{2,3}(\SL_2(\C))$. From this we
get
$$
C_{2E_7(\C)}(E) = \left<E,\Delta_{2,3}(\SL_2(\C))\right>\cong
E\times \SL_2(\C).
$$
Since the preimage of $E$ in $\SL_6(\C)\times \SL_3(\C)$ is non-abelian
it follows from Theorem~\ref{alggrpfacts}\eqref{covering} that $E$ is
non-toral in $\SL_6(\C)\circ_{C_3} \SL_3(\C)$. Now
Theorem~\ref{alggrpfacts}\eqref{centr-toral} shows that $E$ is
non-toral in $2E_7(\C)$ (alternatively one could
also just observe that $C_{2E_7(\C)}(E)$ has rank less than $7$).
It follows that $E$ is conjugate to $E^3_{2E_7}$ in $2E_7(\C)$,
so $C_{2E_7(\C)}(E^3_{2E_7}) = E^3_{2E_7}\times
\SL_2(\C)$. Hence $Z(C_{2E_7(\C)}(E^3_{2E_7})) =
E^3_{2E_7}\times Z(2E_7(\C))$ since the center of $2E_7(\C)$ has
order $2$.

To compute the centralizer of $E^4_{2E_7}$ we note that
$C_{2E_7(\C)}(z)$ has centralizer type $E_6T_1$, and that the
$E_6$-component corresponds to the subgroup $3E_6(\C)\subseteq 2E_7(\C)$. A
computation shows that the $T_1$-component is given by $\T_1 =
\{h(t^2,t^3,t^4,t^6,t^5,t^4,t^3)|t\in \C^\times \}$, and thus we get
$C_{2E_7(\C)}(z) = 3E_6(\C) \circ_{\left<z \right>}
\T_1$. Theorem~\ref{nontorals:3E6-3} now shows that
$C_{2E_7(\C)}(E^4_{2E_7}) = C_{3E_6(\C)}(E^4_{2E_7})
\circ_{\left<z\right>} \T_1 = E^4_{2E_7}\circ_{\left<z\right>} \T_1$.
\end{proof}

%%%%%%%%%%%%%%%%%%%%%%%%%%%%%%%%%%%%%%%%%%%%%%%%%%%%%%%%%%%%%%%%%%%%%%%%%%%

\section{Calculation of the obstruction groups} \label{obstructionssection}
In this section we show that the existence and uniqueness obstructions
to lifting our diagram in the homotopy category to a diagram in the
category of spaces identically vanish. More precisely, we will show
the following theorem.

\begin{thm}[Obstruction Vanishing Theorem] \label{obstructionsmainthm}
Suppose that $X$ is any of the following $p$-compact groups
$(F_4)\threecom$, $(E_6)\threecom$, $(E_7)\threecom$, $(E_8)\threecom$,
$(E_8)\fivecom$ or $\PU(n)\pcom$ (any $p$), or suppose that $p$ is odd
and $X$ is connected with $H^*(BX;\Z_p)$ a polynomial algebra. Then 
$$\lim_{\A(X)}\hspace{-5pt}{}^i \pi_j(B\cZ(\cC_X(-))) = 0,\:\:\mbox{for all $i,j$}.$$
\end{thm}

(See Theorem~\ref{polyinvthm} for an explanation of why exactly these
$p$-compact groups need attention.) Note that for the purpose of
Theorem~\ref{ndetthm} we only need to calculate the above groups for
$j=1,2$ and $i=j$ or $i=j+1$.

We prove the theorem by filtering the functor $F_j =
\pi_j(B\cZ(\cC_X(-)))$, and showing that all filtration quotients
vanish (with a small twist for $\PU(2)\twocom$). First we show that the
quotient functor of $F_j$ concentrated on the toral elementary abelian
$p$-subgroups has vanishing limits, using a Mackey functor argument
which first appeared in \cite{DW93}. This takes care of the case where
$H^*(BX;\Z_p)$ is a polynomial algebra since in this
case all subgroups are toral by Lemma~\ref{torsionfreetoral}. For the
exceptional compact connected Lie groups we then
continue and filter the non-toral part of the functor such that the filtration quotients are
concentrated on only one non-toral subgroup, and use a formula of
Oliver \cite{bob:steinberg} to show that the higher limits of these
subquotient functors all vanish. For $\PU(n)$ we use a variant of this
technique by suitably grouping the non-toral subgroups and using a
combination of Oliver's formula and the Mackey functor argument we
used for the toral part. We divide the proof up in three subsections corresponding to the toral part, the non-toral part for the exceptional groups, and the non-toral part for the projective unitary groups.

We use the notation $\St_G$ to denote the Steinberg module over $\Z_p$
of a finite group of Lie type $G$ of characteristic $p$, defined as
the top homology group with $\Z_p$ coefficients of the Tits building
of $G$ (see e.g., \cite{humphreys87}). In the special case $G=\GL(E)$
we also write $\St(E)$ for the Steinberg module.

\subsection{The toral part}
Define a quotient functor $F^{\operatorname{tor}}_j$ of $F_j$ by
setting $F^{\operatorname{tor}}_j(V) = F_j(V)$ if $V$ is toral and
$F^{\operatorname{tor}}_j(V) = 0$ if $V$ is non-toral.  Let
$\A^{\operatorname{tor}}(X)$ denote the full subcategory of $\A(X)$
consisting of toral subgroups.  From the chain complex defining higher
limits (see e.g., \cite[App.~II, 3.3]{GZ67}) it follows that
$$
\lim_{\A(X)}\hspace{-5pt}{}^* F^{\operatorname{tor}}_j  \cong \lim_{\A^{\operatorname{tor}}(X)}\hspace{-10pt}{}^* F^{\operatorname{tor}}_j.
$$
In order to use a Mackey functor argument on the right-hand side we
need to see that the functor $F^{\operatorname{tor}}_j$ is indeed a Mackey functor. This is accomplished by the following lemma, whose assumptions are always satisfied for $p$ odd by \cite[Thm.~7.5]{dw:center}.

\begin{lemma} \label{centralizercenter}
Fix a connected $p$-compact group $X$ and let $\dT$ be the discrete
approximation to a maximal torus $T$ in $X$. Let $V\subseteq\dT$ be a non-trivial elementary abelian $p$-subgroup.

If $\dT^{W_{\cC_X(V)_1}}$ is a discrete approximation to $\cZ(\cC_X(V)_1)$
then ${\dT}^{W_{\cC_X(V)}}$ is a discrete approximation to
$\cZ(\cC_X(V))$. In particular in this case $\pi_1(B\cZ(\cC_X(V))) =
H^1(W_{\cC_X(V)};L_X)$ and $\pi_2(B\cZ(\cC_X(V))) = (L_X)^{W_{\cC_X(V)}}$, where $L_X = \pi_1(T)$.
\end{lemma}

\begin{rem} \label{sowierdness}
For a connected $p$-compact group $X$ and $p$ odd, the fixed point set
$\dT^{W_X}$ always equals a discrete approximation to the center of $X$
by \cite[Thm.~7.5]{dw:center}. If $X$ is the $\F_p$-completion of a
compact connected Lie group then this is likewise the case for $p=2$
unless $X$ contains a direct factor isomorphic to $\SO(2n+1)\twocom$,
by \cite[Thm.~1.6]{matthey02}.
\end{rem}

\begin{proof}[Proof of Lemma~\ref{centralizercenter}]
Set $Y = \cC_X(V)$ and $\pi = \pi_0(Y)$ for short.
Since $V$ is toral, $\dT$ is in a canonical way a discrete
approximation to a maximal torus in $Y$.

First observe that the center of $Y$ has discrete approximation in
$\dT$. Indeed, otherwise there would by \cite[Thm.~6.4]{dw:center}
exist a central homomorphism $f: \Z/p^n \to \dN_{p,Y}$ with image not
in $\dT$, which would produce a homomorphism $f': \Z/p^n \to
\dN_{p,X}$ commuting with $\dT$ but not in $\dT$, which contradicts
the fact that $T$ is self-centralizing in $X$ by
\cite[Prop.~9.1]{DW94}, since $X$ is connected.

Suppose that $\dT^{W_{\cC_X(V)_1}}$ is a discrete approximation to
$\cZ(Y_1)$ and set $C = \dT^{W_{\cC_X(V)}}$. We want to show that $C$ is
central in $Y$.
Let $f:BC \to BY_1$ be the natural inclusion. We have an obvious
diagram with horizontal maps fibrations
$$
\xymatrix{
\map(BC,BY_1)_{\{f\}} \ar[r] \ar[d] & \map(BC,BY)_f \ar[r] \ar[d] & \map(BC,B\pi)_0 \ar@{=}[d] \\
BY_1   \ar[r] & BY \ar[r] & B\pi
}
$$
where $\{f\}$ denotes the set of homotopy classes of maps $BC \to
BY_1$ generated by $f$ under the $\pi$-action on $BY_1$. If we can
show that $\{f\}$ consists of just $f$ then it follows from the
five-lemma that the middle vertical map is a homotopy equivalence,
since our assumption implies that $C$ is central in $Y_1$.

To see that the action is trivial consider the following diagram:
$$
\xymatrix{
 &  B\N_{Y_1} \ar[r]^-{\tilde g} \ar[d] & B\N_{Y_1} \ar[d] \\
BC \ar[r]^-f \ar[ur]^-{\tilde f} & BY_1 \ar[r]^-g  &  BY_1
}
$$
where $\tilde f$ is the natural inclusion of $BC$ in $B\N_{Y_1}$, $g$
is an element in $\Aut(BY_1)$ induced by an element in $\pi$ and
$\tilde g$ is the corresponding self-map of $B\N_{Y_1}$ defined via
Lemma~\ref{themap}. However by the definition of $C$, the composite
$\tilde g \tilde f$ is homotopic to $\tilde f$ for all $g$ induced by an
element in $\pi$, so $f$ is homotopic to $gf$ as well. Hence we have
shown that $C$ is central in $Y$ and since the center of $Y$ has
discrete approximation in $\dT$ it is obviously the largest subgroup
with this property. So $C$ is a discrete approximation to $\cZ Y$ as
wanted.

The last statement about the homotopy groups now follows easily using
the long exact sequence in group cohomology.
\end{proof}

\begin{rem}
The above lemma should be compared to Lemma~\ref{ndetcomponent} and
Remark~\ref{centralizerassumptionremark} which have slightly different
assumptions and conclusions.
\end{rem}

We are now ready for the proposition, essentially contained in \cite[\S 8]{DW93}, which will take care of the toral part. See \cite[Def.~7.3]{dw:center} and
Remark~\ref{singularsetremark} for the definition of the singular set $\sigma(s)$ of a reflection $s$, and note that the assumptions are always satisfied for $p$ odd by \cite[Thm.~7.5]{dw:center} and the definition of $\sigma(s)$.

\begin{prop} \label{toralobstructionslemma}
Let $X$ be a connected $p$-compact group, and assume that for each
non-trivial toral elementary abelian $p$-subgroup $V \in \A(X)$ the
fixed point set $\dT^{W_{\cC_X(V)_1}}$ is a discrete approximation to
$\cZ(\cC_X(V)_1)$. Then
$$
\lim_{\A(X)}\hspace{-5pt}{}^{i} F^{\operatorname{tor}}_j = \begin{cases}
H^{2-j}(W_X;L_X) & \text{if $i=0$ and $j=1$, $2$} \\
0 & \text{otherwise.}
\end{cases}
$$
where $H^{2-j}(W_X;L_X) \cong \pi_j(B\cZ(X))$ if $\dT^W$ is a discrete
approximation to $\cZ(X)$.

In particular if for all reflections
$s \in W_X$ the singular set $\sigma(s)$ equals the fixed point set
$\dT^{\langle s \rangle}$ then the assumptions above are satisfied and
$\lim_{\A(X)}{}^{i} F^{\operatorname{tor}}_j = \pi_j(B\cZ(X))$ if $i=0$
and zero otherwise. 
\end{prop}

\begin{proof}
The first part of the proof consists of a translation of \cite[\S
8]{DW93} into the current notation.
By \cite[Prop.~3.4]{dw:split} all morphisms in $\A(X)$ between toral
subgroups $V \to X$ and $V' \to X$ are induced by inclusions and
action by elements of $W_X$.  Hence we can identify
$\A^{\operatorname{tor}}(X)$, up to equivalence of
categories with a category which has objects non-trivial subgroups of
${}_p\dT \cong (\Z/p)^r$ (where $r$ is the rank of $T$) and morphisms
the homomorphisms between subgroups induced by inclusions and action by
$W_X$. Also, by \cite[Thm.~7.6(1)]{dw:center}, $W_{\cC_X(V)}$ consists of the
elements in $W_X$ which pointwise fix $V$.
Hence Lemma~\ref{centralizercenter} shows that the functor
$F_2^{\operatorname{tor}}$ on $\A^{\operatorname{tor}}(X)$ is isomorphic to the
functor $\alpha^{0}_{\Gamma,M}$ on $\A_\Gamma$ from \cite[\S 8]{DW93},
where $\Gamma = W_X$ and $M = L_X$. Likewise $F_1^{\operatorname{tor}}$
is isomorphic to
$\alpha^1_{\Gamma,M}$. (Note that there is the difference in formulation from
\cite[\S 8]{DW93} that $M$ is a $\Z_p\Gamma$-module rather than an
$\F_p\Gamma$-module, but the proof is identical.) Therefore \cite[\S
8]{DW93} (which is a Mackey functor argument, which can also be
deduced from \cite{DW92} or \cite{JM92}) implies the first part of the
lemma about obstruction groups.

To see the last part about the singular set recall that for an abelian
subgroup $A \subseteq \dT$,
$$
\bigcap_{\substack{\mbox{ reflections $s\in W_X$}\\ \mbox{ such that
      $A \subseteq \sigma(s)$}}}\sigma(s)
$$
is a discrete approximation to $\cZ(\cC_X(A)_1)$ by
\cite[Thms.~7.5 and 7.6]{dw:center}. Hence if $\sigma(s) = \dT^{\langle s
  \rangle}$ for all reflections $s\in W_X$, then the assumption of the first
part is obviously satisfied.
\end{proof}

\begin{rem} \label{singularsetremark}
Let $G$ be a compact connected Lie group with maximal torus $T$, and
let $\alpha$ be a root of $G$ relative to $T$ with corresponding
reflection $s_{\alpha}$. In this case the singular set
$\sigma(s_{\alpha})$ is just the discrete approximation of the kernel
$U_{\alpha}$ of $\alpha$ on $T$. To see this note that by \cite[\S 4,
no.\ 5]{bo9} the reflection $s_{\alpha}$ lifts to an element
$n_{\alpha}$ (denoted by $\nu(\theta)$ in \cite{bo9}) which satisfies
$n_{\alpha}^2=\exp(\alpha^{\vee}/2)$; the statement now
follows cf.\ \cite[Pf.~of~Prop.~3.1(ii)]{matthey02}.
\end{rem}

Explicit calculations \cite[Prop.~3.1(ii)]{matthey02} (see also
\cite{DW00}, \cite[Prop.~3.2(vi)]{JMO95}, and \cite[\S 4]{osse97})
show that for a compact connected Lie group $G$,
$\sigma(s)$ in fact always equals
$\dT^{\langle s\rangle}$ except when $G$ contains a direct factor
isomorphic to $\SO(2n+1)$, $p=2$ and $s$ is a reflection corresponding
to a short root.  Combining this with
Proposition~\ref{toralobstructionslemma}, now gives the following
calculation of the toral part of the obstruction groups, whose full
strength at $p=2$ we will however not use here.

\begin{cor} \label{strongcor}
Let $G$ be a compact connected Lie group with no direct factors
isomorphic to $\SO(2n+1)$ when $p=2$. Set $X = G \pcom$. Then
$$
\lim_{\A(X)}\hspace{-5pt}{}^{i} F^{\operatorname{tor}}_j =
\begin{cases}
\pi_j(B\cZ(X)) & \text{if $i=0$} \\
0 & \text{otherwise.}
\end{cases}
$$\QED
\end{cor}

\begin{proof}[Proof of {Theorem~\ref{obstructionsmainthm}} when $H^*(BX;\Z_p)$ is polynomial, $p$ odd]
If $H^*(BX;\Z_p)$ is a polynomial algebra concentrated in even degrees
then all elementary abelian $p$-subgroups are toral by
Lemma~\ref{torsionfreetoral}, so $F = F^{\operatorname{tor}}$. Since
$p$ is odd the assumptions of Proposition~\ref{toralobstructionslemma} are
satisfied and Theorem~\ref{obstructionsmainthm} follows.
\end{proof}

\subsection{The non-toral part for the exceptional groups}
In this subsection we prove Theorem~\ref{obstructionsmainthm} when $X$
is the $\F_p$-completion of one of the exceptional groups and $p$ is
odd. Let $F^E_j$ denote the subquotient functor of $F_j$ concentrated
on a non-toral elementary abelian $p$-subgroup $E$. By Oliver's formula
\cite[Prop.~4]{bob:steinberg}
$$
\lim_{\A(X)}\hspace{-5pt}{}^iF_j^E =
\begin{cases}
\Hom_W(\St(E), F_j(E)) & \text{if $i = \rk E -1$}\\
0 & \text{otherwise.}
\end{cases}
$$
We now embark on proving some lemmas which will be used to show that
these obstructions groups identically vanish. (For the use in
Theorem~\ref{ndetthm} we actually only need this when $E$ has rank at
most four.)

Since $\cZ(\cC_X(E))$ is the
$\F_p$-completion of an abelian compact Lie group (cf.\ \cite[Thm.~1.1]{dw:center}), $F_j =0$ unless $j=1,2$. The
following lemma reduces the problem of showing that the obstruction
groups vanish to showing that $\Hom_{W(E)}(\St(E),{}_p\dZ(\cC_X(E))) = 0$,
where ${}_p\dZ(\cC_X(E))$ is the finite group of elements of order $p$
in the discrete approximation $\dZ(\cC_X(E))$.

\begin{lemma} \label{reduction}
Let $A$ be an abelian compact Lie group, and let ${}_pA$ and $A^p$
denote the kernel and the image of the $p$th power map on $A$ (using
multiplicative notation). Let $P$ be a finitely generated projective
$\Z_pW$-module for a finite group $W$, and assume that $A$ has a
module action of $W$.
Then $\Hom_W(P,{}_pA)=0$ if and only if $\Hom_W(P,\pi_1(A) \otimes
\Z_p) = \Hom_W(P,\pi_0(A) \otimes \Z_p) = 0$.
\end{lemma}

\begin{proof}
The long exact sequence of homotopy groups associated  to the exact
sequence of groups $1 \to A^p \to A \to A/A^p \to 1$ shows that the
inclusion $A^p \into A$ induces an isomorphism $\pi_1(A^p) \xrightarrow{\cong}
\pi_1(A)$ and an injection $\pi_0(A^p) \into \pi_0(A)$.

Hence the exact sequence $1 \to {}_pA \to A \xrightarrow{p} A^p \to 1$
produces the following diagram, where the row, as well as the sequence
going through $\pi_i(A)$ instead of $\pi_i(A^p)$, is exact. 
$$
\xymatrix{
\pi_1(A) \ar[r]^-{p} \ar[dr]^-{p} & \pi_1(A^p)
  \ar[d]^-{\cong} \ar[r] & \pi_0({}_pA) \ar[r] & \pi_0(A) \ar[r]^-{p}
  \ar[dr]^-{p} & \pi_0(A^p) \ar@{^{(}->}[d] \\
& \pi_1(A) & & & \pi_0(A) &
}
$$
Apply the exact functor $\Hom_W(P,-\otimes \Z_p)$ to this diagram. The
lemma now follows from Nakayama's lemma, using that $\pi_0(A)$ is
finite and $\pi_1(A)$ is finitely generated. 
\end{proof}

The following elementary observation is so useful that it is worth
stating explicitly.

\begin{lemma}\label{minusone}
Suppose that $W$ is a subgroup of $\GL(E)$, $p$ odd, such that $-1 \in W$.
Then $\Hom_W(\St(E),E) = 0$.
\end{lemma}

\begin{proof}
Set $Z = \langle -1 \rangle$. Since $Z$ acts trivially on $\St(E)$ we
have
$$\Hom_W(\St(E),E) \subseteq \Hom_Z(\St(E),E) = \Hom_Z(\St(E),E^Z) = 0.$$
\end{proof}

We also need the following lemma, which follows from a theorem of S.~Smith \cite{smith82}.

\begin{lemma}\label{smithlemma}
Let $G$ be a finite group of Lie type of characteristic $p$, and let
$P$ be a parabolic subgroup of $G$ with corresponding unipotent
radical $U$ and Levi subgroup $L \cong P/U$. Suppose that $W$ is a
subgroup with $U \subseteq W \subseteq P$, and let $M$ be an
$\F_pW$-module. 
\begin{enumerate}
\item \label{smith1}
If $U$ acts trivially on $M$, then $\Hom_W(\St_G,M) = \Hom_{W/U}(\St_L,M)$.
\item \label{smith2}
If $\St_L \otimes \F_p$ is irreducible as an $\F_pW/U$-module
and if $M$ has a finite filtration as an $\F_pW$-module, with
filtration quotients of $\F_p$-dimension strictly less than
$\operatorname{rank}_{\Z_p}\St_L$ then $\Hom_{W}(\St_G,M) = 0$.
\end{enumerate}
\end{lemma}

\begin{proof}
Since $U$ acts trivially on $M$, $\Hom_W(\St_G,M) =
\Hom_{W/U}((\St_G)_U,M)$ where $(-)_U$ denotes coinvariants. But since
the Steinberg module is self-dual, as is clear from its definition as
a homology module, $(\St_G)_U \cong (\St_G)^U$. Now Smith's theorem
\cite{smith82} says that $(\St_G)^U \cong \St_L$, which proves the
first part of the lemma.

For the second part, we can assume that the filtration quotients are
simple $\F_pW$-modules. Since $U \subseteq O_p(W)$, $U$ acts trivially
on any irreducible $\F_pW$-module, by elementary representation
theory. Hence the second part follows from the first together with a
dimension consideration.
\end{proof}

The above lemma is usually used in conjunction with the following
obvious observation.

\begin{lemma}\label{boundrank}
Let $E$ be a non-toral elementary abelian $p$-subgroup of a compact
Lie group $G$. Then the $\F_p$-dimension of ${}_pZC_G(E)$ is at most 
equal to the maximal dimension of a non-toral elementary abelian
$p$-subgroup of $G$, and $E$ is a $W(E)$-submodule of ${}_pZC_G(E)$.\QED
\end{lemma}

The last lemma we shall need is a concrete calculation. First we need
the following.

\begin{lemma} \label{steinbergcomp}
Let $\Gamma$ be a subgroup of $\GL_n(\F_p)$ and let $M$ be a
$\Z_{(p)}\Gamma$-module. Then
$$
\Hom_{\Z_{(p)}\Gamma}(\St_{\GL_n(\F_p)},M) \cong \sum_{I\subseteq
  \{1, \ldots, n-1\}} (-1)^{|I|} \bigoplus_{g\in \Gamma\setminus \GL_n(\F_p)/P_I} M^{\Gamma
  \cap g P_I g^{-1}},
$$
where $P_I$ is the parabolic subgroup of $\GL_n(\F_p)$ corresponding
to the subset $I$.
\end{lemma}

\begin{proof}
This follows easily from the fact that $\St_{\GL_n(\F_p)}$ is
isomorphic to $\sum_{I\subseteq \{1, \ldots, n-1\}}
(-1)^{|I|} 1_{P_I}^{\GL_n(\F_p)}$ as $\Z_{(p)}\GL_n(\F_p)$-modules
\cite[Cor.~1.2]{KM86} combined with Frobenius reciprocity and the
double coset formula (see e.g.\ \cite[Prop.~3.3.1(ii) and Cor.~3.3.5(iv)]{benson91}).
\end{proof}

\begin{lemma} \label{calclemma}
Let $E$ be a rank 4 elementary abelian $3$-group, and let $W =
\SL_3(\F_3) \times 1 \subseteq \GL(E)$. Then $\Hom_W(\St(E),E) = 0$.
\end{lemma}

\begin{proof}
This is most easily checked by implementing the formula from
Lemma~\ref{steinbergcomp} on a computer, e.g., using {\sc Magma}
\cite{magma1}. However in this case the calculation is sufficiently
small to be redone by hand with some effort. Alternatively
one can use some ad hoc Lie theoretic arguments. (We are grateful to
A.\ Kleschev and H.\ H.\ Andersen for sketching a couple of such
arguments to us---however, since these arguments are rather involved
compared to the size of the calculation at hand we will not provide
them here.)
\end{proof}

Before we start going through the exceptional groups, we need to
introduce a bit of notation. For an $\F_p$-vector space $E = \langle
e_1, \ldots, e_n\rangle$, we let $E_{ij\ldots}$ denote the subspace
generated by $e_i$, $e_j$, $\ldots$. Likewise we let $P_{ij\ldots}$
(resp.\ $U_{ij\ldots}$) denote the parabolic subgroup (resp.\
its unipotent radical) of $\GL(E)$ corresponding to the simple roots
$\alpha_i, \alpha_j,\ldots$ in the standard basis and notation. For
example in $\GL_3(\F_p)$, $U_{2}$ is the subgroup
$$
\begin{bmatrix}
1 & \ast & \ast \\
0 & 1 & 0 \\
0 & 0 & 1
\end{bmatrix}
$$

\begin{proof}[Proof of Theorem~\ref{obstructionsmainthm} when $X = G\pcom$, $(G,p) = (E_8,5)$, $(F_4,3)$, $(E_6,3)$, $(E_7,3)$, or $(E_8,3)$]

By Lemma~\ref{reduction} it is enough to see that $\Hom_{W(E)}(\St(E),
{}_pZC_G(E)) = 0$ for all non-toral elementary abelian
$p$-subgroups $E$ of $G$. We proceed case-by-case.

{\noindent $(E_8,5)$ and $(F_4,3)$:} By \cite[Lem.~10.3 and
Thm.~7.4]{griess} $G$ has, up to conjugacy, exactly one non-toral elementary
abelian $p$-subgroup $E$, which has rank $3$, Weyl group $W=\SL(E)$, and
(since $E$ is necessarily maximal) $E = {}_pZC_G(E)$. Since $\St(E)$
is an irreducible $\SL(E)$-module of dimension $p^3$ we have
$\Hom_W(\St(E),E)=0$.

{\noindent  $(E_7,3)$:} By Theorem~\ref{nontorals:2E7-3} $E_7$ has,
up to conjugacy, two non-toral elementary abelian $3$-subgroups
$E^3_{2E_7}$ and $E^4_{2E_7}$ of rank $3$ and $4$ respectively. Since
$W(E^3_{2E_7}) = \SL_3(\F_3)$  a dimension consideration as above
gives $\Hom_{W(E^3_{2E_7})}(\St(E^3_{2E_7}),{}_3ZC_G(E^3_{2E_7}))
=0$. For $E^4_{2E_7}$  (whose Weyl group is listed in
Theorem~\ref{nontorals:2E7-3}) we use Lemma~\ref{smithlemma}\eqref{smith2},
taking $U= U_{23}$, which immediately gives that also
$\Hom_{W(E^4_{2E_7})}(\St(E^4_{2E_7}),{}_3ZC_G(E^4_{2E_7}))= 0$.

{\noindent $(E_6,3)$:} By Theorem~\ref{nontorals:E6-3} $E_6$ has, up to
conjugacy, eight non-toral elementary abelian $3$-subgroups all of rank at most
$4$. 
For $E = E^{2b}_{E_6}$, $E^{3b}_{E_6}$, $E^{3c}_{E_6}$ and
$E^{4a}_{E_6}$ we have $\Hom_W(\St(E),{}_3ZC_G(E)) = 0$ by
Lemma~\ref{smithlemma}\eqref{smith2} taking $U = 1$, $1$, $U_2$, and
$U_1$ respectively (note that we do not need to know ${}_3ZC_G(E)$
exactly since the rough bound from Lemma~\ref{boundrank} will do.)
For $E = E^{2a}_{E_6}$, $E^{3a}_{E_6}$, $E^{3d}_{E_6}$, and $E^{4b}_{E_6}$ we
use that by Theorem~\ref{nontorals:E6-3} $E = {}_3ZC_G(E)$ in these
cases (a fact we did not need above), and also $-1 \in W(E)$ so
Lemma~\ref{minusone} applies to show that $\Hom_W(\St(E),{}_3ZC_G(E))=0$.

{\noindent $(E_8,3)$:} By Theorem~\ref{nontorals:E8-3} $E_8$ has
seven conjugacy classes of non-toral elementary abelian $3$-subgroups.
If $E = E^{3a}_{E_8}$, $E^{3b}_{E_8}$, $E^{4b}_{E_8}$, or $E^{4c}_{E_8}$ then
$\Hom_W(\St(E),{}_3ZC_G(E))=0$ by Lemma~\ref{smithlemma}\eqref{smith2}
taking $U = U_1$, $1$, $U_2$, and $U_{23}$ respectively (note again
that we do not need to know ${}_3ZC_G(E)$ exactly).
Now consider $E = E^{4a}_{E_8}$. By
Theorem~\ref{nontorals:E8-3} we have ${}_3ZC_G(E) = E$, and by
Lemma~\ref{calclemma} $\Hom_{W(E)}(\St(E),E)=0$.
Next, suppose that $E = E^{5a}_{E_8}$. Then $E$ has the $W$-invariant
subspace $E_1$ and $U = U_{234}$ acts trivially on $E_1$ and $E/E_1$. By
Lemma~\ref{smithlemma}\eqref{smith1} we get $\Hom_W(\St(E),E_1) =
\Hom_{W/U}(\St(E_1)\otimes \St(E_{2345}),E_1)$. Since $W/U$ contains
$\GL_1(\F_3)\times 1$ and $\Hom_{\GL_1(\F_3)}(\St(E_1),E_1) = 0$ by
Lemma~\ref{minusone} we get $\Hom_W(\St(E),E_1) =
0$. Lemma~\ref{smithlemma}\eqref{smith1} also shows that
$\Hom_W(\St(E),E/E_1) = \Hom_{W/U}(\St(E_1)\otimes
\St(E_{2345}),E/E_1)$, and since
$W/U$ contains $1\times \SL_3(\F_3)\times 1$ we get
$\Hom_W(\St(E),E/E_1) = 0$ by Lemma~\ref{calclemma}.
Thus $\Hom_W(\St(E),E) = 0$ as desired.
Finally, let $E = E^{5b}_{E_8}$. The subspace $E_{1234}$ is
$W$-invariant and $U = U_{123}$ acts trivially on $E_{1234}$ and
$E/E_{1234}$. By Lemma~\ref{smithlemma}\eqref{smith1} we have
$\Hom_W(\St(E),M) = \Hom_{W/U}(\St(E_{1234})\otimes \St(E_5),M)$ for
$M=E_{1234}$ and $M=E/E_{1234}$ and since $W/U$
contains $\Sp_4(\F_3)\times 1$ it suffices to prove that
$\Hom_{\Sp_4(\F_3)}(\St(E_{1234}),E_{1234}) = 0$ and
$\Hom_{\Sp_4(\F_3)}(\St(E_{1234}),\F_3) = 0$. Since $-1\in
\Sp_4(\F_3)$, the first claim follows from Lemma~\ref{minusone}.
The second claim follows from \cite{henningbrev} or a direct computer
calculation based on Lemma~\ref{steinbergcomp}. Thus $\Hom_W(\St(E),E)
= 0$ in this case as well. This exhausts the list.
\end{proof}

\subsection{The non-toral part for the projective unitary groups}

We now embark on proving Theorem~\ref{obstructionsmainthm} for $X =
\PU(n)\pcom$. We will throughout this subsection use the notation for
elementary abelian $p$-subgroups of $X$ introduced in
Section~\ref{punsection}.

We first treat the toral case directly (see also Corollary~\ref{strongcor}).

\begin{lemma} \label{puntoralobs}
Let $X = \PU(n)\pcom$. Then
$$
\lim_{\A(X)}\hspace{-5pt}{}^i F^{\operatorname{tor}}_j =
\begin{cases}
\Z/2 & \text{if $n=p=2$, $i=0$ and $j=1$} \\
0 & \text{otherwise.}
\end{cases}
$$
\end{lemma}

\begin{proof}
If $n \neq 2$ then it is immediate to check that
$\dT^{\langle s\rangle}$ is $p$-divisible for an arbitrary
reflection $s \in W_X$, so $\sigma(s) = \dT^{\langle s \rangle}$ by
the definition of $\sigma(s)$. Hence if $n \neq 2$ or $p$ odd the
lemma follows by Proposition~\ref{toralobstructionslemma}.

Now suppose that $X=\PU(2)\twocom$. Since for the non-trivial
elementary abelian $2$-subgroup $V
\subseteq \dT$ we have $W_{X}(V)_1=1$ and $\cC_{X}(V)_1
\cong T$, the first part of Proposition~\ref{toralobstructionslemma} still
applies to finish the proof also in this case. 
\end{proof}

We next record the following general lemma, which is obvious from the
K\"{u}nneth formula.

\begin{lemma} \label{joinlemma}
Suppose $\D_1$ and $\D_2$ are two categories with only finitely many
morphisms. Let $C\D_i$ be `the cone on $\D_i$' i.e., the category
constructed from $\D_i$ by adding an initial object $e$ to $\D_i$, and
let $\D_1 \star \D_2 = C\D_1 \times C\D_2 - \{(e,e)\}$, `the join
of $\D_1$ and $\D_2$' (see \cite[\S 1]{quillen78}). If $F_i: C\D_i
\to \Zppmod$, $i=1,2$ are functors then
$$C^*(C\D_1 \times C\D_2, \D_1 \star \D_2; F_1\otimes F_2) \cong
C^*(C\D_1,\D_1;F_1) \otimes C^*(C\D_2,\D_2;F_2).$$
In particular if one of the chain complexes has torsion free homology
or if everything is defined over $\F_p$ then
$$H^*(C\D_1 \times C\D_2, \D_1 \star \D_2; F_1\otimes F_2) \cong
H^*(C\D_1,\D_1;F_1) \otimes H^*(C\D_2,\D_2;F_2).$$\QED
\end{lemma}

The following result gives that certain filtration quotients have
(almost) vanishing cohomology.

\begin{thm} \label{punnontoralobs} 
Set $X = \PU(n)\pcom$ and fix an integer $r>0$. Let $F^r_j: \A(X) \to
\Zppmod$ denote the functor on objects given by $F^r_j(E) =
\pi_j(B\cZ\cC_{X}(E))$ if $E$ is of the form $\bar \Gamma_r \times \bar A$
(in the notation of Section~\ref{punsection}) and zero
otherwise. Writing $n = p^r k$ we have
\begin{eqnarray*}
\lim_{\A(X)}\hspace{-5pt}{}^iF^r_j =
\begin{cases}
\Z/2 & \text{if $j=i=k=r=1$ and $p=2$} \\
0 & \text{otherwise.}
\end{cases}
\end{eqnarray*}
\end{thm}

\begin{proof}
Define a functor $\tilde F^r_j$ by $\tilde F^r_j(E) =
\pi_j(BZC_{\PU(k)}(\bar A)\pcom)$ if $E = \bar \Gamma_r \times \bar A$
for a fixed $r$ and zero otherwise. This is a subfunctor of $F^r_j$
via the identification $\PU(k) \cong C_{\PU(n)}(\bar \Gamma_r)_1$. Set
$\tilde {\tilde F}^r_j = F^r_j/\tilde F^r_j$ and observe that this is
the trivial functor unless $j=1$ where it is given by $\tilde{\tilde
  F}^r_j(E) = \bar \Gamma_r$ if $E$ is of the form $\bar \Gamma_r
\times \bar A$ and zero otherwise.

Consider the category 
$$\D=  \A^e(X)_{\subseteq \bar \Gamma_r} \times \A^e(\PU(k)\pcom) - \{(e,e)\}$$
where the superscript $e$ means that we do not exclude the trivial subgroup. 
We have a natural inclusion of categories
$\iota: \D \to \A(X)$ on objects given by
$(V,\bar A) \mapsto V \times \bar A$.

{\em Step 1:}
We claim that this map induces an isomorphism
$$ \lim_{\A(X)}\hspace{-5pt}{}^*F^r_j \to \lim_{\D}{}^*F^r_j.$$
By filtering the functor and using Nakayama's lemma it is enough to
show this for $\tilde F^r_j \otimes \F_p$ and $\tilde {\tilde F}^r_j
\otimes \F_p$. We can furthermore replace these functors by the
subquotient functors which are only concentrated on one subgroup $\bar
\Gamma_r \times \bar A$.

Consider first such a subquotient of $\tilde F^r_j \otimes \F_p$. In
this case the formula of Oliver \cite[Prop.~4]{bob:steinberg} together
with Lemma~\ref{joinlemma} shows that the higher limits on both sides
are only non-zero in a single degree, where the map identifies with
the restriction map
\begin{multline*}
\Hom_{W_{X}(\bar \Gamma_r \times \bar A)}(\St(\bar \Gamma_r \times
\bar A), \pi_j(BZC_{\PU(k)}(\bar A)) \otimes \F_p) \longrightarrow \\
\Hom_{\Sp(\bar \Gamma_r) \times W_{\PU(k)}(\bar A)}(\St(\bar \Gamma_r)
\otimes \St(\bar A),\pi_j(BZC_{\PU(k)}(\bar A)) \otimes \F_p).
\end{multline*}
Let $U$ be the subgroup of elements in $W_{\PU(n)}(\bar \Gamma_r \times \bar
A)$ which act as the identity on $\bar A$ and $(\bar \Gamma_r \times \bar A)/\bar A$. Then $U$ acts trivially on $\pi_j(BZC_{\PU(k)}(\bar A)) \otimes \F_p$. Furthermore, by the
theorem of Smith \cite{smith82} (and self-duality of the
Steinberg module) $\St(\bar \Gamma_r \times \bar A)_U \cong
\St(\bar \Gamma_r \times \bar A)^U \cong \St(\bar \Gamma_r) \otimes
\St(\bar A)$, where $(-)_U$ and $(-)^U$ denote coinvariants and
invariants respectively. Hence this map is an isomorphism. The case of
a subquotient of $\tilde {\tilde F}^r_j \otimes \F_p$ is completely
analogous. This shows the isomorphism.

{\em Step 2:}
We now proceed to calculate the higher limits over $\D$, which we do
by calculating the limits of $\tilde F^r_j$ and $\tilde {\tilde F}^r_j$.
We first consider $\tilde{\tilde F}^r_j$. We have already remarked
that only $\tilde{\tilde F}^r_1 \neq 0$. Furthermore if $k \neq 1$
then
$H^*(C\A^{\operatorname{tor}}(\PU(k)\pcom),\A^{\operatorname{tor}}(\PU(k)\pcom);\Z_p)
= 0$ by \cite[\S 8]{DW93} so $\lim^*_\D \tilde {\tilde F}^r_1 = 0$ by
Lemma~\ref{joinlemma}. If $k=1$ we get $\lim^i_{\D} \tilde{\tilde F}^r_1
\cong \Hom_{\Sp(\bar \Gamma_r)}(\St(\bar \Gamma_r),\bar \Gamma_r)$ if
$i=2r-1$ and zero otherwise, by
\cite[Prop.~4]{bob:steinberg}.

Now consider $\tilde F^r_j$. By Lemma~\ref{joinlemma}
$$
\lim_\D{}^i \tilde F^r_j = \Hom_{\Sp(\bar \Gamma_r)}(\St(\bar
\Gamma_r),\Z_p) \otimes \lim_{\A^{\operatorname{tor}}(\PU(k)\pcom)}\hspace{-20pt}{}^{i-2r}
\pi_j(B\cZ\cC_{\PU(k)\pcom}(\bar A)).
$$
By Lemma~\ref{puntoralobs}
$\lim^{i-2r}_{\A^{\operatorname{tor}}(\PU(k)\pcom)}
\pi_j(B\cZ\cC_{\PU(k)\pcom}(\bar A)) = 0$ unless $p=k=2$, $j=1$ and $i-2r
= 1$ where we get $\lim^{1}_{\A^{\operatorname{tor}}(\PU(2)\twocom)}
\pi_j(B\cZ\cC_{\PU(2)\twocom}(\bar A)) = \Z/2$.

By an argument of H.\ H.\ Andersen and C.\ Stroppel \cite{henningbrev} we
have 
$$\Hom_{\Sp(\bar \Gamma_r)}(\St(\bar \Gamma_r),\Z_p) = 0$$
for all $r$ and $p$.

To sum up we get 
$$\lim_{\A(X)}\hspace{-5pt}{}^iF^r_j \cong \lim_{\D}{}^i \tilde{\tilde F}^r_j \cong \Hom_{\Sp(\bar \Gamma_r)}(\St(\bar \Gamma_r),\bar \Gamma_r)$$
if $j=k=1$ and $i=2r-1$ and zero otherwise. 

However the same argument of H.\ H.\ Andersen and C.\ Stroppel
\cite{henningbrev} shows that
$$\Hom_{\Sp(\bar \Gamma_r)}(\St(\bar \Gamma_r),\bar \Gamma_r) = 0$$
unless $r=1$ and $p=2$ where it equals $\F_2$. (Note that this is
obvious if $p$ is odd by Lemma~\ref{minusone}.) This shows the wanted
formula.
\end{proof}

\begin{rem}
Note that slightly non-trivial statements from \cite{henningbrev} are
only used above for $p=2$ and furthermore become trivial when $r=1$,
where $\Sp(\bar \Gamma_1) \cong \SL(\bar \Gamma_1)$, and this is in
fact the only case which involves obstruction groups in the range
needed for the proof of (the $p=2$ version of) Theorem~\ref{ndetthm}.
\end{rem}

\begin{rem}
It is in fact possible to give a short proof of Smith's theorem in the
special case used above, using the geometric definition of the
Steinberg module $\St(E)$ via flags.
\end{rem}

\begin{proof}[Proof of Theorem~\ref{obstructionsmainthm} for $X = \PU(n)\pcom$] 
Theorem~\ref{punnontoralobs} and Lemma~\ref{puntoralobs} directly show the
conclusion unless $n=2$ and $p=2$, so assume we are in this case.
From the definition it follows that $\lim^0_{\A(X)}F_1 =
0$. Now, consider the exact sequence of functors $0 \to F^1_1 \to F_1
\to F^{\operatorname{tor}}_1 \to 0$. The long exact sequence of higher
limits starts out as
$$ 0 \to 0 \to 0 \to
\lim_{\A(X)}\hspace{-5pt}{}^0F^{\operatorname{tor}}_1 \to
\lim_{\A(X)}\hspace{-5pt}{}^1 F_1^1 \to
\lim_{\A(X)}\hspace{-5pt}{}^1F_1 \to 0 \to \cdots.
$$ 
So, since $\lim^0_{\A(X)}F^{\operatorname{tor}}_1 \cong \Z/2 \cong
\lim^1_{\A(X)}F^1_1$ we get $\lim^i_{\A(X)}F_1 =0$ for $i>0$ as
well. This concludes the proof of this last case of
Theorem~\ref{obstructionsmainthm}.
\end{proof}

%%%%%%%%%%%%%%%%%%%%%%%%%%%%%%%%%%%%%%%%%%%%%%%%%%%%%%%%%%%%%%%%%%%%%%%%%%%

\section{Consequences of the main theorem} \label{consequencessection}
In this section we prove the theorems listed in the introduction which
are consequences of the main theorem.

\begin{proof}[Proof of Theorem~\ref{splittheory}]
The theorem follows directly from Theorem~\ref{classthm} together with
the classification of finite $\Z_p$-reflection groups
(Theorem~\ref{lattice-split}), noting that by the proof of
Theorem~\ref{classthm} (and Theorem~\ref{classicalborel}) all
exotic $p$-compact groups have torsion free $\Z_p$-cohomology.
\end{proof}

\begin{proof}[Proof of Theorem~\ref{homogeneous}]
By \cite[Thm.~1.4]{MN94} $X$ is isomorphic to a $p$-compact group of
the form $(X' \times T'')/A$, where $X'$ is a simply connected
$p$-compact group, $T''$ is a $p$-compact torus, and $A$ is a finite
central subgroup of the product. Hence we have $X/T \cong X'/T'$,
where $T$ and $T'$ are maximal tori of $X$ and $X'$ respectively. So
we can without restriction assume that $X$ is simply connected.

For compact connected Lie groups the statement of this theorem is the
celebrated result of Bott \cite[Thm.~A]{bott56}. Hence by
Theorem~\ref{splittheory} it is enough to prove the theorem when $X$
is an exotic $p$-compact group.
In that case $H^*(BX;\Z_p)$ is a polynomial algebra with generators in
even degrees, and the number of generators equals the rank of $X$ (by the
proof of Theorem~\ref{ndetthm}). The same is true over $\F_p$, and
since $H^*(BT;\F_p)$ is finitely generated over $H^*(BX;\F_p)$ by
\cite[Prop.~9.11]{DW94}, $H^*(BX;\F_p) \to H^*(BT;\F_p)$ is
injective by a Krull dimension consideration. But since they
are both polynomial algebras it follows by e.g., \cite[\S
11]{DW98elementary} that $H^*(BT;\F_p)$ is in fact free over
$H^*(BX;\F_p)$. Hence the Eilenberg-Moore spectral sequence of the
fibration $X/T \to BT \to BX$ collapses and
$$H^*(X/T;\F_p) \cong \F_p
\otimes_{H^*(BX;\F_p)}  H^*(BT;\F_p).$$
In particular $H^*(X/T;\F_p)$
is concentrated in even degrees so the rank equals the Euler
characteristic $\chi(X/T)$ which again equals $|W_X|$ by
\cite[Prop.~9.5]{DW94}. By the long exact sequence in cohomology and
Nakayama's lemma, $H^*(X/T;\Z_p)$ is a free $\Z_p$-module
of rank $|W_X|$ as wanted.
\end{proof}

\begin{rem}
Let $H^*_{\Q_p}(\cdot) = H^*(\cdot;\Z_p) \otimes \Q$.
For any connected $p$-compact group $X$ the natural map $X/T \to BT$
induces an isomorphism
$$H^*_{\Q_p}(X/T) \cong \Q_p \otimes_{H^*_{\Q_p}(BX)} H^*_{\Q_p}(BT)$$
since the Eilenberg-Moore spectral sequence of the fibration $X/T \to
BT \to BX$ collapses by \cite[Prop.~9.7]{DW94} and \cite[\S
11]{DW98elementary}. It then follows from \cite{chevalley55} that the
natural $W_X$-action on $H^*(X/T;\Z_p) \otimes \Q_p = H^*_{\Q_p}(X/T)$
is isomorphic to the regular representation of $W_X$ when ignoring the
grading. Just as for compact connected Lie groups this is not true over $\Z_p$ when 
$p \mid |W_X|$ (cf.\ e.g.\ \cite[p.~221]{smith95}).
\end{rem}

\begin{proof}[Proof of Theorem~\ref{peterweyl}]
By Theorem~\ref{splittheory} it is enough to prove the statement in
the case where $X$ is the $\F_p$-completion of a compact connected Lie
group and the case where $X$ is an exotic $p$-compact group
separately. The case where $X$ is the $\F_p$-completion of a compact
connected Lie group of course follows directly from the classical Peter-Weyl
theorem (cf.\ e.g.\ \cite[Thm.~III.4.1]{BT85}), so we can concentrate
on the case where $X$ is exotic. If $p$ does not divide the order of
the Weyl group the statement is also obvious: The inclusion
$\dT \to \U(r)$ induces a map $\dT \semi W \to
\U(r|W|)$ whose $\F_p$-completion is a monomorphism. The
remaining cases have been shown to have faithful representations by
Castellana: If $(W,L)$ is in family $2a$ then this is carried out
in \cite[Thm.~E]{castellana2a} and if $(W,L)$ is one of the pairs
$(G_{12},p=3)$, $(G_{29},p=5)$, $(G_{31},p=5)$, or $(G_{34},p=7)$ this
is carried out in \cite{castellanaexotic}.
\end{proof}

We now turn to Theorem~\ref{fundamentalgroup} which in fact follows
easily from the classification. But let us first state the part which
one can see by elementary means. (See also \cite[Cor.~5.6]{MN94} and
\cite[Lem.~9.3]{DW98elementary}.)

\begin{prop}[{cf. \cite[Lem.~6.11]{dw:split} and \cite[Lem.~9.3]{DW98elementary}}] \label{pioneprop}
Let $X$ be a connected $p$-compact group. Then the natural composite
map
$$(L_X)_W \cong H_0(W;H_2(BT;\Z_p)) \to H_2(BX;\Z_p) \cong \pi_1(X)$$
induced by the inclusion $T \to X$ is surjective with finite kernel.

In particular if $(L_X)_W$ is torsion free then it is an isomorphism. 
\end{prop}

\begin{proof}
By \cite[Thm.~1.4]{MN94} $X$ is isomorphic to a $p$-compact group of
the form $(X' \times T'')/A$, where $X'$ is a simply connected
$p$-compact group, $T''$ is a $p$-compact torus, and $A$ is a finite
central subgroup of the product. Since the center of a connected
$p$-compact group is contained in a maximal torus by
\cite[Thm.~7.5]{dw:center} we can assume that $A$ is a subgroup of $T'
\times T''$, where $T'$ is maximal torus for $X'$, and hence $(T'
\times T'')/A$ is a maximal torus for $X$.
 Therefore we get the following diagram of fibration sequences
$$
\xymatrix{
BA \ar[r] \ar@{=}[d] & BT' \times BT'' \ar[d] \ar[r] & B((T' \times T'')/A) \ar[d] \\
BA \ar[r] & BX' \times BT'' \ar[r] & BX
}
$$ 
The long exact sequence of homotopy groups and the five-lemma now
show that $\pi_2(B((T \times T'')/A)) \to \pi_2(BX)$ is surjective
which is the first statement in the proposition. To see that the
kernel is finite note that by \cite[Thm.~9.7(iii)]{DW94} $H_{\Q_p}^2(BX)
\to H^2_{\Q_p}(BT)^W$ is an isomorphism, which by dualizing to
homology shows the claim.

That we get an isomorphism when $(L_X)_W$ is torsion free is obvious
from the general statement.
\end{proof}

\begin{rem}\label{genpioneremark}
One easily shows that the image of the differential $d_3: H_3(W;\Z_p)
\to H_0(W;H_2(BT;\Z_p))$ in the Serre spectral sequence for the
fibration $BT \to B\N_X \to BW$ is always in the kernel of the
surjective map of Proposition~\ref{pioneprop}. By standard group
cohomology (cf.\ \cite{CV66}) the image of this differential identifies
with the image of the map given by capping with the $k$-invariant
$\gamma \in H^3(W;H_2(BT;\Z_p))$ of the extension. If one knew that
the double coset formula held for $p$-compact groups (more precisely
that $H^*(B\N;\Z_p) \xrightarrow{\operatorname{tr}^*} H^*(BX;\Z_p) \xrightarrow{\res} H^*(BT;\Z_p)$ is
the restriction map, cf.\ \cite[Ex.~VI.4]{feshbach79}) then it would
easily follow that this image is in fact equal to the kernel of the
map in Proposition~\ref{pioneprop}, which would give a conceptual proof of the
formula for the fundamental group.
Note that by a result of Tits
\cite{tits66} (see also \cite{DW01,neumann99,kksa:thesis}) the
extension class $\gamma$ is always of order $2$ for compact connected
Lie groups. The next proposition gives the complete answer in the Lie
case.
\end{rem}

\begin{prop}\label{genpioneremark2}
Let $G$ be a compact connected Lie group. Then the map $\pi_1(T)_W \to
\pi_1(G)$ is surjective with kernel $(\Z/2)^s$, where $s$ is the
number of direct factors of $G$ isomorphic to a symplectic group
$\Sp(n)$, $n\geq1$.
\end{prop}

\begin{proof}
That the map is surjective follows as in the $p$-compact case, so we
just have to identify the kernel.
By \cite[Thm.~1.6]{matthey02}, for any compact connected Lie group
$G$, $T^W = Z(G) \oplus (\Z/2)^s$, where $s$ is the number of direct
factors of $G$ isomorphic to an odd special orthogonal group
$\SO(2n+1)$, $n\geq 1$.

Consider the dual group $G^\vee$ of $G$ obtained as the taking the
dual root diagram (see \cite[\S 4, no.\ 8]{bo9}). Then $G^\vee$ has
fundamental group isomorphic to $\widehat{Z(G)}$, where the hat
denotes the Poincar{\'e} dual group (see \cite[\S 4, no.\ 9]{bo9}). Likewise
$\widehat{(L_{G{\,}\check{}})_W}$ is canonically isomorphic to
$T^W$. Since duality is an involution on the set of compact connected
Lie groups which sends direct summands to direct summands and $\SO(2n+1)$ to
$\Sp(n)$ the claim about the fundamental group follows directly from
the dual result about the center.
\end{proof}

\begin{proof}[Proof of Theorem~\ref{fundamentalgroup}]
By Theorem~\ref{lattice-split} $(L_X)_{W_X} = 0$ for all exotic
$p$-compact groups $X$, so Proposition~\ref{pioneprop} shows the formula
in this case. By Theorem~\ref{splittheory} we are hence reduced to
showing the formula for $X$ of the form $G\pcom$ for some compact connected
Lie group $G$. In this case the formula is well known and
easy. Namely it follows from Remark~\ref{genpioneremark} that the
kernel of $(L_X)_{W_X} \to \pi_1(X)$ is an elementary abelian
$2$-group. Alternatively the same conclusion follows from the formula
for the fundamental group of a compact connected Lie group (see \cite[\S 4,
no.\ 6, Prop.\ 11]{bo9} or \cite[Thm.~5.47]{adams69}, noting that in the
notation of \cite{adams69} $(1-\varphi_r)\gamma_r =
2\gamma_r$).
\end{proof}

We now start to prove Theorems~\ref{borelelementaryabelian} and
\ref{borelelementaryabelianII}.

\begin{lemma} \label{toralnproperty}
Suppose $X$ and $X'$ are two connected $p$-compact groups with the same
maximal torus normalizer $\N$. Then all elementary abelian
$p$-subgroups of $X$ are toral if and only if all elementary abelian
$p$-subgroups of $X'$ are toral.

Furthermore, if for all toral elementary abelian $p$-subgroups $V \to
X$ the centralizer $\cC_X(V)$ is connected then all elementary abelian
$p$-subgroups in $X$ are toral.
\end{lemma}

\begin{proof}
Suppose that $X$ has a non-toral elementary abelian
$p$-subgroup $V \to X$. We can assume that it is minimal, in the
sense that any elementary abelian $p$-subgroup of smaller rank is
toral. Write $V = V' \oplus V''$, where $V'$ has rank one. We can
assume that $V \to X$ factors through $\N$ (indeed through $\N_p$ by
\cite[Prop.~2.14]{dw:center}) and that the restriction to $V''$
factors through $T$ (by the minimality of $V$).

We want to show that the resulting map $V \to \N \to X'$ is also
non-toral, by proving that the adjoint $V' \to \cC_{X'}(V'')$ does not
factor through the identity component $\cC_{X'}(V'')_1$, which would
be the case if $V \to X'$ was toral. In detail, proceed as follows:
Let $\N''$ denote the maximal torus
normalizer in $\cC_X(V'')_1$, which by
\cite[Thm.~7.6(2)]{dw:center} can be described in terms of $V''$ and
$\N$.  The adjoint map $V' \to \cC_\N(V'')$ cannot factor through
$\N''$ since otherwise $V' \to \cC_X(V'')$ would factor through
  $\cC_X(V'')_1$ and hence be toral in
$\cC_X(V'')$ by
\cite[Prop.~5.6]{DW94}, contradicting that $V$ is
assumed to be non-toral.
Note that $\dN''$ is normal in $C_{\dN}(V'')$ and
$C_{\dN}(V'')/\dN'' \cong \pi_0(\cC_X(V'')) \cong
\pi_0(\cC_{X'}(V''))$ (see \cite[Rem.~2.11]{dw:center}). Hence $V' \to
\pi_0(\cC_{X'}(V''))$ is non-trivial so $V \to \N \to X'$ is non-toral in $X'$
  as desired.
The last part of the lemma is clear from the proof of the first part.
\end{proof}

\begin{rem}
Despite the above lemma it is not a priori clear how to determine
whether a $p$-compact group $X$ has the property that all elementary
abelian $p$-subgroups are toral just by
looking at $\N_X$ (but see \cite[Thm.~2.28]{steinberg75} for the Lie
case). However, by a case-by-case analysis
(Theorem~\ref{borelelementaryabelian}), this is the case if and only
if all toral elementary abelian $p$-subgroups have connected
centralizers.
\end{rem} 

\begin{rem}
Note that by Lannes theory \cite[Thm.~0.4]{lannes92} the property that
every elementary abelian $p$-subgroup of $X$ is toral is equivalent to
$H^*(BX;\F_p) \to H^*(BT;\F_p)^{W_X}$ being an $F$-isomorphism. (See
also Theorem~\ref{classicalborel} and Remark~\ref{classicalborelextra}.)
\end{rem}

We state the following well known lemma for easy reference.

\begin{lemma}\label{torsionfreetoral}
Suppose that $X$ is a connected $p$-compact group such that
$H^*(BX;\Z_p)$ is a polynomial algebra with generators concentrated in
even degrees. Then all elementary abelian $p$-subgroups of $X$ are toral.
\end{lemma}

\begin{proof}
For every elementary abelian $p$-subgroup $\nu: E \to X$,
$H^*(B\cC_X(\nu);\F_p)$ is a polynomial algebra with generators
concentrated in even degrees by \cite[Thm.~1.3]{DW98} (note that
Lannes $T$-functor preserves objects concentrated in even degrees by
\cite[Prop.~2.1.3]{lannes92}). In particular $\cC_X(\nu)$ is connected
so by Lemma~\ref{toralnproperty} all elementary abelian $p$-subgroups
are toral. (Alternatively one can use
Remark~\ref{classicalborelextra}.)
\end{proof}

\begin{proof}[Proof of Theorem~\ref{borelelementaryabelian}]
First note that the implications $\eqref{borel1}\Rightarrow
\eqref{borel2}$ and $\eqref{borel1}\Rightarrow \eqref{borel4}$ follow
from Theorem~\ref{classicalborel}. The implication \eqref{borel2}
$\Rightarrow$ \eqref{borel3} follows easily from
Theorem~\ref{nakajima}. Namely for all toral elementary abelian
$p$-subgroups $V \to X$, Theorem~\ref{nakajima} implies that
$W_{\cC_X(V)}$ is a reflection group, so by \cite[Thm.~7.6]{dw:center}
$\cC_X(V)$ is connected, using the assumption that $p$ is odd. But this
implies that all elementary abelian $p$-subgroups are toral by
Lemma~\ref{toralnproperty}.

We now prove the implication \eqref{borel3} $\Rightarrow$
\eqref{borel1}. First note that by Theorem~\ref{lattice-split} and
\cite[Thm.~1.4]{dw:split} we can write $X \cong X' \times X''$ where
$X'$ has Weyl group $(W_G, L_G \otimes \Z_p)$, for some compact
connected Lie group $G$, and $(W_{X''},L_{X''})$ is a product of exotic
finite $\Z_p$-reflection groups. Furthermore, since the normalizer of
a connected $p$-compact group is split for $p$ odd by
\cite[Thm.~1.2]{kksa:thesis} we have $\N_{G\pcom} \cong \N_{X'}$. Since by
Lemma~\ref{toralnproperty} the property of having all elementary
abelian $p$-subgroups toral is a property which only depends on $\N$
we conclude that $G$ has this property as well. But this
implies that $G$ has torsion free $\Z_p$-cohomology by
\cite[Thm.~B]{borel61} (see also \cite[Thm.~2.28]{steinberg75}). In
the exotic case, we know by the proof of Theorem~\ref{ndetthm} that
we can find a $p$-compact group $\tilde X''$ which has the same
maximal torus normalizer as $X''$ and which has torsion free
$\Z_p$-cohomology. Hence we have found a $p$-compact group $G\pcom
\times \tilde X''$ which has the same maximal torus normalizer as $X$
and has torsion free $\Z_p$-cohomology.
Since by Theorem~\ref{ndetthm} a $p$-compact group is determined by
its maximal torus normalizer we conclude that $X$ in fact has torsion free
$\Z_p$-cohomology. (Alternatively, one can appeal to
Remark~\ref{polyvspoly} which shows that the property of having
torsion free $\Z_p$-cohomology only depends on $\N$.)

Finally we prove the implication $\eqref{borel4} \Rightarrow
\eqref{borel1}$, where we seem to need the full strength of
Theorem~\ref{classthm}. Note that by Theorem~\ref{splittheory} we can
write $X \cong G\pcom \times X'$ where $G$ is a compact connected Lie
group and $X'$ has torsion free $\Z_p$-cohomology. Likewise if
$BG\pcom$ has torsion free $\Z_p$-cohomology then $G\pcom$ has torsion
free $\Z_p$-cohomology by \cite[p.~93]{borel67}.
\end{proof}

\begin{rem} \label{torsionremark}
We make some remarks on Theorem~\ref{borelelementaryabelian}.
Since the implication $\eqref{borel1} \Rightarrow \eqref{borel3}$
follows from Lemma~\ref{torsionfreetoral}, we see that
the implications $\eqref{borel1} \Rightarrow \eqref{borel2}$,
\eqref{borel1} $\Rightarrow$ \eqref{borel4}, and $\eqref{borel1}
\Rightarrow \eqref{borel3}$ follow by general arguments. In the case
of compact connected Lie groups the implication $\eqref{borel3} \Rightarrow
\eqref{borel2}$ has a general proof, by combining
\cite[Thm.~2.28]{steinberg75} with \cite{demazure73}, and likewise
\eqref{borel4} $\Rightarrow$ \eqref{borel1} has a general proof by
\cite[p.~93]{borel67}. We do not know non-case-by-case proofs of
these implications for $p$-compact groups. (The implication
$\eqref{borel4} \Rightarrow \eqref{borel1}$ is stated in
\cite[Thm.~4.2]{MN98} but the proof is incorrect.) The remaining
implications do not seem to have general proofs even for compact
connected Lie groups. See also \cite[\S 4]{steinberg75}, Remark~\ref{polyvspoly}, and Theorem~\ref{classicalborel}.
\end{rem}

\begin{proof}[Proof of Theorem~\ref{borelelementaryabelianII}]
By adjointness (cf.\ Construction~\ref{adjoint}) it is obvious that $(2) \Rightarrow (3)$ since a rank one elementary abelian $p$-subgroup of a connected $p$-compact group is toral by \cite[Prop.~5.6]{DW94}. We prove the implications $(3) \Rightarrow (2) \Rightarrow (1)$ by reducing them to theorems for Lie groups via the classification of finite $\Z_p$-reflection groups. For the proof that $(1)$ implies $(2)$ or $(3)$ we also have to rely on Theorem~\ref{fundamentalgroup}.

By Theorem~\ref{lattice-split} and \cite[Thm.~1.4]{dw:split} we can write $X
\cong X' \times X''$ where $X'$ has Weyl group $(W_G, L_G \otimes
\Z_p)$, for some compact connected Lie group $G$, and the Weyl group
$(W'',L'')$ of $X''$ is a product of exotic $\Z_p$-reflection
groups. Note that the each of the properties $(1)$, $(2)$ and $(3)$ holds for $X$ if and only if it holds for $X'$ and $X''$.

By Theorem~\ref{lattice-split} and Proposition~\ref{pioneprop} we have $\pi_1(X'')=0$ so $X''$ satisfies $(1)$. By Theorems~\ref{polyinvthm} and \ref{nakajima} combined with \cite[Thm.~7.6]{dw:center}, $(2)$ also holds for $X''$. So we can without restriction assume that $X \cong X'$.
Furthermore we have $\N_{G\pcom} \cong \N_{X}$, by \cite[Thm.~1.2]{kksa:thesis} since $p$ is odd. By \cite[Thm.~7.6]{dw:center}, $(2)$ is a property which only depends on the maximal torus normalizer. By Lemma~\ref{toralnproperty} and its proof, this is also the case for $(3)$. Since $(2) \Leftrightarrow (3)$ for compact connected Lie groups by \cite[Thm.~2.27]{steinberg75} we see that $(2)$  and $(3)$ are equivalent for $p$-compact groups as well. If $X$ satisfies $(2)$, then so does $G\pcom$ by the above and hence $\pi_1(G\pcom)$ is torsion free by \cite[Thm.~2.27]{steinberg75}. Combining Propositions~\ref{genpioneremark2} and \ref{pioneprop} then shows that $\pi_1(X)$ is torsion free which proves
$(2) \Rightarrow (1)$. Finally to see $(1) \Rightarrow (2)$, observe that, using Theorem~\ref{fundamentalgroup}, $\pi_1(X) = \pi_1(G\pcom)$, which reduces us to the Lie case, where it again follows from  \cite[Thm.~2.27]{steinberg75}.
\end{proof}

\begin{rem}
It is not hard to see that the conjectural classification
for $p=2$ implies that Theorem~\ref{borelelementaryabelianII}
and the implications \eqref{borel1} $\Leftrightarrow$
\eqref{borel4} $\Leftrightarrow$ \eqref{borel3} $\Rightarrow$
\eqref{borel2} in Theorem~\ref{borelelementaryabelian} also hold true
for $p=2$. However, in Theorem~\ref{borelelementaryabelian},
\eqref{borel2} is not equivalent to the other conditions. To see this
observe that $\Z_2[L_{\SO(2n+1)} \otimes \Z_2]^{W_{\SO(2n+1)}}$ is a
polynomial algebra (because this holds for $\Sp(n)$, which has the
same Weyl group) even though $\SO(2n+1)$ has $2$-torsion.
\end{rem}

\begin{rem}\label{polyvspoly}
Notbohm states his classification of connected $p$-compact groups with
$\Z_p[L]^W$ a polynomial algebra in the setup of spaces $BX$ with
polynomial cohomology (cf.\ \cite{notbohm98,notbohm99}). This means that his
uniqueness statement is a priori only uniqueness among $p$-compact
groups with torsion free $\Z_p$-cohomology (cf.\
Theorem~\ref{classicalborel}). We will here briefly sketch a direct
but case-by-case way (following a line of argument given in a special
case in \cite[Pf.~of~Thm.~5.3]{MN98}) to show that for a $p$-compact group
the property of having torsion free $\Z_p$-cohomology depends only on
$(W,L)$, which allows us to remove the extra assumption.

Assume that $X$ is a connected $p$-compact group, $p$ odd, such that
$\Z_p[L_X]^{W_X}$ is a polynomial algebra. We want to show that
$H^*(BX;\Z_p)$ is a polynomial algebra as well. By
Theorem~\ref{polyinvthm}\eqref{cover-lattice}, $(L_X)_{W_X}$ is torsion
free, so $\pi_1(X) = (L_X)_{W_X}$ by Proposition~\ref{pioneprop}. By
the Serre and Eilenberg-Moore spectral sequences $H^*(BX;\Z_p)$ is a polynomial
algebra if and only if $H^*(B(X\onecov);\Z_p)$ is a polynomial
algebra. Furthermore by construction $L_{X\onecov} = SL_X$ and by
Theorem~\ref{polyinvthm}\eqref{cover-lattice}
$\Z_p[L_{X\onecov}]^{W_{X}}$ is also a polynomial algebra, so we can
without loss of generality assume that $X$ is simply connected. By
\cite[Thm.~1.4 and Rem.~1.6]{dw:split} we can furthermore assume that
$X$ is a simple $p$-compact group.

By \cite[Thm.~1.2]{kksa:thesis} $\dN_X = \dT_{X} \semi W_X$. Using
Theorems~\ref{lattice-split} and \ref{polyinvthm} we first show
that the cohomology of $\dN_X$ is detected by elementary abelian
$p$-subgroups. More precisely we show that in each case there is a
compact connected Lie group $H$ such that $\dN_X$ contains a subgroup
isomorphic to $\dN_{H\pcom}$ with index prime to $p$ having the
required property. When $p\nmid |W|$ we take $H$ to be a torus, and if
$(W_X,L_X)$ is in family $1$ or family $2a$ we take $H=\SU(n)$ and
$H=\U(n)$ respectively. If $(W_X,L_X)$ is one of the exotic
$\Z_p$-reflection groups $(G_{12},p=3)$, $(G_{29},p=5)$,
$(G_{31},p=5)$, or $(G_{34},p=7)$ we take $H=\SU(p)$, cf.\ the proof
of Theorem~\ref{nakajima}. The only remaining cases are the ones where
$(W_X,L_X)=(W_G,L_G\otimes \Z_p)$ for one of the following pairs
$(G,p)$: $(G_2,p=3)$, $(3E_6,p=5)$, $(2E_7,p=5)$, $(2E_7,p=7)$, and
$(E_8,p=7)$. In these cases we can by \cite[Prop.~6.11]{JMO92} take
$H=\SU(3)$, $\SU(2)\times_{C_2} \SU(6)$, $\SU(8)/C_2$, $\SU(8)/C_2$
and $\SU(9)/C_3$ respectively. Since both $\dN_{\U(n)\pcom}$ and
$\dN_{\SU(n)\pcom}$ have cohomology which is detected by elementary
abelian $p$-subgroups by \cite[Prop.~3.4]{quillen71} (for
$\N_{\U(n)\pcom}$; $\N_{\SU(n)\pcom}$ follows from this, cf.\
\cite[Lem.~12.6]{notbohm94}) we see that in all cases the cohomology
of $\dN_{H\pcom}$ is detected by elementary abelian $p$-subgroups.
Hence, by a transfer argument, the mod $p$ cohomology of $BX$ is
detected by elementary abelian $p$-subgroups.

Next, we want to show that all elementary abelian $p$-subgroups of $X$
factor through a maximal torus. By Lemma~\ref{toralnproperty}
we just have to show that we can find {\em some} $p$-compact group
$X'$ with the same maximal torus normalizer which has this
property. If $(W_X,L_X)$ is of Lie type this follows by combining Theorem~\ref{polyinvthm}\eqref{lattice-simplycon} with Borel's theorem
\cite[Thm.~B]{borel61}. If $(W_X,L_X)$ is exotic this is also true since
we know (by Theorem~\ref{polynomialrealizationthm} or Notbohm's work
\cite{notbohm99}) that there exists a $p$-compact group with Weyl group
$(W_X,L_X)$ and classifying space having polynomial $\Z_p$-cohomology
algebra.

The fact that all elementary abelian $p$-subgroups of $X$ are toral
combined with the fact that the cohomology is detected by elementary
abelian $p$-subgroups implies that the mod $p$ cohomology of $BX$ is
concentrated in even degrees. Hence $H^*(BX;\Z_p)$ is torsion free as
wanted.
\end{rem}

\begin{proof}[Proof of Theorem~\ref{Tconj}]
Let $X$ be a connected finite loop space with maximal torus $i: T \to X$.
Note that $(X/T)\pcom \simeq X\pcom/T\pcom$ by the fiber lemma
\cite[II.5.1]{bk}, and consequently, by the definition of Euler
characteristic, $\chi(X/T) = \chi(X\pcom/T\pcom)$. Hence $T\pcom \to
X\pcom$ will be a maximal torus for the $p$-compact group
$X\pcom$, for all primes $p$.

For our connected finite loop space $X$, define $W_X(T)$ to be the set
of conjugacy classes of self-equivalences $\varphi$ of $T$ such that
$i$ and $i\varphi$ are conjugate. We obviously have an injective
homomorphism $W_X \to W_{X\pcom}$ for all primes $p$ and we now want
to see that this map is surjective as well, so that we can naturally
identify $(W_X, \pi_1(T) \otimes \Z_p)$ with
$(W_{X\pcom},L_{X\pcom})$.
 
First note that by \cite[Pf.~of~Thm.~9.7]{DW94} we can view
$W_{X\pcom}$ as the Galois group of the extension of polynomial
algebras $H^*_{\Q_p}(BX) \to H^*_{\Q_p}(BT)$. But, since $BX$ has
finitely many cells in each dimension and since $BX$ is nilpotent, we
can identify
$$
\xymatrix{
H^*(BX;\Q) \otimes_\Q \Q_p \ar[r]  \ar@{-}[d]^-\cong & H^*(BT;\Q)
\otimes_\Q \Q_p  \ar@{-}[d]^-{\cong}\\
H^*_{\Q_p}(BX) \ar[r]  & H^*_{\Q_p}(BT)
}
$$
so the extensions $H^*(BX;\Q) \to H^*(BT;\Q)$ and $H^*_{\Q_p}(BX) \to
H^*_{\Q_p}(BT)$ have canonically isomorphic Galois groups. Hence any
element in $W_{X\pcom}$ lifts to a canonical element in the Galois
group of the extension $H^*(BX;\Q) \to H^*(BT;\Q)$. However, since
$BX_\Q$ and $BT_\Q$ are products of Eilenberg-Mac~Lane spaces
(cf.\ e.g.\ \cite[Ch.~V, \S 4, Prop.~6]{serre53}), this Galois group just
identifies with the self-equivalences $BT_\Q \to BT_\Q$ over $Bi_\Q: BT_\Q \to BX_\Q$,
where as usual $Bi_\Q$ has been replaced by an equivalent fibration. Hence any element in $W_{X\pcom}$ gives rise to a compatible
family of self-equivalences of $BT_\Q$ and $BT\lcom$,
for all primes $l$. So by the
arithmetic square \cite[VI.8.1]{bk}, we get a self-equivalence of $BT$
over $BX$, i.e., an
element in $W_X$. The constructed element is a lift of the element in
$W_{X\pcom}$ we started with, so the map $W_X \to W_{X\pcom}$ is
surjective as well.

Likewise, the argument above showed that $W_X$ is the Galois group of
the extension $H^*(BX;\Q) \to H^*(BT;\Q)$, so we have an isomorphism
$$H^*(BX;\Q) \xrightarrow{\cong} H^*(BT;\Q)^{W_X}.$$
Since $H^*(BX;\Q)$ is a polynomial algebra, $(W_X,\pi_1(T))$ is a $\Z$-reflection group
by the Shephard-Todd-Chevalley theorem (see \cite[Thm.~7.2.1]{benson93}). Hence, by
Theorem~\ref{lattice-split}, $(W_X,\pi_1(T))$ is the Weyl group of
some compact connected Lie group $G$.

For each $p$ we have an extension class $\gamma_p \in H^3(W_X;\pi_1(T)
\otimes \Z_p)$ corresponding to the fibration sequence $BT\pcom \to
B\N_{X\pcom} \to BW_X$. Since $H^3(W_X;\pi_1(T))$ is a finite abelian
group, and hence given as a sum of its $p$-primary parts, these
extension classes identify with a unique extension class
$\gamma \in H^3(W_X;\pi_1(T))$. We define the loop space $\N_X$ to be
the loop space of the total space in the fibration sequence $BT \to
B\N_X \to BW_X$ with the canonical action of $W_X$ on $BT$ and
extension class $\gamma$. Since the fiber-wise $\F_p$-completion of
$B\N_X$ with respect to this defining fibration identifies with
$B\N_{X\pcom}$, the arithmetic square produces a canonical morphism
$\N_X \to X$. ($\N_X$ is, quite naturally, called the maximal torus
normalizer of the finite loop space $X$ \cite[Def.~1.3]{MN98}.)

By \cite{kksa:thesis} (see also \cite{neumann99}, \cite{DW01}) the
extension classes defining $T \to \N_X \to W_X$ and $T \to N_G(T) \to
W_G(T)$ are both $2$-torsion. Let $\widetilde{BN}$ denote the
fiber-wise $\Z[\frac12]$-localization of the total space of the
fibration $BT \to BN_G(T) \to BW_G(T)$ or equivalently the
corresponding fibration with $B\N_X$. We hence have embeddings
$$
\xymatrix{
 & {\widetilde{BN}} \ar[dr] \ar[dl]\\
BX[\frac12] & & BG[\frac12]
}
$$
By the arithmetic square \cite[VI.8.1]{bk}, the following square is a pullback
$$
\xymatrix{
BX[\frac{1}{2}] \ar[r] \ar[d] & {\prod_{p \neq 2}} {BX\pcom} \ar[d]\\ 
BX_{\Q} \ar[r]& (\prod_{p \neq 2} BX\pcom)_\Q
}
$$
and similarly for $BG$.
By Theorem~\ref{ndetthm} we can construct unique maps between
$\F_p$-completions under $\widetilde{BN}$, and we obviously also have a
unique map between the rationalizations. By construction (as maps
under $\widetilde{BN}$) these maps agree on the rationalization of the
product of the $\F_p$-completions, so by the arithmetic square we get an
induced map $BX[\frac{1}{2}] \to BG[\frac{1}{2}]$, which by
construction is an $\F_p$-equivalence for all primes $p$. Since both
spaces are one-connected this implies that the map is a homotopy
equivalence.
\end{proof}

%%%%%%%%%%%%%%%%%%%%%%%%%%%%%%%%%%%%%%%%%%%%%%%%%%%%%%%%%%%%%%%%%%%%%%%%%%%

\section{Appendix: The classification of finite $\Z_p$-reflection
groups} \label{latticesection}

The purpose of this appendix is to give a short proof of the classification of finite
$\Z_p$-reflection groups (Theorem~\ref{lattice-split}), simplifying work of Notbohm
\cite{notbohm96,notbohm99lattice}, and as a by-product extending his results to all primes. 
We likewise explain how this classification relates to the
classification of finite $\Z$-reflection groups and the classification
of compact Lie groups (Theorem~\ref{liegroup-lattice}).

We start by recalling some definitions. Let $R$ be an integral domain
with field of fractions $K$. An {\em $R$-reflection group} is a pair
$(W,L)$ where $L$ is a finitely generated free $R$-module, and $W$ is
a subgroup of $\Aut(L)$ generated by elements $\alpha$ such
that $1-\alpha$ has rank one viewed as a matrix over $K$. Two finite
$R$-reflection groups $(W,L)$ and $(W',L')$ are called
{\em isomorphic}, if we can find an $R$-linear isomorphism $\varphi: L
\rightarrow L'$ such that the group $\varphi W \varphi^{-1}$ equals $W'$. A
finite $R$-reflection group $(W,L)$ is said to be {\em irreducible} if
the corresponding representation of $W$ on $L \otimes_R K$ is
irreducible. If $R$ has characteristic zero we define the
{\em character field} of an $R$-reflection group $(W,L)$ as the field
extension of $\Q$ generated by the values of the character of the
representation $W\hookrightarrow \Aut(L)$. For $R=\Z_p$ or $\Q_p$ we define
an {\em exotic} $R$-reflection group to be a finite irreducible
$R$-reflection group with character field strictly larger than $\Q$.

The classification of finite $\Z_p$-reflection groups is based on the
work of Clark-Ewing \cite{CE74} and Dwyer-Miller-Wilkerson \cite{DMW92},
which is again based on the classification of finite $\C$-reflection
groups by Shephard-Todd \cite{ST54} (see also \cite{cohen76}). The
result of Clark-Ewing and Dwyer-Miller-Wilkerson is that there is a bijection between finite $\Q_p$-reflection groups and finite
$\C$-reflection groups whose character field embeds in
$\Q_p$ (for details see \cite[Prop.~5.4, Prop.~5.5 and Pf.~of~Thm.~1.5]{DMW92}). The classification of finite complex reflection groups by Shephard-Todd
\cite{ST54} is as follows: Up to isomorphism, the irreducible finite complex reflection groups fall into $3$ infinite families
and $34$ sporadic cases. We follow the notation of Shephard-Todd and label the three infinite families as $1$, $2$ and $3$ and the sporadic cases as $G_i$, $4\leq i\leq 37$. Moreover any finite complex reflection group can be written as a direct product of irreducible finite complex reflection groups,
cf.\ \cite[Rem.~2.3]{DW00} (in fact this holds over any field of
characteristic $0$).

It is convenient to split family $2$ further
depending on the character field. The associated complex reflection
group is the group $G(m,r,n)$ (where $m$, $r$ and $n$ are integers
with $m,n\geq 2$, $r\geq 1$, $r\!\mid\!m$ and $(m,r,n)\neq (2,2,2)$) from
\cite[p.~277]{ST54} which consists of monomial $n\times n$-matrices
such that the non-zero entries are $m$th roots of
unity and the product of the non-zero entries is an $(m/r)$th root of
unity. Thus $G(m,r,n)$ is the semidirect product of its subgroup
$A(m,r,n)$ of diagonal matrices with the subgroup of permutation matrices.
Let $\zeta_m=e^{2\pi i/m}$. For $n\geq 3$ or $n=2$ and $r\neq m$ the
character field of $G(m,r,n)$ equals $\Q(\zeta_m)$, and for $n=2$ and
$r=m$ it equals $\Q(\zeta_m+\zeta_m^{-1})$ (see
\cite[p.~432--433]{CE74}). Following \cite{CE74} these are denoted as family $2a$ and family $2b$ respectively.

A complete list of the irreducible finite complex reflection groups,
their character fields and the primes for which these embed in $\Q_p$
can be found in \cite[p.~165]{kane88}, \cite[Table~7.1]{benson93}, or
\cite[Table~1]{kksa:thesis}.

If $(W,V)$ is a finite $\Q_p$-reflection group, then by
\cite[Prop.~23.16]{CR81} we can find a (non-unique) finitely generated
$\Z_pW$-submodule $L\subseteq V$ with $L\otimes \Q=V$. Thus any finite
$\Q_p$-reflection group may be obtained from a finite
$\Z_p$-reflection group by extension of scalars, but in general there
are several non-isomorphic $\Z_p$-reflection groups which give rise to
the same $\Q_p$-reflection group. The following result extends
\cite[Thm.~1.5 and Prop.~1.6]{notbohm99lattice} to all primes. (See also the addendum
Theorem~\ref{liegroup-lattice} for an elaboration.)

\begin{thm}[The classification of finite $\Z_p$-reflection groups]
\label{lattice-split}
Let $(W,L)$ be a finite $\Z_p$-reflection group. Then there exists a
decomposition
$$
(W,L) = (W_1 \times W_2, L_1 \oplus L_2)
$$
where $(W_1,L_1) \cong (W_G, L_G \otimes \Z_p)$, for some (non-unique)
compact connected Lie group $G$ with Weyl group $W_G$ and
integral lattice $L_G$, and
$(W_2,L_2)$ is a (up to permutation unique) direct product of exotic
$\Z_p$-reflection groups.

The canonical map $(W,L) \mapsto (W,L \otimes \Q)$ gives a one-to-one
correspondence between exotic $\Z_p$-reflection groups up to
isomorphism and exotic $\Q_p$-reflection groups up to isomorphism.

If $(W,L)$ is any exotic $\Z_p$-reflection group, then $L\otimes \F_p$
is an irreducible $\F_pW$-module, and in particular we have $(L\otimes
\Z/p^\infty)^W=0$ and $H_0(W;L)=0$.
\end{thm}

\begin{rem}
For odd primes $p$ the last two statements say by definition that any exotic
$\Z_p$-reflection group is respectively {\em center-free} and
{\em simply connected}, cf.\ \cite{notbohm96}.

Note also that \cite[Thm.~1.5]{notbohm99lattice} imposes the
unnecessarily strong condition that the invariant ring $\Z_p[L]^W$ is
a polynomial algebra, but this condition is not actually used in
\cite{notbohm99lattice}.
\end{rem}

Before the proof of Theorem~\ref{lattice-split} we need a couple of lemmas. First recall the following elementary fact about elements of finite order
in $\GL_n(\Z_p)$.

\begin{lemma} \label{modpreduction}
Let $G\subseteq \GL_n(\Z_p)$ be a finite subgroup. Then the
mod $p$ reduction $G\hookrightarrow \GL_n(\Z_p) \rightarrow
\GL_n(\F_p)$ is injective if $p$ is odd. For $p=2$ the kernel of the
composition is an elementary abelian $2$-subgroup of rank at most
$n^2$. In particular the kernel is contained in $O_2(G)$, the largest
normal $2$-subgroup of $G$.
\end{lemma}

\begin{proof}
It is easy to see directly that any non-trivial finite order element
in $\GL_n(\Z_p)$ has non-trivial reduction mod $p$ if $p$ is odd
(cf.\ \cite[Pf.~of~Lem.~10.7.1]{smith95}). For $p=2$ the same argument
shows that this is true if we reduce mod $4$. The result now follows.
\end{proof}

\begin{lemma} \label{exoticlattice}
For any exotic $\Q_p$-reflection group
$(W,V)$ there exists a finitely generated $\Z_pW$-submodule $L\subseteq
V$ with $L\otimes \Q=V$, such that $L\otimes \F_p$ is an irreducible
$\F_pW$-module.
\end{lemma}

\begin{proof}
Assume first that $p\nmid |W|$. By \cite[Prop.~23.16]{CR81} we can
find a finitely generated $\Z_pW$-submodule $L\subseteq V$ with
$L\otimes \Q=V$. It follows from \cite[75.6 and 76.15]{CR62} that
$L\otimes \F_p$ is automatically an irreducible $\F_pW$-module.

Assume now that $W$ has order divisible by $p$. From
Clark-Ewing's list we see that the only exotic $\Q_p$-reflection
groups satisfying this condition are the groups $G(m,r,n)$ from family
$2a$ and the groups $G_{12}$ for $p=3$, $G_{24}$ for $p=2$,
$G_{29}$ and $G_{31}$ for $p=5$ or $G_{34}$ for $p=7$.

In case $W=G(m,r,n)$ from family $2a$ we get the extra conditions $m\geq 3$,
$p\equiv 1 \pmod{m}$ and $p\leq n$. Note in particular that $n\geq
3$. The description above directly gives a representation with entries
in $\Z_p$ since the multiplicative group of $\Z_p$ contains the
$(p-1)$th roots of unity. Let $L=(\Z_p)^n$ be the natural
$\Z_pW$-module, i.e., the set of columns with entries in $\Z_p$.
Assume that $0\neq M\subseteq L\otimes \F_p$ is a
$\F_pW$-submodule of $L\otimes \F_p$. Choose $x\in M$ with $x\neq
0$ and let $\theta\in \F_p$ be a primitive $m$th root of unity. Since
$W$ contains the permutation matrices and the diagonal matrix
$\diag(\theta,\theta^{-1},1,\ldots,1)$ we see that $M$ contains an
element of the form $x'=(x_1,x_2,0,\ldots,0)$ with
$x_1\neq 0$. Since $n\geq 3$, $W$ also contains the diagonal matrix
$\diag(\theta,1,\theta^{-1},1,\ldots,1)$ and hence $M$ contains
$((1-\theta) x_1,0,\ldots,0)$. As $\theta\neq 1$ and $W$
contains all permutation matrices we conclude that $M=L\otimes \F_p$,
proving the claim for the groups from family $2a$.

Next consider $W=G_{12}$ at $p=3$. Since $W$ is isomorphic to
$\GL_2(\F_3)$, Lemma~\ref{modpreduction} shows that for any
finitely generated $\Z_3W$-submodule $L\subseteq (\Q_3)^2$ of rank $2$, we may
identify $L\otimes \F_3$ with the natural $\F_3\GL_2(\F_3)$-module. In
particular $L\otimes \F_3$ is an irreducible $\F_3W$-module.

For $W=G_{24}$ at $p=2$ we have $W\cong \Z/2 \times
\GL_3(\F_2)$. Hence Lemma~\ref{modpreduction} shows that for any
finitely generated $\Z_2W$-submodule $L\subseteq (\Q_2)^3$ of rank $3$, we may
identify $L\otimes \F_2$ with the $\F_2(\Z/2\times \GL_3(\F_2))$-module
where $\Z/2$ acts trivially and $\GL_3(\F_2)$ acts naturally. In
particular $L\otimes \F_2$ is an irreducible $\F_2W$-module.

Next consider the groups $G_{29}$ and $G_{31}$ at $p=5$. Since
$G_{29}$ is contained in $G_{31}$ it suffices to show the result for
$W=G_{29}$. The representation in \cite[p.~298]{ST54} is defined over
$\Z[\frac{1}{2},i]$ and hence we get a representation over $\Z_5$ by
mapping $i$ to a primitive $4$th root of unity in $\Z_5$. Let
$L=(\Z_5)^4$ be the natural $\Z_5W$-module. There are $40$ reflections
in $G_{29}$: The $24$ reflections in the hyperplanes of the form
$x_j-i^\alpha x_k=0$, $j\neq k$ and the $16$ reflections in the hyperplanes
of the form $\sum_{j=1}^4 i^{\alpha_j} x_j =0$ with $\sum_{j=1}^4
\alpha_j \equiv 0 \pmod{4}$. In particular $G_{29}$ contains the
reflections in the hyperplanes $x_j-x_k=0$ and thus $G_{29}$ contains
all permutation matrices. The product of the reflections in the
hyperplanes $x_1-i x_2=0$ and $x_1-x_2=0$ equals the diagonal matrix
$\diag(i,-i,1,1)$ and thus this element is also contained in
$G_{29}$. Now the same argument used in the case of the groups from
family $2a$ shows that $L\otimes \F_5$ is an irreducible $\F_5W$-module.

The argument for the group $W=G_{34}$ at $p=7$ is similar. The
representation given in \cite[p.~298]{ST54} is defined over
$\Z[\frac{1}{3},\omega]$, $\omega=\zeta_3$ and hence we get a
representation over $\Z_7$ by mapping $\omega$ to a primitive $3$rd
root of unity in $\Z_7$. Let $L=(\Z_7)^6$ be the natural
$\Z_7W$-module. There are $126$ reflections in $G_{34}$: The $45$
reflections in the hyperplanes of the form $x_j-\omega^\alpha x_k=0$,
$j\neq k$ and the $81$ reflections in the hyperplanes of the form
$\sum_{j=1}^6 \omega^{\alpha_j} x_j =0$ with $\sum_{j=1}^6 \alpha_j
\equiv 0 \pmod{3}$. In particular $G_{34}$ contains all permutation
matrices. The product of the reflections in the hyperplanes
$x_1-\omega x_2=0$ and $x_1-x_2=0$ equals the diagonal matrix
$\diag(\omega,\omega^2,1,1,1,1)$ and thus this element is also
contained in $G_{34}$. As above we then see that $L\otimes \F_7$ is an
irreducible $\F_7W$-module.
\end{proof}

\begin{proof}[Proof of Theorem~\ref{lattice-split}]
Assume first that $(W,V)$ is an exotic $\Q_p$-reflection group. Lemma~\ref{exoticlattice} shows that there exists a finitely generated
$\Z_pW$-submodule $L\subseteq V$ with $L\otimes \Q=V$, such that
$L\otimes \F_p$ is an irreducible $\F_pW$-module. It then follows from
\cite[15.2, Ex. 3]{serre77} that $L$ is unique up to a homothety (i.e. up to scaling by a unit in $\Q_p$). This gives the bijection between exotic $\Z_p$-reflection groups and exotic $\Q_p$-reflection groups.

Since $L\otimes \F_p$ is an irreducible $\F_pW$-module we also
conclude that $(L\otimes \F_p)^W=0$ and $H_0(W;L\otimes \F_p)=0$.
Hence we get $(L\otimes \Z/p^\infty)^W=0$ as claimed. We also see that
multiplication by $p$ is surjective on $H_0(W;L)$ and from this we
obtain $H_0(W;L)=0$ by Nakayama's lemma. This proves the part of
the theorem pertaining to exotic $\Z_p$-reflection groups.

Now consider a finite $\Z_p$-reflection group $(W,L)$ such that
we have a direct sum decomposition $L\otimes \Q = V_1 \oplus V_2$ as
$\Q_pW$-modules. Let $W_1$ (resp.\ $W_2$) be the subgroup of $W$ which
fixes $V_2$ (resp.\ $V_1$) pointwise. It it easy to see
(cf.\ \cite[Lem.~6.3]{dw:split}) that $(W_i,V_i)$ is a
$\Q_p$-reflection group and that we get the decomposition
$(W,L\otimes \Q) = (W_1 \times W_2, V_1 \oplus V_2)$.

We now claim that if $(W_2,V_2)$ is an exotic $\Q_p$-reflection group,
then we have the decomposition $(W,L) = (W_1 \times W_2, L_1
\oplus L_2)$ with $L_i=L\cap V_i$. Let $\alpha: L_1\oplus L_2
\longrightarrow L$ be the addition map. As in \cite[Pf.~of~Thm.~1.5]{dw:split} it suffices to prove that $\alpha\otimes \Z/p^\infty :
(L_1\otimes \Z/p^\infty) \oplus (L_2\otimes \Z/p^\infty)
\longrightarrow L\otimes \Z/p^\infty$ is injective. Assume that
$(x_1,x_2)$ is in the kernel of $\alpha\otimes \Z/p^\infty$, $x_i\in
L_i\otimes \Z/p^\infty$. Thus $x_1+x_2=0$. If $s\in W_2$ is a
reflection we have by definition $s\cdot x_1=x_1$ and hence $s$ also
fixes $x_2=-x_1$. Since $W_2$ is generated by reflections we get
$x_2\in (L_2\otimes \Z/p^\infty)^W$ and hence $x_2=0$ by the results
already proved for exotic $\Z_p$-reflection groups. Hence $x_1=0$ as
well, and thus $\alpha\otimes \Z/p^\infty$ is injective proving the
claim.

Since any finite $\Q_p$-reflection group may be decomposed into a (up
to permutation unique) product of finite irreducible ones, we see by
using the claim repeatedly that any finite $\Z_p$-reflection group
$(W,L)$ may be decomposed as a product $(W,L) \cong (W_1 \times W_2,
L_1 \oplus L_2)$ where $(W_1,L_1)$ is a $\Z_p$-reflection group with
character field equal to $\Q$ and $(W_2,L_2)$ is as in the theorem.

To finish the proof we thus need to show that for any finite
$\Z_p$-reflection group $(W,L)$ with character field equal to $\Q$ we
may find a compact connected Lie group $G$ such that $(W,L)$ is isomorphic to
$(W_G,L_G \otimes \Z_p)$. We start by reducing the problem to finite
$\Z$-reflection groups. The representation $W\rightarrow
\GL(L\otimes \Q)$ is a reflection representation and hence has Schur
index $1$ by \cite[Cor.~p.~429]{CE74}. Thus this representation
is equivalent to a representation defined over $\Q$. Hence
\cite[Cor.~30.10]{CR81} applied to $R=\Z_{(p)}$ shows
that there exists a (unique) finitely generated $\Z_{(p)}W$-submodule
$L'\subseteq L$ with $L'\otimes_{\Z_{(p)}} \Z_p=L$. Now \cite[Cor.~23.14]{CR81}
applied to $R=\Z$ shows that $L'$ contains a (non-unique)
finitely generated ${\Z}W$-submodule $L''\subseteq L'$ with $L'=L'' \otimes
\Z_{(p)}$. We conclude in particular that $(W,L) \cong (W,L''\otimes
\Z_p)$.

We finish the proof by showing that there exists a (non-unique)
compact connected Lie group $G$ whose Weyl group $(W_G,L_G)$ is
isomorphic to $(W,L'')$. For each reflection $s\in W$ the group
$\{x\in L'' \mid s(x)=-x\}$ is an infinite cyclic group with two
generators which we label $\pm \alpha_s$. Let $\Phi = \{ \pm \alpha_s
\mid s \mbox{ is a reflection in } W \}$ and $L''_0=(L'')^W$. It then
follows (cf.\ \cite[p.~85]{osse97}) that $(L'',L''_0,\Phi)$ is a
reduced root diagram whose associated $\Z$-reflection group equals
$(W,L'')$ (see \cite[\S 4, no.\ 8]{bo9} for definitions). From the
classification of compact connected Lie groups (\cite[\S 4, no.\ 9,
Prop.\ 16]{bo9}) it then follows that there exists a compact connected
Lie group $G$ whose root diagram equals $(L'',L''_0,\Phi)$. In
particular $(W_G,L_G)$ is isomorphic to $(W,L'')$ and we are done.
\end{proof}

We now analyze concretely when two
compact connected Lie groups give rise to the same $p$-compact group. For a
compact connected Lie group $G$, let $G\onecov$ denote the universal
cover of $G$. Furthermore let $H$ be the direct product of the
identity component of the center, $Z(G)_1$, with the universal cover
of the derived group of $G$. We have a canonical covering homomorphism
$\varphi: H \to G$ with finite kernel
(cf.\ \cite[\S 1, no.\ 4, Prop.\ 4]{bo9}). If $p$ is a prime number, we let
$\covpp{G}$ denote the covering of $G$ corresponding to the subgroup of
$\pi_1(G)$ given as the preimage of the Sylow $p$-subgroup of
$\pi_1(G)/\varphi(\pi_1(H))$, and let $K$ denote the kernel of $H \to
\covpp{G}$. Write $H = R_1 \times \cdots \times R_n \times S \times T'$ where
each $R_i$ is a special unitary group, $S$ is a simply connected
compact Lie group which contains no
direct factors isomorphic to a special unitary group, and $T'$ is a torus.
Suppose that $K_1$ and $K_2$ are finite central $p$-subgroups of
$H$. We say that $K_1$ and $K_2$ are {\em $p$-equivalent}
subgroups of $H$ if there exists
integers $k_1, \ldots, k_n, k$ prime to $p$ such that the
homomorphism $\Psi =\psi^{k_1} \times \cdots \times
\psi^{k_n} \times 1 \times \alpha_k: T_{R_1} \times \cdots T_{R_n}
\times T_S \times T' \to T_{R_1}
\times \cdots T_{R_n} \times T_S \times T'$ induces an isomorphism from $K_1$
onto $K_2$, where $T_{R_i}$ is a maximal torus of $R_i$, $\psi^l$
is the $l$th power map, and $\alpha_k: T' \to T'$ is a homomorphism
which with respect to some splitting $T' = S^1 \times \cdots \times
S^1$ has the form $1\times\ldots\times 1\times\psi^k$.
More generally we say that $H/K$
and $H'/K'$ are {\em $p$-equivalent} if there
exists an isomorphism between $H$ and
$H'$ such that image of $K$ in $H'$ is $p$-equivalent to
$K'$. This terminology is justified by the following theorem.

\begin{thm}[Addendum to Theorems~\ref{classthm} and \ref{lattice-split}]
\label{liegroup-lattice}
Let $G$ and $G'$ be two compact connected Lie groups and $p$ a prime
number. Then
\begin{enumerate}
\item \label{integrallie}
$(W_G,L_G)$ and $(W_{G'},L_{G'})$ are isomorphic if and only if
$G$ is isomorphic to $G'$ up to the substitution of direct factors
isomorphic to $\Sp(n)$ with direct factors isomorphic to $\SO(2n+1)$.
\item \label{twolie}
$(W_G,L_G\otimes \Z_2)$ and $(W_{G'},L_{G'}\otimes
\Z_2)$ are isomorphic if and only if $\covtwop{G}$ and $\covtwop{G'}$ are
$2$-equivalent up to the substitution of direct factors isomorphic to
$\Sp(n)$ with direct factors isomorphic to $\SO(2n+1)$. Moreover the
following conditions are equivalent:
\begin{enumerate}
\item \label{a2} $(W_G,L_G\otimes \Z_2,L_{G\onecov}\otimes \Z_2)$ and
$(W_{G'},L_{G'}\otimes \Z_2,L_{G'\onecov}\otimes \Z_2)$ are
isomorphic.
\item \label{b2}
$\covtwop{G}$ is $2$-equivalent to $\covtwop{G'}$.
\item \label{c2}
$(BG)\twocom \simeq (BG')\twocom$.
\end{enumerate}
\item\label{oddlie}
For $p$ odd the following conditions are equivalent:
\begin{enumerate}
\item \label{aodd}
$(W_G,L_G\otimes \Z_p)$ and $(W_{G'},L_{G'}\otimes \Z_p)$ are isomorphic.
\item \label{bodd}
$\covpp{G}$ and $\covpp{G'}$ are $p$-equivalent up to the substitution of
direct factors isomorphic to $\Sp(n)$ with direct factors isomorphic
to $\Spin(2n+1)$.
\item \label{codd}
$(BG)\pcom \simeq (BG')\pcom$.
\end{enumerate}
\end{enumerate}
\end{thm}

Note that two simple compact Lie groups $G$ and $G'$ of the form
$\covpp{\cdot}$ are $p$-equivalent if and only if they are
isomorphic. The next two examples shows how this fails in general.

\begin{example} 
Let $\zeta =e^{2\pi i/p}$ and $G =
\SU(p) \times \SU(p)$, and let $\Delta_1$ be the central subgroup of $G$
generated by $(\zeta I,\zeta I)$ and $\Delta_2$ the central subgroup
generated by $(\zeta I,\zeta^2 I)$. Then $G/\Delta_1$ and $G/\Delta_2$
are non-isomorphic as Lie groups if $p \geq 5$, but they are
$p$-equivalent; similar examples can be constructed for $p=2,3$.
\end{example}

\begin{example}
Let $p \geq 5$, $\zeta = e^{2\pi i/p}$, and  $G = \SU(p) \times \SU(p)
\times \SU(p) \times S^1 \times S^1$. Consider the subgroups
$$
\Delta_1 = \left< (\zeta I, I, \zeta I,\zeta,1),(I, \zeta I, \zeta
  I,1,\zeta)\right>,\ 
\Delta_2 = \langle(\zeta I, I,\zeta I,\zeta,1),(I,\zeta I,\zeta
I,1,\zeta^2)\rangle.
$$
The quotients $G/\Delta_1$ and $G/\Delta_2$ are again $p$-equivalent but
not isomorphic. One can check that in this example $\alpha_k$ cannot
be chosen to be the identity. Similar examples can be constructed for
$p=2,3$.
\end{example}

Before proving Theorem~\ref{liegroup-lattice}, we need the following lemma.

\begin{lemma} \label{autosnotsun}
Let $G$ be a simple simply connected compact Lie group not isomorphic
to a special unitary group. Suppose that $\varphi \in N_{\GL(L_G
  \otimes \Z_p)}(W_G)$. Then there exists $\sigma \in \Aut(G)$ such
that $\varphi$ and $\sigma$ induces the same automorphism of $\dZ(G)$.
\end{lemma}

\begin{proof}
If $G = \Spin(2n+1)$, $\Sp(n)$, $2E_7$, $E_8$, $F_4$ or $G_2$, the claim is
obvious, since there are no non-trivial automorphisms of the center.
For $G= \Spin(2n)$, $n\geq 4$, it follows from
Theorem~\ref{outerautos} that $N_{\GL(L_G \otimes \Z_p)}(W_G)$ is
generated by the scalars $\Z_p^\times$, $W$ and the automorphisms of
the Dynkin diagram. Thus the only potential problem occurs for $p=2$
and $n$ odd where $\dZ(G) \cong \Z/4$. However it follows by
\cite[Planche~IV(XI)]{bourbaki-lie4:6} that the  non-trivial
automorphism of the center is induced by the non-trivial graph
automorphism. Finally, by \cite[Planche~V(XI)]{bourbaki-lie4:6},
the same argument works in the case $G=3E_6$, $p=3$ where $\dZ(G)
\cong \Z/3$.
\end{proof}

\begin{proof}[Proof of Theorem~\ref{liegroup-lattice}]
By \cite[\S 4]{osse97} or \cite[Prop.~3.2(vi)]{JMO95} we can recover
the root datum of a compact connected Lie group from its integral
lattice up to substitution of direct factors isomorphic to $\Sp(n)$ with
direct factors isomorphic to $\SO(2n+1)$. Part \eqref{integrallie} now
follows.

Suppose that $G$ and $G'$ are $p$-equivalent. Then by assumption the
associated covering groups $H$ and $H'$ are isomorphic, in such a way
that the image of $K$ in $H'$ differs from $K'$ by a self-map of
$T_{H'}$ of the type $\Psi$. The map $\dT_H \to \dT_{H'}$ induced by
the composite of the isomorphism with the self-map $\Psi$ is hence an
isomorphism sending $K$ to $K'$. Hence we get an induced isomorphism
$(W_G,L_G \otimes \Z_p) \to (W_{G'},L_{G'} \otimes \Z_p)$.

Conversely, suppose that $G = H/K$ and $G' = H'/K'$ do not contain any
direct factors isomorphic to $\Sp(n)$, and that $(W_G,L_G \otimes
\Z_p)$ is isomorphic to $(W_{G'},L_{G'} \otimes \Z_p)$.
By Proposition~\ref{genpioneremark2} the fundamental group
of $G$ equals the coinvariants $(L_G)_W$ and hence $L_H = (L_G)^W \oplus
SL_G$. This shows that $(W,L_H\otimes \Z_p)$ can be reconstructed from
$(W,L_G\otimes \Z_p)$.  By the classification of simply connected
compact Lie groups we can for $p=2$ reconstruct $H$, up to
isomorphism, from $(W,L_H\otimes \Z_p)$. For $p$ odd the only
ambiguity arises from direct factors isomorphic to $\Sp(n)$ or
$\Spin(2n+1)$, but in this case by the assumption on $G$,
$H$ cannot contain any direct factors isomorphic to $\Sp(n)$, and we
conclude that in all cases we can reconstruct $H$, up to isomorphism,
from $(W,L_G\otimes \Z_p)$. Hence we can without loss of generality
assume that $H = H'$.
Note that $K$ is the cokernel of the inclusion $L_H\otimes
\Z_p\rightarrow L_G\otimes \Z_p$, so we can also recover the inclusion
$K\subseteq L_H\otimes \Z/p^\infty$. In other words an isomorphism
between $(W, L_G \otimes \Z_p)$ and $(W,L_{G'} \otimes \Z_p)$ induces
a self-equivalence of $\dT_H$ taking $K$ to $K'$. We have to see that
if there exists such a self-equivalence then there exists an
automorphism of $H$, followed by a map of the type $\Psi$, which also
takes $K$ to $K'$, since this will show that $H/K$ and $H/K'$ are
$p$-equivalent. Write $H = U \times T'$, where $U= R_1 \times \cdots
\times R_n \times S$, in the notation introduced before the
theorem. By Proposition~\ref{algprod} any self-map of $(W_H,L_H
\otimes \Z_p)$ induces a self-map $\varphi$ of  $\dT_U \times \dT'$
which is of the form $\varphi = \varphi_1 \times \varphi_2$. A priori
$\varphi_2$ is an element in $\Aut(\dT') \cong \GL( L_{T'} \otimes \Z_p)$,
but since $K_1$ and $K_2$ are finite we can without loss of generality
replace it by a matrix with integer coefficients and determinant prime
to $p$.
By the `elementary divisor theorem' in linear algebra over a Euclidean
domain (cf.\ e.g., \cite[Thm.~12.4.3]{artin91}) we can find a basis for
$L_{T'}$ in which $\varphi_2$ is given by $AD$ where $A$ is a product
of elementary matrices over $\Z$ and $D$ is
a diagonal matrix with determinant prime to $p$. Since we are only
interested in the effect over $\Z/p^s$ for some fixed large $s$ and
$D$ consists of units modulo $p^s$ we can furthermore change
$\varphi_2$ and $A$ such that $D$ can be assumed to be of the form
$\diag(1,\ldots,1,k)$, with $k$ prime to $p$, since if $x$ is a unit
in $\Z/p^s$ then the matrix $\diag(x,x^{-1})$ can be written as a
product of elementary matrices by a straightforward calculation (done
in \cite[40.25]{CR87}). This shows that we can put $\varphi_2$ on the
correct form. For the map $\varphi_1$ this follows directly from
Proposition~\ref{algprod} together with Lemma~\ref{autosnotsun} and
Theorem~\ref{outerautos}. Hence we have seen that for any $\varphi$ we
can find a map of the stated form which takes $K$ to $K'$.

The above analysis directly shows the first claim in \eqref{twolie} as
well as \eqref{aodd} $\Leftrightarrow$ \eqref{bodd}. From the first
claim in \eqref{twolie}, \eqref{a2} $\Leftrightarrow$ \eqref{b2}
follows, since $\Sp(n)$ and $\SO(2n+1)$ have different
$\Z_2$-reflection data $(W_G,L_G\otimes \Z_2,L_{G\onecov}\otimes
\Z_2)$.  The implications $\eqref{b2} \Rightarrow \eqref{c2}$ follows
from the existence of unstable Adams operations on $\SU(n)$, first
constructed by Sullivan \cite[p.~5.51]{sullivan70}, realizing $\Psi$ on the
level of $\F_p$-completed classifying spaces (or for an overkill, use
Theorem~\ref{classthm} directly). The implication $\eqref{bodd}
\Rightarrow \eqref{codd}$ also follows from this together with the
fact that $B\SO(2n+1)\pcom$ is homotopy equivalent to $B\Sp(n)\pcom$ for $p>2$,
as originally proved by Friedlander \cite[Thm.~2.1]{friedlander75}
(or again as a very special case of Theorem~\ref{classthm}). The remaining
implications  $\eqref{c2} \Rightarrow \eqref{a2}$ and $\eqref{codd}
\Rightarrow \eqref{aodd}$ follow directly from the fundamental
properties of the Weyl group of a $p$-compact group.
\end{proof}

%%%%%%%%%%%%%%%%%%%%%%%%%%%%%%%%%%%%%%%%%%%%%%%%%%%%%%%%%%%%%%%%%%%%%%%%%%%

\section{Appendix: Invariant rings of finite $\Z_p$-reflection
groups, $p$ odd (following Notbohm)} \label{invariantsection}

The purpose of this appendix is to recall Notbohm's determination
\cite{notbohm99lattice} of finite $\Z_p$-reflection groups $(W,L)$,
$p$ odd, such that the invariant ring $\Z_p[L]^W$ is a polynomial
algebra. 

Before stating it let us however for easy reference recall the
following `classical' characterizations of being a `$p$-torsion free'
$p$-compact group, which has a proof by general arguments
which we will sketch below.

\begin{thm}\label{classicalborel}
Let $X$ be a connected $p$-compact group with maximal torus $T$ and
Weyl group $W_X$. The following statements are equivalent:
\begin{enumerate}
\item \label{f}
$H^*(X;\Z_p)$ is torsion free. 
\item \label{g}
$H^*(X;\Z_p)$ is an exterior algebra over $\Z_p$ with generators in
odd degrees (or equivalently with $\F_p$ instead of $\Z_p$).
\item \label{c}
$H^*(BX;\Z_p)$ is a polynomial algebra over $\Z_p$ with generators in
even degree (or equivalently with $\F_p$ instead of $\Z_p$).
\item \label{b}
$H^*(BX;\Z_p)$ is a polynomial algebra and $H^*(BX;\Z_p) \xrightarrow{\cong}
H^*(BT;\Z_p)^{W_X}$.
\end{enumerate}
\end{thm}

We now give Notbohm's classification. The first part (which is a
general argument reducing to the simply connected case) is
\cite[Thm.~1.3]{notbohm99lattice} and the second (which is a
case-by-case argument in the simply connected case) is
a slight extension of \cite[Thm.~1.4]{notbohm99lattice}.
For the benefit of the reader we give a
streamlined proof of the second part. Recall that for a finite
$\Z_p$-reflection group $(W,L)$ we define $SL$ to be the submodule of
$L$ generated by elements of the form $(1-w)x$ with $w \in W$ and $x
\in L$. We call $(W,L)$ {\em simply connected} if $L \cong SL'$ for some
$\Z_pW$-lattice $L'$. (Note that for $p$ odd this is equivalent to
$SL=L$ since $S^2L'=SL'$, cf.\ the discussion of $\Z_p$-reflection data
in the introduction.)

\begin{thm}[Finite $\Z_p$-reflection groups with polynomial
  invariants, $p$ odd] \label{polyinvthm}
Let $p$ be an odd prime and $(W,L)$ a finite $\Z_p$-reflection
group. Then we have the following statements:
\begin{enumerate}
\item \label{cover-lattice}
$\Z_p[L]^W$ is a polynomial algebra if and only if $\Z_p[SL]^W$ is a
polynomial algebra and the group of coinvariants $L_W$ is torsion
free.
\item \label{lattice-simplycon}
Suppose $(W,L)$ is irreducible and simply connected. The following
conditions are equivalent:
\begin{enumerate}
\item
$\Z_p[L]^W$ is a polynomial algebra.
\item
$\F_p[L\otimes \F_p]^W$ is a polynomial algebra.
\item
$(W,L)$ is not isomorphic to $(W_G,L_G\otimes\Z_p)$ for the following
pairs $(G,p)$: $(F_4,3)$, $(3E_6,3)$, $(2E_7,3)$, $(E_8,3)$ and $(E_8,5)$.
\end{enumerate}
\end{enumerate}
In particular, if $X$ is an exotic $p$-compact group then
$\Z_p[L_X]^{W_X}$ is a polynomial algebra and if $(W,L) =
(W_G,L_G\otimes \Z_p)$ for a compact connected Lie group $G$ then
$\Z_p[L_G\otimes \Z_p]^{W_G}$ is a polynomial algebra if and only if
$H^*(G;\Z_p)$ is torsion free.
\end{thm}

\begin{proof}[Sketch of proof of Theorem~\ref{classicalborel}]
The equivalence of \eqref{f}, \eqref{g}, and \eqref{c} is proved by
old $H$-space and loop space arguments which we first very briefly
sketch. By a Bockstein spectral sequence argument (cf.\ e.g., \cite[\S
11-2]{kane88}) $H^*(X;\Z_p)$ torsion free if and only if $H^*(X;\Z_p)$
is an exterior algebra on odd dimensional generators so \eqref{f} is
equivalent to \eqref{g}. This is again equivalent to
$H^*(BX;\Z_p)$ being a polynomial algebra on even dimensional generators
(using the Eilenberg-Moore and the cobar spectral sequence; see e.g.,
\cite[\S 7-4]{kane88}), so \eqref{g} is equivalent to
\eqref{c}. 

That \eqref{b} implies \eqref{c} is obvious. The fact that
\eqref{f}--\eqref{c} also imply \eqref{b} requires more machinery and
is probably first found in \cite[Thm.~2.11]{DMW92}---we
quickly sketch an argument. We want to show that the map $r:
H^*(BX;\Z_p) \to H^*(BT;\Z_p)^W$ is an isomorphism. By
\cite[Thm.~9.7(iii)]{DW94}
\begin{equation} \label{rationaliso}
H^*(BX;\Z_p) \otimes \Q \xrightarrow{\cong} H^*(BT;\Z_p)^W \otimes \Q.
\end{equation}
This implies by comparing Krull dimensions that the number of
polynomial generators equals the rank of $T$. Since $H^*(BT;\F_p)$ is
finitely generated over $H^*(BX;\F_p)$ by \cite[Prop.~9.11]{DW94} it
follows by comparing Krull dimensions again that $H^*(BX;\F_p) \to
H^*(BT;\F_p)$ is injective. Hence $H^*(BX;\Z_p) \to H^*(BT;\Z_p)$ has
to be injective by Nakayama's lemma. Likewise $r$ has to be
surjective: By \eqref{rationaliso} the cokernel of $r$ has to be
$p$-torsion. Since the reduction mod $p$ of $r$ is still injective (as
seen above) the cokernel of $r$ has to be $p$-torsion free as well
(since $\Tor(\coker(r),\F_p) = 0$).
\end{proof}

\begin{rem} \label{classicalborelextra}
If $p$ is odd then $\F_p$-coefficients can also be used in
Theorem~\ref{classicalborel}\eqref{b} by a Galois theory argument
using Lemma~\ref{modpreduction}. For $p=2$, this is not true as can
be seen by taking $X=\SU(2)\twocom$.
See \cite{DW:abelian} for a version
for $p=2$.
\end{rem}

\begin{rem} \label{degree-remark}
If $(W,L)$ is a finite $\Z_p$-reflection group then
$\Z_p[L]^W$ is a polynomial algebra if and only if $\F_p[L\otimes
\F_p]^W$ is a polynomial algebra and the canonical monomorphism
$\Z_p[L]^W\otimes \F_p\longrightarrow \F_p[L\otimes\F_p]^W$ is an isomorphism
as shown in \cite[Lem.~2.3]{notbohm99lattice}. Note that this can be
reformulated as saying that
$\Z_p[L]^W$ is a polynomial algebra if and only if
$\F_p[L\otimes\F_p]^W$ is a polynomial algebra with generators in the
same degrees as the generators of $\Q_p[L\otimes\Q]^W$, since $
\dim_{\Q_p} (\Q_p[L\otimes \Q]^W)_n = \dim_{\F_p} (\Z_p[L]^W\otimes
\F_p)_n \leq \dim_{\F_p} (\F_p[L\otimes \F_p]^W)_n
$ for any $n$.
\end{rem}

\begin{rem} \label{PU3-example}
The finite $\Z_3$-reflection group $(W,L) =
(W_{\PU(3)},L_{\PU(3)}\otimes \Z_3)$ does not have invariant ring a
polynomial ring (e.g., since $L_W \cong \Z/3$ is not torsion
free). However a short calculation shows that $\F_3[L\otimes\F_3]^W$
{\em is} a polynomial ring with generators in degrees $1$ and $6$ (as
opposed to the degrees over $\Q_3$ which are $2$ and $3$). (See also
\cite[Rem.~5.3]{DMW92}.)
It turns out that this example is essentially the only one
since it can be proved that if $(W,L)$ is a finite $\Z_p$-reflection
group, $p$ odd, such that $\F_p[L]^W$ is a polynomial algebra, then
$\Z_p[L]^W$ is also a polynomial algebra unless $p=3$ and $(W,L)$
contains $(W_{\PU(3)},L_{\PU(3)}\otimes \Z_3)$ as a direct factor. We omit
the proof which is an extension of the technique used in the examples
in Section $7$ in a preprint version of \cite{DW98}, which can at the
time of writing be found on Wilkerson's homepage.
\end{rem}

\begin{lemma} \label{nonmodular}
Assume that $L$ is a finitely generated free $\Z_p$-module and
that $W$ is a finite subgroup of $\GL(L)$. If $p\nmid |W|$ and
$\F_p[L\otimes \F_p]^W$ is a polynomial algebra, then $\Z_p[L]^W$
is also a polynomial algebra.
\end{lemma}

\begin{proof}
By assumption we have the averaging homomorphisms
$\Z_p[L]\longrightarrow \Z_p[L]^W$ and
$\F_p[L\otimes\F_p]\longrightarrow \F_p[L\otimes\F_p]^W$ given by
$f\mapsto \frac{1}{|W|} \sum_{w\in W} w\cdot f$. These are
obviously surjective and hence the commutative diagram
$$
\xymatrix{
\Z_p[L] \ar[r] \ar[d] & \Z_p[L]^W \ar[d] \\
\F_p[L\otimes\F_p] \ar[r] & \F_p[L\otimes\F_p]^W
}
$$
shows that the reduction homomorphism $\Z_p[L]^W \rightarrow
\F_p[L\otimes\F_p]^W$ is surjective. The result now follows easily
from Nakayama's lemma (cf.\ \cite[Lem.~2.3]{notbohm99lattice}).
\end{proof}

\begin{proof}[Proof of Theorem~\ref{polyinvthm}]
Part~\eqref{cover-lattice} is contained in
\cite[Thm.~1.3]{notbohm99lattice}. To prove
part~\eqref{lattice-simplycon} note that by Notbohm \cite{notbohm96}
(see also \cite[Thm.~1.2(iii)]{notbohm99lattice} and
Theorem~\ref{lattice-split}), there is a unique finite irreducible
simply connected $\Z_p$-reflection group for each group on the
Clark-Ewing list. We now go through the list, verifying the result in
each case.

If $p\nmid |W|$ the invariant ring $\F_p[L\otimes\F_p]^W$ is
a polynomial algebra by the Shephard-Todd-Chevalley theorem
(\cite[Thm.~7.2.1]{benson93} or \cite[Thm.~7.4.1]{smith95}), and thus
Lemma~\ref{nonmodular} shows that $\Z_p[L]^W$ is a polynomial algebra.

Next, assume that $(W,L)$ is an exotic $\Z_p$-reflection group.
If $(W,L)$ belongs to family $2$,
the representing matrices with respect to the standard basis are
monomial and so $\Z_p[L]^W$ is a polynomial algebra by
\cite[Thm.~2.4]{nakajima79}.

An inspection of the Clark-Ewing list now shows that only $4$ exotic
cases remain, namely $(G_{12},p=3)$, $(G_{29},p=5)$, $(G_{31},p=5)$
and $(G_{34},p=7)$. In the first case we have $G_{12} \cong
\GL_2(\F_3)$ and Lemma~\ref{modpreduction} shows that the action on
$L\otimes\F_3=(\F_3)^2$ is the canonical one. The invariant ring
$\F_3[L\otimes \F_3]^{\GL_2(\F_3)}$ was computed by Dickson
\cite{dickson11}. In the remaining $3$ cases the mod $p$
invariant ring was calculated by Xu \cite{xu94,xuthesis} using computer, see
also Kemper-Malle \cite[Prop.~6.1]{km97}. The conclusion
of these computations is that in all $4$ cases the invariant ring
$\F_p[L\otimes \F_p]^W$ is a polynomial algebra with generators in the same
degrees as the generators of $\Q_p[L\otimes \Q]^W$. By
Remark~\ref{degree-remark} we then see that $\Z_p[L]^W$ is a
polynomial algebra in these cases.

The only remaining cases are the finite simply connected $\Z_p$-reflection
groups which are not exotic. Since $p$ is odd and $\pi_1(G)$ and
$(L_G)_{W_G}$ only differ by an elementary abelian $2$-group
(cf.\ the proof of Theorem~\ref{fundamentalgroup} and
Proposition~\ref{genpioneremark2}), we may assume that
$(W,L)=(W_G,L_G\otimes \Z_p)$ for some simply connected compact Lie group
$G$. In this case Demazure
\cite{demazure73} shows that if $p$ is not a torsion prime for the
root system associated to $G$, then the invariant rings
$\Z_p[L_G\otimes\Z_p]^{W_G}$ and $\F_p[L_G\otimes\F_p]^{W_G}$ are
polynomial algebras.

By the calculation of torsion primes for the simple root
systems, \cite[\S 7]{demazure73}, the excluded pairs $(G,p)$
in the last part of the theorem are exactly the cases where the root
system of $G$ has $p$-torsion. In these cases Kemper-Malle
\cite[Prop.~6.1 and Pf.~of~Thm.~8.5]{km97} shows that
$\F_p[L_G\otimes\F_p]^{W_G}$ is not a polynomial algebra. Hence in these
cases $\Z_p[L_G\otimes\Z_p]^{W_G}$ is not a polynomial algebra by
\cite[Lem.~2.3(i)]{notbohm99lattice}. This proves the second claim.

Finally, let $G$ be a compact connected Lie group with Weyl group $W=W_G$
and integral lattice $L=L_G$. We now prove that $\Z_p[L\otimes
\Z_p]^W$ is a polynomial algebra if and only if $H^*(G;\Z_p)$ is
torsion free. (See also \cite[Prop.~1.11]{notbohm99}.) One direction
follows from Theorem~\ref{classicalborel}, so assume now that
$\Z_p[L\otimes \Z_p]^W$ is a polynomial algebra. From
Theorem~\ref{polyinvthm}\eqref{cover-lattice} we see that
$\Z_p[S(L\otimes\Z_p)]^W$ is a
polynomial algebra and that $(L\otimes\Z_p)_W$ is torsion free. Since
$p$ is odd, we have $(L\otimes\Z_p)_W=\pi_1(G)\otimes\Z_p$ and
$S(L\otimes\Z_p)=L_{G\onecov}\otimes \Z_p$, cf.\ the proofs of
Theorem~\ref{fundamentalgroup} and Proposition~\ref{genpioneremark2}. From
the above we conclude that $H^*(G\onecov;\Z_p)$ is
torsion free. Since $\pi_1(G)$ has no $p$-torsion, it now
follows easily from the Serre spectral sequence that $H^*(G;\Z_p)$ is
torsion free.
\end{proof}

\begin{rem} \label{notbohm-elaboration}
Let $p$ be an odd prime and $(W,L)$ a finite $\Z_p$-reflection
group. We claim that the following conditions are equivalent:
\begin{enumerate}
\item \label{ZpLpoly}
$\Z_p[L]^W$ is a polynomial algebra.
\item \label{FpLpoly+LWtorsfree}
$\F_p[L\otimes \F_p]^W$ is a polynomial algebra and $L_W$ is torsion
free.
\item \label{FpSLpoly+LWtorsfree}
$\F_p[SL\otimes \F_p]^W$ is a polynomial algebra and $L_W$ is torsion
free.
\end{enumerate}
Indeed we have $\eqref{ZpLpoly}\Leftrightarrow
\eqref{FpSLpoly+LWtorsfree}$ by Theorem~\ref{polyinvthm} since
$(W,SL)$ can be decomposed as a direct product of finite irreducible simply
connected $\Z_p$-reflection groups by \cite[Thm.~1.4]{notbohm96}. The
implication $\eqref{ZpLpoly}\Rightarrow \eqref{FpLpoly+LWtorsfree}$
follows from \cite[Thm.~1.3 and Lem.~2.3]{notbohm99lattice}. Finally
$\eqref{FpLpoly+LWtorsfree}\Rightarrow \eqref{FpSLpoly+LWtorsfree}$
follows from \cite[Prop.~4.1]{nakajima79} as $L_W$ torsion free
implies that $SL\otimes \F_p\rightarrow L\otimes \F_p$ is injective.
\end{rem}

%%%%%%%%%%%%%%%%%%%%%%%%%%%%%%%%%%%%%%%%%%%%%%%%%%%%%%%%%%%%%%%%%%%%%%%%%%%

\section{Appendix: Outer automorphisms of finite
$\Z_p$-reflection groups} \label{autosection}

Theorem~\ref{classthm} states that the outer automorphism group of a connected
$p$-compact group $X$, $p$ odd, equals $N_{\GL(L_X)}(W_X)/W_X$, which
makes it useful to have a complete case-by-case calculation of this
group. The purpose of this appendix is to provide such a calculation
based on results of Brou{\'e}-Malle-Michel \cite[Prop.~3.13]{BMM99} over the complex numbers. Calculations in the
case where $W$ is one of the exotic groups from family $2a$ were given
in \cite[\S 6]{notbohm98} (where the non-standard notation $G(q,r;n)$ for
$G(q,q/r,n)$ is used).

Theorem~\ref{lattice-split} and Proposition~\ref{algprod}
reduces the calculation of $N_{\GL(L)}(W)/W$ to the case where $(W,L)$
is exotic or $(W,L)=(W_G,L_G\otimes \Z_p)$ for some compact connected
Lie group $G$. In the second case we can write $G=H/K$ where $H$ is a
direct product of a torus and a simply connected compact Lie group and
$K$ is a finite central subgroup of $H$, and it is easy to use
coverings to compute $N_{\GL(L_G\otimes \Z_p)}(W_G)/W_G$ from
$N_{\GL(L_H\otimes \Z_p)}(W_H)/W_H$. By Proposition~\ref{algprod} this
again reduces to the case where $H$ is simple.

We can hence restrict to the case
where $(W,L)$ is exotic or $(W,L)=(W_G,L_G\otimes \Z_p)$
for some simple simply connected compact Lie group $G$. 
For the statement of our result in these cases (which will take place
in the theorem below as well as in the following elaborations), we fix the
realizations $G(m,r,n)$ of the groups from family $2$ as described in
Section~\ref{latticesection}. Moreover we also fix the realizations of
the complex reflection groups $G_i$, $4\leq i\leq 37$, to be the ones
described in \cite{ST54}. Let $\mu_n$ denote the group of
$n$th roots of unity.
If $G$ is a simply connected compact Lie group, its integral lattice
$L_G$ equals the coroot lattice. Hence the automorphism group $\Gamma$
of the Dynkin diagram of $G$ can be considered as a subgroup of
$N_{\GL(L_G)}(W_G)$, cf.\ \cite[\S 12.2]{humphreys78}.
For $G=\Spin(5)$, $F_4$ or $G_2$, there is an automorphism
$\varphi_l$ of $L_G\otimes \Z[1/\sqrt{l}]$ of order $2$ (here $l=2$
for $\Spin(5)$ and $F_4$ and $l=3$ for $G_2$), see
\cite[p.\ 182--183]{BMM99} or \cite[p.\ 217]{carter:simple} for details.

\begin{thm}[Outer automorphisms of finite $\Z_p$-reflection
  groups] \label{outerautos}
Let $(W,L)$ be a finite irreducible simply connected $\Z_p$-reflection
group, i.e., $(W,L)$ is exotic or of
the form $(W_G,L_G\otimes \Z_p)$ for a simple simply connected compact
Lie group $G$.
Let $(W,V)$ be the associated complex reflection
group. Then $N_{\GL(V)}(W)= \left<W,\C^{\times}\right>$ and hence
$N_{\GL(L)}(W)/W=\Z_p^{\times}/Z(W)$ and $N_{\GL(L)}(W)/\Z_p^{\times}
W=1$ except in the following cases:
\begin{enumerate}
\item \label{fam2exotic}
$W=G(m,r,n)$ is exotic and belongs to family $2$, $(m,r,n)\neq (4,2,2),
(3,3,3)$:
$N_{\GL(V)}(W)= \left<G(m,1,n),\C^{\times}\right>$ and
$N_{\GL(L)}(W)/\Z_p^{\times} W = C_{\gcd(r,n)}$,
cf.\ \ref{normalizer:family2}.
\item \label{G422}
$W=G(4,2,2)$:
$N_{\GL(V)}(W)= \left<G_8,\C^{\times}\right>$ and
$N_{\GL(L)}(W)/\Z_p^{\times} W = \Sigma_3$,
cf.\ \ref{normalizer:G422}.
\item \label{G333}
$W=G(3,3,3)$:
$N_{\GL(V)}(W)= \left<G_{26},\C^{\times}\right>$ and
$N_{\GL(L)}(W)/\Z_p^{\times} W = A_4$,
cf.\ \ref{normalizer:G333}.
\item \label{G5}
$W=G_5$:
$N_{\GL(V)}(W)= \left<G_{14},\C^{\times}\right>$ and
$N_{\GL(L)}(W)/\Z_p^{\times} W = C_2$,
cf.\ \ref{normalizer:G5}.
\item \label{G7}
$W=G_7$:
$N_{\GL(V)}(W)= \left<G_{10},\C^{\times}\right>$ and
$N_{\GL(L)}(W)/\Z_p^{\times} W = C_2$,
cf.\ \ref{normalizer:G7}.
\item
$(W,L)=(W_G,L_G\otimes \Z_p)$ for $G=\Spin(4n)$, $n\geq 2$:
$N_{\GL(V)}(W)= \left<W,\C^{\times},\Gamma\right>$ and
$N_{\GL(L)}(W)/\Z_p^{\times} W \cong \Gamma$, cf.\ \ref{normalizer:typeD}.
\item
$(W,L)=(W_G,L_G\otimes \Z_p)$ for $G=\Spin(5)$, $F_4$ or $G_2$:
$N_{\GL(V)}(W)= \left<W,\C^{\times},\varphi_l\right>$. Moreover
$N_{\GL(L)}(W)/\Z_p^{\times} W=1$ for $p=l$ and
$N_{\GL(L)}(W)/\Z_p^{\times} W=C_2$ for $p\neq l$, cf.\ \ref{normalizer:twist}.
\end{enumerate}
\end{thm}

\begin{lemma} \label{normalizer:basechange}
Let $K\subseteq K'$ be fields of characteristic zero, and $W\subseteq
\GL_n(K)$ an irreducible reflection group. Then $N_{\GL_n(K')}(W) =
\left<N_{\GL_n(K)}(W),{K'}^{\times}\right>$.
\end{lemma}

\begin{proof}
The inclusion `$\supseteq$' is clear, so suppose $g\in
N_{\GL_n(K')}(W)$. Consider the system of equations $X w= gwg^{-1} X$,
$w\in W$ where $X$ is an $n\times n$-matrix. Over $K'$ this has the
solution $X=g$. By \cite[Lem.~2.10]{DW00}, the representation
$W\rightarrow \GL_n(K')$ is irreducible, so the solution space is
the $1$-dimensional space spanned by $g$. Since the coefficients lie
in $K$, the solution space over $K$ is $1$-dimensional as well, so we
can write $g=\lambda g_1$ with $\lambda\in K'$ and $g_1\in M_n(K)$. As
$g\neq 0$ we get $\lambda\neq 0$ and $g_1\in N_{\GL_n(K)}(W)$.
\end{proof}

We can now start the proof of Theorem~\ref{outerautos}. The results on
$N_{\GL(V)}(W)$ follow directly from \cite[Prop.~3.13]{BMM99} except
when $W$ belongs to family $2$ or $W=G_{28}$. The
structure of $N_{\GL(V)}(W)$ in the cases \eqref{fam2exotic},
\eqref{G422} and \eqref{G333} also follows from
\cite[Prop.~3.13]{BMM99} noting that $\left<G(4,1,2),G_6,\C^{\times}\right>=
\left<G_8,\C^{\times}\right>$ and $G(3,1,3)\subseteq G_{26}$.

Now assume that $W$ does not belong to family $2$ and $W\neq G_5,
G_7, G_{28}$. Let $n$ denote the rank of $W$ and $K$ the field extension of
$\Q$ generated by the entries of the matrices representing $W$. Our
assumption ensures that $N_{\GL(V)}(W)=
\left<W,\C^{\times}\right>$. Since $W$ is a reflection group it has
Schur index $1$ and we can assume that $K$ equals the
character field of $W$, cf.\ \cite[Cor.~p.~429]{CE74}.
Then $N_{\GL_n(K)}(W)=
\left<W,K^{\times}\right>$ and Lemma~\ref{normalizer:basechange} now
shows that $N_{\GL_n(\Q_p)}(W)= \left<W,\Q_p^{\times}\right>$. Hence
we get $N_{\GL(L)}(W)= \left<W,\Z_p^{\times}\right>$ and since $W$ is
irreducible we have $W\cap \Z_p^{\times}=Z(W)$,
cf.\ \cite[Lem.~2.9]{DW00}.

This proves Theorem~\ref{outerautos} in case $W$ does not
belong to family $2$ and $W\neq G_5, G_7, G_{28}$. In the cases
\eqref{fam2exotic}, \eqref{G422}, \eqref{G333}, \eqref{G5} and
\eqref{G7} we only need to find the structure of $N_{\GL(L)}(W)$. This
is done in Elaborations~\ref{normalizer:family2}, \ref{normalizer:G422},
\ref{normalizer:G333}, \ref{normalizer:G5} and \ref{normalizer:G7}
below.

This leaves the cases where $(W,L) = (W_G,L_G\otimes \Z_p)$ for
some simple simply connected compact Lie group $G$ such that $W_G$ belongs
to family $2$ and $W_G\neq G_{28}$, i.e.,
$G=\Spin(2n+1)$ for $n\geq 2$, $\Sp(n)$ for $n\geq 3$, $\Spin(2n)$ for
$n\geq 4$, $G_2$ and $F_4$. In the first two cases $W_G$ equals
to $G(2,1,n)$ and hence $N_{\GL(V)}(W_G) =
\left<W_G,\C^{\times}\right>$ when $n\geq 3$ by
\cite[Prop.~3.13]{BMM99}. As above this proves
Theorem~\ref{outerautos} in these cases. The cases $G=\Spin(2n)$, $n\geq
4$ is dealt with in Elaboration~\ref{normalizer:typeD}, and the cases
$G=\Spin(5)$, $G_2$ and $F_4$ are handled in
Elaboration~\ref{normalizer:twist}.

To treat the dihedral group $G(m,m,2)$ from family $2$ we need the
following auxiliary result.
\begin{lemma} \label{dihedral:unitlemma}
Assume that $m\geq 3$ and $p\equiv \pm 1 \pmod{m}$ so that
$\zeta_m+\zeta_m^{-1}\in \Z_p$. Then $2+\zeta_m+\zeta_m^{-1}$ is a
unit in $\Z_p$.
\end{lemma}

\begin{proof}
It suffices to prove that the norm
$N_{\Q(\zeta_m+\zeta_m^{-1})/\Q}(2+\zeta_m+\zeta_m^{-1})$ is not
divisible by $p$. Since its square equals the norm
$N_{\Q(\zeta_m)/\Q}(2+\zeta_m+\zeta_m^{-1})$ it is enough to see that
this norm is not divisible by $p$. In $\Q(\zeta_m)$ we have
$2+\zeta_m+\zeta_m^{-1}=(1+\zeta_m)^2/\zeta_m$ and since $\zeta_m$ is
a unit it is enough to see that $N_{\Q(\zeta_m)/\Q}(1+\zeta_m)$ is not
divisible by $p$. By definition
$$
N_{\Q(\zeta_m)/\Q}(1+\zeta_m)=
\prod_{\substack{0\leq k\leq m\\ \gcd(k,m)=1}} (1+\zeta_m^k)=
(-1)^{\phi(m)} \prod_{\substack{0\leq k\leq m\\ \gcd(k,m)=1}} (-1-\zeta_m^k)=
\Phi_m(-1).
$$
The claim now follows from \cite[Lem.~2.9]{washington97}.
\end{proof}

\begin{con}[Family $2$, generic case] \label{normalizer:family2}
Let $W=G(m,r,n)$ from family $2$ and let $p$ be a prime number such
that $W$ is an exotic $\Z_p$-reflection group. Thus if $n\geq 3$ or $n=2$
and $r<m$ we have $m\geq 3$ and $p\equiv 1 \pmod{m}$, and for $n=2$
and $m=r$ we have $m\geq 5, m\neq 6$ and $p\equiv \pm 1
\pmod{m}$. Assume moreover that $(m,r,n)\neq (4,2,2), (3,3,3)$ (these
two cases are dealt with in Elaborations~\ref{normalizer:G422} and
\ref{normalizer:G333} below).

Assume first that $p\equiv 1 \pmod{m}$. The realizations of the groups
$G(m,r,n)$ and $G(m,1,n)$ from above are both defined over the ring
$\Z[\zeta_m]$ which embeds in
$\Z_p$. Lemma~\ref{normalizer:basechange} shows that
$N_{\GL_n(\Z_p)}(W) = \left<G(m,1,n), \Z_p^{\times}\right>$
whence the natural homomorphism
$(A(m,1,n)/A(m,r,n)) \times \Z_p^{\times} \longrightarrow
N_{\GL_n(\Z_p)}(W)/W$ is surjective. The kernel is the cyclic group
generated by the element
$([\zeta_m I_n],\zeta_m^{-1})$ (here $[\zeta_m I_n]\in
A(m,1,n)/A(m,r,n)$ denotes the coset of $\zeta_m I_n$) and thus
$N_{\GL_n(\Z_p)}(W)/W = (A(m,1,n)/A(m,r,n)) \circ_{C_m}
\Z_p^{\times}$. Note that $A(m,1,n)/A(m,r,n)$ is cyclic of order $r$
generated by the element $x=[\diag(1,\ldots,1,\zeta_m)]$ and that
$[\zeta_m I_n]=x^n$.

If the assumption $p\equiv 1 \pmod{m}$ is not satisfied, then
$W=G(m,m,2)$ is the dihedral group of order $2m$ with $m\geq 5, m\neq
6$ and $p\equiv -1 \pmod{m}$. Conjugating the realization of
$G(m,m,2)$ from above with the element
$$
g=
\begin{bmatrix}
1 & -\zeta_m^{-1} \\
1 & -\zeta_m
\end{bmatrix}
$$
gives a realization $G(m,m,2)^g$ defined over the character field
$\Q(\zeta_m+\zeta_m^{-1})$. Note that if $m$ is odd, then
$N_{\GL_2(\C)}(G(m,m,2))= \left<G(m,1,2),\C^{\times}\right>=
\left<G(m,m,2),\C^{\times}\right>$ and hence
$N_{\GL_2(\Z_p)}(G(m,m,2)^g)/G(m,m,2)^g=\Z_p^{\times}$,
so we may assume that $m$ is even. Since $G(m,1,2)$ is generated by
$G(m,m,2)$ and $\diag(1,\zeta_m)$ we find
$$
N_{\GL_2(\Z_p)}(G(m,m,2)^g)=\left<G(m,m,2)^g,
\begin{bmatrix}
1 & 1\\
-1 & 1+\zeta_m+\zeta_m^{-1}
\end{bmatrix}
,\Q_p^{\times}\right> \cap \GL_2(\Z_p)
$$
using Lemma~\ref{normalizer:basechange}. From Lemma~\ref{dihedral:unitlemma}
we see that the above matrix is invertible over $\Z_p$ and hence
$$
N_{\GL_2(\Z_p)}(G(m,m,2)^g)=\left<G(m,m,2)^g,
\begin{bmatrix}
1 & 1\\
-1 & 1+\zeta_m+\zeta_m^{-1}
\end{bmatrix},\Z_p^{\times}\right>
$$
Thus the homomorphism $\Z \times (\Z_p^{\times}/\mu_2) \longrightarrow
N_{\GL_2(\Z_p)}(G(m,m,2)^g)/G(m,m,2)^g$ which maps $(k,[\lambda])$ to
the coset of
$\lambda
\begin{bmatrix}
1 & 1\\
-1 & 1+\zeta_m+\zeta_m^{-1}
\end{bmatrix}^k$
is surjective. The kernel is easily seen to be the infinite cyclic
group generated by the element $(-2,[1+\zeta_m+\zeta_m^{-1}])$ and
thus we get
$N_{\GL_2(\Z_p)}(G(m,m,2)^g)/G(m,m,2)^g \cong \Z\circ_{\Z}
(\Z_p^{\times}/\mu_2)$. It is easily checked that
$[2+\zeta_m+\zeta_m^{-1}]$ has a square root in $\Z_p^{\times}/\mu_2$
if and only if either $m\equiv 0 \pmod{4}$ or $m\equiv 2 \pmod{4}$ and
$p\equiv -1 \pmod{2m}$. In this case we have
$N_{\GL_2(\Z_p)}(G(m,m,2)^g)/G(m,m,2)^g \cong C_2\times
(\Z_p^{\times}/\mu_2)$.
\end{con}

\begin{con}[$G(4,2,2)$] \label{normalizer:G422}
The realization of the group $G(4,2,2)$ from above and the realization
of the group $G_8$ from \cite[Table~II]{ST54} are both defined over
their common character field $\Q(i)$. Thus the relevant primes $p$ are
the ones satisfying $p\equiv 1 \pmod{4}$. More precisely the
representations are defined over $\Z[\frac{1}{2}, i]$ and as this ring
embeds in $\Z_p$ for all $p$ as above, we get
$N_{\GL_2(\Z_p)}(G(4,2,2)) = \left<G_8, \Z_p^{\times}\right>$ by
Lemma~\ref{normalizer:basechange}. It is easily checked that
$G_8=\left<G(4,2,2),H\right>$, where $H$ is the group of order $24$
generated by the elements
$$
\begin{bmatrix}
0 & i \\
1 & 0
\end{bmatrix},
\frac{1+i}{2}
\begin{bmatrix}
1 & 1 \\
i & -i
\end{bmatrix}
$$
Since $G(4,2,2)\cap \left<H,\Z_p^{\times}\right> = Z(H) = \mu_4$
we conclude that $N_{\GL_2(\Z_p)}(G(4,2,2))/G(4,2,2) \cong (H/Z(H))\times
(\Z_p^{\times}/\mu_4) \cong \Sigma_3\times (\Z_p^{\times}/\mu_4)$.
\end{con}

\begin{con}[$G(3,3,3)$] \label{normalizer:G333}
The realization of the group $G(3,3,3)$ from above and the realization
of the group $G_{26}$ from \cite[p.~297]{ST54} are both defined over
their common character field $\Q(\omega)$ where
$\omega=e^{2\pi i/3}$. Thus the relevant primes $p$ are the ones
satisfying $p\equiv 1 \pmod{3}$. More precisely the representations
are defined over $\Z[\frac{1}{3},\omega]$ and as this ring embeds in
$\Z_p$ for all $p$ as above, we see that $N_{\GL_3(\Z_p)}(G(3,3,3)) =
\left<G_{26}, \Z_p^{\times}\right>$ using
Lemma~\ref{normalizer:basechange}. It is easily checked that $G_{26}$
is the semidirect product of $G(3,3,3)$ with the group $H\cong
\SL_2(\F_3)$ generated by the elements
$$
R_1=
\begin{bmatrix}
1 & 0 & 0\\
0 & 1 & 0\\
0 & 0 & \omega^2
\end{bmatrix},
R_2=\frac{1}{\sqrt{-3}}
\begin{bmatrix}
\omega & \omega^2 & \omega^2\\
\omega^2 & \omega & \omega^2\\
\omega^2 & \omega^2 & \omega
\end{bmatrix}
$$
The center of $H$ is generated by the element
$$
z=
\begin{bmatrix}
0 & -1 & 0\\
-1 & 0 & 0\\
0 & 0 & -1
\end{bmatrix}
$$
and $G(3,3,3)\cap \left<H,\Z_p^{\times}\right> =
\left<-z,\mu_3\right>$. Thus
$$
N_{\GL_3(\Z_p)}(G(3,3,3))/G(3,3,3) \cong
H\circ_{C_2} (\Z_p^{\times}/\mu_3) \cong \SL_2(\F_3)\circ_{C_2}
(\Z_p^{\times}/\mu_3)
$$
where the central product is given by identifying $z\in H$ with
$\left[-1\right]\in \Z_p^{\times}/\mu_3$.
\end{con}

\begin{con}[$G_5$] \label{normalizer:G5}
The realization of the group $G_5$ from \cite[Table~I]{ST54} is
defined over the field $\Q(\zeta_{12})$, but the character
field is $\Q(\omega)$ and thus the relevant primes $p$ are the ones
satisfying $p\equiv 1 \pmod{3}$. Conjugation by the matrix
$$
g=
\begin{bmatrix}
2 & \sqrt{3}-1\\
(\sqrt{3}-1)(1-i) & i-1
\end{bmatrix}
$$
gives a realization defined over $\Z[\frac{1}{3},\omega]$ which embeds
in $\Z_p$ for all $p$ as above. Its easily checked that
$G_{14}$ is generated by $G_5$ and the reflection
$$
S=\frac{1}{\sqrt{2}}
\begin{bmatrix}
-1 & i\\
-i & 1
\end{bmatrix}
$$
By Lemma~\ref{normalizer:basechange} we then get:
$$
N_{\GL_2(\Z_p)}(G_5^g)=\left<G_5^g,
\begin{bmatrix}
0 & 1\\
-2\omega & 0
\end{bmatrix},\Z_p^{\times}\right>
$$
and thus the homomorphism $\Z \times (\Z_p^{\times}/\mu_6) \longrightarrow
N_{\GL_2(\Z_p)}(G_5^g)/G_5^g$ which maps $(k,[\lambda])$ to the coset of
$\lambda
\begin{bmatrix}
0 & 1\\
-2\omega & 0
\end{bmatrix}^k$
is surjective. The kernel is easily seen to be the infinite cyclic
group generated by the element $(-2,[2])$ and we get
$N_{\GL_2(\Z_p)}(G_5^g)/G_5^g \cong \Z\circ_{\Z}
(\Z_p^{\times}/\mu_6)$. It is easy to check that the element $[2]$ has a square
root in $\Z_p^{\times}/\mu_6$ if and only if $p\equiv 1,7,19
\pmod{24}$ (that is unless $p\equiv 13 \pmod{24}$). In this case we
get the simpler description $N_{\GL_2(\Z_p)}(G_5^g)/G_5^g \cong
C_2\times (\Z_p^{\times}/\mu_6)$.
\end{con}

\begin{con}[$G_7$] \label{normalizer:G7}
The realizations of the groups $G_7$ and $G_{10}$ given in
\cite[Tables~I and II]{ST54} are both defined over their common character
field $\Q(\zeta_{12})$. Thus the relevant primes $p$ are the ones
satisfying $p\equiv 1 \pmod{12}$. More precisely the representations
are defined over $\Z[\frac{1}{2},\zeta_{12}]$ and as this ring
embeds in $\Z_p$ for all $p$ as above, we get
$N_{\GL_2(\Z_p)}(G_7) = \left<G_{10}, \Z_p^{\times}\right>$ by
Lemma~\ref{normalizer:basechange}. It is easily checked that
$G_{10}=\left<G_7,C_4\right>$, where $C_4$ is the cyclic group
generated by $\begin{bmatrix}
1 & 0 \\
0 & i
\end{bmatrix}$. Since $G_7 \cap (C_4\times \Z_p^{\times}) = C_2\times
\mu_{12}$ we conclude that $N_{\GL_2(\Z_p)}(G_7)/G_7 \cong C_2 \times
(\Z_p^{\times}/\mu_{12})$.
\end{con}

\begin{con}[$\Spin(2n)$, $n\geq 4$] \label{normalizer:typeD}
The group $G=\Spin(2n)$, $n\geq 4$ has Weyl group $G(2,2,n)$ and by
\cite[Prop.~3.13]{BMM99} $N_{\GL(V)}(W)$ equals
$\left<G(2,1,n),\C^{\times}\right>$ for $n\geq 5$. For $n$ odd we
have $\left<G(2,1,n),\C^{\times}\right> = \left<W_G,\C^{\times}\right>$,
and as above this proves Theorem~\ref{outerautos} in these cases.

Now assume that $n\geq 4$ is even. The automorphism of the Dynkin
diagram which exchanges
$\alpha_{n-1}$ and $\alpha_n$ equals $\diag(1,\ldots,1,-1)$ and hence 
$N_{\GL(V)}(W)$ equals $\left<W_G,\C^{\times},\Gamma\right>$ for
$n\geq 6$. For $n=4$ this also holds since $N_{\GL(V)}(W) =
\left<G(2,1,n),W(F_4),\C^{\times}\right> =
\left<W(F_4),\C^{\times}\right> =
\left<W_G,\C^{\times},\Gamma\right>$
by \cite[Prop.~3.13]{BMM99} and a direct computation.
Since $\Gamma\subseteq \GL(L_G)$,
Lemma~\ref{normalizer:basechange} shows that
$N_{\GL(L_G\otimes \Z_p)}(W_G) =
\left<W_G,\Z_p^{\times},\Gamma\right>$ in all cases and hence
$N_{\GL(L_G\otimes \Z_p)}(W_G)/\Z_p^{\times}W_G \cong \Gamma$ since
$\Gamma \cap \Z_p^{\times}W_G = 1$.
\end{con}

\begin{con}[$\Spin(5)$, $F_4$ and $G_2$] \label{normalizer:twist}
For $G=F_4$ the first claim follows directly from
\cite[Prop.~3.13]{BMM99}.
For $G=\Spin(5)$, $W(G)$ is
conjugate to $G(4,4,2)$ in $\GL_2(\C)$ and \cite[Prop.~3.13]{BMM99}
shows $N_{\GL_2(\C)}(G(4,4,2)) = \left<G(4,1,2),\C^{\times}\right>$.
Similarly, for $G=G_2$, $W(G)$ is conjugate to $G(6,6,2)$ in
$\GL_2(\C)$ and \cite[Prop.~3.13]{BMM99} shows
$N_{\GL_2(\C)}(G(6,6,2)) = \left<G(6,1,2),\C^{\times}\right>$.
From this it is easy to check the first claim in these cases.
Thus $N_{\GL(V)}(W) = \left<W,\C^{\times},\sqrt{l} \varphi_l\right>$
in all cases. By construction $\sqrt{l} \varphi_l$ stabilizes $L_G$,
so Lemma~\ref{normalizer:basechange} shows that
$N_{\GL(L\otimes \Z_p)}(W) = \left<W,\Z_p^{\times},\sqrt{l}
  \varphi_l\right>$. Since $\varphi_l^2=1$ and $\sqrt{l} \varphi_l$
has determinant a power of $l$ the remaining claims follows.
\end{con}
 
\bibliographystyle{plain}
\bibliography{poddclassification}

\begin{thebibliography}{100}

\bibitem{adams69}
J.~F. Adams.
\newblock {\em Lectures on {L}ie groups}.
\newblock W. A. Benjamin, Inc., New York-Amsterdam, 1969.

\bibitem{AM76}
J.~F. Adams and Z.~Mahmud.
\newblock Maps between classifying spaces.
\newblock {\em Inv. Math.}, 35:1--41, 1976.

\bibitem{AW80}
J.~F. Adams and C.~W. Wilkerson.
\newblock Finite {$H$}-spaces and algebras over the {S}teenrod algebra.
\newblock {\em Ann. of Math. (2)}, 111(1):95--143, 1980.
\newblock (Erratum: {A}nn. of {M}ath. (2) {\bf 113} (1981), no. 3, 621--622).

\bibitem{aguade89}
J.~Aguad{\'e}.
\newblock Constructing modular classifying spaces.
\newblock {\em Israel J. Math.}, 66(1-3):23--40, 1989.

\bibitem{henningbrev}
H.~H. Andersen.
\newblock Letter to {K}. {A}ndersen and {J}. {G}rodal.
\newblock October 2002.

\bibitem{kksa:thesis}
K.~K.~S. Andersen.
\newblock The normalizer splitting conjecture for $p$-compact groups.
\newblock {\em Fund. Math.}, 161(1-2):1--16, 1999.
\newblock Algebraic topology (Kazimierz Dolny, 1997).

\bibitem{ABGP03}
K.~K.~S. Andersen, T.~Bauer, J.~Grodal, and E.~K. Pedersen.
\newblock A finite loop space not rationally equivalent to a compact {L}ie
  group.
\newblock {\em Invent. Math.}, 157(1):1--10, 2004.

\bibitem{artin91}
M.~Artin.
\newblock {\em Algebra}.
\newblock Prentice Hall Inc., Englewood Cliffs, NJ, 1991.

\bibitem{benson91}
D.~J. Benson.
\newblock {\em Representations and cohomology. {I}}, volume~30 of {\em
  Cambridge Studies in Advanced Mathematics}.
\newblock Cambridge University Press, Cambridge, 1991.
\newblock Basic representation theory of finite groups and associative
  algebras.

\bibitem{benson93}
D.~J. Benson.
\newblock {\em Polynomial invariants of finite groups}.
\newblock Cambridge University Press, Cambridge, 1993.

\bibitem{borel61}
A.~Borel.
\newblock Sous-groupes commutatifs et torsion des groupes de {L}ie compacts
  connexes.
\newblock {\em T\^ohoku Math. J. (2)}, 13:216--240, 1961.

\bibitem{borel67}
A.~Borel.
\newblock {\em Topics in the homology theory of fibre bundles}, volume~36 of
  {\em Lecture Notes in Math.}
\newblock Springer-Verlag, Berlin, 1967.

\bibitem{boreloeuvresII}
A.~Borel.
\newblock {\em {\OE}uvres: collected papers. {V}ol. {I}{I}}.
\newblock Springer-Verlag, Berlin, 1983.
\newblock 1959--1968.

\bibitem{borel91}
A.~Borel.
\newblock {\em Linear algebraic groups}.
\newblock Springer-Verlag, New York, second edition, 1991.

\bibitem{BT65}
A.~Borel and J.~Tits.
\newblock Groupes r\'eductifs.
\newblock {\em Inst. Hautes \'Etudes Sci. Publ. Math.}, (27):55--150, 1965.

\bibitem{magma1}
W.~Bosma, J.~Cannon, and C.~Playoust.
\newblock The {M}agma algebra system. {I}. {T}he user language.
\newblock {\em J. Symbolic Comput.}, 24(3-4):235--265, 1997.
\newblock Computational algebra and number theory (London, 1993). (See also the
  Magma homepage at http://magma.maths.usyd.edu.au/magma/).

\bibitem{bott56}
R.~Bott.
\newblock An application of the {M}orse theory to the topology of {L}ie-groups.
\newblock {\em Bull. Soc. Math. France}, 84:251--281, 1956.

\bibitem{bourbaki-lie4:6}
N.~Bourbaki.
\newblock {\em \'{E}l\'ements de math\'ematique. {F}asc. {X}{X}{X}{I}{V}.
  {G}roupes et alg\`ebres de {L}ie. {C}hapitre {I}{V}: {G}roupes de {C}oxeter
  et syst\`emes de {T}its. {C}hapitre {V}: {G}roupes engendr\'es par des
  r\'eflexions. {C}hapitre {V}{I}: {S}yst\`emes de racines}.
\newblock Hermann, Paris, 1968.
\newblock Actualit\'es Scientifiques et Industrielles, No. 1337.

\bibitem{bourbaki-lie7:8}
N.~Bourbaki.
\newblock {\em \'{E}l\'ements de math\'ematique}.
\newblock Hermann, Paris, 1975.
\newblock Fasc. XXXVIII: Groupes et alg\`ebres de Lie. Chapitre VII:
  Sous-alg\`ebres de Cartan, \'el\'ements r\'eguliers. Chapitre VIII:
  Alg\`ebres de Lie semi-simples d\'eploy\'ees, Actualit\'es Scientifiques et
  Industrielles, No. 1364.

\bibitem{bo9}
N.~Bourbaki.
\newblock {\em \'{E}l\'ements de math\'ematique: Groupes et alg\`ebres de
  {L}ie}.
\newblock Masson, Paris, 1982.
\newblock Chapitre 9. Groupes de Lie r\'eels compacts. [Chapter 9. Compact real
  Lie groups].

\bibitem{bousfield75}
A.~K. Bousfield.
\newblock The localization of spaces with respect to homology.
\newblock {\em Topology}, 14:133--150, 1975.

\bibitem{bk}
A.~K. Bousfield and D.~M. Kan.
\newblock {\em Homotopy limits, completions and localizations}.
\newblock Springer-Verlag, Berlin, 1972.
\newblock Lecture Notes in Mathematics, Vol. 304.

\bibitem{BT85}
T.~Br{\"o}cker and T.~tom Dieck.
\newblock {\em Representations of compact {L}ie groups}, volume~98 of {\em
  Graduate Texts in Mathematics}.
\newblock Springer-Verlag, New York, 1985.

\bibitem{BV98}
C.~Broto and A.~Viruel.
\newblock Homotopy uniqueness of ${B}{P}{U}(3)$.
\newblock In {\em Group representations: cohomology, group actions and topology
  (Seattle, WA, 1996)}, pages 85--93. Amer. Math. Soc., Providence, RI, 1998.

\bibitem{BV99}
C.~Broto and A.~Viruel.
\newblock Projective unitary groups are totally {$N$}-determined {$p$}-compact
  groups.
\newblock {\em Math. Proc. Cambridge Philos. Soc.}, 136(1):75--88, 2004.

\bibitem{BMM99}
M.~Brou{\'e}, G.~Malle, and J.~Michel.
\newblock Towards spetses. {I}.
\newblock {\em Transform. Groups}, 4(2-3):157--218, 1999.
\newblock Dedicated to the memory of Claude Chevalley.

\bibitem{carter:simple}
R.~W. Carter.
\newblock {\em Simple groups of {L}ie type}.
\newblock John Wiley \& Sons, London-New York-Sydney, 1972.
\newblock Pure and Applied Mathematics, Vol. 28.

\bibitem{carter:rep}
R.~W. Carter.
\newblock {\em Finite groups of {L}ie type. Conjugacy classes and complex
  characters}.
\newblock John Wiley \& Sons Inc., New York, 1985.
\newblock A Wiley-Interscience Publication.

\bibitem{castellanaexotic}
N.~Castellana.
\newblock The homotopic adjoint representation for the exotic $p$-compact
  groups.
\newblock Preprint; available from http://mat.uab.es/$\sim$natalia.

\bibitem{castellana2a}
N.~Castellana.
\newblock On the $p$-compact groups corresponding to the reflection groups
  ${G}(q,r,n)$.
\newblock Trans. Amer. Math. Soc. (to appear). Available from
  http://mat.uab.es/$\sim$natalia.

\bibitem{CV66}
L.~S. Charlap and A.~T. Vasquez.
\newblock The cohomology of group extensions.
\newblock {\em Trans. Amer. Math. Soc.}, 124:24--40, 1966.

\bibitem{chevalley55}
C.~Chevalley.
\newblock Invariants of finite groups generated by reflections.
\newblock {\em Amer. J. Math.}, 77:778--782, 1955.

\bibitem{chevalley:3}
C.~Chevalley.
\newblock {\em Th\'eorie des groupes de {L}ie. {T}ome {I}{I}{I}.
  {T}h\'eor\`emes g\'en\'eraux sur les alg\`ebres de {L}ie}.
\newblock Hermann \& Cie, Paris, 1955.

\bibitem{CE74}
A.~Clark and J.~Ewing.
\newblock The realization of polynomial algebras as cohomology rings.
\newblock {\em Pacific J. Math.}, 50:425--434, 1974.

\bibitem{cohen76}
A.~M. Cohen.
\newblock Finite complex reflection groups.
\newblock {\em Ann. Sci. \'Ecole Norm. Sup. (4)}, 9(3):379--436, 1976.
\newblock (Erratum: {A}nn. {S}ci. {\'E}cole {N}orm. {S}up. (4), 11, (1978), no.
  4, 613).

\bibitem{CoGr}
A.~M. Cohen and R.~L. Griess, Jr.
\newblock On finite simple subgroups of the complex {L}ie group of type
  ${E}_8$.
\newblock In {\em The Arcata Conference on Representations of Finite Groups
  (Arcata, Calif., 1986)}, pages 367--405. Amer. Math. Soc., Providence, RI,
  1987.

\bibitem{CoWa}
A.~M. Cohen and D.~B. Wales.
\newblock Finite subgroups of ${F}_4(\mathbb{C})$ and ${E}_6(\mathbb{C})$.
\newblock {\em Proc. London Math. Soc. (3)}, 74(1):105--150, 1997.

\bibitem{CCNPW85}
J.~H. Conway, R.~T. Curtis, S.~P. Norton, R.~A. Parker, and R.~A. Wilson.
\newblock {\em Atlas of finite groups}.
\newblock Oxford University Press, Eynsham, 1985.
\newblock Maximal subgroups and ordinary characters for simple groups, With
  computational assistance from J. G. Thackray.

\bibitem{CR62}
C.~W. Curtis and I.~Reiner.
\newblock {\em Representation theory of finite groups and associative
  algebras}.
\newblock Pure and Applied Mathematics, Vol. XI. Interscience Publishers, a
  division of John Wiley \& Sons, New York-London, 1962.

\bibitem{CR81}
C.~W. Curtis and I.~Reiner.
\newblock {\em Methods of representation theory. {V}ol. {I}}.
\newblock John Wiley \& Sons Inc., New York, 1981.
\newblock With applications to finite groups and orders, Pure and Applied
  Mathematics, A Wiley-Interscience Publication.

\bibitem{CR87}
C.~W. Curtis and I.~Reiner.
\newblock {\em Methods of representation theory. {V}ol. {II}}.
\newblock Pure and Applied Mathematics. John Wiley \& Sons Inc., New York,
  1987.
\newblock With applications to finite groups and orders, A Wiley-Interscience
  Publication.

\bibitem{CWW74}
M.~Curtis, A.~Wiederhold, and B.~Williams.
\newblock Normalizers of maximal tori.
\newblock In {\em Localization in group theory and homotopy theory, and related
  topics (Sympos., Battelle Seattle Res. Center, Seattle, Wash., 1974)}, pages
  31--47. Lecture Notes in Math., Vol. 418. Springer, Berlin, 1974.

\bibitem{demazure73}
M.~Demazure.
\newblock Invariants sym\'etriques entiers des groupes de {W}eyl et torsion.
\newblock {\em Invent. Math.}, 21:287--301, 1973.

\bibitem{dickson11}
L.~E. Dickson.
\newblock A fundamental system of invariants of the general modular linear
  group with a solution of the form problem.
\newblock {\em Trans. Amer. Math. Soc.}, 12(1):75--98, 1911.

\bibitem{DDK80equivariant}
E.~Dror, W.~G. Dwyer, and D.~M. Kan.
\newblock Equivariant maps which are self homotopy equivalences.
\newblock {\em Proc. Amer. Math. Soc.}, 80(4):670--672, 1980.

\bibitem{DZ87}
W.~Dwyer and A.~Zabrodsky.
\newblock Maps between classifying spaces.
\newblock In {\em Algebraic topology, Barcelona, 1986}, pages 106--119.
  Springer, Berlin, 1987.

\bibitem{dwyer97}
W.~G. Dwyer.
\newblock Homology decompositions for classifying spaces of finite groups.
\newblock {\em Topology}, 36(4):783--804, 1997.

\bibitem{dwyer98}
W.~G. Dwyer.
\newblock Lie groups and $p$-compact groups.
\newblock In {\em Proceedings of the International Congress of Mathematicians,
  Extra Vol. II (Berlin, 1998)}, pages 433--442 (electronic), 1998.

\bibitem{DK92}
W.~G. Dwyer and D.~M. Kan.
\newblock Centric maps and realization of diagrams in the homotopy category.
\newblock {\em Proc. Amer. Math. Soc.}, 114(2):575--584, 1992.

\bibitem{DKS89}
W.~G. Dwyer, D.~M. Kan, and J.~H. Smith.
\newblock Towers of fibrations and homotopical wreath products.
\newblock {\em J. Pure Appl. Algebra}, 56(1):9--28, 1989.

\bibitem{DMW86}
W.~G. Dwyer, H.~R. Miller, and C.~W. Wilkerson.
\newblock The homotopic uniqueness of ${B}{S}\sp 3$.
\newblock In {\em Algebraic topology, Barcelona, 1986}, pages 90--105.
  Springer, Berlin, 1987.

\bibitem{DMW92}
W.~G. Dwyer, H.~R. Miller, and C.~W. Wilkerson.
\newblock Homotopical uniqueness of classifying spaces.
\newblock {\em Topology}, 31(1):29--45, 1992.

\bibitem{DW:abelian}
W.~G. Dwyer and C.~W. Wilkerson.
\newblock $p$-compact groups with abelian {W}eyl groups.
\newblock Unpublished manuscript (available from http://hopf.math.purdue.edu).

\bibitem{DW92}
W.~G. Dwyer and C.~W. Wilkerson.
\newblock A cohomology decomposition theorem.
\newblock {\em Topology}, 31(2):433--443, 1992.

\bibitem{DW93}
W.~G. Dwyer and C.~W. Wilkerson.
\newblock A new finite loop space at the prime two.
\newblock {\em J. Amer. Math. Soc.}, 6(1):37--64, 1993.

\bibitem{DW94}
W.~G. Dwyer and C.~W. Wilkerson.
\newblock Homotopy fixed-point methods for {L}ie groups and finite loop spaces.
\newblock {\em Ann. of Math. (2)}, 139(2):395--442, 1994.

\bibitem{dw:center}
W.~G. Dwyer and C.~W. Wilkerson.
\newblock The center of a $p$-compact group.
\newblock In {\em The \v Cech centennial (Boston, MA, 1993)}, pages 119--157.
  Amer. Math. Soc., Providence, RI, 1995.

\bibitem{dw:split}
W.~G. Dwyer and C.~W. Wilkerson.
\newblock Product splittings for $p$-compact groups.
\newblock {\em Fund. Math.}, 147(3):279--300, 1995.

\bibitem{DW96}
W.~G. Dwyer and C.~W. Wilkerson.
\newblock Diagrams up to cohomology.
\newblock {\em Trans. Amer. Math. Soc.}, 348(5):1863--1883, 1996.

\bibitem{DW98elementary}
W.~G. Dwyer and C.~W. Wilkerson.
\newblock The elementary geometric structure of compact {L}ie groups.
\newblock {\em Bull. London Math. Soc.}, 30(4):337--364, 1998.

\bibitem{DW98}
W.~G. Dwyer and C.~W. Wilkerson.
\newblock K\"ahler differentials, the ${T}$-functor, and a theorem of
  {S}teinberg.
\newblock {\em Trans. Amer. Math. Soc.}, 350(12):4919--4930, 1998.

\bibitem{DW00}
W.~G. Dwyer and C.~W. Wilkerson.
\newblock Centers and {C}oxeter elements.
\newblock In {\em Homotopy methods in algebraic topology (Boulder, CO, 1999)},
  volume 271 of {\em Contemp. Math.}, pages 53--75. Amer. Math. Soc.,
  Providence, RI, 2001.

\bibitem{DW01}
William~G. Dwyer and Clarence~W. Wilkerson, Jr.
\newblock Cartan involutions and normalizers of maximal tori.
\newblock {\em Publ. Res. Inst. Math. Sci.}, 40(2):295--309, 2004.

\bibitem{feshbach79}
M.~Feshbach.
\newblock The transfer and compact {L}ie groups.
\newblock {\em Trans. Amer. Math. Soc.}, 251:139--169, 1979.

\bibitem{freudenthal}
H.~Freudenthal.
\newblock Beziehungen der $\mathfrak{E}_7$ und $\mathfrak{E}_8$ zur
  {O}ktavenebene. {VIII}.
\newblock {\em Nederl. Akad. Wetensch. Proc. Ser. A 66 = Indag. Math.},
  21:447--465, 1959.

\bibitem{friedlander75}
E.~M. Friedlander.
\newblock Exceptional isogenies and the classifying spaces of simple {L}ie
  groups.
\newblock {\em Ann. Math. (2)}, 101:510--520, 1975.

\bibitem{FM88}
E.~M. Friedlander and G.~Mislin.
\newblock Conjugacy classes of finite solvable subgroups in {L}ie groups.
\newblock {\em Ann. Sci. \'Ecole Norm. Sup. (4)}, 21(2):179--191, 1988.

\bibitem{GZ67}
P.~Gabriel and M.~Zisman.
\newblock {\em Calculus of fractions and homotopy theory}.
\newblock Springer-Verlag New York, Inc., New York, 1967.

\bibitem{gorenstein}
D.~Gorenstein.
\newblock {\em Finite groups}.
\newblock Harper \& Row Publishers, New York, 1968.

\bibitem{griess}
R.~L. Griess, Jr.
\newblock Elementary abelian $p$-subgroups of algebraic groups.
\newblock {\em Geom. Dedicata}, 39(3):253--305, 1991.

\bibitem{GR98}
R.~L. Griess, Jr. and A.~J.~E. Ryba.
\newblock Embeddings of {${\rm PGL}\sb 2(31)$} and {${\rm SL}\sb 2(32)$} in
  {$E\sb 8(\mathbb C)$}.
\newblock {\em Duke Math. J.}, 94(1):181--211, 1998.
\newblock With appendices by M. Larsen and J.-P. Serre.

\bibitem{grodal02}
J.~Grodal.
\newblock Higher limits via subgroup complexes.
\newblock {\em Ann. of Math. (2)}, 155(2):405--457, 2002.

\bibitem{hammerli02}
J.-F. H{\"a}mmerli.
\newblock The outer automorphism group of normalizers of maximal tori in
  connected compact {L}ie groups.
\newblock {\em J. Lie Theory}, 12(2):357--368, 2002.

\bibitem{hochschild}
G.~Hochschild.
\newblock {\em The structure of {L}ie groups}.
\newblock Holden-Day Inc., San Francisco, 1965.

\bibitem{hopf41}
H.~Hopf.
\newblock \"{U}ber die {T}opologie der {G}ruppen-{M}annigfaltigkeiten und ihre
  {V}erallgemeinerungen.
\newblock {\em Ann. of Math. (2)}, 42:22--52, 1941.

\bibitem{humphreys75}
J.~E. Humphreys.
\newblock {\em Linear algebraic groups}.
\newblock Springer-Verlag, New York, 1975.
\newblock Graduate Texts in Mathematics, No. 21.

\bibitem{humphreys78}
J.~E. Humphreys.
\newblock {\em Introduction to {L}ie algebras and representation theory},
  volume~9 of {\em Graduate Texts in Mathematics}.
\newblock Springer-Verlag, New York, 1978.
\newblock Second printing, revised.

\bibitem{humphreys87}
J.~E. Humphreys.
\newblock The {S}teinberg representation.
\newblock {\em Bull. Amer. Math. Soc. (N.S.)}, 16(2):247--263, 1987.

\bibitem{humphreys90}
J.~E. Humphreys.
\newblock {\em Reflection groups and {C}oxeter groups}.
\newblock Cambridge University Press, Cambridge, 1990.

\bibitem{huppert67}
B.~Huppert.
\newblock {\em Endliche {G}ruppen. {I}}.
\newblock Die Grundlehren der Mathematischen Wissenschaften, Band 134.
  Springer-Verlag, Berlin, 1967.

\bibitem{JM92}
S.~Jackowski and J.~McClure.
\newblock Homotopy decomposition of classifying spaces via elementary abelian
  subgroups.
\newblock {\em Topology}, 31(1):113--132, 1992.

\bibitem{JMO92}
S.~Jackowski, J.~McClure, and B.~Oliver.
\newblock Homotopy classification of self-maps of ${B}{G}$ via ${G}$-actions.
  {I}+{I}{I}.
\newblock {\em Ann. of Math. (2)}, 135(2):183--226,227--270, 1992.

\bibitem{JMO95}
S.~Jackowski, J.~McClure, and B.~Oliver.
\newblock Self-homotopy equivalences of classifying spaces of compact connected
  {L}ie groups.
\newblock {\em Fund. Math.}, 147(2):99--126, 1995.

\bibitem{JO96}
S.~Jackowski and B.~Oliver.
\newblock Vector bundles over classifying spaces of compact {L}ie groups.
\newblock {\em Acta Math.}, 176(1):109--143, 1996.

\bibitem{kane88}
R.~M. Kane.
\newblock {\em The homology of {H}opf spaces}.
\newblock North-Holland Publishing Co., Amsterdam, 1988.

\bibitem{km97}
G.~Kemper and G.~Malle.
\newblock The finite irreducible linear groups with polynomial ring of
  invariants.
\newblock {\em Transform. Groups}, 2(1):57--89, 1997.

\bibitem{KM86}
N.~J. Kuhn and S.~A. Mitchell.
\newblock The multiplicity of the {S}teinberg representation of {${\rm GL}\sb
  n{\bf F}\sb q$} in the symmetric algebra.
\newblock {\em Proc. Amer. Math. Soc.}, 96(1):1--6, 1986.

\bibitem{lannes92}
J.~Lannes.
\newblock Sur les espaces fonctionnels dont la source est le classifiant d'un
  $p$-groupe ab\'elien \'el\'ementaire.
\newblock {\em Inst. Hautes \'Etudes Sci. Publ. Math.}, (75):135--244, 1992.
\newblock With an appendix by Michel Zisman.

\bibitem{lannes96}
J.~Lannes.
\newblock Th\'eorie homotopique des groupes de {L}ie (d'apr\`es {W}. {G}.
  {D}wyer et {C}. {W}. {W}ilkerson).
\newblock {\em Ast\'erisque}, (227):Exp.\ No.\ 776, 3, 21--45, 1995.
\newblock S\'eminaire Bourbaki, Vol.\ 1993/94.

\bibitem{Liebeck-Seitz}
M.~W. Liebeck and G.~M. Seitz.
\newblock On finite subgroups of exceptional algebraic groups.
\newblock {\em J. Reine Angew. Math.}, 515:25--72, 1999.

\bibitem{lusztig75}
G.~Lusztig.
\newblock On the discrete series representations of the classical groups over
  finite fields.
\newblock In {\em Proceedings of the International Congress of Mathematicians
  (Vancouver, B. C., 1974),Vol. 1}, pages 465--470. Canad. Math. Congress,
  Montreal, Que., 1975.

\bibitem{matthey02}
M.~Matthey.
\newblock Normalizers of maximal tori and cohomology of {W}eyl groups.
\newblock Preprint; available from
  http://www.math.uni-muenster.de/inst/sfb/about/publ/heft203.ps, March 2002.

\bibitem{may67}
J.~P. May.
\newblock {\em Simplicial objects in algebraic topology}.
\newblock Van Nostrand Mathematical Studies, No. 11. D. Van Nostrand Co., Inc.,
  Princeton, N.J.-Toronto, Ont.-London, 1967.

\bibitem{miller84}
H.~Miller.
\newblock The {S}ullivan conjecture on maps from classifying spaces.
\newblock {\em Ann. of Math. (2)}, 120(1):39--87, 1984.
\newblock (Erratum: {Ann.} of {Math.} 121 (1985), no. 3, 605--609).

\bibitem{moller95}
J.~M. M{\o}ller.
\newblock Homotopy {L}ie groups.
\newblock {\em Bull. Amer. Math. Soc. (N.S.)}, 32(4):413--428, 1995.

\bibitem{moller96}
J.~M. M{\o}ller.
\newblock Rational isomorphisms of $p$-compact groups.
\newblock {\em Topology}, 35(1):201--225, 1996.

\bibitem{moller97puppreprint}
J.~M. M{\o}ller.
\newblock {$PU(p)$} as a $p$-compact group.
\newblock unpublished manuscript, April 1997.

\bibitem{moller98}
J.~M. M{\o}ller.
\newblock Deterministic $p$-compact groups.
\newblock In {\em Stable and unstable homotopy (Toronto, ON, 1996)}, pages
  255--278. Amer. Math. Soc., Providence, RI, 1998.

\bibitem{moller99}
J.~M. M{\o}ller.
\newblock Normalizers of maximal tori.
\newblock {\em Math. Z.}, 231(1):51--74, 1999.

\bibitem{MN94}
J.~M. M{\o}ller and D.~Notbohm.
\newblock Centers and finite coverings of finite loop spaces.
\newblock {\em J. Reine Angew. Math.}, 456:99--133, 1994.

\bibitem{MN98}
J.~M. M{\o}ller and D.~Notbohm.
\newblock Connected finite loop spaces with maximal tori.
\newblock {\em Trans. Amer. Math. Soc.}, 350(9):3483--3504, 1998.

\bibitem{nakajima79}
H.~Nakajima.
\newblock Invariants of finite groups generated by pseudoreflections in
  positive characteristic.
\newblock {\em Tsukuba J. Math.}, 3(1):109--122, 1979.

\bibitem{neumann99}
F.~Neumann.
\newblock A theorem of {T}its, normalizers of maximal tori and fibrewise
  {B}ousfield-{K}an completions.
\newblock {\em Publ. Res. Inst. Math. Sci.}, 35(5):711--723, 1999.

\bibitem{ns02}
M.~D. Neusel and L.~Smith.
\newblock {\em Invariant theory of finite groups}.
\newblock American Mathematical Society, Providence, RI, 2002.

\bibitem{notbohm94}
D.~Notbohm.
\newblock Homotopy uniqueness of classifying spaces of compact connected {L}ie
  groups at primes dividing the order of the {W}eyl group.
\newblock {\em Topology}, 33(2):271--330, 1994.

\bibitem{notbohm95}
D.~Notbohm.
\newblock On the ``classifying space'' functor for compact {L}ie groups.
\newblock {\em J. London Math. Soc. (2)}, 52(1):185--198, 1995.

\bibitem{notbohm96}
D.~Notbohm.
\newblock $p$-adic lattices of pseudo reflection groups.
\newblock In {\em Algebraic topology: new trends in localization and
  periodicity (Sant Feliu de Gu\'\i xols, 1994)}, pages 337--352. Birkh\"auser,
  Basel, 1996.

\bibitem{notbohm98}
D.~Notbohm.
\newblock Topological realization of a family of pseudoreflection groups.
\newblock {\em Fund. Math.}, 155(1):1--31, 1998.

\bibitem{notbohm99lattice}
D.~Notbohm.
\newblock For which pseudo-reflection groups are the $p$-adic polynomial
  invariants again a polynomial algebra?
\newblock {\em J. Algebra}, 214(2):553--570, 1999.
\newblock (Erratum: {J}. {A}lgebra {\bf 218} (1999), no. 2, 286--287).

\bibitem{notbohm99}
D.~Notbohm.
\newblock Spaces with polynomial mod-$p$ cohomology.
\newblock {\em Math. Proc. Cambridge Philos. Soc.}, 126(2):277--292, 1999.

\bibitem{notbohm00}
D.~Notbohm.
\newblock Unstable splittings of classifying spaces of $p$-compact groups.
\newblock {\em Q. J. Math.}, 51(2):237--266, 2000.

\bibitem{notbohm00di4preprint}
D.~Notbohm.
\newblock On the 2-compact group {$DI(4)$}.
\newblock {\em J. Reine Angew. Math.}, 555:163--185, 2003.

\bibitem{bob:steinberg}
B.~Oliver.
\newblock Higher limits via {S}teinberg representations.
\newblock {\em Comm. Algebra}, 22(4):1381--1393, 1994.

\bibitem{oliver94}
B.~Oliver.
\newblock {$p$}-stubborn subgroups of classical compact {L}ie groups.
\newblock {\em J. Pure Appl. Algebra}, 92(1):55--78, 1994.

\bibitem{onishchik-vinberg}
A.~L. Onishchik and {\`E}.~B. Vinberg.
\newblock {\em Lie groups and algebraic groups}.
\newblock Springer-Verlag, Berlin, 1990.
\newblock Translated from the Russian and with a preface by D. A. Leites.

\bibitem{osse97}
A.~Osse.
\newblock $\lambda$-structures and representation rings of compact connected
  {L}ie groups.
\newblock {\em J. Pure Appl. Algebra}, 121(1):69--93, 1997.

\bibitem{quillen71}
D.~Quillen.
\newblock The {A}dams conjecture.
\newblock {\em Topology}, 10:67--80, 1971.

\bibitem{quillen72}
D.~Quillen.
\newblock On the cohomology and ${K}$-theory of the general linear groups over
  a finite field.
\newblock {\em Ann. of Math. (2)}, 96:552--586, 1972.

\bibitem{quillen78}
D.~Quillen.
\newblock Homotopy properties of the poset of nontrivial {$p$}-subgroups of a
  group.
\newblock {\em Adv. in Math.}, 28(2):101--128, 1978.

\bibitem{RS74}
D.~Rector and J.~Stasheff.
\newblock Lie groups from a homotopy point of view.
\newblock In {\em Localization in group theory and homotopy theory, and related
  topics (Sympos., Battelle Seattle Res. Center, Seattle, Wash., 1974)}, pages
  121--131. Lecture Notes in Math., Vol. 418. Springer, Berlin, 1974.

\bibitem{rector71}
D.~L. Rector.
\newblock Subgroups of finite dimensional topological groups.
\newblock {\em J. Pure Appl. Algebra}, 1(3):253--273, 1971.

\bibitem{serre99}
J.-P. Serre.
\newblock Letter to J. M. M{\o}ller.
\newblock Sept.\ 24, 1999.

\bibitem{serre51}
J.-P. Serre.
\newblock Homologie singuli\`ere des espaces fibr\'es. {A}pplications.
\newblock {\em Ann. of Math. (2)}, 54:425--505, 1951.

\bibitem{serre53}
J.-P. Serre.
\newblock Groupes d'homotopie et classes de groupes ab\'eliens.
\newblock {\em Ann. of Math. (2)}, 58:258--294, 1953.

\bibitem{serre77}
J.-P. Serre.
\newblock {\em Linear representations of finite groups}.
\newblock Springer-Verlag, New York, 1977.
\newblock Translated from the second French edition by Leonard L. Scott,
  Graduate Texts in Mathematics, Vol. 42.

\bibitem{ST54}
G.~C. Shephard and J.~A. Todd.
\newblock Finite unitary reflection groups.
\newblock {\em Canadian J. Math.}, 6:274--304, 1954.

\bibitem{smith95}
L.~Smith.
\newblock {\em Polynomial invariants of finite groups}.
\newblock A K Peters Ltd., Wellesley, MA, 1995.

\bibitem{smith82}
S.~D. Smith.
\newblock Irreducible modules and parabolic subgroups.
\newblock {\em J. Algebra}, 75(1):286--289, 1982.

\bibitem{spr:linalggrp}
T.~A. Springer.
\newblock {\em Linear algebraic groups}.
\newblock Birkh\"auser Boston Inc., Boston, MA, second edition, 1998.

\bibitem{SS70}
T.~A. Springer and R.~Steinberg.
\newblock Conjugacy classes.
\newblock In {\em Seminar on Algebraic Groups and Related Finite Groups (The
  Institute for Advanced Study, Princeton, N.J., 1968/69)}, pages 167--266.
  Springer, Berlin, 1970.

\bibitem{steinberg68}
R.~Steinberg.
\newblock {\em Endomorphisms of linear algebraic groups}.
\newblock Memoirs of the American Mathematical Society, No. 80. American
  Mathematical Society, Providence, R.I., 1968.

\bibitem{steinberg75}
R.~Steinberg.
\newblock Torsion in reductive groups.
\newblock {\em Advances in Math.}, 15:63--92, 1975.

\bibitem{steinberg99}
R.~Steinberg.
\newblock The isomorphism and isogeny theorems for reductive algebraic groups.
\newblock {\em J. Algebra}, 216(1):366--383, 1999.

\bibitem{sullivan70}
D.~Sullivan.
\newblock Geometric topology, {P}art {I}, {L}ocalization, periodicity and
  {G}alois symmetry.
\newblock Mimeographed notes, MIT, 1970.

\bibitem{tits66}
J.~Tits.
\newblock Normalisateurs de tores. {I}. {G}roupes de {C}oxeter \'etendus.
\newblock {\em J. Algebra}, 4:96--116, 1966.

\bibitem{viruel00}
A.~Viruel.
\newblock ${E}\sb 8$ is a totally ${N}$-determined 5-compact group.
\newblock In {\em Une d\'egustation topologique [Topological morsels]: homotopy
  theory in the Swiss Alps (Arolla, 1999)}, pages 223--231. Amer. Math. Soc.,
  Providence, RI, 2000.

\bibitem{viruel01}
A.~Viruel.
\newblock Mod 3 homotopy uniqueness of ${B}{F}\sb 4$.
\newblock {\em J. Math. Kyoto Univ.}, 41(4):769--793, 2001.

\bibitem{washington97}
L.~C. Washington.
\newblock {\em Introduction to cyclotomic fields}, volume~83 of {\em Graduate
  Texts in Mathematics}.
\newblock Springer-Verlag, New York, second edition, 1997.

\bibitem{wells71}
C.~Wells.
\newblock Automorphisms of group extensions.
\newblock {\em Trans. Amer. Math. Soc.}, 155:189--194, 1971.
\newblock (Erratum: {T}rans. {A}mer. {M}ath. {S}oc. {\bf 172} (1972), 507).

\bibitem{wilkerson74}
C.~Wilkerson.
\newblock Rational maximal tori.
\newblock {\em J. Pure Appl. Algebra}, 4:261--272, 1974.

\bibitem{wilkerson74selfmaps}
C.~Wilkerson.
\newblock Self-maps of classifying spaces.
\newblock In {\em Localization in group theory and homotopy theory, and related
  topics (Sympos., Battelle Seattle Res. Center, Seattle, Wash., 1974)}, pages
  150--157. Lecture Notes in Math., Vol. 418. Springer, Berlin, 1974.

\bibitem{winter66}
D.~J. Winter.
\newblock On automorphisms of algebraic groups.
\newblock {\em Bull. Amer. Math. Soc.}, 72:706--708, 1966.

\bibitem{wojtkowiak87}
Z.~Wojtkowiak.
\newblock On maps from ${\rm ho}\varinjlim {F}$ to ${Z}$.
\newblock In {\em Algebraic topology, Barcelona, 1986}, pages 227--236.
  Springer, Berlin, 1987.

\bibitem{xu94}
C.~Xu.
\newblock Computing invariant polynomials of $p$-adic reflection groups.
\newblock In {\em Mathematics of Computation 1943--1993: a half-century of
  computational mathematics (Vancouver, BC, 1993)}, pages 599--602. Amer. Math.
  Soc., Providence, RI, 1994.
\newblock Proceedings of Symposia in Applied Mathematics, 48.

\bibitem{xuthesis}
C.~Xu.
\newblock {\em The existence and uniqueness of simply connected $p$-compact
  groups with {W}eyl groups {$W$} such that the order of {$W$} is not divisible
  by the square of $p$}.
\newblock PhD thesis, Purdue University, 1994.

\bibitem{zabrodsky84}
A.~Zabrodsky.
\newblock On the realization of invariant subgroups of $\pi \sb\ast ({X})$.
\newblock {\em Trans. Amer. Math. Soc.}, 285(2):467--496, 1984.

\end{thebibliography}

\end{document}